\documentclass[11pt]{article}
\usepackage{amssymb,amsbsy,latexsym,amsmath,bbm,epsfig,psfrag,amsthm,mathrsfs,stmaryrd}
\usepackage[swedish,english]{babel}
\usepackage{graphicx,tikz,esint}
\usepackage[T1]{fontenc}
\usepackage[latin1]{inputenc}
\usepackage{amsfonts}
\usepackage{amsmath}
\usepackage{amssymb}
\usepackage{caption}
\usepackage{subcaption}
\usepackage{epsf}
\usepackage{graphics}
\usepackage{graphicx}
\usepackage{latexsym}
\usepackage{psfrag}
\usepackage{relsize}
\usepackage{verbatim}
\usepackage{hyperref}
\usepackage{pgfplots}
\usepackage{float}
\usetikzlibrary{positioning}
\usetikzlibrary{arrows}
\usetikzlibrary{decorations.pathreplacing}

\usepackage{float}
\usetikzlibrary{positioning}
\usetikzlibrary{arrows}
\usetikzlibrary{decorations.pathreplacing}

\theoremstyle{plain}

\newtheorem{Thm}{Theorem}
\numberwithin{Thm}{section}
\newtheorem{Lem}[Thm]{Lemma}
\newtheorem{Prop}[Thm]{Proposition}
\newtheorem{Cor}[Thm]{Corollary}
\newtheorem{Con}[Thm]{Conjecture}
\newtheorem{Def}[Thm]{Definition}
\newtheorem{rem}[Thm]{\bf{Remark}}

\theoremstyle{definition}
\newtheorem{ex}{Example}
\numberwithin{ex}{section}

\theoremstyle{remark}

\numberwithin{equation}{section}

\newcommand{\dv}{\partial}

\newcommand{\Om}{\Omega}
\newcommand{\dbar}{\partial_{\bar{z}}}

\newcommand{\eps}{\varepsilon}

\newcommand{\R}{{\mathbb R}}
\newcommand{\Pro}{{\mathbb P}}

\newcommand{\SSS}{{\mathbb S}}
\newcommand{\C}{{\mathbb C}}
\newcommand{\Z}{{\mathbb Z}}

\newcommand{\x}{\mathbf{x}}

\newcommand{\e}{\eta}

\newcommand{\Di}{\mathbb{D}}
\newcommand{\Sii}{\mathbb{S}}
\newcommand{\Lq}{\mathscr{L}}

\newcommand{\LL}{\mathcal{L}}

\newcommand{\Af}{\mathscr{A}_{N}(\Om,h_0)}

\newcommand{\UU}{\mathcal{U}}
\newcommand{\DD}{\mathcal{D}}
\newcommand{\TT}{\mathscr{T}}

\newcommand{\Gg}{\mathscr{G}}
\newcommand{\Pp}{\mathscr{P}}
\newcommand{\Qq}{\mathscr{Q}}

\begin{document}

\title{Dimer Models and Conformal Structures}

\author{\small Kari Astala 
\\ \small University of Helsinki \and \small Erik Duse 
\\ \small Kungliga Tekniska H\"ogskolan \and \small Istv\'an Prause 
\\ \small \AA bo Akademi University \and \small Xiao Zhong 
\\ \small University of Helsinki}

\date{}
\maketitle

\begin{abstract}
Dimer models have been the focus of intense research efforts over the last years. Our paper 
grew out of an effort to develop new methods to study minimizers or the asymptotic \emph{height functions} of general dimer models and the geometry of their \emph{frozen boundaries}. 

We prove a \emph{complete classification} of the regularity of  minimizers  and frozen boundaries for \emph{all dimer models} for a natural class of polygonal 
domains much studied in numerical simulations and elsewhere. In particular, we show that the frozen boundaries are always algebraic curves.
Our classification also implies that the Pokrovsky-Talapov law holds for all dimer models at a generic point on the frozen boundary and in addition shows a very strong {local rigidity} of dimer models which can be interpreted as a  \emph{geometric universality} result. Indeed, we  prove a converse result, showing that any geometric situation for any dimer model is, in the simply connected case,  realised  already by the lozenge model. 

To achieve these goals we develope a new study on the \emph{boundary regularity} for a class of Monge-Ampère equations in \emph{non-strictly convex} domains, of independent interest, as well as a new  approach to minimality for a general dimer functional.  In the context of polygonal domains, we give the first general 
results for the existence of \emph{gas domains} for minimizers. 
\end{abstract}

\subsection*{Keywords } 
Dimer models, frozen boundary, liquid region, height function, Pokrovsky-Talapov law, Monge-Amp\`ere equation, Beltrami equation, conformal structure.

\tableofcontents

\bigskip

\section{Introduction}\label{sect:introduction}

\subsection{Nontechnical overview of the paper and the main results.}\label{sect:introduction1}

Dimer models are two dimensional statistical mechanical models of perfect matchings on planar bipartite graphs, also known as random tiling models. They form the so called \emph{determinantal point processes} and are therefore sometimes also referred to as \emph{free fermion models}. In recent years they have been the focus of intense research in particular in the field of \emph{integrable probability}, drawing on connections to algebraic combinatorics and representation theory. One of the most striking features of these models seen in numerical simulations (see for example Figures \ref{fig.first}, \ref{figQuasifrozen1}) are the geometrically induced phase transitions when the systems become large, giving rise to \emph{frozen, liquid and gas phases} also known as \emph{frozen, rough and smooth phases}. Therefore much effort is spent on understanding various asymptotic properties as the systems become infinitely large, a highly challenging task as these statistical models are \emph{highly correlated.} This has been done primarily using steepest descent methods in complex analysis, and also by \emph{Riemann-Hilbert methods.} While this approach has been very successful for special classes of domains and probability measures it has not been able to treat more general cases. 

Another approach, going back to W. Thurston \cite{Th90} is based on associating a Lipschitz function, known as \emph{height function}, to each perfect matching configuration. Using concentration of measure and large deviation principles it has been shown (see \cite{CKP,Kuchumov,Tassy20}) that almost surely every random height function converges to a deterministic height function. Furthermore, this limiting height function is characterized as the unique solution of a {\it variational problem} \eqref{ELbasic} on class of Lipschitz functions with gradient constraints.  In more detail this is described in Subsection \ref{sect:introduction2} of this introduction. This characterization has the advantage of applying to any domain and any dimer model and is the starting point of the present work. 

Using the variational approach we give a \emph{complete characterization} of the regularity of  the frozen boundaries, as described in  Definition \ref{LiqFrozen},  for all dimer models having a periodic weight structure on the edges of the bipartite graph. In particular we will show in this generality that the boundaries  are all algebraic curves, classify their singularities and  show that the only singularities that can occur for \emph{any dimer model} are either first order cusps or first order tacnodes. This is a very strong and surprising rigidity result for general dimer models and as far as we know was not anticipated in the literature. Indeed, this result can be seen as a \emph{geometric universality result} for frozen boundaries of dimer models. For other further aspects of universality, see Theorem \ref{Second.thm} below. Moreover, when the curves bound a simply connected domain we give a complete characterization of their rational parametrizations.

 These results are stated in Subsection \ref{sect:introduction3} of this introduction and proven in Section \ref{propersec}.
 Furthermore, we also achieve a complete classification and give a detailed description of the regularity of the height function at the frozen boundaries, showing that the local rigidity for dimer models extends to height functions as well. In the generic case its {\it gradient is Hölder-$1/2$} continuous up to the frozen boundary, a phenomenon known in the physics literature as the Pokrovsky-Talapov law. This is the first time the law is proven rigorously in full generality for an entire class of models. This is important also since the local regularity of the height function is conjecturally assumed to determine the local \emph{universal stochastic scaling limits} at frozen boundaries. Indeed, our classification was essential in the recent work \cite{AH21,Huang21} of Aggarwal and Huang, proving Airy process universality for \emph{lozenge models} using the results of the present work.

 In addition, we also prove the first general existence result for gas regions within weighted dimer models. Finally, the regularity results lead in a natural way to a  conjecture, see Section \ref{conjectures}, about the possible universal random scaling limits that can occur at the frozen boundary of a dimer model.

In numerical simulations one typically studies random tilings on polygonal domains
related to the particular dimer model. Furthermore, on each such domain one has very special type of boundary values for the associated height function. To capture these features we introduce in Section \ref{sect:terminology} the notion of \emph{natural domains} and \emph{natural boundary values}. These notions will be essential also as they will allow us, using the work \cite{DeS10} of D. De Silva and O. Savin, to prove that in these settings the gradient of the height function is a proper map, the analytic property equivalent for the boundary of the liquid domain to be frozen. 

Naturally, an essential ingredient of the variational problem \eqref{ELbasic} are the \emph{surface tensions} in the dimer functional. As we will show it is the special properties of the surface tensions that are responsible for many of the features of dimer models. It is a striking result in \cite{KOS} that all surface tensions associated to dimer models solve the Dirichlet boundary value problem for a Monge-Ampère equation on a compact convex polygon. Because the convex polygons are not uniformly convex,
  this takes us outside the standard theory of Monge-Ampère equations. For that end we prove new results regarding, in particular, boundary behaviour of solutions of Monge-Ampère equations, which are also interesting in their  own right. These results are stated in Subsection \ref{sect:introduction5}  and proven in Section \ref{MAcomplex.structure}. 

Our line of attack of the variational problem in this paper is to reduce the Euler-Lagrange equation on liquid domains to a first order system, the so called Beltrami equation familiar from the theory of quasiconformal mappings.
The Beltrami equation has appeared earlier in connection  with 
some dimer models, see 
e.g. \cite{KeHoney} or \cite{RES17}, in describing the complex structure of the liquid domain. In our work, however, we take the Beltrami equation as the fundamental tool to  determine the various  geometric properties of the liquid domains.  This  also 
allows us to approach general dimer models. 
It is the Beltrami equations that endow each liquid domain with an \emph{intrinsic conformal structure}. In addition, in our setting  these  equations are \emph{nonlinear degenerate elliptic systems} and therefore fall outside the standard theory of Beltrami equations. This is again due to the special properties of the surface tensions, but on the other hand,  they give rise to many beautiful, specific and surprising features for these associated  Beltrami equations. The new properties developed and the  related results are stated in
 Subsection \ref{sect:introduction4}  and proven in Sections \ref{B+hodo} - \ref{propersec}. Using a hodograph transform we show how these equations can be linearized and this will lead to the important notion of \emph{teleomorphic maps} that will allow one to parametrize all algebraic curves that arise as frozen boundaries. We also believe that these results are of independent interest in the theory of degenerate elliptic PDE systems in the plane. 

Due to the boundary behaviour of the surface tension of any  dimer model, the associated energy functionals (see Section \ref{sect:terminology}\,)
 \emph{are not differentiable}. This causes 
problems, in particular in finding an effective way to characterize minimizers. One consequence in particular is that compared to other variational problems with gradient constraints, in the case of \eqref{ELbasic} the gradient constraint cannot be reduced to a \emph{double obstacle constraint}. To circumvent this issue we give in Section \ref{section:minimality} a sufficient condition, {\it the frozen star ray property,} for the minimality. Using this sufficient condition we prove a surprising universality result for dimer models, namely that for any dimer model \emph{every liquid domain with a frozen boundary} is the liquid domain (with frozen boundary) also for the standard lozenge model. This result is stated in Subsection \ref{sect:introduction6} 
of the introduction and proven in Section \ref{section:minimality}.

  Finally, two other works relevant for us and studying the properties of  the variational problem for dimer models are the work of Kenyon and Okounkov in \cite{KeOk07} and D. De Silva and O. Savin in \cite{DeS10}.  
 In \cite{DeS10}, a very important partial (local) $C^1$-regularity result is proven, of great importance to our work. This result however does not give a classification of the boundary regularity of height functions but applies to a very general class of variational problems in the plane, beyond dimer models, and does not use properties specific to dimer models. 
In the work \cite{KeOk07}  
on the other hand, 
the authors showed for 
 \emph{lozenge tilings} that a family of arctic Jordan curves for this model are (special) algebraic curves. The method of \cite{KeOk07} applies only for the lozenges  where the  spectral curve associated to the model $P(z,w) = z + w -1$ is linear; for a discussion on this c.f.  \cite[p. 271]{KeOk07}.

 In this paper we cover all dimer models, with spectral curve of any degree. As a particular example, we cover now the arctic curves of the important case of the classical  \emph{domino tilings}, for which the corresponding spectral curve is quadratic. Furthermore, our method allows for dimer models with gas phases. 
 In contrast, the proof  of \cite{KeOk07}  based on the variational approach 
  constructs minimizers by algebraic geometry methods. Unfortunately, the argument of \cite{KeOk07} is  incomplete and at this point it is not clear how to fix it. The issue is that in  \cite{KeOk07} one  implicitly assumes that the gradient constraint of the present variational problem can be reduced to a double obstacle constraint, which need not be true. 
  This is explained in more detail  in Subsection \ref{sect:introduction8}, with references to explicit counterexamples.

\subsection{Introduction to the variational problem for dimer models.}\label{sect:introduction2}

\begin{figure}[H]
\hspace{4em}
\includegraphics[scale=0.45]{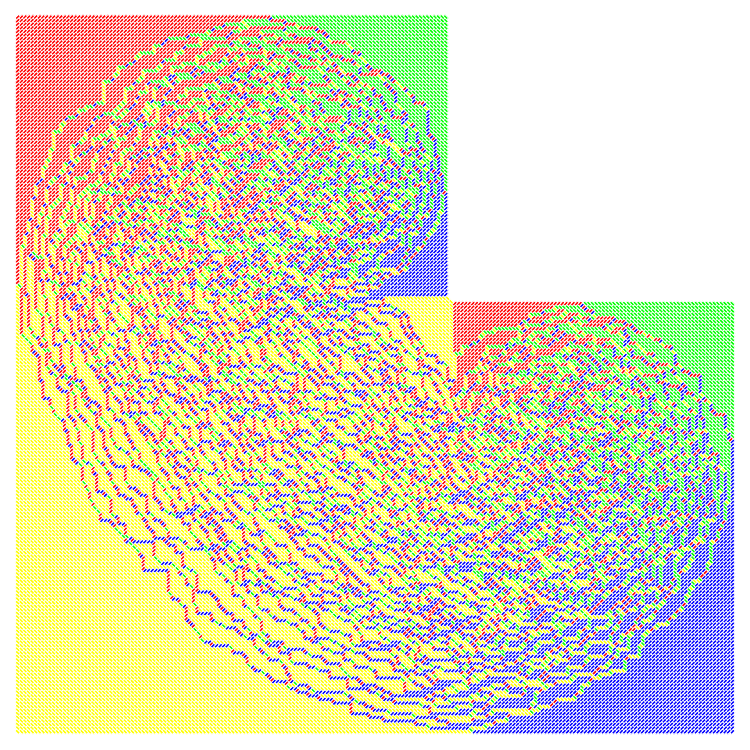}
\hspace{3em}
\includegraphics[scale=0.43]{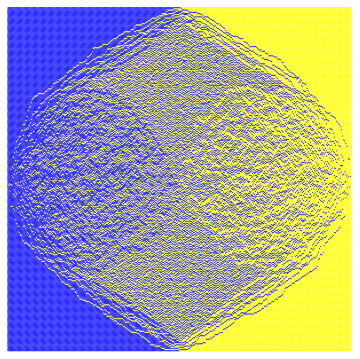}
\caption{ \, Left: A uniform domino tiling.
Right: Weighted domino tiling with gas phases. Image courtesy of S. Chhita and T. Berggren.}
\label{fig.first}
\end{figure}

The theme of this paper are the  geometric properties associated to general (bipartite) dimer models and their asymptotic limit configurations. 
Particular special cases, widely studied in  literature, include for instance the random tilings by lozenges, and weighted or unweighted domino tilings. 
 Our purpose is to present a new  systematic  and uniform approach to  the geometry of their arctic boundaries, such as those in Figure \ref{fig.first}. This will show, for instance, that the geometry of simply connected limit domains with such boundaries is universal and independent of the particular dimer model.

{ \em Dimer models} are certain two dimensional random lattice models 
in statistical physics. 
More precisely,  consider perfect matchings on finite subgraphs $G=(V,E)$ of an infinite bipartite and doubly periodic planar graph. Here $V$ denotes the vertex set of $G$ and $E$ denotes its edge set. A perfect matching of a graph $G$ is a choice of edges  that covers of all the vertices of the graph exactly once. On the set of all perfect matchings $\mathscr{M}_G$ of a graph $G$ one then associates a Gibbs probability measure $\mu$ defined by
\begin{align*}
\mu(M)=\frac{1}{Z}\prod_{e\in M}w(e), \qquad M \in \mathscr{M}_G,
\end{align*} 
where $w:E\to \R_{>0}$ is a 
weight function which is periodic on the set of edges $E$, and where 
\begin{align*}
Z=\sum_{M\in \mathscr{M}_G}\prod_{e\in M}w(e)
\end{align*} 
is a normalisation constant called the partition function, defining the probability distribution of the matchings on $G$. 

 Thurston \cite{Th90}   associated  to each perfect matching a Lipschitz function, called the {\it discrete height function}, 
  which completely describes the combinatorics of the matching, see  e.g. \cite{CKP}, \cite{CEP}, \cite{KOS}, or also Example \ref{microdominoes}. Therefore, one can also view random perfect matchings as random Lipschitz functions.
  For more on the combinatorics of dimer models we refer the reader e.g. to \cite{G19}, \cite{Ke}, \cite{KOS} and \cite{Th90}.

Under a limiting boundary condition, Cohn, Kenyon and Propp \cite{CKP} 
show, see also \cite{Kuchumov}, 
that when one rescales the size of the graph so that the edge lengths tend to zero while the macroscopic size of the  subgraph $G = G_n$ remains fixed,  
 the discrete height function converges almost surely to a deterministic  asymptotic height function $h$.
 In this setting one observes  \cite{G19},  \cite{JohHav},  \cite{Ke} fascinating phenomena for the limiting random surfaces, the graphs of the respective height functions:  For suitable polygonal domains, the surfaces present  ordered and disordered - or frozen and liquid -  parts, c.f. Figure \ref{fig.first}. For a precise definition of these notions see Definition \ref{LiqFrozen},  and for an intuitive discussion c.f. Remark \ref{add.this} with Subsection \ref{subsubsect:frozen}.
 
 The study of these asymptotic random surfaces, or {\em limit shapes}, started in the 90's with the Aztec diamond and the regular hexagonal lozenges tilings  \cite{CLP},  \cite{EKLP1}, \cite{EKLP2}, \cite{JPS}, \cite{Th90}. In particular it is shown in \cite{CKP,KOS} that the asymptotic height function is a solution to the following \emph{variational problem}
\begin{align} \label{ELbasic}
\inf_{u\in \Af}\int_{\Om}\sigma(\nabla u)dx,
\end{align}
where $\Om$ is a Lipschitz domain which is the ``limit'' of the sequence of graphs $G_n$, and $\Af$ is the class of admissible  Lipschitz functions with appropriate boundary value $h_0 = h_{| \partial  \Omega}$; for a definition  see below or Section \ref{sect:terminology}. The energy function $\sigma$ in the integrand is called the {\it surface tension}. Furthermore in \cite{KOS} the surface tensions $\sigma$ that arise in dimer models are completely characterized as follows. To each dimer model there is an associated 
\emph{Kasteleyn matrix} $K(\zeta,\omega)$ 
over a fundamental domain of the periodic bipartite graph (where the 
periodicity is with respect to the edge weights of the graph) depending on two complex parameters; for an overview of the Kasteleyn theory see e.g. \cite{Ke}. 
From here one defines the \emph{the spectral curve} $P(\zeta,\omega)=\det(K(\zeta,\omega))$ of the dimer model. From the spectral curve one defines the \emph{Ronkin function}
\begin{align*}
R(x,y)=\frac{1}{(2\pi i)}\int_{\mathbb{T}^2}\log P(e^x\zeta,e^y\omega) \frac{d\zeta}{\zeta}\frac{d\omega}{\omega}.
\end{align*}
Taking the \emph{Legendre transform} of the Ronkin function gives the surface tension $\sigma$. In particular the domain of $\sigma$ (as a proper convex function) is the \emph{Newton polygon} of $P$. It is furthermore shown in \cite{KOS} that every spectral curve of a dimer model is a so called \emph{Harnack curve}. These have several equivalent characterizations, all related to some maximality property of Harnack curves. The characterization most useful to us is that their surface tension solves a certain type of boundary value problem for a 
 \emph{Monge-Ampère equation},  namely
\begin{equation}\label{MAharnack}
\det\big(D^2\sigma\big)=\pi^2 \quad \text{in }\ N^\circ\setminus \Z^2, \text{ \; with $\sigma$ piecewise affine on the boundary $\dv N$}.
\end{equation}
 Here  $N = N(P)$ denotes the \emph{Newton polygon} of the spectral curve $P$ and $N^\circ$ its interior. In particular $N$ is a compact convex polygon. For the present work, our starting point  is the variational principle \eqref{ELbasic}, with  surface tensions slightly more general than in \eqref{MAharnack}. The key condition is that the Hessian determinant is \emph{constant} -- we derive from this property detailed information of the frozen and liquid regions.
Thus we do not use  spectral curves directly,  their role appears rather only via the
Monge-Amp\`ere equation \eqref{MAharnack}.
It is convenient to assume  the surface tension to be
defined in an arbitrary compact convex polygon in $N \subset \R^2$. Let  $\mathscr P=\{p_1,...,p_k\}$ be the vertices of $N$ and  
in addition,  allow the existence of (an arbitrary) finite subset $\mathscr G=\{ q_1,...,q_l\} \subset N^\circ$ of {\it gas points}.
For simplicity of presentation we normalise by $1/\pi^2$,  in comparison to \eqref{MAharnack}.
Then the surface tension  $\sigma: N\to \R$ in this generalised setting  is  the function bounded and convex in $N$ and solving the following Monge-Amp\`ere equation in the Aleksandrov sense
\begin{equation}\label{Pst}
\begin{cases}
\det\big(D^2\sigma\big)=1+\sum_{j=1}^l c_j\delta_{\{q_j\}} \quad &\text{in }\ N^\circ,\\
\sigma= L &\text{on }\ \partial N,
\end{cases}
\end{equation}
where $L$ is continuous and piecewise affine on $\partial N$, see Subsection (\ref{subsect:generalized}) for the discussions of the above equation.
Moreover, $c_j > 0$ are positive  and $\delta_{\{ q_j\}}$ are  Dirac masses at points
$q_j, \, j=1,2,...,l$. 
Finally, we let $\mathscr Q =  \{p_{k+1},...,p_{k+m}\} \subset \partial N$  denote the possible discontinuities of the derivative of the boundary map $L$, and call these the {\it quasi-frozen points}. For an explanation  of this terminology see Subsection \ref{subsubsect:frozen}. Because the domain $N$ is not \emph{strictly convex}, in this setting the Dirichlet problem for the Monge-Ampère equation \eqref{Pst} has in fact  not been studied before in the literature, and the standard techniques used for the Monge-Ampère do not apply \cite{F}. We therefore introduce new methods in Section \ref{MAcomplex.structure} to give a characterization of the boundary behaviour of its solutions, i.e. the boundary behaviour of surface tensions associated to general dimer models.  

For the lozenges model, $N = N_{_{Lo}}$ is the triangle with corners $\{ (0,0), (0,1), (1,0) \} $ and in this case \cite{CKP}
discovers an explicit expression for the surface tension $\sigma = \sigma_{_{Lo}}$, c.f. Subsection \ref{subsect:generalized}. In fact, all surface tensions $\sigma$  from \eqref{Pst} are very singular on the boundary of their respective domains $N = N_\sigma$, see Theorem \ref{sigma1} below. Thus 
$\sigma$ has no other  convex extension but  $\sigma_{\big | \R^2 \setminus N} \equiv \infty$. Accordingly, in  the variational problem \eqref{ELbasic} it is necessary to allow only the admissible functions,  Lipschitz functions on $\Omega$ with gradient constrained by $\nabla u(z) \in N$ a.e. and with boundary value $u_{\big | \partial \Omega} = h_0$.  It is this class of  functions that we denote by $\mathscr{A}_N(\Omega,h_0)$. However, the singular behaviour of $\sigma$ makes the functional \eqref{ELbasic} \emph{nondifferentiable}, and it is this fact that makes the variational problem for dimer models so challenging.

The above gradient constraint is, of course, consistent with the fact that for $N$ a Newton polygon, the (discrete) gradient of the discrete height function take values in the corners and quasifrozen points of $N$, and hence  their weak limit, the asymptotic height function, satisfies almost everywhere $\nabla h(z) \in N$. 

 A convenient and general notion, covering the liquid and frozen parts  of the  limiting random surfaces of  dimer models, comes as follows;  see also Remark \ref{add.this} below and  Subsection \ref{subsubsect:frozen}.

 \begin{Def}\label{LiqFrozen}
  If  $h$ is the minimizer of \eqref{ELbasic} as above, then the \emph{liquid region} of $h$ is the open subset
\begin{equation}\label{def:liquiddomain2}
\LL:=\{z\in \Omega: \text{$h$ is differentiable at $z$, with }\,\nabla h(z)\in N^{\circ}\setminus \mathscr{G}\}.
\end{equation} 
Further, we say a subset $F \subset \partial \LL$ is \emph{frozen} if  
\begin{equation}\label{frozen.def}
\nabla h(z) \to \partial N \cup \mathscr G  \qquad {as} \quad z \to F, \; z \in\LL.
\end{equation} 
 \end{Def} 
 
 Here (an in the sequel) \eqref{frozen.def} is a short hand notation meaning that  ${\rm dist}(\nabla h(z_j), \partial N \cup \Gg) \to 0$ for every sequence $z_j \in \LL$ with ${\rm dist}(z_j, F) \to 0$.
 Note also  that with Definition \ref{LiqFrozen}    we can have either a part or all of the boundary of $\LL$ frozen, even if it is the latter case that  interests us most. There we also have the direct but surprisingly useful characterisation,
\begin{equation} \label{frozen bdry}
 \partial \LL \; \mbox{ is frozen} \quad   \Leftrightarrow \quad  \nabla h : \LL \to N^\circ \setminus \Gg \mbox{\, is a proper map.}
\end{equation}
 Here we recall that a map $f$ on a domain $\mathcal U$ is proper, if $f^{-1}(K)$ is compact in $\mathcal U$ for every compact set $K$ in the image of $f$. This has the equivalent characterization that $f$ maps boundary sequences to boundary sequences, which then connects \eqref{frozen.def} and \eqref{frozen bdry}.
 
 For  \eqref{def:liquiddomain2} - \eqref{frozen bdry}  to make sense,  and for $\LL$ to be indeed open, we apply the work \cite{DeS10} of De Silva and Savin, who developed the basic regularity theory  for the variational problem \eqref{ELbasic}, for  a general $\sigma$ bounded and strictly convex in $N^\circ \setminus \mathscr G$.
 We will recall their fundamental work in more detail in Section \ref{sect:terminology}. Their Theorem \ref{thm:DeS22} below  implies that  inside the liquid region $\LL$ the minimizer $h$ is  $C^\infty$-smooth, and  thus satisfies 
the Euler-Lagrange equation
\begin{equation}
\label{eq:EL34}
\text{div}\, \big(\nabla \sigma(\nabla h)\big)= 0 \quad \text{in } \LL.
\end{equation} 

This gives also yet another view on the frozen boundary, as  the subset of $\overline{\Om}$ where the Euler-Lagrange equation \eqref{eq:EL34} degenerates. 
Further, we show in Subsection \ref{subsect:topologyFiniteConn} that for any minimizer of \eqref{ELbasic}, the liquid domains with frozen boundary have finitely many components, with each  component  finitely connected.

\begin{rem} \label{add.this}
A glance at the computer simulations, such as Figures \ref{fig.first} or \ref{figQuasifrozen1}, shows (at least intuitively)  that the liquid region corresponds to the part where fluctuations can occur, while  the frozen region arises from the part of the graph  
which becomes completely ordered in the limit. The latter fact is perhaps not so obvious from the above definition, but will become clear by the analysis  below; see Theorem \ref{thm:main2}  and \eqref{general}. Note also that for convenience we 
 include the boundaries of the gas phases to the frozen boundary.
For more exact and detailed descriptions of these correspondences see e.g. \cite{JohHav}.
\end{rem}

\subsection{Results regarding the classification of frozen boundaries and the regularity of height functions: Algebraic curves, classification of singularities and the Pokrovsky-Talapov law.}\label{sect:introduction3}

With the notions introduced in Subsection \ref{sect:introduction2} we can state our first  general results, which are also interesting  to compare with  simulations as e.g. in Figure  \ref{fig.first}.

\begin{Thm} \label{First.thm}
Suppose $ \LL  \subset \C$ is a bounded finitely connected domain, and  $h$ 
 a solution to \eqref{eq:EL34}  in $\LL$, where $\sigma$ satisfies \eqref{Pst}. If $\partial \LL$ is frozen 
 for $h$, i.e. if  $\, \nabla h : \LL \to N^\circ \setminus \mathscr{G}$ is proper, then
 \medskip

a)   $\partial \LL $ is the real locus of an algebraic curve  (minus the isolated points of  the  curve). If  $\LL$ is

\quad simply  connected, then  $\partial  \LL  = R(\Sii^1)$,   the image of the unit circle  under a rational map $R(z)$. 
\medskip

b)  There are at most finitely many singularities $\{\zeta_j\}_{j=1}^m  \subset \partial \LL $ on the boundary, and they are all 

\quad either first order (inward) cusps or tacnodes. 
\vspace{.3cm}

c)  At every  $\zeta \in \partial \LL$, excluding the cusps and tacnodes,  the boundary $\partial \mathcal{L}$ is locally strictly convex: 
 
\hspace{.3cm}  $ B(\zeta,\varepsilon) \cap \LL$ is  strictly convex for all $\varepsilon > 0$ small enough. 

\noindent In particular, the above holds whenever  $\LL$ is liquid for any dimer model,  with or without gas,  with the boundary  $\partial \LL$ frozen.
\end{Thm} 

Moreover, in Theorem \ref{thm:h:proper} and Theorem \ref{thm:h:proper2} we show that oriented polygonal boundary conditions, natural for the given dimer model, produce liquid domains with frozen boundaries, thus  
 the conclusions of Theorem 1.2 hold in this setting. In fact, we will reach a complete classification, in terms of their rational parametrisation, of all simply connected domains $\LL$ that have frozen boundary in the sense of \eqref{frozen bdry} - \eqref{eq:EL34}, see Theorem \ref{characterize}. 

\begin{rem}
 The above Theorem 
 describes the ``cloud curve'' geometry of the frozen boundary. Cloud curves have been introduced in \cite{KeOk07} in connection with the lozenge model. 
Related works are \cite{Bu16,Boutillier-Li} where the authors study fluctuations in domino tilings and in the dimer model on the square-hexagon lattice, respectively. The boundary conditions they considered fit into our framework but are more restrictive -- in both  works the frozen boundaries are found to be cloud curves. This general phenomenon is explained by 
Remark \ref{rem:cloud-curves}.
\end{rem}

A main issue in \cite{KeOk07} is the lack of proof of boundary continuity in general for solutions to the complex Burger equation, see \cite[Sect.2.4 and Prop.2]{KeOk07}. The authors can prove this only by an explicit construction using tools from algebraic geometry, in the special case of lozenges with uniform weights and only on simply connected polygonal domains with $3d$ sides changing cyclicly. Therefore, for general dimer models on general polygonal domains, in this paper we have chosen a different approach with emphasis on complex analytic PDE-methods and calculus of variations. With this approach, for instance, the required boundary continuity  is quickly obtained in any dimer model. We can also give a rigorous derivation of the famous Pokrovsky-Talapov law for general dimer models, c.f  Theorem \ref{thm:main2} and \eqref{Pokrovsky-Talapov} below. \medskip

\begin{Thm}[Pokrovsky-Talapov law]
\label{thm:main2}
Let $h$ be a solution to the Euler-Lagrange equation \eqref{eq:EL34} in a bounded domain $\LL$, with $\sigma$ as in \eqref{Pst}.  Suppose  that $\partial \LL$ is frozen,
i.e.  $ \nabla h : \LL \to N^\circ \setminus \Gg$ is a proper map.   Then

\smallskip

a) There is a finite  set $\{z_j\}_{j=1}^n\subset \partial \LL $, such that for every $z_0\in \partial \LL \setminus \{z_j\}_{j=1}^n$ we have 
\begin{equation} \label{limitp}
 \lim_{z\to z_0, \, z\in \LL} \nabla h(z)=p_0, \quad {\rm where } \; \; p_0 \in \mathscr{P}\bigcup\mathscr{Q}\bigcup\mathscr{G}. 
 \end{equation}

b)  If   \eqref{limitp}  holds for $p_0 \in \mathscr{P}$ and $z_0 \in  \partial \LL$, then there is \, $p\in N^\circ\setminus \mathscr{G} \,$ such that, first,
the  vector

 \hspace{.3cm} $p-p_0$ is normal  
  to  $\partial \LL$ at  $z_0$, 
 and second, that for the arguments
\begin{equation} \label{dir.amoeba}
 \lim_{z\to z_0,z\in\LL}\arg \big(\nabla h(z)-p_0 \big)= \arg(p-p_0).
 \end{equation}
\vspace{-.4cm}

 c)  Outside the above finite set of singularities the minimizer $h \in C^{1,1/2}$, with no better H\"older 
 
   \quad exponent  
  at any point of  $\partial \LL$ { \,{\rm (The Pokrovsky-Talapov law) }.}
 \end{Thm}

In other words, Theorem \ref{thm:main2} tells us that the gradient of the height function is continuous up to the boundary from within the liquid region, except for a finite number of points on the frozen boundary. Here \eqref{dir.amoeba}  gives  the
  {\it limiting direction}  of $\nabla h(z)$ as $z \to z_0$ inside $\LL$, c.f. Figure  \ref{fig.amoeba}.   Finally, outside a finite set of points the gradient $\nabla h$ is Hölder continuous up to the boundary with Hölder coefficient $\frac{1}{2}$.

The Pokrovsky-Talapov law is a general physical law describing transition between commensurate and incommensurate phases in a large class of statistical mechanical models. In particular, for crystal surfaces it predicts a generic exponent of $3/2$ at the transition between a smooth and a rough phase, (e.g. between liquid and frozen or liquid and gas phases). In our case in view of the representation formula \eqref{general}, to be explained below, the relevant quantity is the tiling density which is controlled by the gradient of the height function. Here in fact,   if $z_0 \in \partial \LL$ is a regular boundary point with $\nabla h(z) \to p_0$ as $z \to z_0$,
then  c) of Theorem \ref{thm:main2} can be expressed more precisely as follows: if $n_{z_0}$ is the inner normal to $\LL$ at $z_0$, then
\begin{equation} 
\label{Pokrovsky-Talapov}
h(z_0 + \delta \, n_{z_0}) - h(z_0) -  \delta \,  \langle p_0, n_{z_0} \rangle \,  \simeq  \, \delta^{3/2}, \qquad {\rm as} \; \delta \to 0^+,
\end{equation}
see  Subsection \ref{proofPTalapov}. 

At the points where $\nabla h$ fails to be continuous or fails to be $C^{1/2}$, by Theorem \ref{main.belt.regularity} in Section \ref{sect:introduction4}  one still has control of   $\partial \LL$ and  $\nabla h$  at the singularities, for instance, via an explicit representation of $\nabla h$ described  in \eqref{general}. On the other hand, in addition to the cusps and tacnode singularities arising from the geometry of  $\LL$,  there are also finitely many points on $\partial \LL$ where $\nabla h$ oscillates between two corners of $N$, thus $h$ failing to be $C^1$.

Our approach allows  also an analysis of the gas phases, where $\nabla h \equiv q$ with $q \in \mathscr{G} \subset N^\circ$. There the asymptotic tile correlations  decay exponentially, as opposed to polynomial decay in liquid parts, see \cite[Section 4]{KOS}.  The existence of gas domains has been shown before for some special surface tensions $\sigma$ and domains $\Omega$, see e.g.  \cite{BefChiJoh}, \cite{Bgren}. Here we have the existence  in the general setting.

 \begin{Thm}\label{GAS} Suppose $\Omega \subset \C$ is a bounded and simply connected Lipschitz domain and $h_0$ an admissible boundary value on $\partial \Omega$. Assume that the minimizer $h$ of  \eqref{ELbasic}, for a surface tension $\sigma$ as in \eqref{Pst}, admits a liquid domain $ \LL \subset \Omega$ with $\partial \LL$ frozen.  
 
 If $\sigma$ has gas points  $q \in  \Gg$, then for each $q$ there is (at least one)  non-empty  gas domain $U_q \subset \Omega$. More precisely, 
$U_q$ is open and simply connected, $ \partial U_q $ is one of the components of $\partial \LL $ and
$$
\nabla h \equiv q \;\; {\rm in } \;\; \overline{U_q} \quad {\rm with} \quad \nabla h(z) \to q \quad {\rm as } \quad z \to \partial U_q \quad  {\rm in } \quad \LL.
 $$
  \end{Thm}
  
  The geometry of the gas domains $U_q$ is inherited from $\partial \LL$ and Theorem \ref{First.thm}, for more details see  Theorem \ref{MultiplyConect2}. In particular, unless there are tacnodes,
  $\partial U_q$ is  concave outside the three or more cusps, all directed outwards from $U_q$.

\subsection{Results regarding universality of the Lozenge model.}\label{sect:introduction6}

The PDE-approach with analysis of the complex structures associated to dimer models allows us also  to uncover results  in  different directions, such as the universality in the geometry of dimer models and their frozen boundaries.
Quite unexpectedly, this goes beyond the limit surfaces of discrete (dimer) lattice models, as indicated by the following result, to be proven in Section \ref{subsection:proofoftheoremmain3}.

\begin{Thm}[Universality of frozen boundaries] \label{Second.thm}

Let $\LL \subset \C$ be a bounded Jordan domain and  $\sigma$ any  surface tension as in \eqref{Pst}, with $\mathscr{G} = \emptyset$. 
Suppose  that the Euler-Lagrange equation \eqref{eq:EL34} admits a    solution $ h$  in $\LL$, such  
that   $\, \nabla  h : \LL \to N^\circ$ 
 is a continuous and proper  map. 
 
Then  $\LL$ is liquid with $\partial \LL$ frozen for the Lozenges model.
That is, there exists a polygonal domain $\Omega \supset \LL$ and a piecewise affine boundary value $h_0$ on $\partial \Omega$,   such that if $\, h^\star \in C^1(\Omega)$ is the minimizer  of Lozenges model in $\Omega$, i.e. 
$$\int_\Omega \sigma_{_{Lo}}(  \nabla \, h^\star)  = \inf \left\{ \int_\Omega \sigma_{_{Lo}}(\nabla u):   u_{|{\partial \Omega}} = h_0, \; \nabla u(z) \in N_{_{Lo}} \; \; {\rm for} \; z \in \Om \right\},$$
where $\sigma_{_{Lo}}$ is the surface tension for Lozenges, with domain $ N_{_{Lo}}$, then 
\begin{equation} \label{unilozenges}
\LL \equiv \{ z \in \Om :
  \, \nabla \, h^\star(z)   \in (N_{_{Lo}})^\circ  \}, \quad {\rm and} \quad  \nabla \, h^\star (z) \to \partial N_{_{Lo}}  \;\;  {\rm as}\;  \;  z \to \partial \LL \;\;  {\rm in}\;  \; \LL. 
     \end{equation}
\end{Thm}

 Once \eqref{unilozenges} is shown to hold, then \cite{CKP} and \cite{KOS} prove that $\, h^\star$ is also equal to an asymptotic random lozenges  height function, on a scaling limit of finite subgraphs
converging to $\LL$.  
 In particular, this means that within Jordan domains, all possible geometries of frozen boundaries of  all dimer models occur  already for the Lozenges model !
 
  For multiply connected domains, and in particular for dimer models with gas points, the universality holds {\it locally}, i.e. locally any frozen boundary of any dimer model is locally frozen also for the lozenges model.  
  For details, see Remark \ref{loc.reg}.
  
 \subsection{A short roadmap.} 
We finish this introduction by providing an overview of the logical roadmap of the paper;  see also the next section for an overview of the  tools from Geometric Analysis which we need to develop in order to describe the various properties of  dimer models covered in this work.
Here, in particular, the properness of different (gradient) maps arises as a  very flexible notion in analysing and describing the frozen boundaries.

After describing the basic concepts and notations in Section \ref{sect:terminology}, we represent in Section \ref{MAcomplex.structure}  the surface tension $\sigma$ in harmonic coordinates. This, in turn, leads to a definition of an intrinsic complex structure on liquid regions, described in terms of an appropriate Beltrami equation.
That also allows us  to represent, for any dimer model, all height functions as compositions of harmonic maps and this complex structure. The Euler-Lagrange equation of the variational problem \eqref{ELbasic} then becomes equivalent to a Beltrami equation of the form \eqref{Specific}.

In view of  Definition \ref{LiqFrozen} and our approach to frozen boundaries, this motivates and requires a detailed study of proper solutions of the universal Beltrami equation, which 
we present in Sections \ref{B+hodo} and \ref{propersec}. 
Unpacking these   in terms of the height function lead then  to our results on the variational problem stated in Section \ref{sect:introduction}. In particular, we arrive at algebraicity of frozen boundaries from purely the properness property, as well as at  detailed geometric descriptions of the frozen boundaries.

With a'prior regularity results of the variational problem developed in Section 4 and the method of frozen extensions from Section 8 we also show that in natural domains the minimizers of \eqref{ELbasic} for (oriented) natural boundary values have the property that the boundary of the corresponding  liquid region is frozen, i.e. the solution of the Beltrami equation \eqref{Specific} and its universal counterpart \eqref{belt2} are proper maps. 

Section \ref{section:minimality} is of a complementary nature and provides a sufficient condition for a PDE solution to be a minimizer. This is then used to establish the universality of the lozenge model within all frozen boundaries.

\subsection*{Acknowledgements}
The work of K.A. was supported by ERC grant 834728 and Academy of Finland grant 13316965; part of this work was done when K.A. was visiting Mathematical Sciences Research Institute, Berkeley.   E.D. was supported by the Knut and Alice Wallenberg Foundation KAW grant 2016.0416 and  Academy of Finland grant 12719831, I.P.  by Academy of Finland grants 1266182, 1303765, 13316965 and 365297,  and X.Z. by ERC grant 834728 and Academy of Finland grant 308759.

Our interest in the topic originated in part from a question of Andrei Okounkov related to the analysis of equation \eqref{beltrami2345}. In addition, we thank Rick Kenyon, Stas Smirnov, Kurt Johansson and Konstantin Izyurov for related discussions. 
We want to thank Franc Forstneric for discussion on \cite{FW}, Vadim Gorin for pointing out the crucial reference \cite{DeS10} and Ovidiu Savin for clarifications 
regarding Theorem 1.4 in \cite{DeS10}. 
For discussions and clarifications regarding parametrizations of real algebraic curves we thank Björn Gustafsson, Kathlén Kohn and Kaie Kubjas. Moreover, we thank Tomas Berggren and Maurice Duits for providing us with many figures of simulations of lozenge and domino tilings. We would like to especially thank Sunil Chhita for making many domino tiling simulations for us. Finally, we thank the referees for many valuable comments.

\section{Methods from Geometric Analysis}\label{sect:MGA}

 For the theorems described in the  introduction, we need to apply and develop a multitude of results from geometric analysis and PDE's. In fact, for many of the above theorems, the known and well established methods are not sufficient, and therefore a number of new results from geometric analysis, PDE's and variational calculus need to be established. For the  help of the reader and before the more technical considerations, in this section we give an overview of these tools  necessary for the proofs of results in Theorems   \ref{First.thm}, \ref{thm:main2}, \ref{GAS} and \ref{Second.thm}.

 \subsection{Degenerate Beltrami equations}
\label{sect:introduction4}

 A key point   
 in our approach to general dimer models is that, as  will be shown in Section \ref{MAcomplex.structure}, in simply connected domains $\LL$ there is an explicit one-to-one correspondence between solutions of the Euler-Lagrange equation \eqref{eq:EL34} and those of the specific Beltrami equation \eqref{32kaksi23}, associated to the surface tension $\sigma$, see Theorems \ref{connection} - \ref{converse}. 
 This relation provides a powerful spectrum of tools  to analyse the minimizers,  compared for instance to the
basic methods   from calculus of variations. With these we cover all dimer models with periodic weight structure and, as an example, the boundary continuity of the gradient $\nabla h: \LL \to N^\circ \setminus \Gg$ then quickly follows  even in domains with tacnodes on the boundary, see Theorem \ref{propStoilow35} and Remark \ref{added}.

To explain our approach in more detail, recall that in two dimensions in particular, the classical Leray-Lions equations
 \begin{equation} \label{LionsLeray}
  \text{div}\, \mathcal A(\nabla u) = 0 \qquad {\rm in \; a \; domain \, \; \mathcal U},
\end{equation}
with a monotone and, say,  $C^1$-smooth structure function $\mathcal A$, are closely related  to a complex system of PDE's,  the quasilinear Beltrami equations, see e.g. \cite{ATG}, \cite{BN},  \cite{Cacc},  \cite{IS01}, and for  recent new aspects also  \cite{ACFJK}. 

This holds, of course, also when $u = h$, the minimizer $h$ of \eqref{ELbasic} satisfying the Euler-Lagrange equation \eqref{eq:EL34} in the liquid region. However, the special fact that the surface tension $\sigma$ is a solution to the particular Monge-Ampere equation \eqref{Pst} allows us to develop this much further. For instance, we will use a non-linear and invertible explicit expression,  the so-called {\it Lewy-transform} $L_\sigma$, 
with definition recalled in  \eqref{LewyT}, and  show in Theorem \ref{connection} that for any such  $\sigma$ the composition
$ f := L_\sigma(\nabla h)$
satisfies the  Beltrami equation
\begin{eqnarray}\label{32kaksi23}
f_{\overline{ z }}   = {\mathcal H'_\sigma}( f) f_{ z} ; \qquad  {\mathcal H}_\sigma(w) :=  (I-\nabla \sigma) \circ (I + \nabla \sigma)^{-1}(\overline w), 
\end{eqnarray} 
where ${\mathcal H} = {\mathcal H}_\sigma $ 
is complex analytic! Here analyticity also allows the derivative ${\mathcal H'_\sigma}$ in the first equation above. In addition, see Proposition \ref{Hproper} for details, for any surface tension $\sigma$ the function ${\mathcal H}_\sigma'$
 is a proper map
from its domain $Dom({\mathcal H}_\sigma)$ onto the open unit disc $\Di$ (hence \eqref{32kaksi23} is e.g. not uniformly elliptic).
As we will see, these special properties will lead to strong conclusions. In the sequel we will call  ${\mathcal H}_\sigma $ the  
{\it structure function} associated to the surface tension $\sigma$. 

Quite remarkably, the above relation between $f$ and $h$ has also a converse: In a simply connected domain $\mathcal U$, for every $C^1$-solution  $f$ to \eqref{32kaksi23} the relation $\nabla h(z) = L_\sigma^{-1}\bigl(f(z)\bigr)$ defines a solution to  the Euler-Lagrange equation  \eqref{eq:EL34}. 

All these intimate relations, between solutions to the Euler Lagrange equation  \eqref{eq:EL34} and
 those to the Beltrami equation \eqref{32kaksi23}, ask us to systematically develop the properties of 
 this Beltrami equation. 
 Much of Sections \ref{B+hodo}  and \ref{propersec} \,  is devoted for this purpose, with eye on results having immediate implications for the geometry of liquid domains and frozen boundaries.
 
 In cases where one uses solutions $f: \LL \to Dom({\mathcal H}_\sigma)$ of  Beltrami equations such as in \eqref{32kaksi23} to describe the geometry, structure or regularity of a frozen boundary, there is yet another symmetry one can apply, namely  the {\it conformal invariance} properties of these equations. Namely, given any solution  $f_0$ to  $f_{\overline{ z }}   = {\mathcal H}'_\sigma( f) f_{ z}$, with the chain rule we see that the function $f:=  {\mathcal H}'_\sigma (f_0)$ satisfies the equation 
\begin{eqnarray}\label{beltrami2345}
 f_{\overline{ z }}(z) = f(z)  f_z(z), \qquad z \in \LL.
 \end{eqnarray} 
 We call this the {\it universal Beltrami equation} for the dimer models, first since every solution to a specific equation as in \eqref{32kaksi23}   determines also a solution to \eqref{beltrami2345}. But second, there is also more to this terminology: As we will see in Subsection \ref{sect:introduction7} and in full detail in Section \ref{MAcomplex.structure}, composing ${\mathcal H}'_\sigma$ with a Riemann map makes  \eqref{32kaksi23}  in case of the Lozenges model  equivalent to the universal equation \eqref{beltrami2345}. These two properties are actually the starting points for the universality of the Lozenges geometry described  in Theorem \ref{Second.thm}.

 Additionally, given a minimizer or a solution $h$ to \eqref{eq:EL34}, with the above procedure we can associate to it a solution $f$ of the universal equation, explicitly $f = {\mathcal H}_\sigma' \circ L_\sigma \circ \nabla h$, and now the entire boundary $\partial \LL$ is frozen for $h$ - the most interesting case suggested by 
 simulations -  if and only if   this  $f: \LL \to \Di$
  is a  proper map, c.f. Theorem \ref{eqNonLinB}. 
 Even more, see Corollaries \ref{uniqueSC} and \ref{uniquef},  proper maps $f:\LL \to \Di$ solving \eqref{beltrami2345} in a bounded domain $\LL$ are {\it unique}: given $\LL$ there is  at most one such function.

On the other hand, as mentioned above, the universal Beltrami equation \eqref{beltrami2345} is a degenerate elliptic PDE, and thus many of its interesting properties are  not covered by the standard methods or the existing litterature. Therefore one of the key goals of this work is to provide a deep understanding of these equations. To illustrate their role further, we present  here a selection of results and consequences that we will later prove. We restrict here to simply connected domains, where these attain a particularly transparent form.

 \begin{Thm} \label{New.stuff1} Suppose $\mathcal L \subset \C $ is a bounded and  simply connected domain. Then the following are equivalent:
   \medskip
  
 a) $\mathcal L$ is a liquid domain with frozen boundary for some dimer model with $\Gg = \emptyset$.
 
   \medskip
 b) For some solution $h$ of the Euler-Lagrange equation \eqref{eq:EL34} in $\LL$, where $\sigma: N \to \R$ is a surface tension as in \eqref{Pst} with  $\Gg = \emptyset$, the gradient $\nabla h:\mathcal L \to N^{\circ}$ is a proper map.
 
     \medskip
  c) There is a solution $f: \LL \to \Di$ to the Beltrami equation \,$\partial_{\overline{z}} f(z)   = f(z) \, \partial_z f(z)$ such that $f: \LL \to \Di$ is a proper map.
   \medskip
  
 d) $\mathcal L$ is a liquid domain with frozen boundary for the Lozenges model.
 
   \end{Thm}
   \smallskip
   
   Indeed, Theorem \ref{New.stuff1} is a  consequence of our studies in the subsequent chapters. More precisely, $a) \Rightarrow b)$ follows from Definition \ref{LiqFrozen}, $b) \Rightarrow c)$ is a consequence of Theorem \ref{key.connection}, 
 $c) \Rightarrow d)$ follows from Remark \ref{domH}, Corollary \ref{Dconverse} and Theorem \ref{Second.thm}, while the last implication $d) \Rightarrow a)$ is clear.
 
  Note that above the simple connectivity of $\LL$ requires $\Gg = \emptyset$, c.f. Theorem \ref{GAS}. Also note that in $b)$ we do not ask that $h\in \Af$, as required in Definition \ref{LiqFrozen}.

 Appropriate versions of the above equivalences hold also for multiply connected domains and surface tensions with gas or quasifrozen points, as well as for domains which have a partially frozen boundary.

 We thus see that  the proper solutions to the universal equation \eqref{beltrami2345}  determine in an explicit manner the  properties of all  frozen boundaries. Developing the different aspects of these solutions  we will arrive at the following Theorem. The description here holds for every (multiply or simply  connected) bounded domain and for every dimer model, and gives strong rigidity for the geometry of their limit shapes.
 
 Here and throughout this work, $W^{1,2}_{loc}(\Omega)$ stands for the Sobolev space of measurable functions which together with their distributional derivatives are locally $L^2$-integrable in the domain $\Omega$. For  an overview of the properties of such functions and Sobolev spaces in general, see e.g. \cite[Sections A.1 - A.7]{ATG}.

\begin{Thm} \label{main.belt.regularity} Suppose  $\LL \subset \C$ is a  bounded domain and $f \colon  \LL   \to \Di$    a continuous and proper map, contained in  $W^{1,2}_{loc}(\LL)$. If $f$ is also a solution to universal equation
 \eqref{beltrami2345}, then
  \medskip
  
 a)  $\partial \LL$ is the real locus of an algebraic curve, with   properties a) - c) of Theorem \ref{First.thm}.
  \medskip

b)   The map $f \colon \mathcal{L} \to \mathbb{D}$ is real analytic inside $ \LL  $ and extends continuously up to  $\partial \LL$. 

 \medskip

c) The tangent vectors $\tau(\zeta)$ of $\partial \LL$ and  the boundary values $f(\zeta) \in \partial \Di$ are  related via the identity
\begin{align} \label{tau2}
f(\zeta)=-\tau(\zeta)^2, \qquad \zeta \in \partial \LL \setminus \{{\rm cusps} \}.
\end{align} 
 
 \vspace{-.2cm}

d)  $f \in C^{1/3}( \, \overline{ \LL \, }\, )$. 
\medskip

e) Moreover, \, $f \in  C^{1/2}$ \, in \,  $\overline{ \LL \, }\, \setminus \bigcup_{j=1}^m B(\zeta_j, \varepsilon),$ where $\{ \zeta_j \}$ are the cusps of $\partial \LL$.
\medskip

f) At each cusp singularity of  $ \partial  \LL $,  there is a line $\ell$ transversal to the cusp at $\zeta_j$,  such that  
$$f \in C^{1/3}(\ell \cap B(\zeta_j,\varepsilon) )$$
\quad \ for $\varepsilon>0$ small.

\quad   However,  in the direction $\tau$ of the cusp, $f \in C^{1/2}(\tau \cap B(\zeta_j,\varepsilon) )$  for $\varepsilon>0$ small.
\end{Thm} 

It is this result precisely, and the connection between $f$ and the asymptotic height function
$h$, that we are able to prove Theorem \ref{First.thm} as a direct consequence. This interaction has also several further consequences, discussed in this and the subsequent sections.

This relation  also gives an explicit finite dimensional parametrisation of all simply connected liquid domains with frozen boundary, having a given number of  cusps on $\partial \LL$.  Moreover,  in Theorems \ref{thm:geometry-curve} and \ref{MultiplyConect2} we show that  the number of cusps is explicitly determined by the degree of $f$ on the respective boundary component.

For general domains $\Omega$ and general admissible boundary values $h_0$, it may happen that not all of the boundary of the liquid region $\mathcal{L}$
defined  in
\eqref{def:liquiddomain2} is frozen. However, some part of the boundary may still exhibit frozen phenomena in the sense of \eqref{frozen.def}, and one then speaks of  the {\it frozen part}  
$\partial_F \LL $  of  the boundary $ \partial \LL$. 

In  complete  analogy with Theorems \ref{First.thm} and \ref{main.belt.regularity},  also the frozen parts can be approached with this method.
This gives them strong regularity properties.

\begin{Thm}[Local regularity of frozen boundaries]
\label{thm:localregboundary}
Let $\LL$ be a bounded  domain and $f: \LL \to \Di$ a continuous $W^{1,2}_{loc}$-solution to the equation
\,$ \partial_{\overline{ z }} f(z) = f(z)  \partial_z f(z)$.

Suppose  $\Gamma \subset \dv \LL$ is connected,   $\gamma$ is a smooth crosscut of  $\LL$ and  that $\gamma \cup \Gamma$ forms the boundary of a simply connected domain contained in $\LL$.
Assume further that $|f(z)| \to 1$ as  $ z \to \Gamma$ in $\LL$.

 Then $\Gamma$ is an analytic curve. Moreover,  the classification of singularities in Theorem \ref{First.thm}  carries over to this case.

\end{Thm}
According to   Theorem \ref{thm:localregboundary},  we have $\Gamma = b(-1,1)$ for a function analytic and locally injective on the interval. Note, however, that this
 does not exclude cusps or tacnodes. Also note, that given a limiting height function or a solution $h$ to \eqref{eq:EL34},  via the auxiliary map  $f = {\mathcal H}_\sigma' \circ L_\sigma \circ \nabla h$ the above properties of a locally frozen boundary hold also for $h$.
 \medskip

\subsection{Results regarding properties of surface tensions and the Monge-Ampère equation.}\label{sect:introduction5}

The above result allows one to analyse the liquid domains from inside, up to their frozen boundaries, and lead e.g. to understanding and classification of the boundary geometry. However,  for the universality of the Lozenges geometry, as described in Theorem \ref{Second.thm}, we have to go beyond the boundary and construct polygonal domains $\Omega \supset \LL$ and boundary values on $\partial \Omega$, so that the original boundary $\partial \LL$ is frozen for \eqref{ELbasic} also in this new setting. For this, for instance, 
we need to develop in Section \ref{MAcomplex.structure}  quite detailed and specific properties of  the surface tensions $\sigma$ defined in \eqref{Pst}. 
  The results are also of independent interest for the study of Monge-Amp\`ere equations, considering the specific properties that arise when uniform convexity of solutions and domains is lost on the boundary.  
\begin{figure}[H]
\centering
\begin{tikzpicture}[xscale=1,yscale=1]
\draw[thick,dashed,->,>=stealth] (-2.3,2.3)--(-3,3);
\draw[thick,color=orange,->,>=stealth] (0,0.4)--(0,0.9);
\draw[thick,color=teal,->,>=stealth] (-1,0)--(-0.5,0);
\draw[thick,color=blue,domain=-1.6:1.6]   plot (cosh{\x},sinh{\x});
\draw[thick,color=blue,domain=-1.6:1.6]   plot (-cosh{\x},sinh{\x});
\draw[thick,color=blue,rotate=90,domain=-1.6:1.6]   plot (cosh{\x},sinh{\x});
\draw[thick,color=blue,rotate=-90,domain=-1.6:1.6]   plot (cosh{\x},sinh{\x});

\draw[thick,color=blue] (0.4,-0.1)--(0.4,0.1);
\draw[thick,color=blue] (-0.1,-0.4)--(0.1,-0.4);
\draw[thick,color=blue,domain=0:1,scale=0.4]   plot ({\x},{(1-\x^3)^(1/3)});
\draw[thick,color=blue,domain=0:1,scale=0.4]   plot ({\x},{-(1-\x^3)^(1/3)});
\draw[thick,color=blue,domain=0:1,rotate=90,scale=0.4]   plot ({\x},{(1-\x^3)^(1/3)});
\draw[thick,color=blue,domain=0:1,rotate=-90, scale=0.4]   plot ({\x},{-(1-\x^3)^(1/3)});

\draw (3,0) node {$\mathcal{A}=\nabla \sigma(N^\circ \setminus \{g\})$};
\draw (-3,3.4) node {$\propto\,\, i(p_2-p_3)$};
\draw (-0.75,-0.3) node {\tiny$\nu(z)$};
\filldraw[color=teal] (-1,0)  circle (1.5pt);
\draw (-1.3,0) node {$z$};
\draw (-1.95,0.7) node {$\partial \sigma(p_3)$};
\end{tikzpicture}
\begin{tikzpicture}[xscale=1,yscale=1]
\draw[thick] (2,0)--(0,2)--(-2,0)--(0,-2)--(2,0);
\filldraw[color=blue] (0,0)  circle (1.5pt);
\filldraw[color=blue] (2,0)  circle (1.5pt);
\filldraw[color=blue] (0,2)  circle (1.5pt);
\filldraw[color=blue] (-2,0)  circle (1.5pt);
\filldraw[color=blue] (0,-2)  circle (1.5pt);

\draw[thick,color=orange,->,>=stealth] (0,0)--(0,0.5);
\draw[thick,color=teal,->,>=stealth] (-2,0)--(-1.5,0);
\draw (2.3,0) node {$p_1$};
\draw (0,2.3) node {$p_2$};
\draw (-2.3,0) node {$p_3$};
\draw (0,-2.3) node {$p_4$};
\draw (1,2.3) node {$N$};
\draw (0,-0.5) node {$g$};
\end{tikzpicture}
\caption{The figure on  the right displays the gradient constraint $N$. The figure on the left displays the image of $\nabla \sigma(N^\circ\setminus\mathscr{G})$, also known as the amoeba in the literature. The vector $i(p_2-p_3)$ indicates the direction of one of the asymptotes of the amoeba. Furthermore, the vectors at $p_3$ and $g$ in $N$, c.f. Theorem \ref{sigma1}, show  directional limits that get mapped to the boundary of the amoeba under the map $\nabla \sigma$.}
\label{fig.amoeba}
\end{figure}

  That there is for any positive $c_j$, any gas points $q_j$ and any piecewise affine $L$,  admitting a convex extension
 to $N$,  a unique solution to \eqref{Pst} follows from Theorem 1.1 in \cite{Har}.
 However, neither the domain $N^\circ$ nor the boundary values of $\sigma$ are strictly convex, and in this case the detailed boundary behaviour of $\sigma$ 
 does not seem to have been studied in the existing literature. 
 The question is not easy, however comparison methods and, in particular, developing an analysis with the  Lewy transform, see \eqref{LewyT}, leads to the following first step.
 
 \begin{Thm} \label{sigma1}
 Let $\sigma$ be a  surface tension solving \eqref{Pst} in a polygon $N$, and let $\mathscr P,
  \mathscr Q \subset \partial N$ and $\mathscr G \subset N^\circ$  be the corresponding  corners,  quasifrozen points and gas points in $N$, c.f. \eqref{corners} - \eqref{gas}.
   \smallskip

a) Suppose $J \subset \partial N$ is a closed interval not containing any of the points of $\mathscr P \cup \mathscr Q$. Then 
          \begin{equation}  \label{flat}
    \    |\nabla \sigma(p)| \to \infty \qquad {as} \quad p \to J, \;  p \in N^\circ.     
             \end{equation}

  \smallskip
  
b) If $p_0 \in \mathscr P\cup \mathscr Q\cup\mathscr G$, then for every point $p\in  N^\circ\setminus {\mathscr G}$, the limit
\begin{align}
\label{GeneralizedGradient37}
 \lim_{\tau\to 0^+} \nabla \sigma\big( p_0+\tau(p-p_0)\big) =: \widehat{\nabla} \sigma({p_0}, p-p_0) 
\end{align}
exists and is finite. In addition, for each fixed $p_0 \in \mathscr P\cup \mathscr Q\cup\mathscr G$, the limits   give a  continuous one-to-one correspondence between the directions or arguments $\arg(p-p_0)$, where $p\in  N^\circ\setminus {\mathscr G}$, and the points on an  analytic curve, the boundary of the subdifferential $\partial \sigma(p_0)$. 
  \end{Thm}

 For the notion of a {\it subdifferential} of a convex function see \cite{Sch}, or  Section \ref{subsect:generalized} below. We have, in fact, quite  precise further information of the boundary behaviour of $\sigma$ and its subdifferentials. As an example, in \eqref{GeneralizedGradient37} the vector $p - p_0$ is always  (an outer)
normal to $\partial \sigma(p_0)$ at the boundary point $\widehat{\nabla}\sigma({p_0},p-p_0)$, c.f. Figure \ref{fig.amoeba} above. For details see
 Sections \ref{subsect:generalized} and  \ref{surface}.

  \vspace{.2cm}
 We will often consider  the gradient map 
$\nabla \sigma:  N^\circ\setminus {\mathscr G}\to \R^2$ of the surface tension. 
Since $\sigma$ is strictly convex in $N^\circ\setminus {\mathscr G}$, 
the gradient  is a homeomorphism in its domain. 
Lemma \ref{blowup} and the above discussion show that the complement of  the image domain $\nabla \sigma( N^\circ \setminus \mathscr G)$ in $\R^2$ is precisely the disjoint union   of the subdifferentials 
  $$
  \bigcup_{p \in \mathscr P\cup \mathscr Q\cup\mathscr G} \partial \sigma(p),
  $$
 a finite union of closed convex sets.   
 In the terminology of \cite{KeOk07}, the  {\it amoeba} associated to the surface tension $\sigma$ is the closure 
 ${\overline{\, \nabla \sigma( N^\circ \setminus \mathscr G) \, }} \subset \R^2$. In Figure \ref{fig.amoeba} above, see e.g. \cite{CJ} or the later Example \ref{microdominoes},  on the right we have the gradient constraint $N^\circ\setminus {\mathscr G}$ for  weighted domino tilings, a square punctured at the gas point at the origin, and on the left, the amoeba, the image domain $\nabla \sigma( N^\circ \setminus \{ g \})$.

\subsection{Results regarding representation formulae and regularity for limiting height functions.}\label{sect:introduction7} 

The conformal invariance of the Beltrami equation \eqref{32kaksi23} is useful also for finding convenient representations for the limiting height functions and solutions to the Euler-Lagrange equations \eqref{eq:EL34}. The representations  are particularly  transparent 
 in the case where there are no gas points,  i.e. the gradient constraint $N$ is just the convex hull of its corners and quasifrozen points $\{ p_j\}_{j=1}^k$ on $\partial N$. 
 Indeed, then also $Dom({\mathcal H}_\sigma)$ is simply connected, c.f. Lemma \ref{analyticH}.
 
 Thus given any solution $f_0: \LL \to Dom({\mathcal H}_\sigma)$ to \eqref{32kaksi23},
 we can 
 take the Riemann map $\psi: \Di \to Dom({\mathcal H}_\sigma)$ and define $ \mu_\sigma(z) := {\mathcal H}_\sigma' \circ \psi(z)$ with  $ f(z) := \psi^{-1} \circ f_0(z)$, for $z \in \Di$.
 The arrangement gives us  a solution to 
  \begin{eqnarray} \label{belt478} 
   f_{\overline{ z }}  = \mu_\sigma \bigl(  f \bigr) \, f_z, \qquad  f:\LL \to \Di,
       \end{eqnarray} 
 an equation equivalent to  \eqref{32kaksi23}. Note that here $\mu_\sigma: \Di \to \Di$ is analytic and proper, by the properties of the structure function ${\mathcal H}_\sigma$. Thus, in fact, $\mu_\sigma(z)$ is simply a Blaschke product. For lozenges $\mu_\sigma(z) = z$, while for dominoes $\mu_\sigma(z) = z^2$, see Remark \ref{domH}.

In this setup, 
given now a $C^1$-solution $h:\LL \to \R$ to the Euler-Lagrange equation  \eqref{eq:EL34} with $\sigma$ as in \eqref{Pst},  then by Subsection \ref{sect:introduction4} the relation $f = \psi^{-1} \circ L_\sigma \circ \nabla h$ defines a solution to \eqref{belt478}. 
 It turns out, see Theorem \ref{repone}, that the functions are  related by the explicit identity
 \begin{eqnarray} \label{general}
\nabla h (z) = \sum_{j=1}^k \, p_j \, \omega_{\Di}\bigl( f(z); I_j \bigr), \qquad  z \in \LL,
      \end{eqnarray} 
 where $I_j \subset \partial \Di$ are disjoint open intervals depending  only on the surface tension $\sigma$ (and the Riemann map $\psi(z)$), with closures  covering the unit circle, and where $\omega_{\Di}(\zeta; I_j)$ is the harmonic measure of $I_j$ in the unit disc. 
  
 In addition, if $\partial \LL$ is frozen for $h$ in the (rather weak) sense \eqref{frozen.def}, then $ f: \LL \to \Di$ is proper and extends continuously to $\partial \LL$. It hence follows that $\,\nabla h(z) \to p_j\,$ as $z$ approaches a point on  $f^{-1}(I_j) \subset \partial \LL$. For $j= 1,\dots,k\,$ this family of arcs covers $\partial \LL$ up to a finite set. 
 
 Even further, the representation \eqref{general} leads to a very suggestive  probabilistic interpretation.   
 Namely, in the Aztec diamond the well known work of Cohn, Elkies and Propp   \cite{CEP}  describes the asymptotic tile density of a domino, i.e. the limiting probability for it to occur at a given place of the liquid domain. In that case $\LL$ is the unit disc $\Di$ and $N$ is the convex hull  of the points $\{ \pm 1, \pm i \}$,  see Example \ref{microdominoes}. 
 We will show   in Example \ref{density} that in this setup and for $p_1 = i$ the north corner of $N$, the asymptotic tile density from 
 \cite{CEP} is exactly equal to the term $\omega_{\Di}\bigl( f(z); I_1 \bigr)$ in \eqref{general}, where  $f:\LL \to \Di$ is a proper map solving \eqref{belt478}. Furthermore, this proper map $f:\LL \to \Di$ is  unique up to a conformal isomorphism of $\Di$ preserving $\mu_\sigma$, see Corollary \ref{uniqauto}., and thus for any given  dimer model or surface tension $\sigma$, the terms $\omega_{\Di}\bigl( f(z); I_j \bigr)$ are intrinsic to $\LL$.
 
 This raises the natural {\it conjecture} that, at least in the absence of quasi-frozen and gas phases,  for all dimer models and liquid domains with frozen boundary the asymptotic tile or edge densities are given  by the pull back of the appropriate harmonic measure, under coordinates $f$ satisfying \eqref{belt478}. In case of dominoes this is related to  \cite[Conjecture 13.1]{CKP} of Cohn, Kenyon and Propp. Our conjecture suggests, in particular, a canonical expression for the limiting  tile or edge densities, valid for all dimer models.  
 
   Combined with the classical Stoilow factorization theorem \cite[Theorem 5.5.1]{ATG}, the identity  \eqref{general} allows an even more detailed picture. Namely, we have   $f = B \circ G$, where $B(z)$ is a finite Blaschke product of $\Di$ and $G: \LL \to \Di$  a {\it homeomorphic} solution to $G_{\overline{ z }}  =  \mu_\sigma \bigl(  f \bigr) \, G_z$.
  With this decomposition \eqref{general} gets the form
  \begin{eqnarray} \label{generalB}
\nabla h (z) = \sum_{j=1}^k \, p_j \, \omega_{\Di}\bigl( G(z); B^{-1}I_j \bigr), \qquad  z \in \LL,
      \end{eqnarray} 
where now the arcs $B^{-1}I_j  \subset \partial \Di$ depend on the boundary values of $h$, and the coordinate $G(z)$ gives a complex structure for $\LL$.
In fact,   according to the famous conjecture of Kenyon and Okounkov \cite{KeHoney}, \cite{KeOk07} the fluctuations of the  random height function around its expectation $h(z)$ are given by the Gaussian Free Field with respect to the very same complex structure  which $G(z)$ determines.

  Representations analogous to \eqref{general}  hold also for a general dimer model (or surface tension) where $\sigma$ has gas points or quasi-frozen points, see Remark \ref{myygen} 
  and the identity \eqref{hrepre3456} there. All this gives a direct suggestion that  there  should be a proper probabilistic interpretation  to all terms in \eqref{hrepre3456}, in case the solution is a minimizer to \eqref{ELbasic} and the surface tension $\sigma$ of a dimer model has gas or quasi-frozen points.

  \subsection{Results regarding the variational problem, natural boundary values and frozen extensions.}\label{sect:introduction71}

The above considerations, and Definition \ref{LiqFrozen} in particular, show that the problem of understanding  and finding liquid domains $\LL$ with frozen boundary
  is exactly  a free boundary problem defined by the gradient constraint  $\nabla u \in N^\circ \setminus \Gg$.  Also, by Theorems \ref{thm:DeS12} and \ref{thm:DeS22}
the free boundary inside $\Omega$ is precisely $\partial \LL \cap \Omega$.

The question of frozen boundary is in spirit similar to other variational problems with gradient constraints, for example the elastic-plastic torsion problem, see \cite{CF}, or the question of space-like hypersurfaces in Minkowski space with prescribed boundary values and mean curvature, see \cite{BarSim}.  In the case of the elastic-plastic torsion problem it is possible to reduce the gradient constraint to a \emph{double obstacle problem}, see \cite{BS}, and use the highly developed methods for free boundary problems coming from obstacle constraints \cite{PSU}. However, due to the singularity of $\nabla \sigma$ at $\dv N$ as in Theorem \ref{sigma1}, this reduction is not possible in our case. Similar problems occur in \cite{BarSim}, however with the difference that there $\vert \nabla \sigma\vert$ blows up at every point of the boundary of the gradient constraint. 

This boundary singularity of the surface tension has strong implications on the study of the variational problem \eqref{ELbasic}.
There does exist a unique minimizer, see Proposition 4.5 of \cite{DeS10}, but to study its properties one cannot use 
directly the usual first variations, nor can one have a pointwise characterisation via an Euler-Lagrange variational inequality. The only characterisation of the minimizer $h$ is through the weak  Gâteaux directional derivative inequality
\begin{align} \label{gateaux1}
dI_\sigma[h;u-h]=\int_\Om d\sigma(\nabla h(z), \nabla u(z)-\nabla h(z))dz\geq 0,
\end{align}
 valid only for the $u \in \mathscr{A}_N(\Omega,h_0)$ contained  in the space of admissible functions.
 Unfortunately, the expression $dI_\sigma[h;u-h]$ is \emph{nonlinear} in $u-h$.  Therefore, this characterisation is of limited value, and we are forced to consider alternative approaches. See, however,  Theorem \ref{thm:hstar:minimizer} for a sufficient pointwise criterion for minimality.

 On the other hand, what saves us is  the strong correspondence the variational problem \eqref{ELbasic} has to the  solutions $f$ of the specific Beltrami equation \eqref{32kaksi23}, described in the previous subsections, 
 allowing many features  of $f$ to carry over to the minimizers $h$.

 Returning to \eqref{general}, it is quite remarkable that also this relation has a direct  converse, as shown in  Corollary \ref{Dconverse}. That is,  for the intervals $I_j$ as in \eqref{general}, given  any solution to the Beltrami equation \eqref{belt478}, the identity \eqref{general} defines a  Lipschitz function   $h:\LL \to \R$ which solves the Euler-Lagrange equation \eqref{eq:EL34}, with the given surface tension $\sigma$. 
    
 One can go even deeper into this picture. Namely, the above allows one to extend $h$  as a piecewise affine function to  a polygonal domain $\Omega \supset \LL$, such that  
 $\nabla h(\Omega \setminus \LL) \subset \{ p_j, 1 \leq j \leq k \}$. 
One would then like to show that this extended function, with  new boundary values $h_0$ of the extension on $\partial \Omega$, is the unique minimizer for the variational problem 
  \begin{eqnarray} \label{minimi235}
   \inf \{ \int_{\Omega}\sigma(\nabla u(z))dz : u\in \mathscr{A}_N(\Omega,h_0) \}.
      \end{eqnarray} 
That, however, is  far from trivial, and we can show this only under certain geometric restrictions.  For details and discussion 
see Section  \ref{section:minimality}. 
  \medskip

This takes us to the last  basic 
question in the theme of frozen boundaries. Given a dimer model with a periodic weight structure, for  which polygonal domains  $\Omega$ and which boundary values $h_0$ on 
$\partial \Omega$ does the minimizer in \eqref{minimi235} determine a liquid domain with frozen boundary? It is clear that if no assumption is made, this does not happen.  For example,
one can have $\Omega = \LL$ and no part of $\partial \LL$  frozen. Also, if one does not  assume special properties of the boundary values, even for $C^\infty$-smooth $h_0$  the  liquid domain can be only partially frozen,  with $\partial \LL$ exhibiting extremely complicated behaviour. For explicit examples of this latter phenomena see \cite{Duse15a}. 

In this work we  only consider  piecewise affine boundary values. Indeed, in typical simulations such as in Figure \ref{fig.first} and \ref{figQuasifrozen1}, the limiting minimizer $h$ changes color along a side of the polygon $\Omega$, that is, along the side $h$ is affine but its boundary gradient attains values at two neighbouring corners $p_j$ and $p_{j+1}$ of $N$. In particular, this makes it natural  to study polygonal domains $\Omega = \Omega_N $ for which each side is orthogonal to a side of $N$.   We call such a polygonal $\Omega$ a  {\it natural domain}, see Definition \ref{Def:naturaldomain} and Figure \ref{figGradDomino}. Associated to this one has the notion of {\it natural boundary values} $h_0$. 
On each side of $\Om$, these  are  affine  with the correct gradient, see Definition \ref{def:extremal}. The notions here are for general gradient constraints $N$ as in \eqref{Pst}, but even for triangle $N$ our class of domains and boundary values is strictly larger than the class studied in  \cite{KeOk07}.

  On the other hand, already for lozenges tilings there are simple situations with no liquid domain at all. For instance,  Example \ref{noliquid} presents a  
  piecewise affine (but non-affine) natural boundary value $h_0$, such that 
  there is only one admissible function $u\in \mathscr{A}_N(\Omega,h_0)$ 
  and its gradient takes values in the corners of $N$. However, for natural domains and natural boundary values, such phenomena are the only obstruction for the existence of liquid domains.

\begin{Thm}\label{thm:liquid:affine} Let $\Omega$ be a natural domain and $h_0$  a natural boundary value.  Suppose that $h$ is the minimizer of variational problem (\ref{ELbasic})
among the class $\mathscr{A}_N(\Omega,h_0)$ with (possibly empty) liquid domain $\LL$. 
Then $h$ is (countably) piecewise 
affine in $\Omega\setminus\overline{\LL}$, with gradient having values in $ \mathscr{P}\bigcup\mathscr{Q}\bigcup\mathscr{G}$.
\end{Thm}

  Thus either the liquid domain is non-empty or the minimizer is trivial, piecewise affine in $\Omega$. Therefore, should $\LL$ be empty, Theorem \ref{thm:liquid:affine} reduces the variational problem to a combinatorial problem.
  
  In case the minimizer of (\ref{ELbasic}) is not piecewise affine, the next question is then the structure of the liquid domain $\LL$. 
Here we show  for $N$  a triangle, that  $\partial \LL$ is  frozen 
  whenever $\Omega$ is a natural domain and  $h_0$ is natural. 
  This is a more general model than the Lozenges one: we do allow quasifrozen points $ \mathscr{Q}$ and  also gas points $ \mathscr{G}$ for $N$.

\begin{Thm}\label{thm:h:proper}
Let $\Omega$ be a natural domain, $h_0$  a natural boundary value and $N$ a triangle. Suppose that $h$ is the minimizer of variational problem (\ref{ELbasic})
among the class $\mathscr{A}_N(\Omega,h_0)$. 

Then either $h$ is piecewise affine in $\Omega$ with $\nabla h  \in \mathscr{P}\bigcup\mathscr{Q}\bigcup\mathscr{G}$, or else  there is a liquid domain $\LL \subset \Omega$ and  $\nabla h: \LL\to N^\circ\setminus
\mathscr G$  is a  proper map, i.e. $\partial \LL$ is frozen.

In particular, when $\LL \neq \emptyset$ all properties  from Theorems \ref{First.thm} and  \ref{thm:main2} hold for the minimizer $h$ and the liquid domain $\LL$. 

Similarly, if  $q \in  \Gg$ is a gas point for $\sigma$, then there is a non-empty gas domain $U_q \subset \Omega$ with $\nabla h \equiv q$ in  $U_q$. 
\end{Thm}
In the case where $N$ is not a triangle, we will show in Section \ref{sec:topology} that the conclusions of Theorem \ref{thm:h:proper} hold under   an extra condition on $h_0$, that it should be oriented in the sense of Definition \ref{def:extremalorient}, c.f. Theorem \ref{thm:h:proper2}. For an illustration of this concept in the case of domino tilings see Figure \ref{fig:OrientedCorner}. We conjecture that Theorem \ref{thm:h:proper} holds for all surface tensions $\sigma$ and gradient constraints $N$ from \eqref{Pst}, without assuming that the natural boundary value $h_0$ is oriented, and plan to return to these questions in a future work.  Similarly, we expect that in the setting of Theorem \ref{thm:liquid:affine} the minimizer is always  finitely piecewise affine outside $\LL$.

The issue with  Theorem 
\ref{thm:h:proper} is that due to the singular behaviour of the surface tension $\sigma\,$ e.g. as in \eqref{flat}, there is no general method to analyse the boundary behaviour of the minimizers $h$. To remedy this we introduce the notion of {\it frozen extensions} for the pairs $(\Omega, h_0)$,  see Definition \ref{oldold}.
 Using these with 
  Theorem \ref{thm:DeS22}, due to De Silva and Savin \cite{DeS10}, we prove Theorem \ref{thm:h:proper} in Section \ref{sec:topology}.

  Combining now the previous themes leads us back to Theorem \ref{Second.thm}.   To prove the Theorem and the universal  geometry of frozen boundaries, we first apply the relation \eqref{unisol}  to construct a solution to the Beltrami equation \eqref{beltrami2345}. Next, for the lozenges model the coefficient function  $\mu_{\sigma}(z) \equiv z$ in Equation \eqref{beltrami2345}, so that one can use the procedure described in \eqref{belt478} - \eqref{general}, giving  us  a solution to the Euler-Lagrange equation  $ \,\text{div}\, \big(\nabla \sigma_{_{Lo}}(\nabla \,h^\star)\big)= 0$ in $\LL$. These with Theorem \ref{First.thm} give  detailed information on how to attach on $\partial \LL$  polygonal sides which make a natural domain  $\Omega$, and extend  $h^\star$ from \eqref{general} to have natural boundary values on $\partial \Om$. 

Thus for Theorem \ref{Second.thm} it remains to show that  the extended  height function $h^\star$ is indeed 
the minimizer of the variational problem \eqref{ELbasic}. 
This is a very delicate question and far from obvious, due to the singular features of the surface tension $\sigma$. Our proof of the minimality of $h^\star$ involves a divergence free extension of $\nabla\sigma(\nabla h^\star)$ from $\LL$
to $\Omega$. To show that such an extension is possible we apply Theorem \ref{thm:main2}  with further results from Subsection \ref{complexsec}. Actually, via Theorem \ref{thm:hstar:minimizer} we give a general way to show that a given function is a minimizer. 
We believe that the approach is of its own interests and will have applications in other problems in Calculus of Variations.
Our argument also gives a rigorous proof of the minimality of the 
constructed  height function $h^\star$, without volume constraint,  in  Theorem 2 of 
\cite{KeOk07}, though possibly in a different polygon than the one given there.

\subsection{Relation  to other work and to the Complex Burgers equation.}\label{sect:introduction8}

In their work  \cite{KeOk07}, inspired by the connections of dimer models to algebraic geometry,  Kenyon and Okounkov  composed the partial Legendre transform of the surface tension $\sigma$ from \eqref{MAharnack}
 with the gradient of the asymptotic height function, to express the  Euler-Lagrange equations  \eqref{eq:EL34}
 (with an added volume constraint term) as the system
\begin{equation}\label{CB1}
   \left\{
    \begin{array}{ll}
      P\big(\zeta(x,y),\omega(x,y)\big)=0,\\
     \displaystyle \frac{\zeta_x}{\zeta}+ \frac{\omega_y}{\omega}=c.
      \end{array} \right.
\end{equation}
where $P$ is the spectral curve of the given dimer. In particular, in the case of the lozenge model with uniform probability the spectral curve is 
given by $P(z,w)=z+w-1$. In this case one can solve for $w$, giving for $c=0$ the \emph{complex Burgers equation}

\begin{align*}
\frac{\zeta_x}{\zeta_y}=\frac{\zeta}{1-\zeta}. 
\end{align*}

In \cite{KeOk07} they show that the Jacobian $J(e^{-cx}\zeta, e^{-cy}\omega) \equiv 0$ 
and then claim that this shows there exists an analytic function $Q$ of two variables such that $Q(e^{-cx}\zeta,  e^{-cy}\omega) = 0$. In that case,  the system \eqref{CB1} would reduce to the coupled algebraic equations 
\begin{equation}\label{CB2}
   \left\{
    \begin{array}{ll}
      P\big(\zeta(x,y),\omega(x,y)\big)=0,\\
    Q\big(\e^{-cx}\zeta(x,y),  e^{-cy}\omega(x,y)\big)=0.
      \end{array} \right.
\end{equation}
There are however a few issues with this   approach. 
First,  their argument gives $Q$, but only locally outside the critical points of the functions. A second issue arises via the Remark below, c.f. also \cite[Lecture 10]{G19}.

\begin{rem}\label{KOF}																
For every pair of  functions analytic in a domain $\UU \subset \C$, the Jacobian  $J(f,g) = f_x g_y - g_x f_y \equiv 0$. However, e.g.  Forstneric and  Winkelmann \cite[Theorem 1]{FW} show that there are (many) analytic functions $f$ and $g$ in the unit disc, such that the only analytic function $Q(z,w)$ with $Q(f,g) \equiv 0$ in the disc is $Q \equiv 0$. 
 \end{rem}

And third, in order for this approach to be useful, further properties of $Q$ must be deduced directly from the boundary values in the variational problem. This however seems very difficult and is not done in \cite{KeOk07}. 
 A particular important issue, pointed out  but left  open in  \cite{KeOk07}, is that in order to conclude that $Q$ is an analytic polynomial one needs that the functions $\zeta$ and $\omega$ are continuous up to the boundary, see \cite[Proposition 2]{KeOk07}.

  For the above reasons, 
we instead approach the  frozen boundaries by developing the properties of the special Beltrami equation \eqref{beltrami2345} and its relatives. This, for instance, allows us to show the boundary continuity  even for a general dimer model, see e.g. Theorem \ref{propStoilow35} and Remark \ref{added}, and in this way opens up a uniform  way to understand the geometry of frozen boundaries for all dimer models.

In \cite{KeOk07}, to construct solutions to  \eqref{CB2} and to produce  liquid domains $\LL$  inscribed in a polygonal domain with 3d cyclicly changing tangent directions, they 
 use  tropical algebraic geometry and a deformation 
 with respect to a Lagrange multiplier associated to a volume constraint. 
 On the other hand, this construction  
 only applies to lozenge models and only for a very special class of polygonal domains which in particular are simply connected. In their construction in \cite{KeOk07} they  get a tentative height function $h^\star$, the ansatz for the solution of the variational problem. To show that the constructed function $h^\star$ is the unique minimizer and thus equals the asymptotic height function,
  the authors estimate the directional derivatives  of the energy functional in the direction of each admissible function $g$. To show the minimality, i.e. that these G\^{a}teaux derivatives at $h^\star$ are non-negative, 
in  \cite[Section 4]{KeOk07} it is assumed that  each frozen facet  is either a bottom facet where $g \geq h^\star$  for any admissible function $g$, or a top facet where $g \leq h^\star$ for every such $g$. But that would  mean that the frozen facets
$\mathcal{F}$ are part of the so called coincidence set $ \Lambda$, see Subsections \ref{subsubsect:domains} and  \ref{subsubsect:frozen} for the precise definitions of the sets $\mathcal{F}$ and $\Lambda$.
However, as discussed in   Subsection \ref{subsubsect:frozen}, it is not true that  $\mathcal{F}$ is a subset of $\Lambda$. For a counterexample, see Figure \ref{facet} and the discussion before it; see also \cite[p. 491]{DeS10}, and for a simulation where that fails, see e.g. the concave corner of Figure 4 in \cite{MPet}.
Hence in such situations the argument of  \cite{KeOk07} remains incomplete (we thank Rick Kenyon for helpful discussions regarding this).

The failure of $\mathcal{F}\subset \Lambda$ ultimately depends on the lack of smoothness of the surface tension $\sigma$. Indeed, if $\sigma\in C^1(\overline{N})$, then a modification of the arguments in \cite{BS} proves $\mathcal{F}\subset \Lambda$.  For dimer models, however, we always  have $\sigma \notin C^1(\overline{N})$. 

A further issue is that \cite{KeOk07} constructs functions $h_n$ approximating $h^\star$, which are, in turn, built by patching together $n$ volume constrained solutions $f_j$'s. 
To estimate the energy integrals of the approximants $h_n$ one has to control the discontinuities of $\nabla \sigma(\nabla f_j)$'s along the patches, to make the total contribution negligible. Here \cite{KeOk07} uses the fact that $\| \nabla f_j \|=O(n^{-1/2})$ near the patch. As far as we can see, c.f. \eqref{gradsigma23}, this information alone gives only an $O(1/n)$ estimate for the discontinuities of $\nabla \sigma(\nabla f_j)$'s and may well add up to a non-neglible error term.
For these reasons we have chosen in Section \ref{section:minimality} a different route towards the minimization questions.

 On the other hand,  the complex Burgers equation is closely related to the Beltrami equation used in this paper, even if the basic methods  the complex Burgers equation allows are quite different. For the benefit of the reader we briefly describe the relations between these two equations. Suffice to say that this part requires familiarity with the theory developed later in this paper.  Indeed, one can approach the complex Burgers equation also from the point of view of \eqref{ELsigma}
 -   \eqref{1toista}. For this note first that once one has the stream function \eqref{stream}, then in the notation  \eqref{1toista} the Euler-Lagrange equation \eqref{eq:EL34} is  simply  the relation  \begin{eqnarray}\label{simplee}
   (F_x)_y - (F_y)_x = 0.
    \end{eqnarray} 
  In the special case where   $\sigma = \sigma_0 \,$ is the surface tension \eqref{gradsigma23} - \eqref{gradsigma45} of the lozenges model, then upon substituting the explicit form of $ \sigma_0$    to  \eqref{1toista}, we see that the partial derivatives of $F$ satisfy the special relation $\, e^{-i \frac{\pi}{2} F_y } \,  + \, e^{i \frac{\pi}{2} F_x }\,  \equiv 1 $. This suggests  one to introduce
 \begin{eqnarray}\label{legendre}
 \log \zeta \, = \, -\frac{i \pi}{2} F_y  \, = \, \pi  \, \sigma_1(\nabla h) - i \pi h_y   \quad {\rm and } \quad \log \omega \, =  \, \frac{i \pi}{2 } F_x \, = \, \pi \, \sigma_2(\nabla h) + i \pi h_x
 \end{eqnarray} 
with the function $F$ as in  \eqref{1toista}. These  functions $\zeta$ and $\omega$ are the ones used
in  \cite{KeOk07} for the lozenges model.  The Euler-Lagrange equation  \eqref{simplee} gets the form
 \begin{eqnarray}\label{burgers44}
 \frac{\zeta_x}{\zeta} + \frac{\omega_y}{\omega} = 0,
 \end{eqnarray} 
 where for lozenges  one has  additionally the relation $\zeta + \omega = 1$.
  Combining this with \eqref{burgers44}  leads to the complex Burgers equation for $u = \zeta/(1-\zeta)$. 
Moreover, the  map $\zeta$ takes the liquid domain $\LL$ to the lower half plane (as $\nabla h$ lies in the triangle with corners $\{ 0, 1, i\}$). In fact, if one  composes $\zeta$ with a M\"obius transform and sets
   \begin{equation} \label{beltburgers}
   f := M \circ \zeta, \qquad M(\zeta) = \frac{(1-i)\zeta +i}{(1+i)\zeta -i},
      \end{equation}
then a straightforward  calculation shows that $f(z)$ is a solution to \eqref{belt2} with values in  the unit disc $\Di$. 

For a general dimer model, \cite{KeOk07} uses the associated Harnack  spectral  curve $P(\zeta,\omega) \equiv 0$, 
 a certain Laurent  polynomial. Also in this case there is a natural route  arriving at the Beltrami equation.  Indeed,  
 note that $P(\zeta,\omega)$ determines a Riemann surface $\Sigma\subset \hat \C^2$ and a decomposition $\Sigma=\Sigma^+\cup \Sigma^-$, with the zero locus $\mathcal{Z}_\R(P)=(\dv \Sigma^+)\cap (\dv\Sigma^-)$ as the common boundary. Moreover, an implicit differentiation 
gives
$\,  \dv_\zeta P(\zeta,\omega)\zeta_y+\dv_\omega P(\zeta,\omega)\omega_y=0$, and 
combining this with \eqref{burgers44} leads 
to the complex Burgers equation 
\begin{align} \label{mikhalkin}
\zeta_x-\gamma(\zeta) \zeta_y=0, \qquad \gamma(\zeta) :=\frac{\zeta\dv_\zeta P(\zeta,\omega)}{\omega\dv_\omega P(\zeta,\omega)}\bigg\vert_{\omega=\omega(\zeta)}.
\end{align}
Here by a result of Mikhalkin, the {\it logarithmic Gauss map} $\gamma: \Sigma^+\to \mathbb{H}_{-}$ 
 is a proper holomorphic map, and 
because of the properties of Harnack curves, for any $(\zeta,\omega)\in \Sigma^+$ we can uniquely solve $\omega=\omega(\zeta)$ as a holomorphic function of $\zeta$. 

If now one  expresses \eqref{mikhalkin} in terms of the complex derivatives and chooses a Riemann map   $\psi: \Sigma^+\to \mathcal{D}$  to a circle  domain 
 in $\C$, then  
 $ f(z) :=\psi(\zeta(z))$ solves the Beltrami equation
\begin{align*}
\overline{\dv}f(z)=\mu(f(z))\dv f(z), \qquad \mu:=\frac{\gamma \circ \psi^{-1} + i}{\gamma \circ \psi^{-1} -i}.  \qquad  
\end{align*}

As  an example, for domino tilings the spectral curve is determined by 
   $$
     \zeta + \omega + \zeta \omega = 1.
$$
This leads to $\gamma(\zeta) = \frac{- 2\zeta}{1-\zeta^2}$, where $\zeta$ takes values in the upper half-plane. 
Thus   $\mu(\zeta) =  \left( i \frac{\zeta - i}{\zeta + i} \right)^2$, 
 c.f. also Remark \ref{domH}   below. In particular,  now $\widehat f = (M_1 \circ \zeta)^2\,$ with $M_1(\zeta) = i\frac{\zeta - i}{\zeta +i }\, $ defines a solution  $\widehat f: \LL \to \Di $  to \eqref{belt2}. Thus one can arrive at the Beltrami equation in the liquid domain $\LL$ also from this point of view, and conversely.

 Our approach is not the only way how Beltrami equations can be used to study limit surfaces.
For instance, it turns out \cite{KP1} that {\em any} strictly convex surface tension carries an intrinsic complex structure in which the Euler-Lagrange equation is expressed (for a different $f$) in the form $f_{\overline{ z }}=\mu(f)f_z$. However, the coefficient $\mu$ is analytic - the key property for our approach - only if the determinant of the Hessian of $\sigma$ is constant and we have the dimer models. Hence the properties of the solutions are different from those studied in this paper. The latter point of view is useful in studying limit shapes when the condition \eqref{Pst} fails, as is the case for the five-vertex model, see \cite{KP2}.

\addtocontents{toc}{\vspace{-4pt}}
\section{Terminology and Preliminary Results}\label{sect:terminology}

\subsection{Surface tensions}\label{subsect:generalized}
In the rest of the paper, $N$ is a compact convex polygon in $\R^2$ and $\sigma$ a solution to a Monge-Amp\`ere equation
in $N$. Also, let $L:N\to \R$ be a convex function which is piecewise affine on $\partial N$. Finally, let $q_1,...,q_l \in N^\circ$ with $c_1,...,c_l$  positive numbers.

Then we call a {\it surface tension} any convex function $\sigma: N\to \R$   such that
\begin{equation}\label{Pst2}
\begin{cases}
\det\big(D^2\sigma\big)=1+\sum_{j=1}^l c_j\delta_{\{q_j\}} \quad &\text{in }\ N^\circ;\\
\sigma= L &\text{on }\ \partial N,
\end{cases}
\end{equation}
where $\delta_{\{q_j\}}$ is the Dirac mass at the point $q_j$, $j=1,...,l$.  This setting has three finite sets of special points, which all will have a special role throughout this paper: 
\begin{eqnarray} 
\mathscr P& = & \{p_1,...,p_k\}   \quad \textrm{are the corners of} \; N. \label{corners} \\
\mathscr Q& = & \{p_{k+1},...,p_{k+m}\} \quad \textrm{are the quasi frozen points}, 
\textrm{ points on} \; \partial N \setminus  \mathscr P \; \label{qfrozen}
\\ && \textrm{\hspace{.1cm}  where} \; L \;  \textrm{is not differentiable. } \nonumber \\
\mathscr G& = & \{ q_1,...,q_l\} \quad \textrm{are the gas points, locations of the Dirac masses in \eqref{Pst2}. }\label{gas}
\end{eqnarray}
The terms ``quasi frozen'' points and ``gas points'' are explained by their role in the geometry of
 the height functions, see Subsection \ref{subsubsect:frozen}.  From now on we give for $\mathscr P \cup \mathscr Q = \{p_j: 1 \leq j \leq  k+m\}$ the  cyclic numbering or order induced by $\partial N$, in the counterclockwise direction.

It follows from Theorem 1.1 in \cite{Har}, that for every possible data as in  \eqref{corners} - \eqref{gas},  the equation \eqref{Pst2} admits a unique convex solution. Also, $\sigma \in C^{\infty}(N^\circ \setminus \Gg)$, c.f. \cite{Schulz} or Section \ref{complexsec}. As an example, in the case where $N = N_{_{Lo}}$ is the triangle with corners $\{ 0, 1, i\}$, the boundary value $ L = 0$ and there are no gas or quasifrozen points, then one is studying the lozenge tiling model. By Theorem 8 in \cite{Ke},  the gradient $\nabla \sigma = \nabla \sigma_{_{Lo}} $ of the lozenge surface tension is explicitly given  by
    \begin{equation}  \label{gradsigma23}
    \pi  \nabla \sigma_{_{Lo}}(s,t) = \left(  \log\left(\frac{\sin(\pi s)}{\sin(\pi(t+s))} \right), \, \log\left(\frac{\sin(\pi t)}{\sin(\pi(t+s))} \right) \right), \quad (s,t) \in \text{int}N_{_{Lo}}.
 \end{equation}
After an integration one can write the surface tension function $\sigma_{_{Lo}}:N_{_{Lo}} \to \R$ as  
\begin{align}   \label{gradsigma23b}
\sigma_{_{Lo}}(s,t) =-\frac{1}{\pi^2}(\Lq(\pi s)+\Lq(\pi t)+\Lq(\pi (1-s-t))),  
\end{align}
where 
\begin{align} \label{gradsigma45}
\Lq(\theta)=-\int_{0}^{\theta}\log \vert 2\sin x\vert dx
\end{align}
is the Lobachevsky function.

The surface tension $\sigma$ from \eqref{Pst2} has special features, which are the key to the detailed analysis in several themes studied in this work. We start by proving the first claim in 
Theorem \ref{sigma1}. 

\begin{Lem} \label{blowup}
  Let $\sigma:N\to \R$ be the unique bounded convex function that solves equation (\ref{Pst2}), and
suppose $J \subset \partial N$ is a closed interval not containing any of the points of $\mathscr P \cup \mathscr Q$. Then 
          \begin{equation}  \label{flat2}
    \    |\nabla \sigma(p)| \to \infty \qquad {as} \quad p \to J, \;  p \in N^\circ.     
             \end{equation}
\end{Lem}

\begin{proof}
By rotation and translation we may assume that $0 \in J \subset \R$ 
and that $N^\circ$ lies in the upper half plane. 
For $0 < x_1 \in J$,  let $T$ be the isosceles triangle with corners $\{0, x_1, ix_1\}$, and assume  $x_1$ is so small that $T  \subset N^\circ \setminus  \Gg$. Finally,  adding a linear map and making an affine change of coordinates keeps the form of the Monge-Amp\`ere equation (but may change affinely its data in \eqref{Pst2}). With such a transform we can  assume that $x_1 = 1$.

If $L$ is a linear map such that  $L = \sigma$ on $[0,1]$
 and $L(i) = \sigma(i)$, 
then by convexity $\sigma \leq L$ in all of the triangle $T$.
Furthermore, if $\sigma_0$ is the surface tension from \eqref{gradsigma23} - \eqref{gradsigma45}, then for $\sigma_L := \sigma_0 +L$
 the Hessian
$\det\bigl(D^2(\sigma_L)\bigr)=1$ in $T$. In particular,
by the comparison principle for Monge-Ampère equations, e.g. \cite[Lemma 2.7]{FigDePhi}, it follows that 
\begin{equation} \label{sigmacomp}
\sigma(p) \leq \sigma_L(p), \qquad p\in T.
\end{equation} 

On the other, from \eqref{gradsigma23} we see that $\partial_2 \sigma_L(p) \to -\infty$ 
as $p \to [\delta,1-\delta] \subset J$ and  $0 < \delta < 1$.   Since $\sigma$ and $\sigma_L$ have the same boundary value $L$ on $[0,1]$, with \eqref{sigmacomp} this forces  also $ |\nabla \sigma(p)|  \to \infty$ as $p \to [\delta,1-\delta] \subset J$. Finally, covering $J$ with such subintervals proves the claim.
\end{proof}

On the other hand, at the corners and gas points the surface tension $\sigma$ is more regular.  
 For this, recall  that 
the set of all subgradients of $\sigma$ at  $p_0$,
$$ \partial \sigma(p_0) := \{ \xi \in \mathbb{R}^2 :  \sigma(p)\ge \sigma(p_0)+\langle \xi, p-p_0\rangle, \; \forall \, p \in N \}, 
$$ 
is   called the {\it subdifferential} of $\sigma$ at $p_0$.
It follows  that 
$\partial \sigma(p_0)$ is a closed and convex set at any $p_0\in N$.   If $\sigma$ satisfies  \eqref{Pst2}, then for $p_0\in N^\circ \setminus \Gg$ we have $\partial \sigma(p_0) = \nabla \sigma(p_0)$, while by the previous Lemma, $\partial \sigma(p_0) =  \emptyset $ whenever $p_0 \in \partial N \setminus ( \Pp\cup \Qq)$. 

\smallskip

However, 
at corners, quasifrozen points and gas points, the subdifferential has a rich and smooth structure, c.f. also Figure \ref{fig.amoeba}.

 \begin{Thm}\label{prop:sigma:property3} Suppose $\sigma$ is as  \eqref{Pst2}. Then for gas points  $q_k \in \mathscr{G}$, the boundary  of $\partial \sigma(q_k)$ is an analytic, convex and bounded 
 Jordan curve $\gamma_{q_k}$.

 On the other hand, if $p_j \in \mathscr{P} \bigcup \mathscr{Q}$, the boundary of $\partial \sigma(p_j)$ is an unbounded Jordan arc $\gamma_{p_j}$, analytic except at $\infty$. 
  At endpoints at  $\infty$, the tangents of $\gamma_{p_j}$ are  parallel to $\, i(p_{j-1} - p_{j})$ and to $\, i(p_{j} - p_{j+1})$,  respectively, where $p_{j-1}, p_{j+1} \in \mathscr P\cup \mathscr Q$ are the points  
 neighbouring $p_j$, in the order induced by $\partial N$.
 
 Finally, the boundary of the complete amoeba defined by $\nabla\sigma$, 
 \begin{equation} \label{amoebabdry}
 \partial \left[  \nabla \sigma\left(  N^\circ\setminus {\mathscr G} \right)\right] \; = \; \Bigl( \bigcup_{p_j \in \mathscr P\cup \mathscr Q} \gamma_{p_j} \Bigr)\, \cup \, \bigcup_{q_k \in \mathscr G} \gamma_{q_k},
 \end{equation}
 is a disjoint union of the Jordans arcs $\gamma_{p_j}$ and the boundaries of the gas components $\gamma_{q_k}$.
  \end{Thm} 
\smallskip

Geometrically, the tangential directions $\, i(p_{j-1} - p_{j})$ and $\, i(p_{j} - p_{j+1})$ above are the outer normals  of\,  $N$  on $[p_{j-1}, p_{j}]$
and  $[p_{j}, p_{j+1}]$, respectively.

Proving  Theorem \ref{prop:sigma:property3}   
requires special properties of solutions to the Monge-Amp\`ere equation, developed in Subsection \ref{complexsec}, and hence the proof will be given there, see Remark \ref{proofThm22}. 

In fact, the interaction between $\sigma$ and the geometry of the boundary curves  in Theorem \ref{prop:sigma:property3}  can be further sharpened.
For this, first recall that as stated 
in Theorem \ref{sigma1}  and as will be  proven in  Section \ref{MAcomplex.structure}, for each such special point   $p_0 \in \mathscr P\cup \mathscr Q \cup \mathscr G$,
given any 
$\, p\in  N^\circ\setminus {\mathscr G}$ 
 the  following limit in the direction $p-p_0$,
\begin{align}
\label{GeneralizedGradient3a}
\widehat{\nabla} \sigma({p_0}, p-p_0) :=  \lim_{\tau\to 0^+} \nabla \sigma\big( p_0+\tau(p-p_0)\big),
\end{align}
exists, is finite and lies on the boundary of the subdifferential $\partial \sigma (p_0)$. In addition, every boundary point 
 of $\, \partial \sigma (p_0)$
arises this way. That gives these
subdifferentials also  explicit geometric 
properties as follows.
\smallskip

\begin{Prop}\label{prop:sigma:property4}
For any $p, q\in  N^\circ\setminus {\mathscr G}$ and any $p_0 \in \mathscr P\cup \mathscr Q \cup \mathscr G$ we have 
\begin{equation}\label{amoeba:convex}
\langle\widehat{\nabla}\sigma({p_0},p-p_0)-\widehat{\nabla}\sigma({p_0},q-p_0), p-p_0\rangle \geq 0.
\end{equation}
 In particular, $p - p_0$ is (an outer)
normal to $\partial \sigma(p_0)$ at the boundary point $\widehat{\nabla}\sigma({p_0},p-p_0)$.
\end{Prop}

\begin{proof}  
Once inequality 
\eqref{amoeba:convex} is established, Theorem \ref{prop:sigma:property3} with  convexity of $\partial  \sigma(p_0)$ 
and  smoothness of $\gamma_{p_0}$ implies that $p-p_0$ is normal to the boundary of $\partial \sigma (p_0)$.

For the inequality \eqref{amoeba:convex}, by definition of subdifferentiability we have, for all $p^\ast\in \dv \sigma(p_0)$ and all $p\in N^\circ\setminus {\mathscr G}$,
\begin{align*}
\sigma(p)\geq \sigma(p_0)+\langle p^\ast, p-p_0 \rangle.
\end{align*}
Hence, for any $t\in(0,1]$ and for every $p^\ast\in \dv \sigma(p_0)$, 
\begin{align} \label{gate22}
t \,  \langle p^\ast, p-p_0 \rangle \, \leq  \, \sigma\bigl(p_0+t(p-p_0)\bigr)-\sigma(p_0)  \, \leq  \, 
t\, \langle {\nabla}\sigma\bigl(p_0 + t(p-p_0)\bigr),p-p_0\rangle, 
\end{align}
 where the right hand inequality only uses the convexity of $\sigma$. Combining this with \eqref{GeneralizedGradient3a} and the choice  
 $\, p^\ast = \widehat{\nabla}\sigma(p_0;p-p_0)$ gives\begin{align} \label{gate2}
\lim_{t\to 0^+}\frac{\sigma(p_0+t(p-p_0))-\sigma(p_0)}{t}\le \langle \widehat{\nabla}\sigma(p_0;p-p_0),p-p_0\rangle,
\end{align}
as well as 
\begin{align} \label{defsubdif}
\langle \widehat{\nabla}\sigma(p_0;p-p_0),p-p_0\rangle\geq \langle p^\ast, p-p_0\rangle, \qquad \forall \; p^\ast\in \dv \sigma(p_0).
\end{align}
Since $\widehat{\nabla}\sigma(p_0;q-p_0)\in 
 \dv \sigma(p_0)$ whenever $q\in N^\circ\setminus {\mathscr G}$, one obtains
\begin{align*}
\langle \widehat{\nabla}\sigma(p_0;p-p_0)-\widehat{\nabla}\sigma(p_0;q-p_0),p-p_0\rangle\geq 0,
\end{align*}
which concludes the proof of \eqref{amoeba:convex}.
\end{proof}

With the above theorems, the subgradients $\widehat{\nabla}\sigma({p_0},p-p_0)$ parametrise the boundary arcs $\gamma_{p_0}$, for each $p_0 \in \mathscr P\cup \mathscr Q \cup \mathscr G$.

The study of the minimization  \eqref{ELbasic} will require also other special and quite delicate properties of the subdifferentials of $\sigma$, at the corners and quasifrozen points of $N$, such as the following result. In particular, it tells that  when $\sigma(p_0) = \sigma(p_1)$ at two neighbouring  points, the corresponding curves $\gamma_{p_0}$ and $\gamma_{p_1}$ have 
the same asymptotic line orthogonal to $[p_0, p_1]$. 

 For readers convenience we formulate the result already here, even if its proof can be provided only later, in Corollary \ref{ortoasymptotes}.
  \smallskip

\begin{Prop}\label{prop:sigma:property}
Let $\sigma$ be as  \eqref{Pst2}, and suppose  $p_0, p_1 \in \mathscr P\cup \mathscr Q$ are neighbouring points, in the order induced by $\partial N$. 

Then for any given  \; $\widehat p \in (p_0,p_1) \subset \partial N$,
there exists the limit
 \begin{equation} \label{amoeba:convex3}
\lim_{ N^\circ\setminus {\mathscr G}\, \ni \; p \to \, \widehat p \, }\langle\widehat{\nabla}\sigma({p_0},p-p_0), p -p_0\rangle =  \sigma(\, \widehat p \, ) - \sigma(p_0).
\end{equation}

\end{Prop}

\subsection{Variational problems}\label{subsect:variational}

Let $\Omega\subset \R^2$ be a bounded Lipschitz domain,
and let  $h_0:\partial \Omega\to \R$ be Lipschitz  continuous. 
We say that  $h_0$ is an {\it admissible boundary function} for $N$, if 
the class of the admissible functions
\begin{equation}\label{admissible_function}
\mathscr{A}_N(\Omega,h_0):=\{u\in C^{0,1}(\overline{\Omega}):
\nabla u(z)\in N \text{ for a.e. } z\in \Omega, u=h_0 \text{ on } \partial \Omega\}
\end{equation}
is not empty. For the dimer tiling problems, such a function is also known as a feasible boundary 
height function \cite{KeOk07}.

As shown in \cite{CKP} and \cite{KOS}, the limit height functions $h$ of discrete dimer models are solutions to  the following  variational problem: 
 minimize the functional 
\begin{align}
\label{SFunc}
I_\sigma[u]=\int_{\Omega}\sigma(\nabla u(z))dz
\end{align}
among all  $u\in \mathscr{A}_N(\Omega,h_0)$. From the convexity properties of $\sigma$ it is not hard to 
prove that  the minimization problem 
has a unique minimizer, 
denoted by $h$, 
see Proposition 4.5 of \cite{DeS10}.

Furthermore, for a $\sigma$ as singular as in \eqref{Pst2}, the only practical characterisation for a  minimizer of  \eqref{SFunc} 
is given in terms of the Gâteaux directional derivative. That is,  $h$ is a minimizer in  the class $\mathscr{A}_{N}(\Omega,h_0)$ if and only if 
\begin{equation}
\label{eq:DirectionalDerivative}
d I[h; u-h]=\int_\Omega d\sigma\big(\nabla h(z); \nabla u(z)-\nabla h(z)\big)\, dz\ge 0.
\end{equation} 
for all $u\in \mathscr{A}_N(\Omega,h_0)$. Here  the Gâteaux derivative of $\sigma$ at $p_0 \in N$ in the direction of $p-p_0$ is, by definition,
\begin{align*}
d\sigma(p_0;p-p_0):=\lim_{t\to 0^+}\frac{\sigma\big((1-t)p_0+tp\big)-\sigma(p_0)}{t}. 
\end{align*}
Since $\sigma$ is convex in $N$, we know that the Gâteaux derivative  $d\sigma$ is bounded from above,
\begin{equation}
\label{GatDerivative} 
d\sigma(p_0;p-p_0) \leq \sigma(p) - \sigma(p_0) \leq 2 \max_N |\sigma|, \qquad  \; \forall \; \; p_0, p \in N,
\end{equation} 
  but   in view of Lemma \ref{blowup} it is not  bounded from below. For lower bounds we have the 
inequalities
\begin{align} \label{defsubdif34}
d\sigma(p_0;p-p_0) \geq \langle p^\ast, p-p_0\rangle, \qquad \forall \; p^\ast\in \dv \sigma(p_0).
\end{align} 
 In addition, \eqref{gate2} with Lemma \ref{blowup}			
give the explicit description 
 \[ d\sigma(p_0;p-p_0)=\begin{cases}\langle \nabla \sigma(p_0),p-p_0\rangle \qquad &\text{if } p_0\in N^\circ\setminus \mathscr G;\\
\langle \widehat{\nabla} \sigma(p_0;p-p_0),p-p_0\rangle 
& \text{if } p_0\in \mathscr{P}\cup\mathscr{Q}\cup\mathscr{G};\\
-\infty &\text{if } p_0\in\dv N\backslash (\mathscr{P}\cup \Qq),
\end{cases}
\]
for $p \in N^\circ$.

 \subsubsection{Partial $C^1$-regularity of the minimizers} \label{sub:PRM}

 The basic regularity theory  for minimizers of \eqref{ELbasic}, for  a general $\sigma$ bounded and strictly convex in $N^\circ \setminus \mathscr G$, was developed by De Silva and Savin   \cite{DeS10}.
 We need the following results from their work as our starting point for our much more specific case of dimer models and surface tensions $\sigma$ from \eqref{Pst2}.

\begin{Thm}[De Silva - Savin \cite{DeS10}] 
\label{thm:DeS12} Suppose $h$ is the minimizer of \eqref{ELbasic}, with boundary values $h_0 \in Lip(\partial \Omega)$ on a bounded Lipschitz domain $\Omega$. Then $h\in C^1$ away from the obstacles, i.e. on the set $\{z \in \Omega: m(z) < h(z) < M(z) \}.$  
\end{Thm} Here for an admissible boundary
 function $h_0:\partial\Omega\to \R$, the 
 obstacles are simply 
 the minimum and, respectively, the maximum over all admissible functions in
$\mathscr{A}_N(\Omega,h_0)$,
\begin{equation}\label{def:Mm}
m(z) =\inf\{u(z):u\in \mathscr{A}_N(\Omega,h_0) \}, \quad 
M(z) =\sup\{u(z): u\in \mathscr{A}_N(\Omega,h_0)\}.
\end{equation}
 The function $m$ is called 
the \emph{lower obstacle}, and  $M$ the \emph{upper obstacle}.
Clearly, $m,M\in \mathscr{A}_N(\Omega,h_0)$ with $m(z)\le u(z) \le M(z)$ for all functions $u\in \mathscr{A}_N(\Omega,h_0)$. In particular,  this is true for the minimizer $h$ of the 
variational problem \eqref{SFunc}. 
\smallskip

A second result  from \cite{DeS10}, fundamental e.g. for the very notion  of a liquid domain, is their Theorem 4.1., which makes the convenient Definition \ref{LiqFrozen}  altogether  possible to work with. Here let $\Gamma:N\to \mathbb{S}^2$ be a continuous map such that $\Gamma(\partial N)=\{\xi\}\subset   \mathbb{S}^2$, while $\Gamma$ is a homeomorphism between $N^\circ$ and 
$\mathbb{S}^2 \setminus \{\xi\}$. Here $\mathbb{S}^2$ denotes the two dimensional unit sphere.
 \smallskip

\begin{Thm}[De Silva - Savin  \cite{DeS10}] \label{thm:DeS22} Suppose $h$ is the minimizer of \eqref{ELbasic} as in the previous theorem. Then  
$\Gamma\circ \nabla h: \Om\to  \mathbb{S}^2$ is continuous. 
\end{Thm}

 In \eqref{def:liquiddomain2} it was natural to define the 
liquid region $\LL$ to consist of those points where  the  minimizer $h$ of  \eqref{ELbasic} is  differentiable and the  gradient 
$\nabla h(z) \in N^\circ \setminus \Gg $. A priori such a set $\LL$ is only measurable, however, by Theorem \ref{thm:DeS22}  the set becomes open and the minimizer there is $C^1$-smooth. 
Thus via the  classical regularity theory we have $h \in C^\infty(\LL)$ as well as 
the Euler-Lagrange equation \eqref{eq:EL34} in $\LL$. As another aspect, Theorem \ref{thm:DeS22} gives us tools to control the frozen extensions of solutions to \eqref{eq:EL34},  see Theorem \ref{thm:h:proper3}. That, in turn, is one of the essential steps in proving the universality of the Lozenges  frozen boundaries, c.f. Theorem \ref{Second.thm}.

 It is an interesting open question  if similar regularity results as in Theorems \ref{thm:DeS12} or \ref{thm:DeS22}  hold  in higher dimensions.

   \subsubsection{  Boundary height functions and obstacles}\label{subsubsect:domains}

 From  another point of view, 
  the obstacles $M$ and $m$  given in  \eqref{def:Mm} are  the \emph{McShane} extensions  \cite{Rock}  of the boundary value $h_0$. That is,
  they
coincide on $\Om$ with the following functions which are defined in the whole space $\R^2$,
\begin{equation}\label{def2:Mm}
m(z)=\max_{w\in \partial \Omega}\big(-h_N(w-z)+h_0(w)\big), \quad
M(z)=\min_{w\in\partial\Omega}\big(h_N(z-w)+h_0(w)\big).
\end{equation}
Here $h_N:\R^2\to \R$ is the {\it support function}  \cite{Sch} of the convex polygon $N$ with corners $\{ p_j \}$, 
\begin{equation}\label{supportf} h_N(z)=\sup_{p\in N}\, \langle p, z\rangle = \sup_{1 \leq j \leq m} \, \langle p_j, z\rangle.
\end{equation}

Via a translation in the data we can (and will) often assume that $0 \in N^\circ$. Then from \eqref{supportf} we have the simple, but often useful, observation
\begin{equation}\label{struct.supportfcn} 
h_N(z) = \langle p_n, z \rangle, \qquad {\rm whenever} \; \; \; \langle  p_{n+1} - p_{n}, z \rangle = 0  \quad{\rm and}  \quad  \langle p_n, z \rangle > 0.
\end{equation}

The {\it trivial set}, for given $\Omega$ and boundary function $h_0$,  is the set of points  where the lower and upper obstacles coincide,
\begin{equation}\label{def:trivialset} 
\Omega_T=\{ z\in \overline{\Omega}: m(z)=M(z)\}.
\end{equation}
If  $\Omega = \Omega_T$
there is only one admissible function in  $\mathscr{A}_N(\Omega,h_0)$, 
which is trivially the minimizer of the variational problem \eqref{ELbasic}. 
 This is not a typical case, but can happen even for natural non-trivial boundary values, see Example \ref{noliquid}. Otherwise, the variational problem
is reduced to  studying it in  each connected component of the open set $\Omega\setminus \Omega_T$. We denote by 
\begin{equation}\label{coincide} 
\Lambda_m=\{x\in \overline{\Om}: h(x)=m(x)\} \quad {\rm and } \quad \Lambda_M=\{x\in \overline{\Om}: h(x)=M(x)\}
\end{equation}
the {\it coincidence sets} of $h$ and by $\Lambda := \Lambda_m \cup  \Lambda_M$ the total coincidence set of $h$.

\subsubsection{Free boundary: frozen, liquid and gas domains}\label{subsubsect:frozen}

 De Silva's and Savin's Theorems \ref{thm:DeS12} or \ref{thm:DeS22} 
allow  natural definitions for the liquid and frozen regions of the minimizer $h$ in $\Om$. In view of these theorems, it would have been equivalent to 
require $\LL$ to consist of those points of $\Om$ where $h$ is $C^1$  in a neighbourhood of $z$ and we have $\nabla h (z) \in N^\circ \setminus \Gg$.  
That is actually the approach in \cite{KeOk07}.

Similarly, Definition  \eqref{frozen.def}  of  the frozen boundary allows the minimizer to have three different boundary modes 
(or as you may consider, three different phases),   determined by the corresponding facets in $\Omega \setminus \LL$, c.f. Theorem \ref{thm:main2} and \eqref{limitp} in particular. 
 In analogy of above, we ask these sets of phases to be open,  and call them the \emph{frozen, quasi-frozen} and \emph{gas regions}, respectively. 
\begin{Def}\label{def:domains}
The \emph{frozen region} of $h$ is defined as
\begin{align}\label{def:frozendomain}
\mathcal{F}:= {\rm int} \{z\in \Omega: \text{$h$ is $C^1$ in a neighbourhood of $z$, }\, \nabla h(z)\in \Pp \}. 
\end{align}
The \emph{gas region} of $h$ is defined as
\begin{align}\label{def:gasdomain}
\mathcal{G}= {\rm int}\{z\in \Omega: \text{$h$ is $C^1$ in a neighbourhood of $z$, }\,\nabla h(z)\in \mathscr{G}\}.
\end{align} 
The \emph{quasi-frozen region} of $h$ is defined as
\begin{align}\label{def:qfrozendomain}
\mathcal{Q}= {\rm int}\{ z\in \Om:\text{$h$ is $C^1$ in a neighbourhood of $z$, }\,\nabla h(z)\in \mathscr{Q}\}.
\end{align}
\end{Def}

Clearly, any of these sets can be empty. With the following two simulation examples we illustrate the above concepts.

\begin{ex}[Simulation Lozenge Tiling with Quasi-frozen domain]

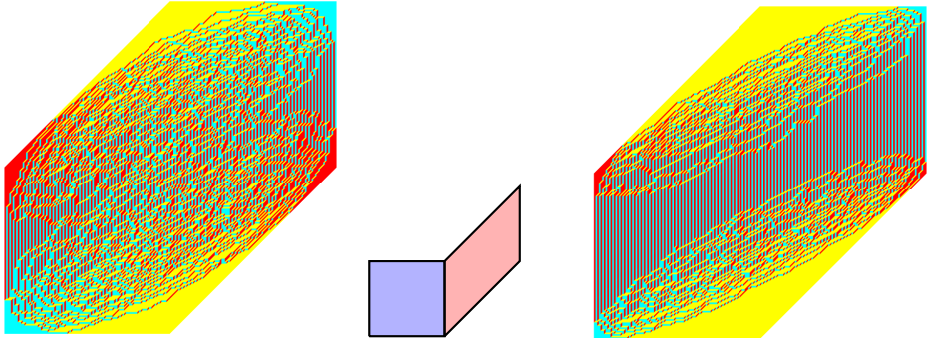
\begin{figure}[H]
\centering{}
\includegraphics[scale=0.4]{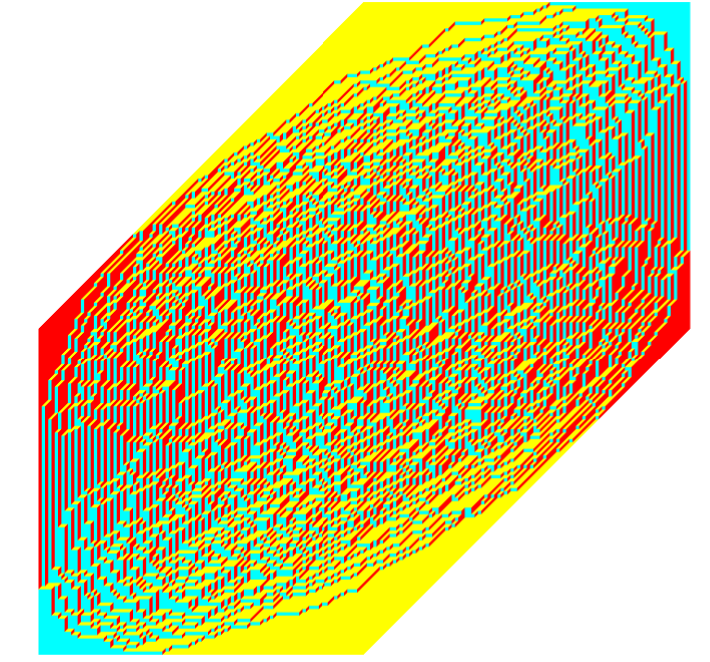}
\hspace{.1em}
\begin{tikzpicture}[xscale=1,yscale=1]
\definecolor{Mycolor}{RGB}{44,221,218}
\filldraw[fill=Mycolor] (1,1)--(2,1)--(2,2)--(1,2)--(1,1);
\filldraw[fill=red](0.3,0.3)--(1,1)--(1,2)--(0.3,1.3)--(0.3,0.3);
\draw[thick] (1,1)--(2,1)--(2,2)--(1,2)--(1,1);
\draw[thick] (0.3,0.3)--(1,1)--(1,2)--(0.3,1.3)--(0.3,0.3);
\end{tikzpicture}
\hspace{1.5em}
\includegraphics[scale=0.4]{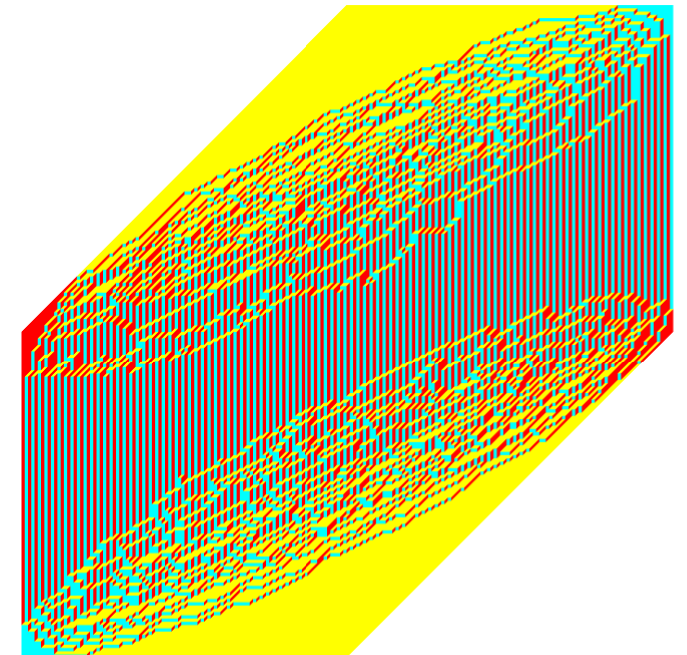}
\caption{Quasifrozen domains with different surface tensions but with the same boundary height function $h_0$. In the middle a 'quasi-particle'.  Image courtesy of M. Duits.}
\label{figQuasifrozen1}
\end{figure}

In Figure \ref{figQuasifrozen1} we have two simulations 
of  lozenge tilings in a hexagonal domain. Both simulations have a liquid domain as well as non-empty frozen and quasifrozen regions. The piecewise affine boundary height function $h_0$ is the same for both,  but 
 the surface tensions have different  boundary values $L$ on $\partial N$. 
The change in $L$ is induced in the limit from the change of periodic weights in the probability measure for the microscopic perfect matching model. 

In  Figure \ref{figQuasifrozen1} on the left, the liquid domain is simply connected with two cusps in the limit. Moreover, if in the simulation we zoom in on the quasi-frozen phase, we see that it consists entirely of the composite tile presented in the middle of the figure, namely the quasifrozen phases are vertical stripes regions which are formed by putting the quasi-particle tiles indicated in the figure  on top of each other.

This explains the name  of the quasi-frozen domain, where a composite tile can be thought of as a quasi-particle. 

At a critical choice of weights, corresponding to a special choice of boundary value for the surface tension $\sigma$, the two cusps of the liquid domain will merge to form a tacnode.  Continuing to change the weights will then separate the two components of the liquid domain, and we obtain the configuration on the right of Figure \ref{figQuasifrozen1}. 
Moreover, the minimizer $h$ will not coincide with either of the obstacles at the quasi-frozen domain, since $\nabla h\in \Qq\subset \dv N\setminus\mathscr{P}$ in the  quasi-frozen domain and the gradient of the obstacles belongs to set $\mathscr{P}$ almost everywhere.
\end{ex}

\vspace{-.1cm}

\begin{figure}[H]
\centering{}
\includegraphics[scale=0.34]{Picture10.png}
\hspace{3em}
\includegraphics[scale=0.075]{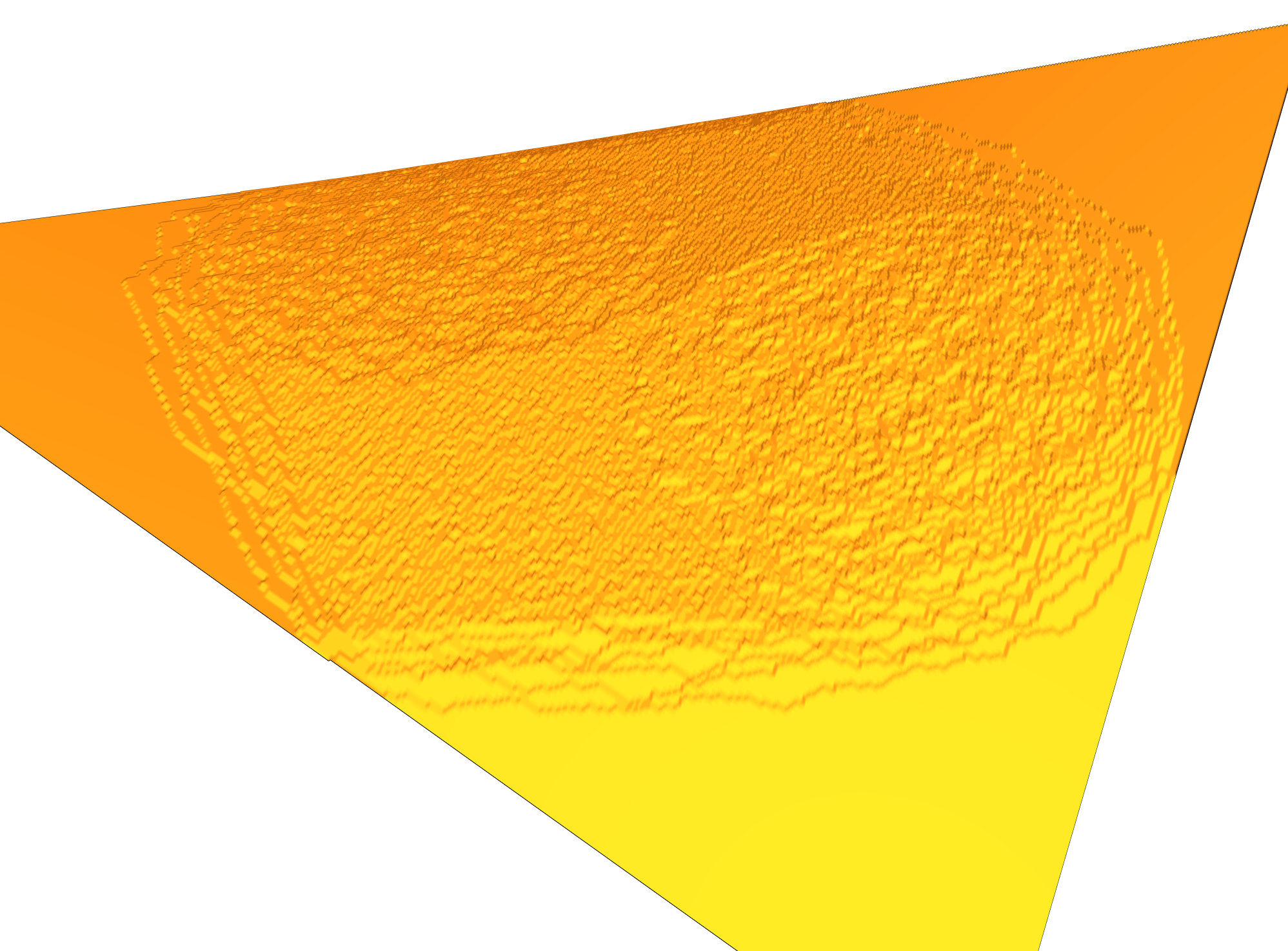}
\caption{Simulation of  random tiling of the Aztec diamond with two gas domains. Graph of the height function pictured from above (left) and from  side (right). Courtesy of Tomas Berggren.}
\label{figGasdomains1}
\end{figure}

\vspace{-.15cm}

\begin{ex}[Simulation of Domino Tiling with Gas Domains] In Figure \ref{figGasdomains1} above we have a simulation of an Aztec diamond, with a weighting that gives rise to a gas region with two different components. In both components the gradient of the height function is constant, but different in the different components. 
\end{ex}

\begin{figure}[H] 
\centering{}
\includegraphics[scale=0.32]{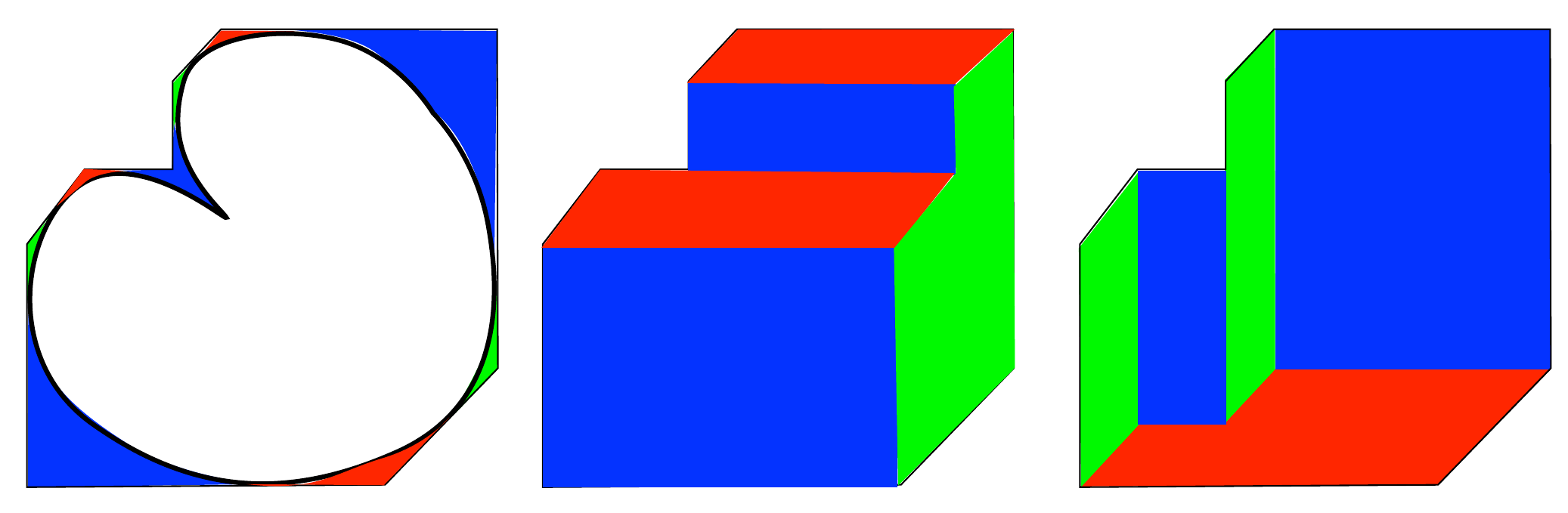} 
\caption{On  left: Cusp, where frozen facet $ \mathcal{F} \not\subset \Lambda$. In middle: $\Lambda_m$. On right: $\Lambda_M$.}
\label{facet}
\end{figure}

\enlargethispage{3\baselineskip}
In some special cases, the frozen set $\mathcal{F}$ is a subset of the coincidence set $\Lambda$ defined as in  \eqref{coincide}. In general, this is not true even in the setting of the lozenges model, 
 as the following example in Figure \ref{facet} illustrates. Here the frozen boundary is a cardioid and 
 Theorem  \ref{thm:hstar:minimizer}  verifies that the minimizer is indeed as is depicted in the left figure below.

One typically thinks of the frozen boundary as a \emph{phase boundary}, i.e. as an interface between the liquid domain and the frozen, quasi-frozen or gas region. However, if the liquid domain has only partially frozen boundary, then many other different type of phenomena can occur, see
e.g. \cite{DeS10} or \cite{Duse15a}. On the issue of partially frozen boundaries  we will in this work only  discuss their regularity, see e.g. Theorem \ref{surprise2}.

\subsubsection{A distinguished class of domains and boundary values}\label{subsubsect:distinguished}

 For any given  dimer model, we next look  for natural candidates of polygonal domains $\Om$ and piecewise affine boundary values on $\partial \Om$ to give rise to frozen phenomena.

In the first basic simulations,  see e.g.  Figure \ref{fig.first}, in  frozen facets the limit height function has gradient $\nabla h = p_j$ lying  in  a corner of $N$. 
Typically two  facets, say, with $\nabla h = p_j$ and $\nabla h = p_{j+1}$,  meet along a side of the polygon $\Omega$  at a point tangential to the liquid domain. 
In addition, 
 the boundary value $h_0$ is affine (thus $C^1$) along each side of $\Omega$. Hence if $S_0 \subset \partial \Omega$ is the side where the above facets meet, we must have  
$S_0 \perp (p_{j+1}-p_j)$. 
 This simple observation leads to the following definition 

  \begin{Def}\label{Def:naturaldomain}
Let $N$ be a closed convex polygon in $\R^2$,  with 
vertices $\mathscr{P}=\{ p_1,...,p_k\}$. We say  that a (possibly non-convex) polygon $\Omega\subset \R^2$,  with vertices 
$\{ z_1,...,z_d\}$,  is a \emph{natural domain} for $N$,
if we can 
 associate to each vertex $z_j$ a vertex $p_n\in \Pp$, such that 
 either
\begin{align}
\label{ND1}
\langle z_{j}-z_{j-1}, p_{n}-p_{n-1}\rangle&=0 \quad {\rm and} \quad
\langle z_{j+1}-z_j, p_{n+1}-p_n\rangle=0,
\end{align}
or 
\begin{align}
\label{ND3}
\langle z_{j}-z_{j-1}, p_{n+1}-p_{n}\rangle&=0  \quad {\rm and} \quad
\langle z_{j+1}-z_j, p_{n}-p_{n-1}\rangle=0.
\end{align}
\end{Def}

We assume above that both sets of vertices $\{ p_1,...,p_k\}$ and $\{ z_1,...,z_d\}$ are given in counterclockwise order induced by the boundaries $\partial N$ and, respectively, $\partial \Omega$. We also set above $p_{k+1}=p_1$ with $p_0=p_{k}$, and similarly,  $z_{d+1}=z_1$ and $z_0=z_d$.  Note that we do not ask for a relation between the corner $p_n$ attached to $z_j$ and those attached to $z_{j\pm 1}$.

The attached corners $p_n$ are needed to construct appropriate  boundary values on $\partial \Om$, see Definition \ref{def:extremal} below.  When  $N$ is a triangle, like for example in the lozenges model, it suffices to require that each side  $[z_j, z_{j+1}] \subset \partial \Om$ is orthogonal to some side of $N$. This already determines unique  corners $p_n$ that satisfy  \eqref{ND1} or \eqref{ND3}. 
Similar property actually holds for a generic polygon $N$, with  geometry already determining the allowed boundary values.
 When  $N$ is a rectangle like for the domino tilings, see the discussion in Example \ref{microdominoes}, again it is equivalent just to ask that each side of $\Omega$ is orthogonal to a side of $N$. This time, however, the choices of $p_n$ are not unique, rather any choice at each corner $z_j$ will do.

When we want to emphasise that a domain is a natural domain for $N$, we write $\Om_N$. 
For the natural domains, the above discussion asks them to be associated with a special class of boundary values,  which we call {\it natural boundary values}.
\begin{Def}\label{def:extremal}
Let $\Omega = \Om_N$ be a natural domain with
$d$ vertices $\{ z_1,...,z_d\}$. We say that an admissible boundary function $h_0:\partial \Omega\to \R$ 
is  \emph{natural}  if for each $j=1,...,d$ we have
\begin{equation} \label{zjpn}
 h_0(z)=\langle p_n, z-z_j\rangle+h_0(z_j)
 \end{equation}
for all $z$ contained on the segments $[z_j,z_{j+1}]$ and $[z_{j-1},z_j]$, where 
$p_n, 1\le n\le k$, is the corner attached to $z_j$,   such that (\ref{ND1})  or (\ref{ND3})   holds.
\end{Def}

Note that an arbitrary natural domain need not admit a natural boundary value function; we leave it to the reader to construct such examples.

One could also consider multiply connected natural domains $\Omega = \Omega_N$ and natural boundary values  on their boundaries,
by removing from the interior of a simply connected natural domain (a finite number of) other  natural domains (associated with the same polygon $N$). The natural boundary 
function is then required to satisfy Definition \ref{def:extremal}  on each boundary component of $\Omega_N$. 

For   lozenge models and triangle gradient constraints $N$,  in a simply connected natural domain $\Om$ the natural
boundary value $h_0$ is uniquely determined by $\dv \Om$, up to an additive constant. However,  if $\Om$ is multiply connected, this additive constant can be different for different components of $\partial \Om$. Thus finding the actually minimal configuration  leads to an additional  problem of minimizing over  the relative height differences between the different boundary components. On the other hand, which ever of these gives the minimal energy, it will have admissible boundary values  and thus our description of its properties  and the corresponding liquid domains still apply.
For simplicity, in the sequel we will consider only simply connected natural domains, as described in Definition \ref{Def:naturaldomain}. 

\begin{ex}[Natural Domains and Boundary Values for Domino Tilings] \label{microdominoes}
It is interesting to observe  how the natural boundary values  arise from the microscopic pictures, 
 for instance in the case of domino tilings. 
 
 There are, of course, many different boundary values  on the microscopic level, allowing small variations for them, that lead 
 in the scaling limit to given affine macroscopic boundary values. However, here we wish just to indicate by an example how the very simplest natural discrete boundary values  lead to the notion described in Definition \ref{def:extremal}. 
 
 We have now four  (orientation taken into account) different tiles, 
 say, yellow, red, green and blue, as in the top part of  Figure \ref{figGradDomino} below. 
  The discrete height function, as defined by Thurston \cite{Th90}, lives on the dual lattice, but an equivalent description \cite{CEP}, \cite{CKP}
  is obtained by assigning the gradient $(0,1)$ to yellow dominos, $(-1,0)$ to red dominos, $(0,-1)$ to green dominos and $(0,1)$ to blue dominos. In the limit when the size of the dominoes goes to  $0$, the gradient of the asymptotic height function $h$ 
 will take values in the convex hull $N$ of the points $\{(0,1),(-1,0),(0,-1),(1,0)\}$. See Figure \ref{figGradDomino}.

\begin{figure}[H]
\centering

\begin{tikzpicture}[xscale=0.4,yscale=0.4]

\filldraw[scale=1.5,xshift=0.7cm,yshift=0.7cm,fill=gray!30!white] (-1,0)--(0,-1)--(1,0)--(0,1)--(-1,0);
\draw[scale=1.5,xshift=0.7cm,yshift=0.7cm,ultra thick] (-1,0)--(0,-1)--(1,0)--(0,1)--(-1,0);
\draw[scale=1.5,xshift=0.7cm,yshift=0.7cm,thick,->,>=stealth]  (0,-0.3)  arc[radius = 3mm, start angle= -90, end angle= 90];
\draw(2.3,2.3) node {$N$};

\filldraw[xshift=-3cm,yshift=0cm,fill=red!80!white] (0,0)--(1,0)--(1,2)--(0,2)--(0,0); 
\draw[xshift=-3cm,yshift=0cm,thick] (0,0)--(1,0)--(1,1)--(0,1)--(0,0)--(1,0)--(1,1);
\filldraw[xshift=-3cm,yshift=0cm,fill=black] (0.5,0.5) circle [radius=0.1cm];
\filldraw[xshift=-3cm,yshift=0cm,fill=white] (0.5,1.5) circle [radius=0.1cm];
\draw[xshift=-3cm,yshift=0cm,->,>=stealth, ultra thick] (1,1)--(-1,1);
\draw[xshift=-3cm,yshift=0cm] (-1,2) node {$p_3$};

\filldraw[xshift=4cm,yshift=0cm,fill=blue!80!white] (0,0)--(1,0)--(1,2)--(0,2)--(0,0); 
\draw[xshift=4cm,yshift=0cm,thick] (0,0)--(1,0)--(1,1)--(0,1)--(0,0)--(1,0)--(1,1);
\filldraw[xshift=4cm,yshift=0cm,fill=white] (0.5,0.5) circle [radius=0.1cm];
\filldraw[xshift=4cm,yshift=0cm,fill=black] (0.5,1.5) circle [radius=0.1cm];
\draw[xshift=4cm,yshift=0cm,->,>=stealth,ultra thick] (0,1)--(2,1);
\draw[xshift=4cm,yshift=0cm] (2,2) node {$p_1$};
\filldraw[xshift=0cm,yshift=-3cm,fill=green] (0,0)--(2,0)--(2,1)--(0,1)--(0,0); 
\draw[xshift=0cm,yshift=-3cm,thick] (0,0)--(2,0)--(2,1)--(0,1)--(0,0)--(1,0)--(1,1);
\filldraw[xshift=0cm,yshift=-3cm,fill=white] (0.5,0.5) circle [radius=0.1cm];
\filldraw[xshift=0cm,yshift=-3cm,fill=black] (1.5,0.5) circle [radius=0.1cm];
\draw[xshift=0cm,yshift=-3cm,->,>=stealth,ultra thick] (1,1)--(1,-1);
\draw[xshift=0cm,yshift=-3cm] (1.5,-1.5) node {$p_3$};


\filldraw[xshift=0cm,yshift=4cm,fill=yellow] (0,0)--(2,0)--(2,1)--(0,1)--(0,0); 
\draw[xshift=0cm,yshift=4cm,thick] (0,0)--(2,0)--(2,1)--(0,1)--(0,0)--(1,0)--(1,1);
\filldraw[xshift=0cm,yshift=4cm,fill=black] (0.5,0.5) circle [radius=0.1cm];
\filldraw[xshift=0cm,yshift=4cm,fill=white] (1.5,0.5) circle [radius=0.1cm];
\draw[xshift=0cm,yshift=4cm,->,>=stealth,ultra thick] (1,0)--(1,2);
\draw[xshift=0cm,yshift=4cm] (2,2) node {$p_2$};

\end{tikzpicture}

\vspace{.4cm}

\begin{tikzpicture}[xscale=0.5,yscale=0.5]

\draw (-3.8,0.8) node {$\circ$};
\draw (-3.6,0.6) node {$\circ$};
\draw (-3.4,0.4) node {$\circ$};
\draw[xshift=6.5cm,yshift=-6.5cm] (-3.8,0.8) node {$\circ$};
\draw[xshift=6.5cm,yshift=-6.5cm] (-3.6,0.6) node {$\circ$};
\draw[xshift=6.5cm,yshift=-6.5cm] (-3.4,0.4) node {$\circ$};

\filldraw[xshift=-1cm,yshift=-1cm,fill=yellow] (-3,0)--(-1,0)--(-1,1)--(-3,1)--(-3,0);
\draw[xshift=-1cm,yshift=-1cm,thick] (-3,0)--(-1,0)--(-1,1)--(-3,1)--(-3,0)--(-2,0)--(-2,1);
\filldraw[xshift=-1cm,yshift=-1cm,fill=black] (-2.5,0.5) circle [radius=0.1cm];
\filldraw[xshift=-1cm,yshift=-1cm,fill=white] (-1.5,0.5) circle [radius=0.1cm];

\filldraw[xshift=0cm,yshift=-2cm,fill=yellow] (-3,0)--(-1,0)--(-1,1)--(-3,1)--(-3,0);
\draw[xshift=0cm,yshift=-2cm,thick] (-3,0)--(-1,0)--(-1,1)--(-3,1)--(-3,0)--(-2,0)--(-2,1);
\filldraw[xshift=0cm,yshift=-2cm,fill=black] (-2.5,0.5) circle [radius=0.1cm];
\filldraw[xshift=0cm,yshift=-2cm,fill=white] (-1.5,0.5) circle [radius=0.1cm];


\filldraw[xshift=2cm,yshift=-4cm,fill=red!80!white] (-3,0)--(-2,0)--(-2,2)--(-3,2)--(-3,0);
\draw[xshift=2cm,yshift=-4cm,thick] (-3,0)--(-2,0)--(-2,1)--(-3,1)--(-3,0)--(-2,0)--(-2,1);
\filldraw[xshift=2cm,yshift=-4cm,fill=black] (-2.5,0.5) circle [radius=0.1cm];
\filldraw[xshift=2cm,yshift=-4cm,fill=white] (-2.5,1.5) circle [radius=0.1cm];

\filldraw[xshift=3cm,yshift=-5cm,fill=red!80!white] (-3,0)--(-2,0)--(-2,2)--(-3,2)--(-3,0);
\draw[xshift=3cm,yshift=-5cm,thick] (-3,0)--(-2,0)--(-2,1)--(-3,1)--(-3,0)--(-2,0)--(-2,1);
\filldraw[xshift=3cm,yshift=-5cm,fill=black] (-2.5,0.5) circle [radius=0.1cm];
\filldraw[xshift=3cm,yshift=-5cm,fill=white] (-2.5,1.5) circle [radius=0.1cm];

\filldraw[xshift=4cm,yshift=-6cm,fill=red!80!white] (-3,0)--(-2,0)--(-2,2)--(-3,2)--(-3,0);
\draw[xshift=4cm,yshift=-6cm,thick] (-3,0)--(-2,0)--(-2,1)--(-3,1)--(-3,0)--(-2,0)--(-2,1);
\filldraw[xshift=4cm,yshift=-6cm,fill=black] (-2.5,0.5) circle [radius=0.1cm];
\filldraw[xshift=4cm,yshift=-6cm,fill=white] (-2.5,1.5) circle [radius=0.1cm];

\draw (4.2,0.8) node {$\bullet$};
\draw (4.4,0.6) node {$\bullet$};
\draw (4.6,0.4) node {$\bullet$};
\draw[xshift=6.5cm,yshift=-6.5cm] (4.2,0.8) node {$\bullet$};
\draw[xshift=6.5cm,yshift=-6.5cm] (4.4,0.6) node {$\bullet$};
\draw[xshift=6.5cm,yshift=-6.5cm] (4.6,0.4) node {$\bullet$};

\filldraw[xshift=-1cm,yshift=-1cm,fill=green] (5,0)--(7,0)--(7,1)--(5,1)--(5,0);
\draw[xshift=-1cm,yshift=-1cm,thick] (5,0)--(7,0)--(7,1)--(5,1)--(5,0)--(6,0)--(6,1);
\filldraw[xshift=-1cm,yshift=-1cm,fill=white] (5.5,0.5) circle [radius=0.1cm];
\filldraw[xshift=-1cm,yshift=-1cm,fill=black] (6.5,0.5) circle [radius=0.1cm];

\filldraw[xshift=0cm,yshift=-2cm,fill=green] (5,0)--(7,0)--(7,1)--(5,1)--(5,0);
\draw[xshift=0cm,yshift=-2cm,thick] (5,0)--(7,0)--(7,1)--(5,1)--(5,0)--(6,0)--(6,1);
\filldraw[xshift=0cm,yshift=-2cm,fill=white] (5.5,0.5) circle [radius=0.1cm];
\filldraw[xshift=0cm,yshift=-2cm,fill=black] (6.5,0.5) circle [radius=0.1cm];


\filldraw[xshift=2cm,yshift=-4cm,fill=blue!80!white] (5,0)--(6,0)--(6,2)--(5,2)--(5,0);
\draw[xshift=2cm,yshift=-4cm,thick] (5,0)--(6,0)--(6,1)--(5,1)--(5,0)--(6,0)--(6,1);
\filldraw[xshift=2cm,yshift=-4cm,fill=white] (5.5,0.5) circle [radius=0.1cm];
\filldraw[xshift=2cm,yshift=-4cm,fill=black] (5.5,1.5) circle [radius=0.1cm];

\filldraw[xshift=3cm,yshift=-5cm,fill=blue!80!white] (5,0)--(6,0)--(6,2)--(5,2)--(5,0);
\draw[xshift=3cm,yshift=-5cm,thick] (5,0)--(6,0)--(6,1)--(5,1)--(5,0)--(6,0)--(6,1);
\filldraw[xshift=3cm,yshift=-5cm,fill=white] (5.5,0.5) circle [radius=0.1cm];
\filldraw[xshift=3cm,yshift=-5cm,fill=black] (5.5,1.5) circle [radius=0.1cm];

\filldraw[xshift=4cm,yshift=-6cm,fill=blue!80!white] (5,0)--(6,0)--(6,2)--(5,2)--(5,0);
\draw[xshift=4cm,yshift=-6cm,thick] (5,0)--(6,0)--(6,1)--(5,1)--(5,0)--(6,0)--(6,1);
\filldraw[xshift=4cm,yshift=-6cm,fill=white] (5.5,0.5) circle [radius=0.1cm];
\filldraw[xshift=4cm,yshift=-6cm,fill=black] (5.5,1.5) circle [radius=0.1cm];

\end{tikzpicture}
\caption{}
\label{figGradDomino}
\end{figure}

For the underlying  bipartite graph (in this case $\Z^2$) all nodes are coloured either white or black.
Furthermore,  in the most common simulations of domino tilings giving rise to frozen phenomena, such as e.g. 
the Aztec diamond   or the simulations in Figure  \ref{fig.first}, one considers domains covered by domino tiles, 
with boundary consisting of stepped sides of the form depicted in the lower part of Figure \ref{figGradDomino}.

In the first case, when the corners of the stepped side has only have white vertices,  the boundary tiles are either yellow or red dominos, see left picture in Figure \ref{figGradDomino}. Moreover, there is a switch from yellow to red dominos only once along such stepped side. In the same way, if the stepped side consists of black vertices, one tiles with green and blue dominos and similarly can  switch between the dominos only once, see Figure \ref{figGradDomino}.

Finally, in the lower Figure \ref{figGradDomino} on left  the affine height functions have  gradients $(0,1)$ and $(-1,0)$ (corresponding to yellow and red dominos, respectively). Note  that the components of these gradients in the ``tangential'' direction $(-1,1)$  are the same. Similar holds between the affine functions 
corresponding to green and blue dominoes, with gradient $(0,-1)$ and $(1,0)$, respectively.
In the scaling limit, the above settings lead  precisely to the natural boundary values of  Definition \ref{def:extremal},
 along each edge.
  Of course, for the entire domain to be natural, one also needs a condition allowing stepped sides to be glued together at corners. 
\end{ex}
  \smallskip
 
 As a last aspect on the natural boundary values, understanding the boundary behaviour of the gradients $\nabla h$ of the minimizers of \eqref{ELbasic}  is of course fundamental for the goals of this paper.  However, the problem here is that due to the singular behaviour of the surface tension $\sigma$, as indicated e.g. in Theorem \ref{sigma1} and studied more  thoroughly in the next section, there appears no general or appropriate  methods that guarantee the boundary continuity of $\nabla h$.
 The only point one can make use of is Theorem \ref{thm:DeS22} due to De Silva and Savin. 
 
 In this situation, even to enable or start the analysis of the minimization problem \eqref{ELbasic},   we introduce the simple but very useful concept of  {\it frozen extensions}. The idea here is that given $z_0 \in \partial \Om$ which, say, is not a corner of $\Omega$, then  the above natural boundary values turn out to  admit an extension to the boundary of a larger domain $\widehat \Omega \supset \Omega$ with $z_0 \in \widehat \Omega$, 
 such that the obstacles  \eqref{def2:Mm} agree on $\widehat \Omega \setminus \Omega$.
 In particular, this means that minimizer for the new boundary value on $\partial \widehat \Omega$ is forced to agree  in $\Omega$ with the old minimizer, 
 while at the same time $z_0$ becomes an interior point for  $\widehat \Omega$. That this is at all possible will be shown in Section \ref{sect:properness}, where we also discuss which corners admit a frozen extension. Here we only present the following definition. 
 
\begin{Def}\label{oldold} Let $\Omega$ be a bounded Lipschitz domain and $h_0$ an admissible boundary value. We say that $\Omega$ admits a {\rm frozen extension} at a point $z_0 \in \partial \Omega$, if there exists a domain $\widehat \Omega \supset \Omega$
with $z_0 \in \widehat \Omega$, and  a boundary value $\widehat h_0$ on $\partial \widehat \Omega$ for which
the upper and lower obstacles $\widehat M (z)$ and $\widehat m(z)$ satisfy
$$ \widehat M (z) = \widehat m(z)\quad {\rm for \; all} \; z \in {\overline {\widehat \Omega}} \setminus \Omega, \quad {\rm with} \quad 
\widehat M (z) = \widehat m(z) = h_0(z) \; {\rm on } \; \partial \Omega.
$$
\end{Def}

\section{Monge-Amp\`ere Equation and the Complex Structure}  \label{MAcomplex.structure}

 In their work  \cite{KeOk07} Kenyon and Okounkov pointed out that the  complex Burgers equation induces a complex structure on  liquid domains of any dimer model. For further details see also \cite{KeHoney}. 
In a sense,  for our study the 
natural complex structure is a starting point as well as a unifying theme. We will develop the many different aspects of the  complex structures  to reach a detailed  understanding 
  of  the geometry of the limits of random surfaces within the general dimer models.  We will see many variants of this point of view in the surface tension, in the minimizers of the variational problem  \eqref{ELbasic}  as well as in the geometry of liquid domains  and the Beltrami equations the domains support.
  
  In addition, what makes the Beltrami equations particularly useful for the study of the frozen boundaries of liquid domains  is that, as we shall see, the solutions to these equations admit continuous extensions to the boundary, a feature not apparent in other methods describing the liquid domains or frozen boundaries.
 
  The purpose of the present Section is, starting from the Euler-Lagrange equation \eqref{ELsigma},  to build up the machinery necessary to describe and understand these questions, including  the geometry of frozen boundaries in general dimer models.
  First, we  develop in Subsections \ref{ELBeltrami} - \ref{universal} the intimate connection between the minimizers of  \eqref{ELbasic} and the appropriate Beltrami equations. Then in Subsection \ref{surface}, we study  the boundary behaviour of the surface tension $\sigma$ and its gradient. This is one of the keys in controlling the Beltrami equation, see Sections \ref{B+hodo} 
   and \ref{propersec}.
    Understanding  the delicate boundary behaviour of $\nabla \sigma$ is essential  also for  the methods  developed in Section  \ref{section:minimality}, for showing that functions $h$ in the admissible class $\mathscr{A}_{N}(\Omega,h_0)$ we construct are actually  minimizers of the variational problem (\ref{ELbasic}).

\subsection{Euler-Lagrange equations in the Beltrami picture}  \label{ELBeltrami}

 Given a surface tension $\sigma$ as in \eqref{Pst2} and a minimizer $h \in  \mathscr{A}_N(\Omega,h_0)$ of the variational integral \eqref{ELbasic}, assume $h$ has a non-empty liquid domain $\LL \subset \Omega$, c.f.  \eqref{def:liquiddomain2}.
 As discussed in Subsection (\ref{sub:PRM}), $\mathcal{L}$ is open. 
 The minimizer $h$ is smooth and satisfies in the liquid domain $\LL$  the Euler-Lagrange equation
   \begin{equation} \label{ELsigma}
  \text{div}\, \big(\nabla \sigma(\nabla h)\big)= 0. 
  \end{equation}
  As is well-known, in two dimensions the Beltrami equation and the quasiconformal mappings give sharp and well adapted  methods to understand linear and non-linear elliptic PDE's such as  \eqref{ELsigma}, see e.g. \cite{ATG}. This is the starting point also in our study of the geometry of limit shapes of dimer models. For reader's convenience we briefly  recall here how this point of view applies to the specific Euler-Lagrange equation \eqref{ELsigma} above. 
  
First, by a $C^{1}$-solution to  \eqref{ELsigma} we mean a function $u \in C^1(\mathcal U)$ such that $\nabla u(z) \in N^\circ \setminus \mathscr{G}$ for all $z \in \mathcal U$ and   \eqref{ELsigma}  holds in the domain ${\mathcal U} \subset \C$.  \,
    If $\, {\mathcal U} $ is  simply connected, the function
    admits an associated conjugate or stream function, defined by $\nabla v = * \nabla \sigma(\nabla u)$, i.e. by the ``Cauchy-Riemann equations''  
\begin{eqnarray} \label{stream}
\;   v_x = - \sigma_y(\nabla u), \quad v_y = \sigma_x(\nabla u).
\end{eqnarray}
This allows one to consider the function
 $F = 2(u + iv) $ 
in ${\mathcal U}$ so that, in particular,
\begin{eqnarray} \label{1toista}
  \left\{ 
 \begin{array}{ll} \hspace{-.1cm} \frac{1}{2}F_x = u_x + i v_x = u_x - i \sigma_y(\nabla u),
  \\ \;  \frac{1}{2}F_y =  u_y + i v_y  =  u_y + i \sigma_x(\nabla u). 
   \end{array} \right.
\end{eqnarray}
 In terms of the complex derivatives $F_{\overline{ z }} = \frac{1}{2}(F_x + i F_y )$ 
  and $F_z = \frac{1}{2}(F_x - i F_y)$ 
   these identities  can be written in the compact form
\begin{eqnarray} \label{1yksi}
F_{\overline{ z }}  = (I - \nabla \sigma)(\nabla u) \quad {\rm and } \quad \overline{\, F_{ z} \,}  = (I + \nabla \sigma)(\nabla u). 
\end{eqnarray}
Combining the two gives us an autonomous Beltrami equation 
\begin{eqnarray}\label{2kaksi2}
F_{\overline{ z }}   = {\mathcal H}( F_{ z} ); \qquad  {\mathcal H}(w) :=  (I-\nabla \sigma) \circ (I + \nabla \sigma)^{-1}(\overline w), 
\end{eqnarray} 
where the non-linear Cayley transform ${\mathcal H} = {\mathcal H}_\sigma$ is called the {\it structure function} of the equation \eqref{ELsigma}. 
Here note that since $\sigma(z)$ is strictly convex on $N^\circ \setminus {\mathscr{G}}$, the transform $\nabla \sigma$ is monotone and thus    ${\mathcal H} $ is well defined.

In multiply connected domains the stream function and 
$F(z)$ 
are only locally
 defined. Hence to use global solutions  in  general domains
 it is advantageous to look for equations for  the complex gradient $f:= F_z$ which is well defined in every domain for any solution to
 \eqref{ELsigma}. 
 
 This applies of course  to a general strictly convex $\sigma$,  where differentiating the corresponding equation  \eqref{2kaksi2}
with respect to $\partial_z$ shows that $f = F_z$ then satisfies the quasilinear Beltrami equation 
   \begin{equation} \label{a}
  f_{\overline{ z }}  = \mu(f) f_z + \nu(f) {\overline{f_z}} \qquad {\rm in} \; \; {\mathcal U},
\end{equation}
where 
$ |\mu(\zeta)| + |\nu(\zeta)| < 1$ for $\zeta \in {\rm {Dom}}({\mathcal H})$, see e.g. \cite[Chapter 16]{ATG}.
\smallskip

However, in our special case of  \eqref{Pst2}   where 
  $\det D^2 \sigma \equiv 1$, we  have the remarkable additional feature that the structure function  ${\mathcal H} = \mathcal H_\sigma \,$  is complex  analytic  !
  \smallskip
 
\begin{Lem} \label{analyticH}
 If $\sigma(z)$ is a convex solution to  the Monge-Amp\'ere equation  \eqref{Pst2}, then the associated 
structure function  ${\mathcal H} = {\mathcal H}_\sigma$ from \eqref{2kaksi2} is complex analytic in its domain,  
    \begin{equation} \label{domain34}
     {\rm Dom} ({\mathcal H}) := \{ w \, : \, \overline{w} \in  \left( I + \nabla \sigma\right)\left(  N^\circ \setminus {\mathscr{G}} \right) \, \} \subset \C.
      \end{equation}
     As a set ${\rm Dom} ({\mathcal H})$ is homeomorphic to 
      $ N^\circ \setminus {\mathscr{G}}$. Moreover, the derivative $|{\mathcal H}'(w)| < 1$ for all $w \in  {\rm Dom}( {\mathcal H})$. 
\end{Lem}

For further important properties of $\mathcal H_\sigma$ see Subsection \ref{Hsection}, and Proposition \ref{Hproper} in particular.

 Note also that the domain ${\rm Dom} ({\mathcal H_\sigma})$ is very closely related to the amoeba $\nabla \sigma(N^\circ \setminus {\mathscr{G}})$. Indeed, there is a one-to-one affine correspondence between the  boundary curves of ${\rm Dom} ({\mathcal H_\sigma})$ and the boundary curves of  $\nabla \sigma(N^\circ \setminus {\mathscr{G}})$, see \eqref{amoeba3} -  \eqref{amoeba34}
and Corollary \ref{subdiff2}. For an illustration of ${\rm Dom} ({\mathcal H_\sigma})$, see Figure \ref{fig.amoeba2}.
\smallskip

\noindent {\it Proof of Lemma \ref{analyticH}.} Via the chain rule, the composition  $a = h \circ b$ of any two $C^1$-functions satisfies
\begin{equation} \label{canalip} a_z {\overline{b_z}} - a_{\overline{ z }}\, {\overline{b_{\overline{ z }}}} = (h_\zeta \circ b) (|b_z|^2 - |b_{\overline{ z }}|^2) \; \quad {\rm and } \quad \; 
 a_{\overline{ z }} \,{{b_z}} - a_{z} {{b_{\overline{ z }}}} = (h_{\overline{ \zeta}} \circ b) (|b_z|^2 - |b_{\overline{ z }}|^2).
\end{equation}
Taking  $h(\zeta) = (I-\nabla \sigma) \circ (I + \nabla \sigma)^{-1}(\zeta)$ and choosing $a(z) =  z - 2\sigma_{\overline{ z }} = z - \nabla \sigma (z)$ with  $b(z)  = z + 2\sigma_{\overline{ z }} $, the left identity in \eqref{canalip} gives 
$$ 1 - \det(D^2 \sigma) =   (h_\zeta \circ b )  \det(D  b).
$$
Here $ \, \det(D  b) = 2 + \Delta \sigma \geq 4$ in $ N^\circ \setminus {\mathscr{G}}$.  
Applying then \eqref{Pst2}  
we  see that ${\mathcal H} = {\mathcal H}_\sigma$ is complex analytic. 

Similarly, the second identity in \eqref{canalip}, with same choices for $a$ and $b$,  gives
$$ |h_{\bar \zeta}|^2 = \frac{|\sigma_{xx} - \sigma_{yy} +2i \sigma_{xy}|^2}{(\Delta \sigma + 2)^2} =  \frac{(\Delta \sigma)^2 - 4}{(\Delta \sigma + 2)^2} \leq 1,
$$
which shows  that ${\mathcal H} = {\mathcal H}_\sigma$ is a contraction. By maximum principle ${\mathcal H}'$ takes values in the open unit disc. Finally, we know that
${\rm Dom} ({\mathcal H})$ is homeomorphic to 
      $ N^\circ \setminus {\mathscr{G}}$, since the map $I+\nabla \sigma$ is a homeomorphism from $N^\circ$ to its image.   \hfill $\Box$

\medskip
 
 In conclusion, we find  a tight connection between the  Euler-Lagrange equation \eqref{ELsigma} and the specific Beltrami equation \eqref{Belt654} below.
 \medskip
 
 \begin{Thm} \label{connection}
 Suppose $\sigma(z)$ is a convex solution to  the Monge-Amp\'ere equation  \eqref{Pst2}, and  $u \in C^{1}({\mathcal U})$ with $\nabla u(z) \in  N^\circ \setminus {\mathscr{G}}$ for $z \in {\mathcal U}$.\;  
  
 If $u$ solves the Euler-Lagrange equation 
 \begin{equation} \label{EL6654}
  \text{div}\, \big(\nabla \sigma(\nabla u)\big)= 0 \qquad {\rm in}  \quad{\mathcal U},
   \end{equation}
 then the identity 
 \begin{equation} \label{identity} \overline{f(z)\,} =  (I + \nabla \sigma)\bigl(\nabla u(z)\bigr), \qquad z \in  {\mathcal U},
  \end{equation}
   defines a solution to the Beltrami equation 
 \begin{equation} \label{Belt654}
  f_{\overline{ z }} = \mathcal H'(f) f_z, 
 \end{equation}
 where $\mathcal H = \mathcal H_\sigma \,$ is the holomorphic structure function associated to $\sigma$. In particular,  $|\mathcal H'(w)| < 1$ for every $w \in Dom(\mathcal H)$.
 \end{Thm}
 \begin{proof} In terms of \eqref{1yksi}, the function $f:= F_z$ is $C^\infty$-smooth by the elliptic regularity theory and by the smoothness of $\nabla \sigma$ in $N^\circ \setminus {\mathscr{G}}$. Moreover,  $\, \overline{f \,} =  (I + \nabla \sigma)\bigl(\nabla u \bigr) $  by definition,  while differentiating the first identity in \eqref{2kaksi2} gives the Beltrami equation \eqref{Belt654}.
 \end{proof} 
 
 Quite remarkably, Theorem \ref{connection} has a direct converse. Namely, for details see Theorem \ref{converse}, given 
 any solution $f(z)$ to  \eqref{Belt654} the identity \eqref{identity}  defines a function $u(z)$, in general domains locally and in simply connected domains globally. Even more: this function $u(z)$ satisfies the Euler-Lagrange equation \eqref{EL6654} whenever $f$ is a solution to \eqref{Belt654}. 
 
 To prove this converse direction we will first need to cover  the basic properties of the so called Lewy transform,
 to be discussed in the next Subsection.

 The following Sections will provide even further ramifications of the above  relations, as well as their connections to the geometry of the surface tension $\sigma$ and its  ``Amoeba'' ${\mathcal A}_\sigma := \nabla \sigma(N^\circ\setminus {\mathscr{G}})$. These  work for any $\sigma$ satisfying $\det D^2 \sigma = 1$ in $N^\circ\setminus {\mathscr{G}}$ and  piecewise affine on $\partial N$, 
 independently of whether or not $\sigma$ arises as  Legendre 
 dual of a Ronkin function  of an algebraic spectral curve.

    \subsection{Surface tension in harmonic coordinates}  \label{complexsec}

  In the interior  $N^\circ \setminus {\mathscr{G}}$ of the gradient constraint  the surface tension $\sigma$ from \eqref{Pst2} is smooth and uniformly convex, thus well under control. However, since with linear boundary values the  strict convexity of  $\sigma(z)$ degenerates on  $\partial N$, the  minimization problem \eqref{ELbasic}, identifying the limiting height functions, presents new and non-standard phenomena  for which appropriate tools will be developed  here and  in the subsequent Sections. 
   
  One of the key points which altogether enables the analysis of the 
  problem \eqref{ELbasic} is the detailed understanding of the boundary behaviour of the surface tensions of different dimer models.  That, on the other hand, one can approach via an analysis based on the complex structure associated to the Monge-Amp\`ere equation, already reflected  in the properties of the non-linear Cayley transform $\mathcal H_\sigma$ in \eqref{2kaksi2}. This point of view turns out to be particularly useful. 
    It will, for instance, also give rise to the representation identities \eqref{general}
 and their extensions to general dimer models, and leads to natural identities involving the different Beltrami equations describing the geometry of the liquid domains.

   More precisely,  we will make use of the fact that in two dimensions  one can express the solutions to \eqref{Pst2} 
     in suitable harmonic coordinates. Different versions of this method go back at least to Lewy, see  \cite[p. 389]{Spiv}. For more recent work see e.g. \cite{FMM}.   
     
  Indeed, as suggested  already  by \eqref{1yksi}, for any surface tension $\sigma(z)$ as in \eqref{Pst2},     it is convenient to  introduce 
 the  following operator, the  {\it Lewy transform}  
     \begin{equation} \label{LewyT}
     L_\sigma(z) := \; \overline{\,z\,} + \overline{\, \nabla \sigma(z)} \; = \;  x -iy+ \sigma_x(z) - i \sigma_y(z), \qquad z=x+iy \in N^\circ \setminus {\mathscr{G}}.
       \end{equation}
  Note that in the above references the Lewy transform is defined without the complex conjugation in \eqref{LewyT}. However, in view of our applications and the  structure of the related Euler-Lagrange equations, 
  see e.g. \eqref{2kaksi2} or Theorem \ref{connection}, Definition \eqref{LewyT} appears  the most flexible for the purposes of this paper.

  Given now a point \,$w   \in  {\rm Dom} ({\mathcal H}) = L_\sigma(N^\circ \setminus {\mathscr{G}})$,
then  \eqref{domain34} gives 
        $ \overline{w } = z + \nabla \sigma (z)$ for some point $z \in N^\circ \setminus {\mathscr{G}}$
     and thus with \eqref{2kaksi2} we have
$  \,  {\mathcal H}(   w  ) = z - \nabla \sigma(z)$.
Adding and arranging  gives 
    \begin{equation} \label{345kolme}
     z = \frac{1}{2} \bigl( \overline{w} +  {\mathcal H}( w )  \, \bigr) \quad \Rightarrow \quad 
     \nabla  \sigma(z) =  \frac{1}{2}\bigl(  \overline{w} -  {\mathcal H}(  w ) \, \bigr), \qquad {\rm for \; \;} w \in  {\rm Dom} ({\mathcal H}), \; z \in N^\circ \setminus {\mathscr{G}}.
   \end{equation}
   
 Here by Lemma \ref{analyticH}  the function ${\mathcal H}(w ) $ is analytic, so that we can express $\nabla \sigma(z)$ in terms of the  {\it harmonic} maps 
  $ \frac{1}{2} \bigl(\,  \overline w \pm   {\mathcal H}(w )   \hspace{.01cm} \bigr)$.  In particular, in these coordinates $\nabla \sigma(z)$ plays the role similar to that of the {\it Hilbert transform}. 
  
  The Lewy transform $L_\sigma: N^\circ \setminus {\mathscr{G}}\to  {\rm Dom} ({\mathcal H_\sigma})$ is a homeomorphism, by Lemma \ref{analyticH} .
  Even more, the identities \eqref{345kolme} and  \eqref{domain34}  tell    that  the inverse of the Lewy transform, 
   \begin{equation} \label{Lewyinverse}
   L_\sigma^{-1}: w \mapsto  \frac{1}{2} \bigl(\,  \overline w +  {\mathcal H}( w ) \,  \hspace{.01cm} \bigr),
   \end{equation} 
   is a harmonic homeomorphism from  ${\rm Dom} ({\mathcal H_\sigma})$ onto $N^\circ \setminus {\mathscr{G}}$.

\begin{figure}[H]
\centering
\begin{tikzpicture}[xscale=1,yscale=1]
\draw[thick,dashed,->,>=stealth] (-2.3,2.3)--(-3,3);
\draw[thick,color=teal,->,>=stealth] (-1.2,0)--(-0.7,0);
\draw[thick,->,>=stealth] (3,3)--(3.2,2.8);
\draw[thick,->,>=stealth] (3,3)--(2.8,3.2);
\draw[thick,color=blue,domain=-1.8:1.8]   plot (cosh{\x}+0.2,sinh{\x}); 
\draw[thick,color=blue,domain=-1.8:1.8 ]   plot (-cosh{\x}-0.2,sinh{\x}); 
\draw[thick,color=blue,rotate=90,domain=-1.8:1.8]   plot (cosh{\x}+0.2,sinh{\x}); 
\draw[thick,color=blue,rotate=-90,domain=-1.8:1.8]   plot (cosh{\x}+0.2,sinh{\x});

\draw (3.7,3.4) node {\small$\vert p_4-p_1\vert$};
\draw [thick,->,>=stealth] (5,-1) to [bend right=20] (3,-1);
\draw (4,-0.5) node {$L_\sigma$};
\draw (2.8,1) node {\small$\overline{p_1} + \overline{\partial \sigma(p_1)}$};
\draw (0,2) node {\small$\overline{p_4} + \overline{\partial \sigma(p_4)}$};
\draw (-3.3,3.4) node {\small$\propto\,\, i(p_2-p_3)$};
\draw (-0.9,-0.3) node {\tiny$\nu(z)$}; 
\draw (0,-0.9) node {$Dom(\mathcal H_\sigma)$};     
\filldraw[color=teal] (-1.2,0)  circle (1.5pt);
\draw (-1.5,0) node {$z$}; 
\end{tikzpicture}
\begin{tikzpicture}[xscale=1,yscale=1]
\draw[thick] (2,0)--(0,2)--(-2,0)--(0,-2)--(2,0);
\draw[thick,color=teal,->,>=stealth] (-2,0)--(-1.5,0);
\draw (2.3,0) node {\small$p_1$};
\draw (0,2.3) node {\small$p_2$};
\draw (-2.3,0) node {\small$p_3$};
\draw (0,-2.3) node {\small$p_4$};
\draw (1,2.3) node {\small$N$};

\filldraw[color=blue] (2,0)  circle (1.5pt);
\filldraw[color=blue] (0,2)  circle (1.5pt);
\filldraw[color=blue] (-2,0)  circle (1.5pt);
\filldraw[color=blue] (0,-2)  circle (1.5pt);

\end{tikzpicture}
\caption{On the right the  Newton polygon $N$ of the uniform domino tilings, and on the left $Dom(\mathcal H_\sigma)$, the image of $N$ under the Lewy transform $L_\sigma$. Boundary components of $Dom(\mathcal H_\sigma)$ and of the amoeba $\nabla \sigma (N^\circ \setminus {\mathscr{G}})$ are in  a one-to-one affine correspondence, c.f. \eqref{LewyT}.\hfill \break
 For details see \eqref{amoeba3} -  \eqref{amoeba34}
and Corollary \ref{subdiff2}.}
\label{fig.amoeba2}
\end{figure}

    As a first consequence, these lead us to the converse of Theorem \ref{connection}.
        
        \begin{Thm} \label{converse}
   Suppose $\sigma(z)$ is a convex function satisfying  the Monge-Amp\'ere equation  \eqref{Pst2}, and let $f:\mathcal U \to Dom(\mathcal H_\sigma)$ be a $C^1$-solution to the Beltrami equation 
      \begin{equation} \label{Belthere}
         f_{\overline{ z }}(z) = {\mathcal H}'_{\sigma} \bigl(f(z)\bigr)\, f_z(z),\qquad z \in \mathcal U,
            \end{equation} 
   in a simply connected domain $\mathcal U \subset \C$.       
   
   Then the vector field $\; L_\sigma^{-1} \circ f\, $ has a potential,
    a $C^1$-function $u$  such that
    $$L_\sigma^{-1} \circ f(z)  = \nabla u(z), \qquad z \in \mathcal U.$$
    Moreover
     \begin{eqnarray} \label{locsol}
 {\rm div} \, \big(\nabla \sigma(\nabla u)\big)= 0 \quad {\rm in} \;\; \mathcal U.
         \end{eqnarray}

        \end{Thm}
        {\it Proof.} In view of  \eqref{Lewyinverse} the complex derivatives of the inverse Lewy transform are
      given by 
      $$\partial_w\,  L_\sigma^{-1}(w) = \frac{1}{2}  {\mathcal H'_\sigma}( w ), \quad {\rm while} \quad  \partial_{\overline w } \,L_\sigma^{-1}(w) = \frac{1}{2}. $$
  Using the chain rule we have
   $$ \partial_z \left[ L_\sigma^{-1} \circ f \right] = \partial_w\,  L_\sigma^{-1}(f)\, \partial_z f + \partial_{\overline w } \,L_\sigma^{-1}(f) {\overline{f_{\overline z}\,}} = \Re e\left[ {\mathcal H}'_{\sigma}(f)\, f_z \right],
   $$
  in particular, 
$${\rm curl} \left[ L_\sigma^{-1} \circ f \right]  \equiv 2 \, \Im m \left( \partial_z \left[ L_\sigma^{-1} \circ f \right] \right) = 0.$$
Since $L_\sigma^{-1} \circ f \,$ is a curl-free vector field in the simply connected domain $\mathcal U$, we have $L_\sigma^{-1} \circ f = \nabla u$ for some $C^1$-function $u(z)$. 

In addition,  from \eqref{345kolme} - \eqref{Lewyinverse} we see that 
  \begin{eqnarray} \label{gradsi}
 \nabla\sigma \circ \nabla u =    \nabla\sigma \circ ( L_\sigma^{-1} \circ f ) =  \frac{1}{2} \left(\, \overline{f \,}  - \mathcal{H_\sigma}(f) \right). 
   \end{eqnarray} 
 Therefore 
  $$  {\rm div} \, \big(\nabla \sigma(\nabla u)\big) \equiv 2 \, \Re e\,  \left[  \partial_z  \big(\nabla \sigma \circ \nabla u)\big) \right] = 
    \, \Re e\,  \left[ \,  \overline{f_{\overline z} } -  \mathcal{H'_\sigma}(f) f_z \right] = 0. \hspace{3cm} \Box$$
    \smallskip
         
         \begin{rem} In general domains $\mathcal U \subset \C$ the solutions $f$ to the Beltrami equation \eqref{Belthere} of course   give locally rise to solutions  of the Euler-Lagrange equations \eqref{locsol}, and globally whenever the vector field $L_\sigma^{-1} \circ f$ has zero monodromy.  Note also, that in case $\sigma$ has gas points, then by Theorem \ref{GAS} the liquid domains with frozen boundary are never simply connected.  
         \end{rem}
         \smallskip

        Next, composing the Lewy transform with a Riemann map of the domain $Dom({\mathcal H_\sigma})$
     reveals  beautiful properties for the transform. In the case where there are no gas points $q \in  {\mathscr{G}}$,
     these are particularly transparent.
     \medskip

   \begin{Thm} \label{harmonic rep} 
   Let $N \subset \R^2$  be a convex polygon and  $\sigma(z)$  a convex solution to $\det  \big(D^2 \sigma\big)=1$
   in $N^\circ$, with $\sigma{\big|_{\partial N}}$ convex and piecewise affine. 
   
   Also, let 
   $\mathscr P \cup \mathscr Q = \{ p_j \}_1^{m}$ consist of 
   the corners of $N$ and of the points of discontinuity of $\sigma{\big|_{\partial N}}$,  with the cyclic (counterclockwise) 
   order  induced by $ \partial N$.
Finally, assume  $ \mathscr  G = \emptyset$.  
   
  If  $\psi: \Di \to {\rm Dom} ({\mathcal H}) $ is a Riemann map, then   $\, U(\zeta) := L_\sigma^{-1} \circ \psi(\zeta) \, $ defines 
  a harmonic homeomorphism  $U = U_\sigma: \Di \to N^\circ$, and it 
   has  the representation   
      \begin{eqnarray} \label{harmonichomeo31}
    U(\zeta) := L_\sigma^{-1} \circ \psi(\zeta) = \sum_{j=1}^m \, p_j \, \omega_\Di(\zeta; I_j), \qquad  \zeta \in \Di,
      \end{eqnarray} 
  where $I_j \subset \partial \Di$ with pairwise disjoint open arcs whose closure covers the unit circle, and where $\omega_\Di(\zeta; I_j)$ is the harmonic measure of $I_j$ in $\Di$.
     \end{Thm}
\begin{proof}    
   Composing  \eqref{345kolme}-\eqref{Lewyinverse} with  the Riemann map $\psi: \Di \to  L_\sigma(N^\circ) = Dom({\mathcal H})$ gives
     \begin{equation} \label{harmonichomeo1}
          U(\zeta) 
          =  \frac{\overline{\, \psi(\zeta) }  +   {\mathcal H}( \psi(\zeta) ) }{2}, 
          \end{equation} 
and  
          \begin{equation} \label{harmonichomeo21}
          \nabla  \sigma \circ  U(\zeta) = \frac{\overline{\, \psi(\zeta) }  -  {\mathcal H}( \psi(\zeta) ) }{2}, \qquad \zeta \in  \mathbb \Di.
\end{equation} 
  
  Here  $U(\zeta)$ is  a harmonic homeomorphism from the unit disc onto $N^\circ$. It is sense reversing and  has  negative Jacobian, since ${\mathcal H}$ is a strict contraction.
We may thus use a theorem of Hengartner and Schober \cite[Theorem 4.3]{HS}, see also \cite[p.35]{Dur}, to the (sense preserving) harmonic homeomorphism $\zeta \mapsto U(\overline{\,\zeta\,})$. Their work implies  that  outside a countable set  $E \subset \partial \Di$ there exist the  unrestricted limits
$$ U(e^{ix}) := \lim_{ \Di  \ni \, \zeta \to e^{ix}} U(\zeta) \in \partial N, \qquad e^{ix} \notin E,
$$
while for the exceptional points $e^{ix} \in E$, the cluster set of $U(\zeta)$ at $e^{ix}$ is a non-degenerate segment $J \subset \partial N$. 
Moreover,  in the complement of $E$ the boundary function $e^{ix} \mapsto U(e^{ix})$ is continuous and sense reversing.
           
           If now $U(e^{ix_0}) \in \partial N \setminus \{p_j\}_{j=1}^m$ for some point $e^{ix_0} \in \partial \Di \setminus E$ outside the exceptional set, then by  continuity at $e^{ix_0} $ and 
 Lemma \ref{blowup}, $$ \vert \nabla  \sigma \circ  U(r e^{ix})\vert  \to \infty\quad {\rm as} \quad r \to 1,$$ whenever $|x-x_0|$ is small with $e^{ix} \notin E$.
  But this is not possible: From \eqref{harmonichomeo1} - \eqref{harmonichomeo21} we see that  $U - \nabla \sigma \circ U$ is analytic so that $- \nabla \sigma \circ U$ is, up to an additive constant,  the harmonic conjugate of $U$. Since $U(\zeta)$ is bounded, by the Theorems of Riesz \cite[Theorem III.2.3]{Garnett} and Fatou \cite[Theorem I.5.3]{Garnett}, the function  $\nabla \sigma \circ U$ has almost everywhere finite radial boundary values.
  
  In conclusion, the unrestricted limits of $U(\zeta)$ at  points outside $E$ 
  are all contained in the finite set $\mathscr{P}\cup \mathscr Q $. 
   As a sense reversing map the boundary function is thus piecewise constant, 
  and has the representation  
 \begin{equation} \label{summa3}
   U(\zeta) =  \sum  p_{k_j} \, \omega(\zeta; I_j) 
    \end{equation}
 where the intervals $I_j \subset \partial \Di$, covering the unit circle, have pairwise disjoint interiors.

   On the other hand,  the cluster set of $U(\zeta)$ at any point in $E$  is a line segment of $\partial N $.  
  Hence each corner of $N$ appears in the sum  \eqref{summa3}.  
  
  For the points $p_j \in \mathscr{Q}$ we need a further argument. 
  For this recall that the harmonic measure of the (counterclockwise oriented) arc $I \subset \partial \Di$ between the points $\eta_1, \eta_2 \in \partial \Di$
  is given by  
  \begin{equation} \label{harmonic345}
    \omega(\zeta; I)  
  = \frac{1}{\pi} \Im m \, \log\left( \frac{\zeta - \eta_2}{\zeta - \eta_1}  \right) + c(I), \qquad \zeta \in \Di.
    \end{equation}
  Next note that by  \eqref{harmonichomeo1} - \eqref{harmonichomeo21} the difference
   $U - \nabla \sigma \circ U$ is analytic so that $- \nabla \sigma \circ U$ is, up to an additive constant,  the harmonic conjugate of $U(\zeta)$.
  Combining this with \eqref{summa3} - \eqref{harmonic345} one sees that if  arcs $I_j$ have the endpoints $\eta_{j-1}$ and $\eta_j$, then
  \begin{equation} \label{conju43}
   \nabla \sigma \circ U(\zeta) = \frac{i}{ \pi}  \sum  p_{k_j} \, \log\left| \frac{\zeta - \eta_{j}}{\zeta - \eta_{j-1}}  \right| + c_0 = 
    \frac{1}{i \pi} \sum  (p_{k_j} - p_{k_{j+1}}) \log \frac{1}{|\zeta-\eta_j|} + c_0,
      \end{equation}
    where  $c_0$ is a constant. 
    
    Next, for a fixed index $j= j_0$,  the term $  i(p_{k_{j_0}} - p_{k_{{j_0}+1}})$ is  the outer normal to $N$ on the side $[p_{k_{j_0}}, p_{k_{{j_0}+1}}] \subset \partial N$. Moreover,   the cluster set of $U(\zeta)$ at $\zeta = \eta_{j_0}$ equals the interval  $ [p_{k_{j_0}}, p_{k_{{j_0}+1}}]$. Therefore the derivative of $\sigma$ in the direction of the tangent  on $[p_{k_{j_0}}, p_{k_{j_0}+1}] \subset \partial N$ 
    is determined from
    \begin{eqnarray} \label{der.tangential}
  |p_{k_{j_0}} - p_{k_{{j_0}+1}}|(\partial_T \sigma) \circ U(\zeta) = \langle p_{k_{j_0}} - p_{k_{{j_0}+1}} , \nabla \sigma \rangle \circ U(\zeta) \qquad  \qquad \nonumber \\
    = \sum_{j \neq j_0}   \langle p_{k_{j_0}} - p_{k_{{j_0}+1}}, \frac{1}{i\pi}(p_{k_j} - p_{k_{j+1}})\rangle \log \frac{1}{|\zeta-\eta_j|} + c_1.
    \end{eqnarray}
   The expression shows that along the side $(p_{k_{j_0}}, p_{k_{{j_0}+1}})$ the tangential derivative of $\sigma$ is continuous.     
    But that means that no quasifrozen point $p \in \mathscr{Q}$ can lie  in such an open interval, rather,  they  all must be among the image points $p_{k_j}$ in the sum \eqref{summa3}. 
    Thus  the representation \eqref{harmonichomeo31} follows.
    \end{proof}

     \begin{rem}  Note that in the notation  of \eqref{harmonichomeo31} with   $I_j = (\eta_{j-1}, \eta_{j})$,   the set of non-tangential limits of $U(\zeta)$  at the endpoint $\eta_j$ is precisely the open interval $(p_j,p_{j+1}) \subset \partial N, \; j= 1,\dots, m$.  

          \end{rem} 
      As we will see later  in Proposition \ref{prop.hrep}, with a  more explicit picture of the geometry of the domain $ L_\sigma(N^\circ) = {\rm Dom} ({\mathcal H})$
      one can as well  represent the inverse Lewy transform directly in terms of  harmonic measures on ${\rm Dom} ({\mathcal H})$.
      This works independently of whether $\sigma$ has gas points or not.

     For later purposes we also formulate as a separate result the above expression \eqref{conju43} of  the gradient of the surface tension in harmonic coordinates. In fact, it appears quite  useful and remarkable that such a presentation exists in explicit coordinates for general surface tensions of Theorem \ref{harmonic rep}. 
        \medskip
        
       \begin{Cor} \label{sigmarep} Suppose 
       $N$,         $\sigma(z)$ and        $\, U(\zeta) = L_\sigma^{-1}  \circ  \psi(\zeta) \, $ with  $\psi: \Di \to {\rm Dom} ({\mathcal H}) = L_\sigma(N^\circ)$ are as in Theorem \ref{harmonic rep}.
 Then
           \begin{eqnarray}  \label{summa4}
              \nabla \sigma \circ U(\zeta) =
    \frac{1}{i \pi} \sum_{j=1}^m  \,  (p_{j} - p_{j+1}) \log \frac{1}{|\zeta-\eta_j|} + c_0,\qquad \zeta \in \Di,
         \end{eqnarray}
         where  the arcs $\{ I_j\}$ are as in Theorem \ref{harmonic rep}, the $\{ \eta_j\}$ are their endpoints and where $c_0$ is a constant.
       \end{Cor}

    \medskip

 In case of the general Monge-Amp\`ere equation  \eqref{Pst2} with gas points, the domain $ {\rm Dom} ({\mathcal H}_\sigma)  = L_\sigma(N^\circ \setminus \mathscr G)$ is multiply connected.  We can still uniformise $ {\rm Dom} ({\mathcal H})$  but now with a conformal map from a circle domain. In the end we arrive at   a  representation similar to \eqref{harmonichomeo31}, for details  
 see Proposition \ref{prop.hrep}.  
       \medskip

     \subsection{Representing height functions and the universal Beltrami equation} \label{universal}

   As a next step, given a solution to 
    \begin{equation} \label{EL654new}
  \text{div}\, \big(\nabla \sigma(\nabla u)\big)= 0 \qquad {\rm in}  \quad \mathcal U,
   \end{equation} 
    it is now natural to combine the identity \eqref{identity} with the representation \eqref{harmonichomeo31} for the inverse Lewy transform. However, the  representation \eqref{harmonichomeo31} is parametrised by the unit disc 
   $\Di$, while the solution $f$ in \eqref{identity}  takes values in $ {\rm Dom} ({\mathcal H_\sigma})$. There are two (basically equivalent) ways to settle this issue,       either use the harmonic measure of $ {\rm Dom} ({\mathcal H})$, as will be done later  in Proposition \ref{prop.hrep}, or as here,
   use the conformal invariance in the Beltrami equation to study solutions $f: \mathcal U \to \Di$ with target the unit disc.

   Namely, if the surface tension $\sigma$ has no gas points, then  $ {\rm Dom} ({\mathcal H_\sigma})$ is simply connected, and we have the conformal parametrisation $\psi: \Di \to  Dom({\mathcal H_\sigma})$ as in  \eqref{harmonichomeo1} - \eqref{harmonichomeo21}. This allows one to define 
     \begin{equation} \label{mu12}
 \mu_\sigma(z) := {\mathcal H_\sigma'} \circ \psi(z), \qquad z \in \Di. 
   \end{equation} 
   By definition $ \mu_\sigma: \Di \to \Di$ is analytic, and we will see later in Proposition \ref{Hproper} that $ {\mathcal H_\sigma'}$ is a proper map. Thus in fact, $ \mu_\sigma$ is a Blaschke product. Note, however, that $ \mu_\sigma$ is uniquely defined only up to a precomposition with a M\"obius transform.
   
    Replacing the complex gradient  $f = F_z$ from \eqref{2kaksi2} by the composition $ f = \psi^{-1} \circ F_z:  \mathcal U \to \Di$, then
via  \eqref{Belt654} and the chain rule, we see that  the Beltrami equation for this new function takes the form
     \begin{equation} \label{beltB1}
   f_{\overline{ z }}(z) = \mu_{\sigma} \bigl(f(z)\bigr)\, f_z(z),\qquad z \in \mathcal U.
  \end{equation}
  This  setting allows an explicit representation for solution to the Euler-Lagrange equation \eqref{EL654new}.
  \medskip
  
       \begin{Thm} \label{repone}  Suppose $\sigma$ is a surface tension as in   \eqref{Pst2}, defined  in a polygon  $N$ with corners and quasifrozen points $\{ p_j\}_1^m $, and assume that there are no gas points, $\Gg = \emptyset$.

 Then if $\,u $ is a $C^1$-solution to the equation $\, {\rm div} \, \big(\nabla \sigma(\nabla u)\big)= 0$ in a bounded domain 
 $ \mathcal U \subset \R^2$,  it has the representation
     \begin{equation} \label{repfirst}
 \nabla u(z) = \sum_{j=1}^m \, p_j \, \omega_\Di(\, f(z) ; I_j), \qquad z \in \mathcal U,
   \end{equation}
   where $f: \mathcal U \to \Di$ is a $C^1$-solution  to the Beltrami equation \eqref{beltB1}.
       \end{Thm}
       \begin{proof} Let first $f_0$ be the function determined by the identity \eqref{identity}, solving the Beltrami equation \eqref{Belt654}. Put briefly, by \eqref{identity} $f_0  = L_\sigma (\nabla u)$ while as discussed above, $f = \psi^{-1} \circ f_0$ solves \eqref{beltB1}, with $\psi: \Di \to Dom(\mathcal H_\sigma)$  the Riemann map. Thus  Theorem \ref{harmonic rep} gives
        \begin{equation} \label{repsecond}
         \nabla u = L_\sigma^{-1} \circ f_0 = U ( \psi^{-1} \circ f_0) =  \sum_{j=1}^m \, p_j \, \omega_\Di(\, f(z) ; I_j),
         \end{equation}
  proving the claim.    
   \end{proof}
   
    By  composing with the Riemann map, one can of course transform any solution of \eqref{beltB1} back to  a solution of the original Beltrami equation \eqref{Belt654}. For instance, doing this transformation in the setting of Theorem \ref{converse} shows that   also Theorem \ref{repone} has direct converse. 
         
     \begin{Cor} \label{Dconverse}
     Suppose $\sigma$ is a surface tension as in Theorem \ref{repone}, and $U = U_\sigma: \Di \to N^\circ$ is the harmonic homeomorphism from \eqref{harmonichomeo31}.
     
    If  $f: \mathcal U \to \Di$ is any $C^1$-solution to the Beltrami equation \eqref{beltB1} in a simply connected domain $\mathcal U  \subset \C$, then the identity
    $$ \nabla u(z) = U_\sigma \circ f(z) =  \sum_{j=1}^m \, p_j \, \omega_\Di(\, f(z) ; I_j), \qquad z \in \mathcal U$$
     defines a $C^1$-solution $u$ to the Euler-Lagrange equation \eqref{EL654new}.
     \end{Cor}

  When the specific properties of a particular dimer model or surface tension $\sigma$ are not needed, for instance in questions such as the study of the geometry of  frozen boundaries within general dimer models,  then  there exists yet {\it a third}  approach or version in describing    the properties of the solutions to the Euler-Lagrange equation \eqref{EL654new} with help of Beltrami equations. 
  Namely, given a solution to
          \begin{equation} \label{Specific}
          \partial_{\overline{ z }} f(z) = \mathcal{H}'_\sigma\bigl( f(z) \bigr) \,  \partial_z f(z), \qquad z \in  \mathcal U, 
            \end{equation}
 call this $ \widehat f$,  using again  the conformal invariance and the chain rule we see that the composition $f(z) :=  \mathcal{H}'_\sigma\bigl( \widehat f(z) \bigr) $
 is a solution to 
      \begin{equation} \label{belt2}
\partial_{\overline{ z }} f(z) = f(z)  \partial_z f(z), \qquad z \in \mathcal U.
  \end{equation}
  
  \begin{Def} We call \eqref{belt2} the Universal Beltrami equation.
  \end{Def}

  The name for equation \eqref{belt2} is a natural one, as any solution to any of the Beltrami equations \eqref{Specific} or \eqref{beltB1} that we study in this paper,
  gives rise to a solution of  \eqref{belt2}.  In particular, this works with all surface tensions $\sigma$ as in \eqref{Pst2}, whether  or not $\sigma$ has gas points. 
  
   However, there is even  more to this terminology. Namely, see Remark \ref{domH},
 for the uniform lozenges model we can take $\mu_{\sigma} (z) \equiv z$ in \eqref{beltB1}. Therefore any solution to the universal equation, and thus in particular any solution to  \eqref{beltB1}, 
 gives via Corollary \eqref{Dconverse} rise at least to a candidate for a height function in the lozenges model, thus suggesting an  approach to Theorem \ref{Second.thm}, the universality of the lozenges geometry. For further details, see also the discussion after  Remark \ref{domH}. 

  Developing this theme in depth and, in particular, showing that one indeed arrives to actual lozenges height functions  requires a variety of different methods to be developed in subsequent sections.

The functions solving \eqref{belt2} are $C^\infty$-smooth and locally quasiregular, but in general not much more can be said of the  solutions  or  the domains $\LL$ supporting them. However,  if the solution happens to  define a proper map $f: \LL \to \Di \,$, so that  the boundary $\partial \LL$ is roughly the locus of points with $|f(z)| = 1$, the situation changes quite dramatically. Then the equation in a sense determines the geometry  of the   liquid domain  $\LL$, see e.g. Section \ref{B+hodo} and \ref{propersec}.

      \subsection{Geometry of the surface tension} \label{surface}

 We have seen in Theorems \ref{connection} and \ref{converse} that there is basically a one-to-one correspondence between solutions of the Euler-Lagrange equation  \eqref{ELsigma} and those of the Beltrami equations  \eqref{Belt654} and \eqref{beltB1}.  
 This connection also gave rise to the Lewy transform and to the  representation theorems for surface tensions $\sigma$, as well as a representation for the solutions to the Euler-Lagrange equation, thus in particular for the limiting height functions.
 
 However, there is one aspect of this interaction that is still missing, namely the relations between the boundary behaviours of  these objects, necessary of course to study and describe the geometry of frozen boundaries. In many respects, the  basis to this is the clear understanding of the geometry and boundary behaviour of surface tension $\sigma$, to be developed in this subsection.
  It will be of independent interested to discuss this within the geometry of a general convex surface tensions and assume of $\sigma$ only  that it satisfies the Monge-Ampere equation  \eqref{Pst2}.   

 With this view, a particular goal of this Subsection is  to prove Theorems \ref{sigma1} and \ref{prop:sigma:property3} and  Proposition \ref{prop:sigma:property}. 
 
 Our starting point is  the representation \eqref{summa4} of $\nabla\sigma$ in the harmonic coordinates as described in Subsection \ref{complexsec}, c.f. also Figure \ref{fig.amoeba}. As recalled at the end of  
 Section  \ref{sect:introduction5}	 
  the complement of the image  $\nabla \sigma (N^\circ \setminus {\mathscr{G}})$ consists of  the disjoint union of the subdifferentials  $\partial \sigma(p)$, where $p \in \mathscr{P} \cup \mathscr{Q} \cup \mathscr{G}$. 
   For general convex maps, the subdifferentials $\partial \sigma(p)$ are convex and closed sets, while as we will see
 for solutions to \eqref{Pst2} the subdifferentials have  quite  specific extra features.
 
\medskip

\begin{Lem}  \label{subdiff} 
Suppose $\sigma$ is a convex solution to \eqref{Pst2}.
Then  for each $p_j \in \mathscr{P} \bigcup \mathscr{Q}$, the boundary of the subdifferential $\partial \sigma(p_j)$ is an unbounded, convex and analytic Jordan arc $\gamma_{p_j}$. 

Moreover, if we give $\mathscr{P} \bigcup \mathscr{Q} = \{ p_j\}$  the order induced by $\partial N$, then  at the  endpoints at  
$\infty$ the tangents of $\gamma_{p_j}$ are  orthogonal, respectively,  to $[p_{j-1}, p_{j}]$ and  to $[p_{j}, p_{j+1}]$.
\end{Lem}

  \begin{proof} 
  
  \noindent $1^\circ$. We study first the surface tensions $\sigma$ without gas points,   $\mathscr{G} = \emptyset$.
In this case one can use 
the harmonic homeomorphism from  \eqref{harmonichomeo31},  the map $U = L_\sigma^{-1} \circ \psi(\zeta) : \Di \to N^\circ$ where  $\psi: \Di \to L_\sigma(N^\circ) = {\rm Dom} ({\mathcal H}) $ is the Riemann map.
  
  With the notation of  \eqref{harmonichomeo31},  the boundary values of $U(\zeta)$ take the interval $I_{j} = (\eta_{j-1}, \eta_j) \subset \partial \Di$ to $p_j$, for each $j=1,\dots, m$. Moreover, \eqref{harmonichomeo1} - \eqref{harmonichomeo21} show that  
      \begin{equation} \label{riemann}
  \overline{\psi(\zeta)} = (\nabla \sigma \circ U)(\zeta) + U(\zeta), \quad \zeta \in \Di.
    \end{equation}
     Since  for any convex function the set of subdifferentials is closed, we thus obtain
    \begin{equation} \label{unifmap} 
  \, \lim_{r\to 1} \overline{ \psi(r\eta) } = (\nabla \sigma \circ U)(\eta) + p_j  \in \partial \sigma(p_j) + p_j, \qquad {\rm for}  \; \eta \in I_{j},
  \end{equation}
  where Corollary \ref{sigmarep} shows the existence of the above limit, 
  for each $\eta$ in the open interval $ I_{j}$. 
  
    On the other hand, since $\sigma$ is convex and by \eqref{Pst2} smooth and strictly convex in $N^\circ\setminus {\mathscr G}$, any 
   boundary point of $\partial \sigma (p_j)$ is obtained as
   a subsequential limit of $ (\nabla \sigma \circ U)(\eta)$ where  $\eta \to I_j$ in $\Di$. 
   
 With the  representation of Corollary \ref{sigmarep} we also see that $\gamma_{p_j}(\eta) := \overline{\psi(\eta)} - p_j$, for $ \eta \in  I_{j}$, is an analytic curve, unbounded at the end points and  a parametrisation of the entire boundary of the set  $\partial \sigma(p_j)$. Since subdifferentials are always closed convex sets, it follows that   $\gamma_{p_j}(I_j)$ is a convex Jordan arc. Moreover, 
  $$            \nabla \sigma \circ U(\zeta) =
    \frac{1}{i \pi} \,  (p_{j} - p_{j+1}) \log \frac{1}{|\zeta-\eta_j|} + c(\zeta), \qquad {\rm as \; }\;  \zeta \to \eta_j,
  $$
with $c(\zeta)$ bounded near $\eta_j$. Hence the boundary of $\partial \sigma(p_j)$ is asymptotically  orthogonal to $p_{j+1} - p_{j}$ when 
$\eta \to \eta_j$ on $I_j$. Similarly, the boundary becomes orthogonal to  $p_{j} - p_{j-1}$ when $\eta \to \eta_{j-1}$ on $I_j$. 

This completes the proof in the case where there are no gas points  $q \in {\mathscr G}$. For later purposes we note the following additional observation: Since $\gamma_{p_j}(\eta) = \nabla \sigma \circ U(\eta) $ on $I_j$, we see via  Corollary \ref{sigmarep}  and \eqref{harmonic345} 
       that $\gamma_{p_j}(\eta)$ admits a complex analytic extension across $I_j$.     Since  the curve $\gamma_{p_j}(I_j)$ is convex, it follows that $\gamma_{p_j}'(\eta) \neq 0$ on the interval $I_j$. Combining this with the Cauchy - Riemann equations and \eqref{riemann} shows that the radial derivative 
    \begin{equation} \label{radial}
    \partial_r U(\eta) \neq 0, \qquad \forall \;  \eta \in I_j.
  \end{equation}

\smallskip

 \noindent $2^\circ$. If $\sigma$ has gas points $q \in {\mathscr G}$   and  $ p_j \in \mathscr P\cup \mathscr Q$ is given,  it is natural to modify the argument above and consider a restriction of $\sigma$ to a neighbourhood of $p_j$. That is, choose $\rho_1\in (p_{j-1},p_{j})$ and $\rho_2\in (p_j,p_{j+1})$ on $\partial N$,
and connect them with a smooth strictly convex arc  $\Gamma \subset {N}^\circ$, intersecting $\partial N$ non-tangentially at the end points $\rho_1$ and $\rho_2$.
 Further, call $M^\circ$ the domain bounded by
$[\rho_1,p_{j}] \cup [p_j, \rho_2] \cup \Gamma$. 

For $\rho_1,\rho_2$ close to $p_j$ we can assume  that $M^\circ $ 
contains no points $q \in {\mathscr G}$. Then $\sigma$ is affine on $[\rho_1,p_{j}]$ and on  $[p_j, \rho_2]$, but not any more on $\Gamma$, where it is strictly convex and smooth. On the other hand, taking the restriction to $M^\circ $ does not change the subdifferential $\partial \sigma(p_j)$.

Moreover, one can still use the Lewy transform restricted to $M^\circ $, which remains to have a harmonic inverse. Thus  if ${\psi}_M: \Di \to L_\sigma(M^\circ)$ is a Riemann map,
we now study the harmonic homeomorphism $U_M = L_\sigma^{-1} \circ {\psi}_M: \Di \to M^\circ$.
The same analysis as in Theorem \ref{harmonic rep}, using the Hengartner-Schober theorem,  shows that we can write
\begin{equation} \label{Udecomp}
 U_M(z) = W_M(z) 
 + p_j \, \omega_\Di(z; I),
\end{equation}
where $I = (\eta_{1}, \eta_2) \subset \partial \Di$  a non-degenerate interval, 
and $W_M(z)$ is  the Poisson integral of the (continuous) boundary values of $U_M$ on the complement $\partial \Di \setminus I$, with 
$\Gamma = {W}_M(\partial \Di \setminus I)$. 

In addition, via \eqref{harmonichomeo1} - \eqref{harmonichomeo21} we observe that  $ U_M(\zeta)  - (\nabla \sigma \circ U_M)(\zeta)$ is analytic, in other words, $i \,\nabla \sigma \circ U_M$ is up to an additive constant  the harmonic conjugate of $U_M(\zeta)$. Since $W_M(\eta)$ vanishes on $I$, with \eqref{Udecomp} we thus see that $\nabla \sigma \circ U_M(\zeta)$ has a radial limit at each $\eta \in I$.

As in step $1^\circ$  this leads to a parametrization 
of the boundary of $\partial \sigma (p_j)$ 
 by $\gamma_{p_j}(\eta) = (\nabla \sigma \circ U_M)(\eta)$ for $\eta \in I$, where  $\gamma_{p_j}$ is an unbounded, convex and analytic Jordan arc.

 Last, to identify the tangent directions of $\gamma_{p_j}$ at end points, via \eqref{harmonichomeo1} - \eqref{harmonichomeo21} we observe that  $ U_M(\zeta)  - (\nabla \sigma \circ U_M)(\zeta)$ is analytic, in other words, $i \,\nabla \sigma \circ U_M$ is up to an additive constant  the harmonic conjugate of $U_M(\zeta)$. 
 Since $W_M(\eta)$ vanishes on $I$ and $W_M(\eta) \to \rho_1$ as  $\eta \to \eta_1$ on $\partial \Di \setminus I$, the conjugate Poisson representation shows that   on the arc $I$ near $\eta_1$ and  for any $\varepsilon > 0$, 
   $$            \nabla \sigma \circ U_M(\eta) =
  \frac{  (p_{j} - \rho_1)+ b_{\varepsilon}(\eta)}{i \pi} \, \log \frac{1}{|\eta-\eta_1|} + C_\varepsilon, \qquad \eta \in I,
  $$
  where $C_\varepsilon$ is a constant and $\|b_{\varepsilon}\|_{\infty} \leq \varepsilon$. 
    
  Therefore again,  the  tangent to the boundary curve $\gamma_{p_j}(\eta)$ becomes orthogonal to  $p_{j} - \rho_1$, i.e. to $p_{j} - p_{j-1}$, when $\eta \to \eta_1$ on $I$. Similarly, the tangent becomes orthogonal  to  $p_{j+1} - p_j$, when $\eta \to \eta_2$ on $I$. 
  
As a last point, a similar argument as above shows that  $  \partial_r U_M(\eta) \neq 0$ for all $ \eta \in I $. 
     \end{proof}

  An analogous description holds also for the subdifferentials at the possible gas points. To see this we need an auxiliary result, useful also elsewhere.

 \begin{Lem} \label{help2}
Suppose $I \subset \SSS^1$ is an open interval, let $0 < \rho < 1$ and assume that $\alpha(z)$, $\beta(z)$ are  bounded analytic functions  in the strip
 $${\mathbb A} = \left\{z: \rho < |z| < 1, \; \frac{z}{|z|} \in I \right\}.$$ Assume further that 
$$ \alpha(z) - \overline{\beta(z)} \to 0  \;\; \mbox{   \, as } \;\; z \to I  \;\; \mbox{   \, in } \;\; {\mathbb A}.
$$
Then $\alpha(z)$ and $\beta(z)$ extend analytically to the double strip ${\mathbb A}_1 = \{z: \rho < |z| < 1/\rho, \; \frac{z}{|z|} \in I  \}$.
\end{Lem}
\begin{proof}
Consider the auxiliary function 
\begin{align} \label{auxiH}
H(z)= \left\{ 
 \begin{array}{ll}  \overline{\beta(1/\overline {z} )},  \quad 1 < |z| < 1/\rho,  \; \frac{z}{|z|} \in I,
  \\  \alpha(z), \qquad \;  \rho < |z| < 1,\; \;   \; \frac{z}{|z|} \in I.
   \end{array} \right.
\end{align}
Then $H(z)$ is analytic and bounded in ${\mathbb A}_1 \setminus I$. Moreover,  by Fatou's theorem \cite[p. 29]{Garnett}, at almost every $e^{i\theta} \in \partial \Di$ we have the non-tangential  limits  $\alpha( e^{i\theta})$ and $\overline{\beta(1/\overline {\,e^{i\theta} \,} )}$.
By our assumption, these  limits  are equal for almost everywhere on $ I$. Then by Morera's theorem $H(z)$
extends analytically across $I$, see e.g. \cite[p. 95]{Garnett}.
\end{proof}

  \begin{Lem}  \label{subdiff12} 
 If  $\sigma$ is a convex solution to \eqref{Pst2}, then for each gas point $q \in \mathscr{G}$, the boundary of $\partial \sigma(q)$ is a bounded, convex and analytic Jordan curve. 
\end{Lem}
  \begin{proof}  
Near each   gas point $q_0 \in {\mathscr G}$ we have $\, \det D^2 \sigma(z) = 1 + c \, \delta_{\{q_0\}}$. Thus  the Jacobian of $\nabla \sigma$ is a positive measure at $q_0$,  with area  $|L_\sigma(\{q_0\})| = 
|\nabla \sigma(\{q_0\})| > 0$. In particular,   given a small disc $\Di_0 = \Di(q_0,\varepsilon)$,  the image of   $\Di_0\setminus \{q_0\}$ under the Lewy transform  is a doubly connected domain,  call it $W_0$, having non-degenerate boundary components. 

Thus $W_0$ is conformally equivalent to an annulus $\mathbb{A}(1,R) = \{z: 1 < |z| < R\}$, and one can now compose the harmonic inverse $L_\sigma^{-1}$ with a Riemann map $\psi:\mathbb{A}(1,R) \to W_0$, giving us the homeomorphism  
  \begin{equation} \label{333three1}
   U := L_\sigma^{-1} \circ \psi : \mathbb{A}(1,R) \to \Di_0\setminus \{q_0\}; \qquad U(\zeta) \to q_0  \quad {\rm as } \quad  |\zeta| \to 1.
  \end{equation}

  In the annulus the harmonic map decomposes as   $U(\zeta) = q_0 + \alpha(\zeta) + {\overline{\, \beta(\zeta) \, }} + C 
  \log|\zeta|, $ for bounded analytic functions $\alpha(\zeta)$ and $\beta(\zeta)$  and for a constant $C$. Moreover, with the boundary condition in \eqref{333three1} we have $ \alpha(\zeta) + {\overline{\, \beta(\zeta) \, }} \to 0$ on the unit circle. Thus the above Lemma \ref{help2} shows that  we can extend $\alpha(\zeta)$ analytically across the unit circle by setting
  $$  \alpha(\zeta) = - {\overline{\, \beta \left(1 / \, {\overline{\zeta}} \,\right)}}, \qquad \frac{1}{R} < |\zeta| \leq 1.
  $$ 
 In brief, via \eqref{345kolme}  and a proper choice for the additive constant in $\alpha(\zeta)$ 
 we have
  \begin{equation} \label{334harmonic}
  U(\zeta) = q_0 + \alpha(\zeta) - \alpha \left(1 / \, {\overline{\zeta}} \,\right) 
  + C \log|\zeta|, \qquad 1 \leq |\zeta| <  R,  
    \end{equation}
  with
   \begin{equation} \label{radialpoint}
   \nabla \sigma \circ U (\zeta) =
 - \alpha(\zeta) - \alpha \left(1 / \, {\overline{\zeta}} \,\right) + C \log|\zeta|.
 \end{equation}

 Next note that 
 $ \partial \sigma(q_0)$ consists of the limits of $ \nabla \sigma \circ U (\zeta)$ as $\vert\zeta\vert\to 1$. Hence
   the boundary of $\partial \sigma(q_0)$ has the parametrisation 
  \begin{equation} \label{unifgas} 
  (\nabla \sigma \circ U) (\eta) = - 2\, \alpha(\eta), \qquad |\eta| = 1,
    \end{equation}
where $\alpha(\zeta)$ is analytic in $\{ 1/R < |\zeta| < R\}$.  
Since $\partial \sigma(q_0)$ is convex with area  $|\partial \sigma(q_0)| = |\nabla \sigma(\{q_0\})| > 0$, the function $\alpha(\zeta)$ is one-to-one on the unit circle.
 \end{proof}

\begin{rem} \label{proofThm22}
Combining Lemmas \ref{subdiff} and \ref{subdiff12} together with Lemma \ref{blowup}
proves Theorem \ref{prop:sigma:property3}.
\end{rem}     

Further, with the approach of Lemmas \ref{subdiff} and \ref{subdiff12} we can give  for any point $p_0 \in \mathscr P\cup \mathscr Q\cup\mathscr G$   a natural geometric interpretation  and parametrisation of the boundary of the subdifferential $\partial \sigma(p_0)$, in terms of  the directional limits of $\nabla \sigma$ at $p_0$,  in  directions $p-p_0$ where $p \in N^\circ$.
We denote this limit by $\widehat{\nabla} \sigma({p_0}, p-p_0)$.

\smallskip

  \begin{Cor} \label{GenGradProp} Let $p_0$ be a point in the set $\mathscr P\cup \mathscr Q\cup\mathscr G$ and $\sigma$ a convex solution to \eqref{Pst2}. 

Then for every interior point $ p \in  N^\circ$, $p \neq p_0$,  the gradient $\nabla \sigma$ has  limit at $p_0$ in direction $p-p_0$,
\begin{align}\label{GeneralizedGradient}
\widehat{\nabla} \sigma({p_0}, p-p_0) :=  \lim_{\tau\to 0^+} \nabla \sigma\big( p_0+\tau(p-p_0)\big) \, \in \, \partial \sigma(p_0). 
\end{align}
Conversely, every point on the boundary of the subdifferential  $\partial \sigma(p_0)$ arises in this manner, as  
$\; \widehat{\nabla} \sigma({p_0}, p-p_0)\, $ with  $ p \in  N^\circ$, $p \neq p_0$.
\end{Cor}

\begin{proof}
 Let us first discuss the case $p_0 = p_j \in \mathscr P\cup \mathscr Q$.

\noindent $1^\circ$.   If there are no points $q \in {\mathscr G}$, we can apply directly Theorem \ref{harmonic rep} and Corollary \ref{sigmarep}. In this case, let $U : \Di \to N^\circ$ be the harmonic homeomorphism in \eqref{harmonichomeo31}. We observe from \eqref{harmonichomeo31} and \eqref{harmonic345} that outside the singularities, the endpoints $\{\eta_j\}$ of the intervals $\{ I_j\}$, the map extends smoothly across the unit circle with radial derivative 
  \begin{equation} \label{harmonic4}
   \partial_r U(\eta) = \frac{1}{\pi} \sum_{j=1}^m (p_j - p_{j+1}) \Im m \left( \frac{\eta}{\eta - \eta_j} \right).
\end{equation}
It follows that the direction of $  \partial_r U(\eta)$ approaches $p_j - p_{j+1}$ when we let $\eta \to \eta_{j}$
  on the arc $(\eta_{j-1}, \eta_j) \subset \partial \Di$. In fact, we have  the quantitative estimate 
    \begin{equation} \label{harmonic.limit}
    \frac{\partial_r U(\eta) }{|\partial_r U(\eta) |} = \frac{p_j-p_{j+1}}{|p_j-p_{j+1}|}\bigl(1 + {\mathcal O}(| \eta - \eta_{j}|)\bigr),
    \qquad \eta \in (\eta_{j-1},\eta_{j}) \subset  \partial \Di.
    \end{equation}

Similarly, the direction of $  \partial_r U(\eta)$ approaches $p_{j} - p_{j-1}$ when $\eta \to \eta_{j-1}$
  on the arc $(\eta_{j-1}, \eta_j)$.
Furthermore, by \eqref{radial}, outside the singularities  $\{\eta_j\}$  the derivative $\partial_r U(\eta)  \neq 0$. 
It hence follows that $ \partial_r U(\eta)$ attains on the arc $(\eta_{j-1},\eta_{j}) \subset  \partial \Di$ every direction between 
$p_j - p_{j+1}$ and $p_{j} - p_{j-1}$.

 On the other hand, $U(\eta)$ maps all of the arc $\, \eta \in (\eta_{j-1}, \eta_j) = I_j$ to the point $p_j$,  with
$$ U(r \eta) = p_j + \partial_r U(\eta) (r-1) +  o_\eta(1-r), \quad {\rm as } \; \; r \uparrow 1.
 $$
 Thus given a direction $p-p_j$  with $p \in N^\circ $,  we can find $\eta \in (\eta_{j-1}, \eta_j)$ such  that 
 $U(r \eta) - p_j $ gets the direction $p-p_j$ as $r \to 1$.  
 But then \eqref{summa4} shows that 
\begin{equation} \label{limits53}
  \lim_{\tau\to 0^+} \nabla \sigma\big( p_j+\tau(p-p_j)\big) =  
\lim_{r \to 1}  \nabla \sigma \circ U(r\eta)
\end{equation}
exists and is finite.
\smallskip

\noindent  $2^\circ$.  If  $\sigma$ has gas points $q \in {\mathscr G}$   and  $p_0 = p_j \in \mathscr P\cup \mathscr Q$,  we modify the argument similarly as in the proof of Lemma \ref{subdiff}. Thus choose  $\rho_1\in (p_{j-1},p_{j})$ and $\rho_2\in (p_j,p_{j+1})$,
 connect them with a smooth strictly convex arc  $\Gamma \subset {N}^\circ$ and call $M^\circ$ the domain bounded by
$[\rho_1,p_{j}] \cup [p_j, \rho_2] \cup \Gamma$. Again we can assume that $M^\circ$ contains no points $g \in {\mathscr G}$. 

Let also the harmonic homeomorphism $U_M: \Di \to M^\circ$ be as in \eqref{Udecomp}, with   $W_M$ the term vanishing on the interval $I  = (\eta_{1}, \eta_2)$, as in
Lemma \ref{subdiff}. Since $W_M(\eta) \to \rho_1$ as $\eta \to \eta_1$ on $I \subset \partial \Di$, estimating the radial derivative of the Poisson representation shows  that when  $\eta \to \eta_{1}, \, \eta \in I $, for any $\varepsilon > 0$
$$ \partial_r U_M(\eta) = \frac{\rho_1 - p_{j}}{\pi} \, \Im m \left( \frac{\eta}{\eta - \eta_1} \right) + \frac{b_{\varepsilon}(\eta)}{|\eta - \eta_1|} + C_\varepsilon, 
$$
 where $C_\varepsilon$ is a constant and $\|b_{\varepsilon}\|_{\infty} \leq \varepsilon$. 
 
 A similar estimate holds for $ \eta \in I $  near $ \eta_{2}$. In addition, from proof of Lemma  \ref{subdiff} we have $  \partial_r U_M(\eta) \neq 0$ for all $ \eta \in I $.
Thus as in the case $1^\circ$  one sees that $ \partial_r U_M(\eta)$ attains on the arc $I$  every direction between 
$p_j - \rho_1$ and $p_{j} - \rho_2$, that is, every direction between  $p_j - p_{j+1}$ and $p_{j} - p_{j-1}$.

In conclusion, again as in case $1^\circ$,  for any given direction $p-p_j$  with $p \in N^\circ $,  we can now find $\eta \in I$ such  that 
 $U_M(r \eta) - p_j $ gets the direction $p-p_j$ as $r \to 1$.  On the other hand, Lemma \ref{subdiff} provided us with the  radial limits of  $\nabla \sigma \circ U_M(r\eta)$, for every $\eta \in I$, and this  completes the proof of step $2^\circ$.

\medskip

\noindent  $3^\circ$. As in the proof of Lemma \ref{subdiff12}, if we have a gas point $q_0 \in N^\circ$, choose a small 
disc $\Di_0 = \Di(q_0,\varepsilon) \subset N^\circ$ and a Riemann map $\psi:\mathbb{A}(1,R) \to L_\sigma (\Di_0\setminus \{q_0\})$.
Composing then the  inverse $L_\sigma^{-1}$ with  $\psi$  gives us the harmonic homeomorphism  
  \begin{equation} \label{333three}
   U := L_\sigma^{-1} \circ \psi : \mathbb{A}(1,R) \to \Di_0\setminus \{q_0\}; \qquad U(\zeta) \to q_0  \quad {\rm as } \quad  |\zeta| \to 1.
  \end{equation}
Moreover, as in the proof of Lemma \ref{subdiff12} we get 
  \begin{equation} \label{334harmonic}
  U(\zeta) = q_0 + \alpha(\zeta) - \alpha \left(1 / \, {\overline{\zeta}} \,\right) 
  + C \log|\zeta|, \qquad 1 \leq |\zeta| <  R,  
    \end{equation}
  and
   \begin{equation} \label{radialpoint}
   \nabla \sigma \circ U (\zeta) =
 - \alpha(\zeta) - \alpha \left(1 / \, {\overline{\zeta}} \,\right) + C \log|\zeta|.
 \end{equation}
 Now  $U(\zeta) $ with its derivatives extends analytically up to the unit circle and  
 therefore as before, the existence of limits \eqref{GeneralizedGradient} is reduced to  showing that 
 $$ \lim_{r \to 1^+}  \nabla \sigma \circ U(r \eta), \qquad |\eta| = 1, 
 $$
 exist for each $\eta \in \partial \Di$. This, on the other hand, is a direct consequence of \eqref{radialpoint}.
 
  Finally, every boundary point  $\xi$ of $\partial \sigma(q_0)$ is given by some subseqential limit of $(\nabla \sigma \circ U)(\eta_j)$, where $\eta_j \to \partial \mathbb{D}$. In order for this sequence to be convergent, $\eta_j$ must converge to some $\eta \in \partial \mathbb{D}$ and then  $\xi$ arises as a suitable limiting value $\hat \nabla \sigma (q_0,p-q_0)$. A similar argument, with $\partial \mathbb{D}$ replaced by $I_j$, applies to points $p_j \in \mathscr{P} \cup \mathscr{Q}$ as well.  Consequently, it follows  that  for each $p_0 \in \mathscr P\cup \mathscr Q\cup\mathscr G$, every point in $\partial \sigma(p_0)$ arises as a limiting value of the  form \eqref{GeneralizedGradient}.
\end{proof}
  \smallskip

A further interesting geometric picture is added to Corollary \ref{GenGradProp}  by the results we showed in Proposition \ref{prop:sigma:property4}. Namely,  in addition to the existence of the limits \eqref{GeneralizedGradient}, 
for any $p \in N^\circ$ and any
$p_0 \in  \mathscr P\cup \mathscr Q \cup G$,  the vector  $p - p_0$ is (an outer)
normal to $\partial \sigma(p_0)$,  at the boundary point $\widehat{\nabla}\sigma({p_0},p-p_0) \in \partial \sigma(p_0)$. 

On the other hand, to prove in Section \ref{section:minimality}  the universality of the geometry of the Lozenges frozen boundaries  within general dimer models we need, c.f. \eqref{Uepsilon3},  one more subtle property of the surface tensions \eqref{Pst2}, stated as follows.

  \begin{Cor} \label{ortoasymptotes}
  Suppose $\sigma$ and $N$ are as in \eqref{Pst2}. Let  $p_j, p_{j+1} \in \mathscr P\cup \mathscr Q$ be neighbouring points, in the order induced by $\partial N$, and assume
  $\, \widehat p  \in (p_{j}, p_{j+1}) \subset \partial N$.
 Then  the following  limits exist and are equal,
 \begin{eqnarray} \label{tangentlimits}
\hspace{-.5cm} \lim_{ N^\circ\setminus {\mathscr G}\, \ni \; p \to \, \widehat p \, }\langle\widehat{\nabla}\sigma({p_j},p-p_j), p-p_j\rangle = 
\lim_{ N^\circ\setminus {\mathscr G}\, \ni \; p \to \, \widehat p \, } \langle\widehat{\nabla}\sigma({p_j},p-p_j),\, \widehat p - p_j\rangle =  \sigma(\,\widehat p\,) - \sigma(p_j).
\end{eqnarray}
\end{Cor}
  \begin{proof}
  When there are no gas points $q \in \mathscr  G$ we use Corollary \ref{sigmarep} and the representation \eqref{summa4}. In that notation,  let us write for brevity
  $$ (\delta_j\sigma )(\zeta) := \frac{1}{\pi} \sum_{k \neq j}  \,  i(p_{k} - p_{k+1}) \log \frac{1}{|\zeta-\eta_k|} + c_0, \qquad \zeta \in \Di.
  $$
  Then $ (\delta_j\sigma )(\zeta)$ has a finite limit as $\zeta \to \eta_j$ in $\overline{\Di}$, while the remaining main term of \eqref{summa4} is orthogonal to the side $[p_j, p_{j+1}]$.  Thus  
  for every $\widehat p \in (p_j, p_{j+1})$ there exist  the finite  limit 
  \begin{equation} \label{tangentlimits2}
   \lim_{\zeta \to \eta_j} \langle \widehat p-p_j  , \nabla \sigma \rangle \circ U(\zeta) = 
   \lim_{\zeta \to \eta_j} \langle (\delta_j\sigma )(\zeta), \widehat p-p_j  \rangle, 
  \end{equation}
  independently of the direction of $\zeta- \eta_j$. As in \eqref{der.tangential} we  see that 
  the above limits give the tangential derivative of $\sigma$ on $[p_j, \widehat p \, ]$. But  $\sigma$ being affine on $(p_j, p_{j+1})$, the above limits are in fact equal to $\sigma(\, \widehat p \,) - \sigma(p_j)$.
  
  With \eqref{limits53} this gives the second identity in \eqref{tangentlimits}. For the first, it is enough to 
  notice that by \eqref{harmonic.limit}, for any $p \in N^\circ$ we have $\widehat{\nabla}\sigma({p_j},p-p_j) = (\nabla \sigma \circ U)(\eta)$ where    $${\rm  dist}(p, [p_j,p_{j+1}]) \leq {\mathcal O} (|\eta - \eta_j|),$$ and that the singularity of $(\nabla \sigma \circ U)(\eta)$ at $\eta_j$ is logarithmic.
 If $\sigma$ has gas point, a similar argument works with $U$ replaced by the harmonic mapping $U_M$ in \eqref{Udecomp}, as in Lemma \ref{subdiff}, step $2^o$. 
   \end{proof}
  \medskip

    \subsubsection{Proof of Theorem \ref{sigma1}.} \label{sigmaregproofs}
    The above results complete the proof of Theorems \ref{sigma1}:    
   Lemma \ref{blowup} proves the first part of Theorem \ref{sigma1}. For the second, let $p_0 \in \mathscr{P} \cup \mathscr{Q} \cup \mathscr{G}$.  Then   Lemmas \ref{subdiff} and \ref{subdiff12}  imply that the boundary of $\partial \sigma(p_0)$ is an analytic curve, while Corollary \ref{GenGradProp} provides us with the limits  $\widehat{\nabla} \sigma({p_0}, p-p_0) :=  \lim_{\tau\to 0^+} \nabla \sigma\big( p_0+\tau(p-p_0)\big)$, and shows that they parametrise the boundary of $\partial \sigma(p_0)$.
   
    Further,  once these limits \eqref{GeneralizedGradient} exist, Proposition \ref{prop:sigma:property4} shows that for any $p\in N^\circ \setminus \Gg$ the vector 
    $p - p_0$ is normal to $\partial \sigma(p_0)$ at the boundary point 
$\widehat{\nabla}\sigma({p_0},p-p_0)$ of the subdifferential. Hence there is a one-to-one continuous correspondence between the limits $\widehat{\nabla}\sigma({p_0},p-p_0)$ and the directions  or arguments $\arg(p-p_0)$. This completes the proof of Theorem \ref{sigma1}. \hfill $\Box$

    \subsection{The analytic structure function  $\mathcal{H}_\sigma$}\label{Hsection}

    With all the results on the surface tension $\sigma$ covered in the previous Subsections, the basic properties of the structure function $ {\mathcal H} =  \mathcal{H}_\sigma$ and its domain $Dom(\mathcal H_\sigma)$  are now easily explained.
  Recall that for any surface tension $\sigma$ from \eqref{Pst2},
   \begin{equation} \label{structf}
{\mathcal H}_\sigma (w) :=  (I-\nabla \sigma) \circ (I + \nabla \sigma)^{-1}(\overline w), \qquad  {\nabla \sigma}(0) = 0,
  \end{equation}
 where the last   condition is achieved  by adding, if necessary,  a suitable linear function to $\sigma$, and where 
   \begin{equation} \label{domH2}
  \text{Dom}({\mathcal H}_\sigma) = L_\sigma(N^\circ \setminus  \mathscr{G}),
 \end{equation}
 in terms of the Lewy transform \eqref{LewyT}.
 
 According to Theorem \ref{prop:sigma:property3}, the boundary of the amoeba $ \nabla \sigma\left(  N^\circ\setminus {\mathscr G} \right)$ is the union 
     \begin{equation} \label{amoeba3}
     \bigcup_{p_0 \in    \mathscr P \cup \mathscr Q \cup \mathscr G} \; \gamma_{p_0},
 \end{equation}
     i.e. the union of the boundary curves of the corresponding subdifferentials $\partial \sigma(p_0)$, for  $p_0$ either a corner, a quasifrozen point or a gas point in $N$. 
     
     The boundary curves of   $\text{Dom}({\mathcal H_\sigma})$ are each a similarity image  of a curve in \eqref{amoeba3}, obtained by a translation and a reflection (complex conjugation) of a curve in that collection. In brief, we set
          \begin{equation} \label{amoeba2}
     \widehat{\,\gamma_{p_0}} := \overline{\,{p_0}\,} + \overline{\, \gamma_{p_0} \,},\qquad {\rm for \; each\; } \; p_0 \in    \mathscr P \cup \mathscr Q \cup \mathscr G,
      \end{equation}
      see also Figure \ref{fig.amoeba2}.
      
     \begin{Cor}  \label{subdiff2}  Suppose the surface tension $\sigma$ and its domain $N$ are as in 
     \eqref{Pst2}.
      Then  the boundary of the domain of the structure function $Dom(\mathcal H_\sigma)$  is the disjoint union of analytic and convex curves,
         \begin{equation} \label{amoeba34}
\partial  \,{\rm{Dom}}({\mathcal H}_\sigma) =    \bigcup_{p_0 \in    \mathscr P \cup \mathscr Q \cup \mathscr G} \; \widehat{ \, \gamma_{p_0} }.
 \end{equation}
    \end{Cor}
  \begin{proof} 
     The Lewy transform $L_\sigma$ is the gradient of the convex function $\widehat{\sigma}(z) = \frac{1}{2}|z|^2 + \sigma(z)$ composed with the complex conjugation. Hence the boundary of $L_\sigma(N^\circ \setminus  \mathscr{G})$ consists of the complex conjugates of the boundaries of the subdifferentials $\partial \widehat{\sigma}(p)$ with $p \in \partial N \cup  \mathscr{G}$. Clearly, 
    \begin{equation} \label{domH3}
 \partial  \widehat{\sigma}(p_0) =  \{p_0\} + \partial \sigma(p_0), \qquad p_0 \in  \mathscr{P} \cup  \mathscr{Q}  \cup  \mathscr{G},
  \end{equation}
 while for all other $p \in \partial N$ we have $\partial   \widehat{\sigma}(p) = \emptyset$. 
   \end{proof}

 \medskip
  
   In view of, for instance, Theorem \ref{connection} it is also important to understand the mapping properties of the derivative of the structure function $\mathcal H_\sigma': Dom(\mathcal H_\sigma) \to \Di$.

 \begin{Prop} \label{Hproper} For a surface tension $\sigma$ from \eqref{Pst2}, the structure function  $ \mathcal{H}_\sigma$ in \eqref{structf} extends analytically across any bounded and connected part of the boundary of $\,\rm{Dom}( \mathcal{H}_\sigma)$. Furthermore, 
 
 i) The derivative $ \mathcal{H}_\sigma' : \rm{Dom}( \mathcal{H}_\sigma) \to \Di$ is an analytic and  proper  map.
 
 ii) The degree $\,{\rm{deg}}(\mathcal{H}'_\sigma) = m-2 + 2\ell$, where  $m = \#(\mathscr{P} \cup \mathscr{Q})$ and $\ell = \# \mathscr{G}$. 
 \end{Prop}
\begin{proof} First, let us choose a uniformisation for the domain $Dom(\mathcal H_\sigma) = L_\sigma(N^\circ \setminus \mathscr G)$. However, if the surface tension $\sigma$ has gas points, the domain  is not simply connected. Hence in this case  we need to uniformize by a circle domain and use the classical Koebe uniformisation theorem \cite{Koebe}, which tells that  every finitely connected planar domain is conformally equivalent to a {\it circle domain} $\DD \subset \C$, i.e. a domain with all boundary components circles or isolated points.

 It is no restriction to assume  that  we have a conformal uniformisation 
     \begin{eqnarray} \label{circled1}
 \psi: \mathcal{D} \to Dom(\mathcal H_\sigma), \qquad     \mathcal{D} = \Di \setminus \bigcup_{k=1}^\ell D(z_k, \delta_k),
     \end{eqnarray} 
where the  subdiscs $D(z_k, \delta_k) \subset \Di$ are pairwise disjoint and have radii $0 < \delta_k < 1-|z_k|$.  Further, we may assume that $\psi$ is unbounded on $\partial \Di$, so that   the inner circles of $ \mathcal{D}$ correspond to gas components and the uniformization  extends to the circles $S_k := \partial D(z_k,\delta_k) $ as a homeomorphism, taking $S_k$ to the analytic curve $    \widehat{\,\gamma_{q_k}} $, where  $q_k \in \Gg$, $k = 1, \dots, \ell$.

     Next, for each $p_0 \in   \mathscr{P} \cup \mathscr{Q} \cup \mathscr{G}$, let $\widehat{\, \gamma_{p_0}}$ be the analytic arc or curve defined in \eqref{amoeba2}, a connected component of the boundary of $Dom(\mathcal H_\sigma) $,  and finally, let 
        \begin{eqnarray} \label{arcs23}
         I_{p_0} = \psi^{-1}( \widehat{\, \gamma_{p_0}} )\subset \partial \Di, \qquad p_0 \in   \mathscr{P} \bigcup \mathscr{Q},  
          \end{eqnarray} 
be the preimages of the unbounded boundary arcs of $Dom(\mathcal H_\sigma) $ under the conformal map $\psi$, arising from the corners and quasi frozen points of $\partial N$. Clearly the arcs $I_{p_0}$ have disjoint interiors with closures covering the whole circle $\partial \Di$.

With the above notation, we again make use of the  harmonic homeomorphism $\, U :=  L_\sigma^{-1} \circ \psi \, $ and, in particular, the identity 
 \begin{eqnarray} \label{HooJaL}
  {\mathcal H}_\sigma\Bigl( \psi(\zeta) \Bigr) = 2 U(\zeta) - \overline {\psi(\zeta)}, \qquad \zeta \in \mathcal{D},
 \end{eqnarray} 
 which follows from \eqref{Lewyinverse}. 

If  $p_0 \in   \mathscr{P} \cup \mathscr{Q}$,  the boundary arc  $\widehat{\, \gamma_{p_0}}$ is analytic, so that
 the uniformisation $\psi: \mathcal{D} \to Dom(\mathcal H_\sigma)$ extends analytically across the arc $ I_{p_0} = \psi^{-1}(\widehat{\, \gamma_{p_0}})\subset \partial \Di$.  
Thus letting  
 $\zeta \to I_{p_0}$ in   \eqref{HooJaL} gives 
 \begin{eqnarray} \label{HooJaL3}
  {\mathcal H}_\sigma \Bigl( \psi(\eta) \Bigr) = 2 p_0 -  \overline {\psi(\eta)}  \qquad \eta \in I_{p_0}  \subset \partial \Di.
 \end{eqnarray} 
As $\overline{\psi(\eta)} = 1/\psi(\eta)$ on $I_{p_0}$ and $\psi(\eta)$  extends analytically across
$I_{p_0}$, by the above identity  $ {\mathcal H}_\sigma$ extends analytically across $\widehat{\, \gamma_{p_0}}$. A similar reasoning shows that 
$$ {\mathcal H}_\sigma \Bigl( \psi(\zeta) \Bigr) = 2 q_k -  \overline {\psi(\zeta)}, \quad {\rm for } \;  \;  \eta \in S_k,
$$
allowing
$ {\mathcal H}_\sigma$ to extend analytically also across the inner circles
$S_k$. Thus the first claim follows.

For the properness of $\mathcal{H}'_\sigma : \text{Dom}( \mathcal{H}_\sigma) \to \Di$, let   $p_0 \in   \mathscr{P} \cup \mathscr{Q}$.  Differentiating 
 \eqref{HooJaL3}  along the arc $ I_{p_0}$  shows that
  \begin{eqnarray} \label{Htangent}
   {\mathcal H}_\sigma' \Bigl( {\psi(\eta)} \Bigr) =  \,  \frac{\,1\,}{ \eta^2 \,}\frac{|\psi'(\eta)|^2}{ \psi'(\eta)^2}, \qquad  \eta \in I_{p_0} \subset \partial \Di,
  \end{eqnarray} 
  where we note that $i \eta \,\psi'(\eta)$ gives  the tangent to $\widehat{\, \gamma_{p_0}}$ at $\psi(\eta)$. An identity analogous  to \eqref{Htangent} holds on the circles $S_k$ as well, with $|\mathcal {H}_\sigma'|  \equiv 1$  on $S_k$.

  Finally, by Lemma \ref{analyticH} the derivative $\mathcal{H}'_\sigma$ is analytic with values in $\Di$, it extends analytically across the finite parts of the boundary of   $\text{Dom}( \mathcal{H}_\sigma)$, and has unit modulus there.  At  the (finitely many) end points of the  arcs $I_{p}  \subset \partial \Di$, for $p \in   \mathscr{P} \cup \mathscr{Q}$, the uniformization $\psi(\zeta) \to \infty$, but  Lemma \ref{subdiff} with \eqref{riemann} and \eqref{HooJaL} show that even then $\mathcal{H}'_\sigma$ has a limit of unit modulus.
  Therefore it follows that $\mathcal{H}'_\sigma: \text{Dom}( \mathcal{H}_\sigma) \to \Di$ is proper.

  To determine the degree of $\mathcal{H}'_\sigma$, note that by Lemma \ref{subdiff},  the change of the tangential argument along 
  $\widehat{\, \gamma_{p_0}}$ is $\alpha(p_0)$, the angle of the polygon $N$ at the corner $p_0$; for $p_0$ quasifrozen $\alpha(p_0) = \pi$. Thus by \eqref{Htangent}, along $\partial \Di$ the total change in argument for $\mathcal{H}_\sigma'$  is $2(m-2)\pi$. Similarly by \eqref{Htangent}, the total change 
  for $\mathcal{H}'_\sigma$ around each ``gas circle'' $S_k$ is $4\pi$. These  together prove the claim ii). 
\end{proof}
\smallskip

Recall also that in \eqref{beltB1} we saw that  it is often convenient to compose the structure function $\mathcal{H}_\sigma$ with a suitable conformal uniformization. In the case where there are no gas points, we can  uniformize with the unit disc and choose a conformal $\psi: \Di \to  \text{Dom}( \mathcal{H}_\sigma)$, which took us to the Beltrami equation
$$
   f_{\overline{ z }}(z) = \mu_{\sigma} \bigl(f(z)\bigr)\, f_z(z),\qquad \quad  {\rm where} \quad  \mu_\sigma \equiv {\mathcal H_\sigma'} \circ \psi  \quad {\rm and} \quad f:  \mathcal U \to \Di, 
$$
an equation equivalent to the original version \eqref{Belt654}. 
Here by Proposition \ref{Hproper} above, $\mu_\sigma: \Di \to \Di$ is an analytic and proper map, hence a Blaschke
product. Thus in this picture the coefficient function $\mu_\sigma$ often allows an explicit expression:
\medskip

\begin{rem}  \label{domH}  For the basic lozenges model, $N$ is a triangle without gas or quasifrozen points, so that by Proposition \ref{Hproper} the degree ${\rm deg}(\mathcal{H}') = 1$ and therefore $\, \mathcal{H}_\sigma':\text{Dom}( \mathcal{H}_\sigma) \to \Di$ is a conformal map;  thus we can simply take $\mu_\sigma(z) \equiv z$.

On the other hand, for the uniform  domino tilings, c.f. Example \ref{microdominoes}, the polygon  $N$ is square with corners $\{\pm 1, \pm i \}$, without gas points, and $\sigma$ has zero boundary values. Thus $\sigma(iz) = \sigma(z)$, by uniqueness of the convex solutions to \eqref{Pst2}. Then  $\nabla \sigma(iz) = i \nabla \sigma(z)$, while
\eqref{structf} gives $\mathcal{H}_\sigma(iz) = -i \mathcal{H}_\sigma(z)$. Also, $\mathcal{H}_\sigma'$ has degree $2$ by Proposition \ref{Hproper}. Hence, with a uniformisation fixing the origin, we have $\, \mu_\sigma(z) = z^2$ for the dominoes.
\end{rem} 
\smallskip

With  Remark \ref{domH} we can now tie together several themes of this Section: Given  a  liquid domain $\LL$ with frozen boundary $\partial \LL$ in {\it some }dimer model, with surface tension $\sigma$, let  $h: \LL \to \R$ be the corresponding limiting height function. Then by Theorem \ref{key.connection} below, 
  \begin{equation} \label{unisol} 
   f(z) := \mathcal{H}'_\sigma \circ L_\sigma \circ \nabla h(z), \qquad z \in \LL,
  \end{equation}
is a proper map $f:\LL \to \Di$ which solves the universal Beltrami equation
$$ f_{\overline z} (z) = f(z) \, f_z(z), \qquad z \in \LL.
$$
Via Remark \ref{domH}, we then can view this as  the Beltrami equation \eqref{beltB1} associated to the Lozenges model. 

Indeed,  let $N_0$ be the gradient constraint of the Lozenges model,  a triangle with corners $p_1, p_2,p_3$, and   $\sigma_0$ the surface tension from \eqref{Pst2}, with ${\sigma_0}_{\big | \partial N_0} = 0$. In case we choose the corners to be $0$, $1$ and $i$, then $\sigma_0$ would be given by \eqref{gradsigma23}. We can then  apply the representation \eqref{harmonichomeo31}  from Theorem \ref{harmonic rep} in the case of the Lozenges, and use Corollary \ref{Dconverse} to see that  the identity
  \begin{equation} \label{Lozerep} 
  \nabla h_0(z) := \sum_{j=1}^3 \, p_j \, \omega_{\Di}\bigl( f(z); I_j \bigr), \qquad  z \in \LL,
\end{equation}
  defines a solution to the corresponding Euler-Lagrange equation  in Lozenges model,
    \begin{equation} \label{double} 
 {\rm div} \, \big(\nabla \sigma_0(\nabla h_0)\big)= 0 \qquad \mbox{ for all }
    \; z \in  \LL.
  \end{equation}
In addition, since by Theorem \ref{key.connection} the map $f:\LL \to \Di$ is a proper map, so is $\nabla h_0 : \LL \to N_0$. 

In other words, the original domain $\LL$ we started with   is a candidate also in the lozenges model, for a liquid domain with frozen boundary. However, to complete the picture, we must show that $\LL$ and its frozen boundary do arise from the minimization problem \eqref{ELbasic} associated with the Lozenges surface tension $\sigma_0$, and also  from some polygonal domain $\Om \supset \LL$ and from  some admissible boundary values on $\partial \Omega$.
   
  For this, note that the representation \eqref{Lozerep} and properness of $\nabla h_0$ allow an extension
  for $h_0$ as a piecewise affine function, with $\nabla h_0(\Omega \setminus \LL) \subset \{ p_1, p_2, p_3\}$,  to  a natural domain $\Omega \supset \LL$. It is not difficult to construct the extension so that it has  natural boundary values  on $\partial \Omega$, in the sense of Definition \ref{def:extremal}.

  One would then like to show that this extended function, still denoted by $h_0$, is the unique minimizer for the variational problem \eqref{ELbasic} for its boundary values in the bigger domain $\Omega$. 
  This, however, is a non-trivial issue. We will  discus it  in detail later in Section  \ref{section:minimality}, see Theorem \ref{Nminimizer}.

\subsection{Representing height functions for general domains and surface tensions}  \label{repsection}

  With a more detailed picture of the domain $L_\sigma(N^\circ \setminus \mathscr G) = {\rm Dom} ({\mathcal H}) $ now at our disposal, let us return to the theme of Theorem \ref{harmonic rep},  but now for $\sigma$  a solution to the general Monge-Amp\`ere equation  \eqref{Pst2} with gas points. To simplify the notation, given a curve ${\mathcal J }$ on the boundary of $ \,{\rm{Dom}}({\mathcal H}_\sigma)$  we write briefly
    \begin{eqnarray} \label{Hdom}
    \omega_{\mathcal H}(\zeta; {\mathcal J }) :=  \omega_\Omega(\zeta; {\mathcal J }), \qquad {\rm where } \quad \Omega = {\rm Dom} ({\mathcal H}_\sigma),
    \end{eqnarray} 
  for the harmonic measure of ${\mathcal J }$ in the domain $\Omega = {\rm Dom} ({\mathcal H}_\sigma)$.
  
  \begin{Prop} \label{prop.hrep}
    Suppose $\sigma$ is a solution to   \eqref{Pst2} in a convex polygon $N$, with corners and quasifrozen points 
    $\mathscr{P} \cup \mathscr{Q} =  \{ p_j\}_1^m $ and gas points $\mathscr{G} =  \{ q_k\}_1^\ell$. 
    For each such point, let also $ \widehat{ \, \gamma_{p_j}}$ and $ \widehat{ \, \gamma_{q_k}}$ be the corresponding curves on the boundary $  \partial  \,{\rm{Dom}}({\mathcal H}_\sigma)$, as in \eqref{amoeba34}.

    Then the Lewy transform
    $L_\sigma: N^\circ \setminus \mathscr G\to  {\rm Dom} ({\mathcal H})$ has  harmonic inverse, with the representation
        \begin{eqnarray} \label{hrep.general}
L_\sigma^{-1}(z) =  \sum_{j=1}^m \, p_j \, \omega_{\mathcal H}\bigl(z; \,\widehat{ \gamma_{p_j}}\bigr) 
+  \sum_{k=1}^\ell q_k \, \omega_{\mathcal H}\bigl(z; \,\widehat{ \gamma_{q_k}}\bigr), \qquad z \in {\rm Dom} ({\mathcal H}).     \end{eqnarray} 
      \end{Prop}

      \begin{proof}        We know from \eqref{Lewyinverse} that $L_\sigma^{-1}$ is harmonic, and from Corollary \ref{subdiff2} that the boundary of ${\rm Dom} ({\mathcal H})$ consists of the unbounded analytic and convex Jordan arcs $ \,\widehat{ \gamma_{p_j} } $ together with the bounded analytic and convex Jordan curves $ \widehat{ \, \gamma_{q_k} }$,   one arc 
      for each $p_j \in \mathscr{P} \cup \mathscr{Q}$, $j = 1,\dots, m,$ as well as  one curve for each gas point $q_k$, $k = 1,\dots, \ell$.
      
Further, as shown in Corollary \ref{subdiff2}, $L_\sigma(z) \to \widehat{ \gamma_{p_j} } $ as $z \to p_j $, and respectively $L_\sigma(z) \to
\widehat{ \, \gamma_{q_k} } $, when $z \to q_k  $. 
But this means that 
 $ L_\sigma^{-1}(z)$ and
             \begin{eqnarray} \label{hrep.general5}
       z \mapsto  
\sum_{j=1}^m \, p_j \, \omega_{\mathcal H}\bigl(z; \,\widehat{ \gamma_{p_j}}\bigr) 
+  \sum_{k=1}^\ell \, q_k \, \omega_{\mathcal H}\bigl(z; \,\widehat{ \gamma_{q_k}}\bigr), \qquad z \in {\rm Dom} ({\mathcal H}),
        \end{eqnarray} 
            are bounded harmonic mappings in ${\rm Dom} ({\mathcal H}_\sigma)$ with the same boundary values. Thus the functions must be equal. 
       \end{proof}  
    
    As in the case without gas points, the identity \eqref{hrep.general} on the Lewy transform leads to a representation of solutions to Euler-Lagrange equation \eqref{eq:EL34}. Indeed, combining Theorem \ref{Belt654} and Proposition \ref{prop.hrep}  gives
        
  \begin{Cor} \label{ELrep.gen}
  Let $\sigma$ be a surface tension as in   \eqref{Pst2}, with corners $\mathscr P$, quasifrozen points $\mathscr Q$
  and gas points $\mathscr G$. 
  
   Assume $u \in C^1({\mathcal U})$ is a solution to the Euler-Lagrange equation  
 \begin{align} \label{nabla234}
  \, {\rm div} \, \big(\nabla \sigma(\nabla u)\big)= 0 \qquad {\rm in} \quad {\mathcal U},
   \end{align}
  in a domain ${\mathcal U} \subset \C$. 
  Then $\nabla u$ has the representation
 \begin{align} \label{nablarep}
\nabla u(z) = \sum_{j=1}^m \, p_j \, \omega_{\mathcal H}\bigl(f(z); \,\widehat{ \gamma_{p_j}}\bigr) 
+  \sum_{k=1}^\ell \, q_k \, \omega_{\mathcal H}\bigl(f(z); \,\widehat{ \gamma_{q_k}}\bigr), \qquad z \in \mathcal U,
 \end{align}
   where $f: \mathcal U \to \text{Dom}(\mathcal{H}_\sigma)$ is a $C^1$-solution  to the Beltrami equation $$  f_{\overline{ z }}(z) = {\mathcal H}'_{\sigma} \bigl(f(z)\bigr)\, f_z(z),$$ and $\, \omega_{\mathcal H}$ is the harmonic measure on $\text{Dom}(\mathcal{H}_\sigma)$.
    \end{Cor}
    
    \begin{proof} For any $C^1$-solution $u$ as in \eqref{nabla234},  \eqref{Belt654} tells that the condition $f = L_\sigma \circ \nabla u$ defines 
    a solution $f$ to the equation $f_{\overline{ z }} = {\mathcal H}'_{\sigma} (f) \, f_z$. The representation \eqref{nablarep} then follows from Proposition \ref{prop.hrep}. 
    \end{proof}
    
  \begin{rem}  \label{myygen}
Above one can of course also take  conformal parametrization  $ \psi: \mathcal{D} \to Dom(\mathcal H_\sigma)$ by a circle domain and argue as in Subsection \ref{universal}, to have the (equivalent) formulation of 
  Corollary \ref{ELrep.gen} representing solutions to \eqref{nabla234} in terms of  the harmonic measure on 
  $\mathcal{D}$,
  \begin{equation} \label{hrepre3456}
\nabla h(z)=\sum_{j=1}^m \, p_j \, \omega_{\mathcal{D}}\bigl(\, f(z); I_j\bigr) +  \sum_{k=1}^\ell \, q_k \,  \omega_{\mathcal{D}}\bigl(\, f(z); S_k\bigr), \qquad z \in \LL.
\end{equation} 
where $f: \mathcal U \to  \mathcal{D}$ is a solution to the Beltrami equation 
   \begin{align} \label{multimyy}
   f_{\overline{ z }} = \mu_{\sigma} (f) \, f_z, 
    \end{align}
 and  $\mu_\sigma := {\mathcal H}'_{\sigma}  \circ \psi$.
\end{rem}

\addtocontents{toc}{\vspace{-4pt}}
\section{Non-linear Beltrami Equation and the Hodograph Transform} \label{B+hodo}

We saw in the previous Section that 
 the minimizers of the integral \eqref{ELbasic} and the limiting height functions of different dimer models, in particular, are tightly connected to the solutions of specific Beltrami equations.  
The following result sums these relations from the point of view of liquid domains with frozen boundaries, c.f. Definition \ref{LiqFrozen}.

   \begin{Thm} \label{key.connection}
 Suppose $\sigma(z)$ is a general surface tension as in \eqref{Pst2}, and
  $h$  is a $C^1$-solution to the Euler-Lagrange equation 
\, $  \rm{div} \big(\nabla \sigma(\nabla h)\big)= 0$\,  in a domain $\LL\subset \C$. 

 Then  the identity 
$ f(z) = \mathcal H'_\sigma \bigl( L_\sigma \circ \nabla h(z)\bigr) \,$ defines a solution to the universal Beltrami equation 
\begin{equation} \label{eqNonLinB}
\partial_{\overline{ z }} f(z) = f(z)  \partial_z f(z),
\end{equation}
where furthermore,  $\partial \LL$ is frozen for $h$ if and only if  $f: \LL \to \Di$ is proper map.
 \end{Thm}
\begin{proof}  That $ f(z) = \mathcal H'_\sigma \bigl( L_\sigma \circ \nabla h(z)\bigr) \,$ is a solution to \eqref{eqNonLinB} follows from
Theorem \ref{connection} and from \eqref{Specific}-\eqref{belt2}.

On the other hand, the Lewy transform is a homeomorpism $L_\sigma: N^\circ \setminus \Gg \to \text{Dom}(\mathcal{H}_\sigma)$, while Proposition \ref{Hproper} tells that $ \mathcal H'_\sigma:  \text{Dom}(\mathcal{H}_\sigma) \to \Di$ is a holomorphic proper map.
Hence the proof follows.
\end{proof}

In fact, in view of Theorems \ref{Second.thm}, \ref{converse} and \ref{key.connection}, the question of finding domains $\LL$ that have a frozen boundary in some dimer model is basically equivalent, at least in the simply connected case, to describing solutions to \eqref{eqNonLinB}  that are proper maps $f:\LL \to \Di$.
\smallskip

\begin{rem} \label{locfrozen}
The characterization of the Theorem \ref{key.connection} works also for locally frozen boundaries:  If $h$ and $f$ are as in the Theorem, then a subset $\Gamma \subset \partial \LL$ is frozen for $h$, in the sense of \eqref{frozen.def},  if and only if 
   \begin{eqnarray} \label{partProper}
    | f (z)| \to 1 \quad {\rm whenever} \quad z \to \Gamma, \;\; z \in \LL.
 \end{eqnarray}
 \end{rem}
 
Thus a key point for the goals of this work is  to  develop and understand the specific features   of   solutions $f:\UU \to \Di$ to  \eqref{eqNonLinB}. In the next Section we then apply them to study the geometry of frozen boundaries. To allow gas points in the surface tensions, we must allow multiply connected domains $\UU$, while for the solution to define a proper map, the domain $\UU$ must necessarily be finitely connected,  c.f. Proposition \ref{prop:FiniteConn}.
Therefore, unless stated otherwise it will always be assumed that \quad  
\medskip

$1^o$ \; The domain $\UU\subset \C$ is bounded and finitely connected.
\medskip

$2^o$ \; The map $f \in W^{1,2}_{loc}(\UU)$  is  continuous and takes values in the open unit disc, $f(\UU) \subset \Di$.
\medskip

Since the equation \eqref{eqNonLinB} is quasilinear, we apply the method of {\it hodograph transformation} for it, and write  \eqref{eqNonLinB}
in the hodograph plane.

 The   hodograph method is most straightforward for simply connected domains $\LL \subset \C$ and solutions to  \eqref{eqNonLinB} that define a proper map
$f:\LL \to \Di$, see   Theorem \ref{propStoilow35}.
In the general case the  method works as follows. Here recall that by  Koebe's  theorem \cite{Koebe}  every finitely connected planar domain is conformally equivalent to a {\it circle domain} $\DD \subset \C$, i.e. a domain with all boundary components circles or isolated points. 

\begin{Thm}
\label{propStoilow}
Every non-constant and continuous solution $f:\UU\to \Di$ to (\ref{eqNonLinB}) is real analytic, and  admits a Stoilow factorization 
\begin{align} \label{factor}
f= b \circ g^{-1},
\end{align}
where  $b(z)$ is analytic in 
$\DD$ with  $\Vert b \Vert_\infty\leq 1$, $\DD$ is a circle domain 
 and where $g:\DD\to \UU$ is a homeomorphism of finite distortion  solving the linear Beltrami equation
\begin{align}
\label{eqLinB3}
\dbar g=- \, b(z) \overline{\, \dv_z g \,}, \qquad z \in \DD.
\end{align}
Moreover, if $\, \UU$ is simply connected we can take $\DD = \Di$, the open unit disc,  and if  also   the map $f:\UU\to \Di$ is proper,  then $b(z) = B(z)$, a finite Blaschke product.
\end{Thm}
\begin{proof}
Since $f(z)$ is continuous with values in the open unit disc it is locally bounded away from the unit circle, i.e. for each relatively compact subdomain $V \subset \UU$ we have 
\begin{align}
\label{locEllip}
 |f(z)| \leq k_V < 1, \qquad z \in V. 
\end{align}
This gives a complex structure to the domain $\UU$, call it ${\mathcal A}$,  by requiring analytic charts to have the form $(\varphi, V)$, where $V$ is open with $\overline{V} \subset \UU$  and  
\begin{align}
\label{charts}
\dbar \varphi (z) = f(z) \,  \dv_z \varphi (z), \qquad z \in V.
\end{align}
Indeed, due to uniform ellipticity within $V$, when domains of two charts intersect the function $\varphi_1 \circ \varphi_2^{-1}$ is analytic, c.f. \cite[p. 179]{ATG}.

As $\UU$ is assumed to be finitely connected, by  Koebe's  uniformisation theorem there is  a conformal homeomorphism $ G: (\UU, {\mathcal A}) \to \DD$ to a  circle domain, equipped with the standard complex structure. Thus 
\begin{align}
\label{Belt}
\dbar G(z) = f(z) \,  \dv_z G(z), \qquad z \in \UU.
\end{align} 
 In particular, the composition $ b:= f \circ G^{-1} : G(\UU) \to \Di$
 is analytic (in the standard sense). 
As  $|f(z)|$ bounded away from $1$ locally in $\UU$, the inverse $g := G^{-1}: \DD \to \UU$ satisfies the equation \eqref{eqLinB3}, see  \cite[p.34]{ATG}. Also, by definition, 
\begin{align}
\label{realanal3}
 f(z) = b \circ g^{-1}(z), \quad
  {\rm  with  \; \; }\| b\|_\infty \leq 1. 
  \end{align} 

Concerning the smoothness of $f$,  our initial assumptions $1^\circ - 2^\circ\, $  
 make it locally quasiregular,
hence it is locally H\"older continuous and by a bootstrap argument, \eqref{eqNonLinB} makes the function $C^{\infty}$-smooth.
 In fact, by  \eqref{eqLinB3} the homeomorphism $g$ is real analytic with non-vanishing derivative, as  its  inverse $G$ is now $C^{\infty}$-smooth. Hence  $G$, and by \eqref{realanal3}  also $f$, are real analytic.
 
Finally, if  $\UU$ is simply connected then (up to a M\"obius transformation) the circle domain is either  $ \DD =  \Di$ or $ \DD =  \C$. As  $b(z)$ is bounded and analytic in $\DD$,   Liouville's theorem shows  the second  case  is not possible. 
Consequently,  if for a simply connected $\UU$ the map $f:\UU \to \Di$ is proper, then in \eqref{realanal3}  $b:\Di \to \Di$ is an analytic and proper map, hence by Fatou's theorem a finite Blaschke product.  
\end{proof}

A  similar   factorization holds also for the  Beltrami equation \eqref{Belt654} associated to a specific surface tension $\sigma$, with structure function ${\mathcal H}_\sigma$. Here and in the sequel, $H^{\infty}(\DD)$ denotes the space of functions bounded and analytic on the domain $\DD$. 

\begin{Cor}\label{cor:propStoilow} Suppose $\UU \subset \C$ is a finitely connected domain. 
Then every
 non-constant and continuous $W^{1,2}_{loc}$-solution $f:\UU\to Dom({\mathcal H}_\sigma)$ to 
\begin{align} \label{multipartic}
  f_{\overline{ z }} = {\mathcal H}_\sigma'(f) f_z,
\end{align}
admits a  factorization 
\begin{align} \label{factor2}
f=\eta\circ g^{-1},
\end{align}
where  $\eta: \mathcal{D} \to Dom({\mathcal H}_\sigma)$ is analytic,  $\mathcal{D}$ is a circle domain 
and where $g:\DD\to \UU$ is a homeomorphism of finite distortion which solves the linear Beltrami equation
\begin{align}
\label{eqLinB3b}
\hspace{1cm} \dbar g=- \, b(z) \overline{\, \dv_z g \,}, \qquad b: ={\mathcal H}_\sigma' \circ \eta  \in H^{\infty}(\DD) \; \; \, {\rm with} \; \;  \Vert b \Vert_\infty\leq 1.
\end{align}
\end{Cor}
\begin{proof} We argue as in Theorem  \ref{propStoilow}, but this time to find a homeomorphic solution $G: \UU \to \DD$ to $\dbar G(z) = ({\mathcal H}_\sigma' \circ  f)(z) \,  \dv_z G(z)$. The classical  Stoilow factorisation \cite[p. 179]{ATG} then gives $f = \eta \circ G$ with $\eta$ analytic, and we see that $g = G^{-1}$ satisfies \eqref{eqLinB3b}.
\end{proof}

\subsection{Teleomorphic maps}

 The factorization \eqref{factor} linearizes the Beltrami equation \eqref{eqNonLinB} and leads one to a study of the equation 
  \eqref{eqLinB3}. Given a domain  $\DD\subset\C$ let 
\begin{align}
\mathcal{B}_{H^\infty}(\DD):=\{b\in H^{\infty}(\DD):\sup_{z\in \DD} \vert b(z)\vert_\infty\leq 1\}.
\end{align}

\begin{Def}
A function $g\in W^{1,2}_{loc}(\DD)$, not necessarily  injective,  which solves the equation 
\begin{align}
\label{eqLinB333}
\dbar g=- \, b(z) \overline{\, \dv_z g \,}, \qquad z \in \DD,
\end{align}
for some $b\in \mathcal{B}_{H^\infty}(\DD)$ is called a \emph{teleomorphic} function on $\DD$. The set of such functions is denoted by $\TT(\DD)$. 
\end{Def}
Note that the set of teleomorphic maps on $\DD$ contains the holomorphic maps $\mathcal{O}(\DD)$, corresponding to the case when $b\equiv 0$ on $\DD$. 

\begin{rem}
The word teleomorphic is taken from mycology.  More precisely, the different life cycles of fungi are called teleomorph and anamorph respectively, where teleomorph is a sexually reproductive life stage, and anamorph is an asexual reproductive life stage. In addition, the whole life cycle is called holomorph. Since we will show in Proposition \ref{LemMapHolo} below that teleomorphic functions are generated by pairs of holomorphic functions, we think the term is appropriate. 
\end{rem}

Teleomorphic maps bear a certain resembles to sense preserving  harmonic maps, as those solve a Beltrami equation of the form
\begin{align*}
\overline{  \dbar u(z) } = \omega(z) \dv_z u(z),
\end{align*}
where $\omega\in H^{\infty}(\DD)$. However, in many respects  the properties of these two classes of functions are quite different.

\begin{Prop}
\label{LemMapHolo}
Assume that $g\in \TT(\DD)$ with the coefficient $b\in \mathcal{B}_{H^\infty}(\DD)$. Then the function
\begin{align}
\label{lemMapTeleoHolo1}
\gamma :=g+b \,\overline{\,g\,},
\end{align}
is analytic. 

Conversely, if  $b\in \mathcal{B}_{H^\infty}(\DD)$ and $\gamma(z)$ is an analytic function  in $\DD$,  then
via the identity
\begin{align}
\label{lemMapHoloTeleo1}
g(z)=\frac{1}{1-\vert b(z)\vert^2}\bigl(\gamma(z) - b(z)\overline{\gamma(z)} \,\bigr), \qquad z \in \DD,
\end{align}
the pair $(b,\gamma)$ determines a solution to  the equation 
\begin{align} \label{teleo34}
 \dbar g=-b(z) \overline{\, \dv_z g \,},
\end{align}
and the relation  \eqref{lemMapTeleoHolo1} holds. 
\end{Prop}

\begin{proof}
This is a direct computation; we get 
\begin{align*}
\dbar \gamma(z)=\dbar [g(z)+b(z)\overline{g(z)}] =\dbar g(z)+b(z)\overline{\dv_z g(z)}= 0.
\end{align*}
Given $\gamma \in \mathcal{O}(\DD)$ as in (\ref{lemMapTeleoHolo1}), we may conversely solve for $g$. Taking the conjugate of (\ref{lemMapTeleoHolo1}) gives
$\overline{\gamma }=\overline{g}+\overline{b}\,g$, and this
 with (\ref{lemMapTeleoHolo1}) leads to a linear system, with unique solution 
 \eqref{lemMapHoloTeleo1}.
 
  Indeed,  by the latter identity any pair $(b,\gamma )$ of holomorphic functions with $b\in \mathcal{B}_{H^\infty}(\DD)$ defines a function $g(z)$. In particular, then $\gamma =g+b \, \overline{g}$ with $\gamma$ holomorphic and thus
 $\dbar g= -b \; \overline{\, \dv_z g \,}$. 
On the other hand, $g$ need not be injective.
\end{proof}

 \begin{rem}     It is often convenient to apply the operation of affine and M\"obius transformations on the teleomorphic functions.
 Via \eqref{factor}, or directly, we  see that affine and M\"obius  transformations operate similarly as well on the Beltrami equation \eqref{eqNonLinB} and the corresponding liquid domains $\LL$.
 \end{rem}

As a first simple application of the hodograph transform we consider the removability of an isolated singularity for solutions to \eqref{eqNonLinB}.
In general, mappings of finite distortion can have isolated singularities, see e.g. \cite[p. 540]{ATG}. However, this does not occur in the present case.

\begin{Cor} \label{remo}
Isolated singularities are removable for solutions to    $ \partial_{\overline{ z }} f(z) = f(z)  \partial_z f(z)$. That is, each continuous $W^{1,2}_{loc}$-solution  $f:\UU \setminus \{z_0\} \to \Di$ extends to a solution in $\UU$. 
\end{Cor}
\begin{proof}  We use the factorisation $f = b \circ g^{-1}$ from Theorem \ref{propStoilow}. Here $g:\DD \to \UU \setminus \{z_0\}$ is a teleomorphic homeomorphism from a circle domain,  with analytic coefficient 
$|b(z)| < 1$ as in \eqref{teleo34}.

If the homeomorphism $g^{-1}$ maps (neighbourhoods of) $z_0$ to  (relative neighbourhoods of) a non-degenerate boundary circle $S(\zeta_k, r)$ of $\DD$, then as in \eqref{lemMapTeleoHolo1} we can consider an auxiliary function
$$\gamma_0(\zeta) = (g(\zeta) - z_0) + b(\zeta) \,{\overline{  (g(\zeta) - z_0)}}, \qquad \zeta \in \DD.
$$ 
The teleomorphic equation \eqref{teleo34} shows that  $\gamma_0$ is analytic, while by our assumptions $\gamma_0(\zeta) \to 0$ as $\zeta \to S(\zeta_k, r)$.
However, an analytic function cannot vanish on non-degenerate boundary arc, unless it vanishes identically. Accordingly, $g^{-1}$ maps $z_0$ to singleton, and thus extends continuously across $z_0$.

Now  the bounded analytic function $b$ extends analytically across the isolated singularity $\zeta_0 = g^{-1}(z_0)$.
Therefore also $f = b \circ g^{-1}$ extends continuously to $\UU$, and solves there  the Beltrami equation $ \partial_{\overline{ z }} f(z) = f(z)  \partial_z f(z)$.
\end{proof}

\subsection{Real analytic extensions} 

Via the hodograph transform,  the teleomorphic maps and their boundary values give  us parametrisation of the frozen boundaries. The maps are real analytic and locally  quasiregular, thus locally well under control in their domain. However, their boundary behaviour is a more subtle issue. Somewhat surprisingly, there we have the best control in the cases where the ellipticity of the equation \eqref{teleo34} degenerates as strongly as possible, see Proposition \ref{surprise1} and Theorem \ref{surprise2}. 

No doubt the most interesting cases are when the degeneration happens on the whole boundary. However, we wish to use our analysis also for partially frozen boundaries, and therefore we formulate accordingly (most of) the results in this and the following subsection. On the other hand, even for dimer models with gas, this approach  allows us often to work with maps on the unit disc and thus makes the descriptions in the general case simpler.

\begin{Prop} \label{surprise1}
Let  $g$ be a bounded $W^{1,2}_{loc}$-solution to 
\begin{align}\label{againeqn11}
\dbar g(z)=-b(z)\, \overline{\dv g(z)}, \quad z\in \Di,
\end{align}
where $b(z)$ is analytic  with $|b(z)| < 1$  in the unit disc $\Di$. 

Suppose 
there is an open arc $I \subset \partial \Di$ on the unit circle, such that for the coefficient function 
 \begin{equation}\label{limit}
 \lim_{z \to w}\vert b(z)\vert =1  \quad {\rm for} \;  \emph{every}  \; \;  w \in I.
 \end{equation}
 Then both  $b(z)$ and  the  auxiliary function $\gamma(z) =g(z)+b(z) \,\overline{\,g(z)\,}$ are analytic on $I$ and  extend meromorphically to the domain 
  \begin{equation}\label{refl.dom}
  \Omega_I := \Di \cup I \, \cup (\overline{ \C} \setminus \overline{ \Di}).
   \end{equation}
 
In addition,   $b(z)$ has no critical points on $I \subset \partial \DD$,
 and we have  the identities
\begin{align} \label{symmetries}
\overline{ b(1/\bar z)} =  1/b(z), \hspace{.6cm} \gamma(z) = b(z) \, \overline{ \gamma(1/\bar z)}, \qquad   z \in \Omega_I.
\end{align}
\end{Prop}
\begin{proof}  
  Since $b(z)$ is bounded and analytic in $\Di$, we may  apply the famous Beurling factorization (cf. \cite[Corollary II.5.7]{Garnett}) which decomposes $b$ into a product of 
  three factors, 
$b(z)=B(z) O(z) S(z)$.

Here $B(z)$ is a (possibly infinite) Blaschke product. The term $ O(z)$ is the outer factor, 
\begin{align*}
O(z)=c\exp\bigg(\frac{1}{2\pi}\int_{-\pi}^{\pi}\frac{e^{i\theta}+z}{e^{i\theta}-z} \,\log|b(e^{i\theta})|\, d\theta\bigg),
\end{align*}
 where $b(e^{i\theta})$ stands for the non-tangential boundary values of the function; these 
  exists for a.e. $\theta$ with $\log|b(e^{i\theta})| \in L^1(\partial \Di)$, by  \cite[Theorem II.5.3]{Garnett}.   Further, $c$ is unimodular constant.
  
   Finally, the third factor is the 
singular function $S(z)$, defined by a measure $\mu$ on $\dv \Di$ which is singular  with respect to the arc length. More precisely, 
\begin{align*}
S(z)=\exp\bigg(\frac{1}{2\pi}\int_{-\pi}^{\pi}\frac{e^{i\theta}+z}{e^{i\theta}-z}d\mu(\theta)\bigg). 
\end{align*}

The factor $S(z)$ has  unimodular radial limits almost everywhere, but at points of the support of $\mu$ the  singular factor 
has radial limit $0$. By \eqref{limit} the measure is hence supported on $\partial \Di \setminus I$. Similarly, the Blaschke product $B(z)$ has no zeros accumulating at any point of 
$I$, and for the outer factor the term $\log|b(e^{i\theta})| = 0$ on $I$. It follows that each of the terms in the Beurling factorisation, and   
hence $b(z)$ itself, extends analytically across the  interval $I \subset \partial \Di$. 

In fact, $b(z)$ extends to a meromorphic function of $\Omega_I$, with no poles or zeros on $I$.
Furthermore, locally an analytic function is a composition of an integer power with a conformal mapping. Since $b(\Di) \subset \Di$ with $b(I) \subset \partial \Di$,  we  see that neither can $b(z)$  have critical points on $I$. Finally, since for $z \in I$ we have $1 = |b(z)|^2 = \overline{ b(1/\bar z)} b(z) $, by unique continuation the first of the identities \eqref{symmetries}
holds throughout $\Omega_I$.

For the auxiliary function $\gamma(z)$ we can make use of the representation 
$$g(z)=\frac{1}{1-\vert b(z)\vert^2}(\gamma(z) - b(z)\overline{\gamma(z)}), \qquad z \in \Di,
$$
and apply  Lemma \ref{help2} 
with the  choice $\alpha(z) = \gamma(z)/b(z)$ and $\beta(z) =  \gamma(z)$. 
Indeed, by \eqref{limit} the term $|b(z)|^2 \to 1$ as $z \to w \in I$ in $\Di$. Since the solution $g(z)$ is bounded by assumption,  
this forces $ \alpha(z) - \overline{\beta(z)} \to 0$  when $z \to w \in I$.

With Lemma \ref{help2} the auxiliary factor
$\gamma(z)$ now extends analytically across the interval $ I$ and with the symmetry \eqref{auxiH}, i.e. $ \overline{\gamma(1/\bar z)} = \gamma(z)/b(z)$, it becomes meromorphic in $\Omega_I$. 
\end{proof}

When the interval $I$ is the entire unit circle, the coefficient function $b(z)$ in \eqref{againeqn11} is just a finite Blaschke product, see Theorem \ref{propStoilow35}. With  Proposition  \ref{surprise1} one sees that  a general  coefficient  $b(z)$ has very similar properties
on intervals where the ellipticity of \eqref{againeqn11} degenerates uniformly, as in \eqref{limit}.
 
To make use of the above phenomena, the following meromorphic function of two complex variables encodes perhaps most efficiently the properties of the teleomorphic function 
$g(z)$ and the other related relevant quantities. 

\begin{Def}\label{2variables}
 If  $g(z)$ is a bounded solution to \eqref{againeqn11}, for $b(z)$ holomorphic with $|b(z)| < 1$ in $\Di $, we set
\begin{equation} \label{bestsym1}
\Phi(z,w) := \frac{\gamma(z)\, b(w) - \gamma(w)\, b(z)}{b(w) - b(z)}, \qquad w \neq z \in \C.
\end{equation} 
 where $\gamma(z) = \gamma_g(z)$ is the holomorphic factor from \eqref{lemMapTeleoHolo1}.
\end{Def}
\smallskip

\begin{rem} \label{rational3}
Under the assumptions of Proposition \ref{surprise1} the above function $\Phi(z,w)$ extends meromorphically to $\Omega_I  \times \Omega_I $.
In particular,  note  that if $$I = \partial \Di$$ is the entire unit circle, then $\Phi(z,w)$ and  $\gamma(z)$ are, in fact,  rational functions of $\C^2$ and $\C$, respectively.
\end{rem}

\medskip
The above approach implies  that   a teleomorphic function in $\Di$  admits a real analytic extension across any interval $I \subset \partial \Di$  where the ellipticity of  \eqref{againeqn11} degenerates, in the sense \eqref{limit}. It is not difficult to see that as a consequence,  analogous results hold on boundary arcs of a general circle domain $\DD$, see Lemma \ref{multirep3}. 

\begin{Thm} \label{surprise2}
Let  $g(z)$ be a bounded  $W^{1,2}_{loc}$-solution to  \eqref{againeqn11} in  the unit disc $\Di$, with the coefficient function $b(z)$ satisfying the assumptions of Proposition \ref{surprise1} for an interval $I \subset \partial \Di$. 

Then with  the function $\Phi(z,w)$ from \eqref{bestsym1}, 
 we have the representation
\begin{equation} \label{gsymmtr}
 g(z) = \Phi(z,1/\bar z), \qquad z \in \Omega_I. 
\end{equation} 
 Thus $g(z)$ extends to a real analytic function on all of $\Omega_I$, satisfying $g(1/\bar z) = g(z)$ there. Moreover, on the interval
\begin{equation} \label{onbdry}
 g_{{|}_I}(z) =  R(z), \qquad z \in I,
\end{equation} 
where $R(z)$ is meromorphic in $\Omega_I$. 

If $I = \partial \Di$ is the entire unit circle, then the boundary value   $R(z)$  is a rational function of $\C$. 
\end{Thm}

\begin{proof} The symmetries \eqref{symmetries} give
$$\Phi(z,1/\bar z) = \frac{\gamma(z)/ \, \overline{b(z)}- b(z) \overline{\gamma(z)}/ \, \overline{b(z)}}{1/ \, \overline{b(z)} - b(z)} = 
\frac{\gamma(z) - b(z)\, \overline{\gamma(z)}}{1- |b(z)|^2} = g(z),
$$
where the last identity is \eqref{lemMapHoloTeleo1}.
  Moreover, \eqref{bestsym1}  shows that on the diagonal $ w= z$, 
\begin{equation} \label{bestsym2}
\Phi(z,z) = \lim_{w \to z} \Phi(z,w)  = \gamma(z) -  \frac{b(z)}{ b'(z)} \,  \gamma'(z) =: R(z), 
\end{equation} 
a meromorphic function of one variable. If $I = \partial \Di$, Remark \ref{rational3} shows that $\Phi$, and hence $R$, is rational.
\end{proof}

We will actually need several slight variants of the above argument. For instance, if \eqref{againeqn11} - \eqref{limit} hold in an annulus $\{z:  1 <|z| < \rho \}$ with $I = \partial \Di$, the natural versions of \eqref{symmetries} and Theorem \ref{surprise2} hold  in the double annulus  $\{z:  1/\rho <|z| < \rho \}$.

 \subsection{Boundary regularity of teleomorphic homeomorphisms} \label{sec.bdryreg}

Our next task is to understand, in as detail as possible, the geometry of the meromorphic boundary function $R(z)$ from Theorem \ref{surprise2}. By precomposing  with suitable analytic functions it is not difficult to find  teleomorphic maps  and boundary values with critical points of arbitrarily high order. However, the situation is different for the maps in \eqref{factor} arising from the hodograph transform, as these are {\it homeomorphic}.  
Therefore, to describe the geometry and boundary regularity of frozen boundaries, as in Theorem \ref{First.thm},
in this section we consider only homeomorphic solutions to \eqref{againeqn11}. The following comes as  part of the proof of the Pokrovsky-Talapov law.

\begin{Prop} \label{bdryreg1}
Suppose $g(z)$ is a bounded  and homeomorphic $W^{1,2}_{loc}$-solution to  \eqref{againeqn11} in  $\Di$,
with coefficient $b(z)$ and the interval $I \subset \partial \Di$ as in Proposition \ref{surprise1}.

 Let further $R(z) = R_g(z)$ be the meromorphic boundary values of  $g(z)$ on $I$, as in \eqref{bestsym2}. If $z_0 \in I$ is not a critical point of $R$, then 
$$ | g(r z_0) - g(z_0)| \simeq C (1-r)^2, \qquad 0 < r < 1.
$$
In directions other than the normal, the derivatives $\partial_\zeta g(z_0) $ are nonzero, tangent to $g(\Di)$ at $g(z_0)$.
\end{Prop} 
\begin{proof} 
$g(z)$ is real analytic in a neighbourhood of $I\subset \partial \Di$,  with radial derivative $$\partial_r g(r e^{i\theta}) = e^{i\theta} g_{z }(r e^{i\theta}) + e^{-i\theta} g_{\overline{ z }}(r e^{i\theta}).$$ We need to show that
$$ \partial_r g(z_0) = 0 \quad {\rm with} \quad \partial_r^2 g(z_0) \neq 0, \qquad {\rm whenever} \quad R'(z_0) \neq 0,\; z_0 \in I.
$$
The first claim is clear since $g(r z_0) = g(z_0/r)$ for  $0 < r < \infty$, by Theorem \ref{surprise2}.
For the other claim,  
 derivate the equation  \eqref{againeqn11} to obtain 
\begin{eqnarray} \label{secondBelt246}
 g_{\bar z \bar z}(z)  = -b(z) \overline{g_{z z}(z)} \quad {\rm and }
 \quad  g_{z \bar z}(z) = -b'(z)\, \overline{g_{z}(z)} - b(z)\, \overline{g_{z \bar z}(z)}.
\end{eqnarray}
Inserting these identities to the expression $ \partial_r^2 g(z) = e^{i2 \theta} g_{z z}(z) + 2 g_{z \bar z}(z) +  e^{-i 2 \theta} g_{\bar z \bar z}(z) $
shows that 
\begin{eqnarray} \label{secondBelt3}
 \partial_r^2 g(z) + b(z) \, \overline{\partial_r^2 g(z) } = (1-|b(z)|^2)  \frac{z}{\bar z} \, g_{z z}(z) - 2 b'(z) \,\overline{g_{z}(z)}, \qquad z \in \Di.
\end{eqnarray} 

By real analyticity \eqref{secondBelt3} extends to $\Di \cup I$. Thus we only need 
to observe from \eqref{gsymmtr} that $g_{z}(z) = \partial_1 \Phi(z,1/\overline{z}) = \partial_1 \Phi(z,z)$ for $z \in I$, while  \eqref{bestsym2} implies 
$$ R'(z) = 2 \partial_1 \Phi(z,z), \qquad z \in \Omega_I.
$$
Hence if  $z_0 \in I$ with $\; R'(z_0) \neq 0$, necessarily $g_{z}(z_0) \neq 0$. As the derivative $b'(z)$ does not vanish on 
$I$, the identity \eqref{secondBelt3} shows that $\partial_r^2 g(z_0) \neq 0$.

For the last claim note that $i z_0 R'(z_0)$ is tangent to $g(\Di)$ at $g(z_0)$.  Thus with \eqref{gsymmtr} and $\zeta_0 = - z_0 e^{i\alpha}$,
differentiating $g(z_0 + t \zeta_0)$  gives 
$\partial_\zeta g(z_0) = - i \sin(\alpha) z_0 R'(z_0)$, which for $0 < |\alpha| < \pi$ is non-zero and parallel to the tangent at $g(z_0) \in \partial \LL$.
\end{proof}
\medskip

 \subsubsection{Critical points on the boundary}
 
 As  the many simulations show, the frozen boundaries typically have cusps. 
 In this subsection we  show that, however, even in the setting of locally frozen boundaries the possible cusps are always simple. 
 For that purpose we need some further analysis of homeomorphic solutions of \eqref{againeqn11}  and their meromorphic boundary parametrisations discovered in Theorem \ref{surprise2}. 
  
 Note  that in general for a  rational map on the unit circle, having an extension to a solution of some Beltrami equation - even a homeomorphic one -  the critical points can easily be of higher order. A simple example of this is   $R(z) =  \frac{1 + 3z^2}{z^3 + 3z}$, which  has a quasiconformal extension to $\Di$.
 
 We thus need to use the special structure of the specific equation \eqref{againeqn11}. For this it is again useful to apply
 the function $\Phi(z,w) $ from Definition \ref{2variables}, with $g(z_0) = \Phi(z_0,z_0)$ for  points $z_0 \in \Sii^1$ on the unit circle.
Let us start with

\begin{Lem} \label{aux1234}
Consider a  bidisc $U = \Di(z_0,\delta) \times  \Di(z_0,\delta) \subset \C^2$ and a function $\Psi(z,w)$ holomorphic in $U$.
Suppose there are analytic functions $\alpha, \beta: \Di(z_0,\delta) \to \C$ of one variable, such that 
  \begin{eqnarray}  \label{beta64}
  [\beta(z) - \beta(w)] \, \Psi(z,w) = \alpha(z) - \alpha(w) \quad {\rm with} \quad \beta'(z)  \neq 0, \qquad \forall \, (z,w) \in U. 
   \end{eqnarray} 
If $n\geq 1$ is the smallest integer with $\, \partial_1^k \partial_2^{n-k} \Psi(z_0,z_0) \neq 0$ for some $0 \leq k \leq n$, then
  \begin{equation} \label{Taylor14}
\Psi(z,w) - \Psi(z_0,z_0) = c_0 \sum_{k=0}^n (z-z_0)^k(w-z_0)^{n-k} + \sum_{s=n+1}^\infty P_s(z-z_0, w-z_0), \qquad z,w \in \Di(z_0,\delta),
\end{equation}
where $ P_s(\zeta, \eta)$ are $s$-homogeneous polynomials and $c_0 \neq 0$.
\end{Lem}
\begin{proof}  By subtracting a constant, we may assume that $\Psi(z_0,z_0) = 0$, this does not change the $\beta(z)$ term in \eqref{beta64}. 
We may also take $z_0 = 0$. Then
  \begin{equation*} \label{TaylorPsi}
   \Psi(z,w) =  \sum_{k=0}^n a_k \, z^k w^{n-k}  + \sum_{s=n+1}^\infty P_s(z, w), \qquad |z|,|w| < \delta,
\end{equation*}
where  $P_s(\zeta, \eta)$ are $s$-homogeneous polynomials, and by assumption some $a_k \neq 0$.

 Take then $z = t\zeta$ and $w = t\eta$, where $t > 0$ and $\zeta, \eta \in \Sii^1$,  and develop $\alpha(z)$ and $\beta(z)$ as Taylor series at $z_0 = 0$.  Letting now $t \to 0$ one observes from  \eqref{beta64} that $\alpha^{(k)}(z_0) = 0$ for $1 \leq k \leq n$. Moreover, the lowest order terms in the same identity give
 $$ (\zeta - \eta)  \sum_{k=0}^n a_k \, \zeta^k \eta^{n-k} =  c_0 (\zeta^{n+1} - \eta^{n+1}), \qquad \forall \zeta, \eta \in \Sii^1,
 $$
 where $c_0= \alpha^{(n+1)}(0)/ \beta'(0)$.
 In particular,  $\alpha^{(n+1)}(z_0) \neq 0$. Since this holds for every $\zeta, \eta  \in \Sii^1$, we must have $a_k = c_0$, for every
 $1 \leq k \leq n$. The claim follows.
 \end{proof}
  \medskip

Since the coefficient function $b(z)$ in Equation \eqref{againeqn11} has no critical points on the unit circle, see Proposition \ref{surprise1}, our function $\Phi(z,w)$ from Definition \ref{2variables} satisfies the requirements of the previous Lemma at any given  $z_0 \in I \subset \Sii^1$, with $ \beta(z) = \frac{1}{b(z)}$ and  $\alpha(z) = \frac{\gamma(z)}{b(z)}$.

We will then make use of the relation \eqref{gsymmtr} between $\Phi(z,w)$ and $g(z)$ in two ways. First, with help of  the homeomorphism $g(z)$ we show that now
only  the cases $n= 1$ and $n = 2$ are possible in \eqref{Taylor14}. Second,
if we have a critical point so that $n=2$, then the relation allows us to describe the exact boundary behaviour of the homeomorphism $g(z)$, see  Corollary \ref{gmap}.

\begin{Thm} \label{simplecuspthm} 
Let $g(z)$ be a bounded  and homeomorphic solution to  \eqref{againeqn11} in  $\Di$, and
assume  that  
 $g$ admits a meromorphic extension $R(z)$ across an interval  $I \subset \partial \Di$, as in Theorem \ref{surprise2}.

Then on the interval $I$,  every critical  point of  $R(z)$   is simple. 
\end{Thm}
\begin{proof} 
Suppose $z_0 \in \partial \Di$ is a critical point of $R(z)$. To show that $R''(z_0) \neq 0$ we use the auxiliary function
$ \Phi(z,w)$ from \eqref{gsymmtr}. Thus 
  \begin{eqnarray} \label{ratio321}
  R(z) =  \Phi(z,z) \quad {\rm for} \quad  z \in \C, \quad {\rm with} \quad  R'(z_0) = 2 \partial_1 \Phi(z_0,z_0) = 0.
   \end{eqnarray} 
   
 Adding a constant, we can assume that $g(z_0) =  \Phi(z_0,z_0) = 0$. Moreover, to simplify the notation we change variables with the Möbius transform
   \begin{eqnarray} \label{mobi} 
   \psi(z) = z_0 \frac{1+iz}{1-iz}, \qquad z \in {\mathbb H}_+ := \{ z \in \C : \Im z > 0 \},
   \end{eqnarray}
 where $ \psi ({\mathbb H}_+) = \Di, \; \psi(0) = z_0$.
  In the new coordinates $ \Phi $ still satisfies the assumptions of Lemma \ref{aux1234}. Thus it
  has the representation
$ \Phi ( \psi(z), \psi(w) ) = \Phi_0(z,w) + \Phi_1(z,w) $,
   \begin{eqnarray} \label{taylor22}
 \Phi_0(z,w) =  c_0 \sum_{k=0}^n z^k w^{n-k}, \qquad  \Phi_1(z,w) =  \sum_{s=n+1}^\infty Q_s(z, w), \qquad z,w \in \Di(0,\varepsilon), 
  \end{eqnarray} 
where $c_0 \neq 0$, the    $Q_s(z, w)$ are  $s-$homogeneous polynomials,  and where by \eqref{gsymmtr},
     \begin{eqnarray} \label{homeo34}
 z \mapsto  g_\psi(z) :=   
  \Phi ( \psi(z), \psi(\overline{ z }) ) \quad \mbox{is a homeomorphism in \;} {\mathbb H}_+ .
 \end{eqnarray}

 Our claim is that \eqref{taylor22} with \eqref{homeo34} forces $n= 2$. Indeed, normalising 
 $c_0 = 1$ we have
   \begin{eqnarray} \label{loworder}   \Phi_0(z, \bar z) =  \sum_{k=0}^n z^k {\bar z}^{n-k} = |z|^n \,  \frac{\sin\bigl((n+1) \arg(z) \bigr)}{\sin\bigl(\arg(z) \bigr)} \, \in \, \R, \qquad  z \in {\mathbb H}_+.
  \end{eqnarray}
In the upper half plane the function $\, \sin\bigl((n+1) \arg(z) \bigr) \, $ vanishes on the rays $\, \arg(z)  = \frac{k\pi}{n+1}, \, k = 1,2,\dots , n$, and changes sign alternatively in between.

 Let us then  consider the following cones, Jordan domains 
 $$\Gamma_k(\varepsilon, \delta) := \{z \in \C:   \left| \arg(z) - \frac{k\pi}{n+1}\right| <  \varepsilon, \; |z| < \delta \} \subset  {\mathbb H}_+, \qquad k = 1,2,\dots , n.
 $$
 The homeomorphsim $g_\psi(z)$ maps the cones $\Gamma(\varepsilon, \delta)$ onto disjoint  Jordan domains, each with $0$ on its boundary.

 Moreover, since   $$ \Phi_0(z, \bar z)  \simeq \,  \pm \,  \varepsilon \,  (-1)^k (n+1)\,   |z|^n \quad \mbox{ on the sides} \quad \arg(z) =  \frac{k\pi}{n+1} \pm \varepsilon \quad \mbox{ of the cone, }$$  and since $| \Phi_1 ( \psi(z), \psi(\bar z) ) | \leq C |z|^{n+1}$, we see that for $\delta > 0$ small, the images of cone sides are Jordan arcs emanating from the origin, one in the left halfplane $\{ z : \Re e \,  z < 0 \}$ and the other in the right half plane  
 $\{ z : \Re e \, z > 0 \}$. In particular, each $g_\gamma\left(  \Gamma_k(\varepsilon, \delta)\right)$ is a Jordan domain containing an interval
 $(0,it], \, t \in \R,$ of the imaginary axis. But since the images of cones are disjoint, there can be at most two such cones in the upper half plane  $ {\mathbb H}_+ $. 
 
 Thus $n =2 $, and a derivation with \eqref{ratio321}-\eqref{taylor22} gives finally
 $ R''(z_0)  \neq 0.$
\end{proof}

 Next, let us use the representation \eqref{Taylor14} to study the behaviour of the mapping $g(z)$ at a cusp point, 
  i.e. at a critical point $z_0$ of the boundary parametrisation $R(z)$ from Theorem \ref{surprise2}.   For an illustration see  Figure \ref{atcusp}   below.  The boundary curve  $\Gamma := g(I)$  has unit tangent vector
$$\tau_{_ \Gamma}(w) := \frac{i \eta R'(\eta)}{|R'(\eta)|}, \qquad  w = R(\eta), \quad  \eta  \in I \setminus \{ z_0\},
$$ 
locally outside the cusp, the critical value $w_0 := R(z_0)$. 
As $z \to z_0$ on $I$, the unit tangent has a well defined  limit  $\tau_{_\Gamma}(w_0) $, the direction of the cusp.

\begin{Cor} \label{gmap}
Let $g(z)$ be a bounded  and homeomorphic solution to  \eqref{againeqn11} in  $\Di$, admitting
 a meromorphic extension $R(z)$ across an interval  $I \subset \partial \Di$, as in Theorem \ref{surprise2}.
  Let  $\Gamma = R(I)$ and $z_0 \in I$ a critical point of $R$ so that  $w_0 := R(z_0) \in \Gamma$ is a (simple) cusp.  
 
 \noindent Then, with $-z_0$ being  the direction of the inner normal of $\partial \Di$ at $z_0$, we have
 
 i)  For $-\pi/2 < \theta < \pi/2; \quad \theta \neq \pm \pi /6$,  the curves $ g(z_0 -  t z_0 e^{\pm i \theta} ), \; 0 < t < \varepsilon$, 
 are asymptotic to 
 
\hspace{.2cm} the line spanned by $\tau_{_\Gamma}(w_0) $, with 
 \begin{equation} \label{critregu}
g(z_0 -  \varepsilon z_0 e^{i \theta} ) - g(z_0) = {\mathcal O}(\varepsilon^2).
\end{equation}

ii) The curves $ g(z_0 -  t z_0 e^{\pm i \pi /6} ), \; 0 < t < \varepsilon$, are both asymptotic to the line orthogonal to 
$\tau_{_\Gamma}(w_0) $.  

\hspace{.2cm}  In addition 
 \begin{equation} \label{maxcomp}
 g(z_0 -  \varepsilon z_0 e^{\pm i \pi /6} ) - g(z_0) = {\mathcal O}(\varepsilon^3).
 \end{equation}
\end{Cor} 
\begin{proof}   
We may take $g(z_0) = 0$ and  the direction of the cusp
  $\tau_{_\Gamma}(w_0) \in \R_+$.  It is again convenient to change coordinates with the M\"obius transform \eqref{mobi}; note that $g_\psi(z)$ satisfies   \eqref{againeqn11} in the upper halfplane, with coefficient $b \circ \psi$.  
  
  Now $z_0 = 0$, and in the upper halfplane the boundary normal lies in the direction of the imaginary axis. Thus in the notation 
 \eqref{homeo34}  we arrive at
 $$ g_\psi (i\varepsilon e^{i\theta}) = c_1  |\varepsilon|^2 \,  \frac{\sin\bigl(3(\theta +  \pi/2) \bigr)}{\sin\bigl(\theta +  \pi/2 \bigr)} + {\mathcal O}(\varepsilon^3), \qquad c_1 < 0.$$
 This proves the first claim \eqref{critregu}. 
 
Concerning the second claim, the asymptotic directions of the curves $g(  i t  e^{\pm i \pi /6} )$, as $t \to 0$, are seen from the proof of Theorem \ref{simplecuspthm}. It remains to analyse the maximal compression 
 in these exceptional angles.  
 For this  we apply on the third order derivatives an analysis similar to 
 \eqref{secondBelt3}.

 Differentiating first   \eqref{secondBelt246} gives for $g = g_\psi$,
   \begin{eqnarray} \label{thirdorder}  
  g_{\bar z \bar z \bar z}(z) & = & - b(z) \overline{g_{z z z}(z)},  
 \quad  g_{z \bar z \bar z }(z)  = - b'(z) \overline{g_{z z}(z)} - b(z) \overline{g_{z z \bar z}(z)} \quad {\rm and }  \nonumber \\  \nonumber \\
   g_{z  z \bar z }(z) &=&  - b''(z) \overline{g_{z}(z)} -2 b'(z)  \overline{g_{z \bar z }(z)}- b(z) \overline{g_{z \bar z \bar z}(z)}.
\end{eqnarray}
Inserting these identities one has  for the directional derivatives
$ \; \partial_\alpha =  
  e^{i\alpha} \partial_z +  e^{-i\alpha} \partial_{\overline{ z }}  , \,
$
 \begin{eqnarray*} \partial_\alpha^3 g(z) &+& b(z)\, \overline{\partial_\alpha^3 g(z)} \\  
&=& e^{i 3 \alpha} (1-|b(z)|^2) g_{z z z}(z) - 3 e^{i\alpha} \left( b''(z) \, \overline{g_{z}(z)} 
+ 2 b'(z) \, \overline{g_{z \bar z}(z)} \right) - 3 e^{-i\alpha} \, b'(z) \, \overline{g_{z  z}(z)}. 
\end{eqnarray*}
At the critical point $z_0 = 0$ the first derivatives of $g$ vanish,   while \eqref{taylor22} - \eqref{loworder} imply for the second derivatives 
 $g_{\bar z  \bar z}(0) = 2 \,  g_{z \bar z}(0) =  \,  g_{z  z}(0) = 2 c_1 \neq 0.$ Combining all these identities  gives
$$\partial_\alpha^3 g(z) +   b(z)\, \overline{\partial_\alpha^3 g(z)}\quad  \to  \; - 6 e^{i\alpha} \, b'(0) \,   \overline{g_{z \bar z}(0)} - 3 e^{-i\alpha} \,  b'(0) \, \overline{g_{z  z}(0)} =  c_2 \sin (\alpha), $$
with $c_2 \neq 0$. Thus, outside the tangential directions, $\partial_\alpha^3 g(z_0) \neq 0$ which completes the proof.
\end{proof}

\begin{figure}[H]
\centering
\begin{tikzpicture}[xscale=1,yscale=1]
\draw[thick,densely dotted,red] (-{sqrt(2)},{sqrt(2)})--({1.2-sqrt(2)},{sqrt(2)+0.1)});
\draw[dashed,red] (-{sqrt(2)},{sqrt(2)})--({1.2-sqrt(2)},{sqrt(2)-0.8)});
\draw[thick] (0,0) circle [radius=2cm];
\draw (0,0)--(-{sqrt(2)},{sqrt(2)});
\draw (-{sqrt(2)},{sqrt(2)})--({1/(sqrt(2))-sqrt(2)},{sqrt(2)+1/(sqrt(2))});
\draw (-{sqrt(2)},{sqrt(2)})--({-1/(sqrt(2))-sqrt(2)},{sqrt(2)-1/(sqrt(2))});
\draw[blue] (-{sqrt(2)},{sqrt(2)})--(-1,0);
\draw[blue] (-{sqrt(2)},{sqrt(2)})--(-0.1,1.1);
\draw [blue,<-,>=stealth] ({-0.6},{sqrt(2)-0.2})--({-0.55},{sqrt(2)-0.209});
\draw [blue,<-,>=stealth] ({-1.15},{sqrt(2)-0.9})--({-1},{sqrt(2)-1.4});
\draw ({-1.15},{sqrt(2)-0.9}) arc (-75:-48:1);
\draw (-0.8,0.4) node{\tiny$\pi/6$};
\draw (-1.6,1.6) node{\tiny$z_0$};
\filldraw (-{sqrt(2)},{sqrt(2)})  circle (0.8pt);
\draw[->,>=stealth] plot [smooth , tension=0.9] coordinates {(3.2,0)(3.8,0.2)(4.4,0)};
\draw (3.8,0.4) node{\tiny$g$};

\draw[blue,xshift=7cm] (0,0)--(0,0.5);
\draw[blue,xshift=7cm] (0,0)--(0,-0.5);
\draw[blue,<-,>=stealth,xshift=7cm] (0,0.4)--(0,0.5);
\draw[blue,<-,>=stealth,xshift=7cm] (0,-0.4)--(0,-0.5);
\draw[blue] plot [smooth , tension=0.9,xshift=7cm] coordinates {(0,0.5)(0.1,1.4)(0.3,1.8)};
\draw[blue] plot [smooth , tension=0.9,xshift=7cm] coordinates {(0,-0.5)(0.1,-1.4)(0.3,-1.8)};
\draw[thick,densely dotted,color=red,domain=0:01.1,xshift=7cm]   plot (\x,{-1.5*\x^2}); 
\draw[dashed,color=red,domain=0:01.1,xshift=7cm]   plot (-\x,{-1.4*\x^2}); 
\draw[thick,domain=0:3.1,xshift=7cm]   plot (\x,{(0.4*\x)^1.5}); 
\draw[thick,domain=0:3.1,xshift=7cm]   plot (\x,{-(0.4*\x)^1.5}); 
\draw[->,>=stealth,xshift=7cm] (0,0)--(-1.5,0);
\draw[xshift=7cm] (-1,0.3) node{\tiny$\tau_\Gamma(w_0)$};
\end{tikzpicture}
\caption{}
\label{atcusp}
\end{figure}

 \addtocontents{toc}{\vspace{-4pt}}
\section{Proper Maps $f(z)$ and the Geometry of the Liquid Domains} \label{propersec}

Once the basic features of general solutions  
 to the   Beltrami equation  \eqref{eqNonLinB}
are established, a next step is to apply these 
 to the study of frozen boundaries. In view of Theorem \ref{key.connection}  this asks us
 to understand solutions   to the   Beltrami equation  that are  proper as maps $f:\LL\to \Di$.
 Indeed, it is this last property that allows a  detailed analysis and classification,
and gives rise to a finite dimensional space of solutions. 

For simply connected liquid domains $\LL$,  Theorems \ref{surprise2} and \ref{propStoilow} apply directly. Indeed, when
$f:\LL\to \Di$ is a proper map, then also the 
 analytic factor 
$b: \Di \to \Di$ in \eqref{factor} is  proper. Hence by Fatou's theorem $b = B$, a finite Blaschke product $B(z)$, 
and \eqref{eqLinB3} gets the form
\begin{align}
\label{eqLinB333}
\dbar g(z)=-B(z)\overline{\dv_z g(z)}, 	\qquad z \in \Di.
\end{align}
Further, the auxiliary analytic function $\gamma :=g+ B \,\overline{\,g\,}$ from \eqref{lemMapTeleoHolo1} is now 
a rational map, by Proposition \ref{surprise1}, and we have the invariance properties   
\begin{equation} \label{oncemore}
\gamma(z) = B(z)  \overline{ \,\gamma(1/ \bar z)\,} \quad {\rm and } \quad g(z) = g(1/ \bar z), \quad z \in \overline \C,
\end{equation}
from \eqref{symmetries} and \eqref{gsymmtr}.

This case 
already presents all essential ideas, even if now   gas points $g \in \Gg$ are not allowed.
Hence we first discuss simply connected domains with some length, and turn to the general multiply connected situation in  Subsections  \ref{subsect:topologyFiniteConn} and \ref{multiply.connected}.

\begin{Thm}
\label{propStoilow35}
Let $\LL \subset \C$ be a bounded simply connected domain, supporting a continuous $W^{1,2}_{loc}$-solution  $f(z)$ to
 \begin{equation} \label{once more1}
  \partial_{\overline{ z }} f(z) = f(z)  \partial_z f(z).
  \end{equation}
If the solution 
is a proper map  $f:\LL\to \Di$, then 

i) The map admits the factorization 
\begin{align} \label{factor333}
f=B\circ g^{-1},
\end{align}
where  $B(z)$ is a finite Blaschke product and $\, g:\Di\to \LL \,$ is a homeomorphic solution to \eqref{eqLinB333}

ii) Moreover, $g$ extends real analytically to ${\overline{\C}}$, with rational boundary values
  $g_{\big | \partial \Di} = R$ on the unit circle. In addition, $g$ is locally injective on $\partial \Di$. 
\end{Thm}

\begin{proof} With unit disc $\Di$  the uniformization domain in Theorem \ref{propStoilow}, 
apply the decomposition \eqref{factor} for the first claim.  The second claim 
follows directly from Theorem \ref{surprise2}, giving
$ g(z) = \Phi(z,1/\bar z)$, where $\Phi(z,w)$  is  rational in $\C^2$. 
On the unit circle  Theorem \ref{simplecuspthm}  shows that the possible critical point of 
$g_{\big| { \partial \Di }}(z) = \Phi(z,z)$ are simple, and hence $g$ 
 is locally injective on $\partial \Di$. 
\end{proof}

\begin{rem} \label{added} The above result implies that the solution $f:\LL\to \Di$ has a continuous  extension to $\overline{ \LL }$ in the sense of prime-ends. Little later, see Theorem \ref{thm:geometry-curve}, we improve conclusion ii) in Theorem \ref{propStoilow35} and show  that in fact $f$ is globally continuous in  $\overline{ \LL }$,
i.e. that  at the possible double points of $\partial \LL$ the one-sided limits of $f$ agree.
\end{rem} 

In fact, Theorem \ref{thm:localboundary}   and Corollary \ref{bdryRegu3} will later establish a quite precise picture of the boundary regularity of the solutions $f$. 
 On the other hand, already the mere boundary continuity is very useful,
see for instance Proposition 2 in \cite{KeOk07}.

As a first quick application of  the factorization in 
Theorem \ref{propStoilow35} 
we have an interpretation of  the asymptotic particle densities in the Aztec diamond  \cite{CEP}.

\begin{ex} \label{density} Random domino tilings of the Aztec diamond 
was  one of  the first   dimer model where  in  simulations the  frozen phenomena were observed.  
The existence of the frozen boundary, the arctic circle,  was then proven for this model by   Jockusch, Propp and Shor \cite{JPS}. 

In their work  \cite{CEP} on the Aztec diamond,
Cohn, Elkies and Propp   identified the asymptotic probability 
for a given domino to occur at a given place. 
 In the notation
of Example \ref{microdominoes}, if we normalise the liquid domain $\LL$ of the Aztec diamond to be the unit disc,
 then  \cite{CEP}  shows that the asymptotic probability for the northbound domino to occur at $(x,y)$ is equal to 
\begin{equation} \label{Elkiesd}
P(x,y) = \frac{1}{2} + \frac{1}{\pi} \arctan \left(     \frac{\sqrt{2} y -1}{\sqrt{1-x^2-y^2}}  \right), \qquad z = x+iy \in \Di.
\end{equation}
Naturally for $|z| \geq 1$, $P(x,y) = 1$ if $y > 1/\sqrt{2}$ and $P(x,y) = 0$ otherwise.

On the other hand, in view of Theorem \ref{repone} and Remark \ref{domH},
 the limiting height function $h$ for the Aztec diamond is given by
\begin{equation} \label{repreex}
 \nabla h(z)=\sum_{k=1}^4 \, p_k \, \omega_{\Di}\bigl( f(z); I_k\bigr), \qquad z \in  \LL = \Di,\quad p_k = i^k,
\end{equation}
 where $f:\Di \to \Di$ is a proper map solving $f_{\overline{ z }}(z) = f(z)^2 \, f_z(z) $, and  intervals $I_k$ all have the same  length, $\pi/2$.  The terms in \eqref{repreex} are intrinsic. In fact from Corollary \ref{uniqueSC} it  follows that the map $f$ is unique up to a choice of sign. We claim that actually  
  the asymptotic tile density  from \eqref{Elkiesd}, 
\begin{equation} \label{tiledens}
P(x,y) = \omega_{\Di}\bigl( f(z); I_1 \bigr), \qquad z = x +iy \in \Di.
\end{equation}
 
Indeed, using Theorem \ref{thm:geometry-curve} below for $\widehat f := f^2$, 
we see that $\deg(f) = 1$ with $f(\eta) = \pm \eta$ on $\partial \Di$. A choice of sign only permutes the $I_k$, so we 
assume $f$ to be the identity on $\partial \Di$. Further, via Theorem \ref{propStoilow35} we have $\widehat f = B \circ g^{-1}$ for a 
homeomorphism $g: \Di \to \Di$, solving  \eqref{eqLinB333} in $\Di$. One can thus take  $B(z) = z^2$, while
 this with  \eqref{lemMapHoloTeleo1} and Proposition \ref{propParaSolutionBlaschke}  gives  $g(z) = 2z(1+|z|^2)^{-1}$. Clearly 
 the choices made above leave \eqref{tiledens} invariant. In particular, 
$$f(z) = g^{-1}(z) = \frac{z}{|z|^2} \left(  1 - \sqrt{1-|z|^2}\right), \qquad z \in \Di. $$
 A simple way to obtain \eqref{tiledens} is now to note that $P \circ g(z) = \omega_{\Di}\bigl( z; I_1 \bigr)$. Namely the functions have the same boundary values on $\Di$, and an elementary derivation shows that $P \circ g$ is harmonic.
  \end{ex}

In \cite{CLP} Cohn, Larsen and Propp worked out the analogous probabilities for the hexagonal lozenges tilings. Again, with an argument similar but more tedious than above, one can show  the asymptotic tile probabilities agree with 
the corresponding expressions $ \omega_{\Di}\bigl( f(z); I_j\bigr)$, now for a solution $f: \LL \to \Di$ to \eqref{eqNonLinB}.  These examples make it very suggestive that, at least in the absence of quasi-frozen and gas phases,  for all dimer models  the   asymptotic edge or particle densities can be obtained as a pull back of the harmonic measure in the correct coordinates, i.e. as a pull back by a proper map solving  \eqref{beltB1}.

\subsection{A characterisation of simply connected domains with  frozen boundary}

It turns out that for every liquid domain $\LL$ with boundary  completely frozen,  $\partial \LL$ is the real locus of an algebraic curve.
 We will show this later in connection with multiply connected domains,  see  Theorem \ref{MultiplyLocus}.  
In the simply connected case  we  have  even a complete classification,  in terms of rational parametrizations of their boundary,  
of  domains $\LL$ with frozen boundary, i.e. domains supporting  proper maps $f:\LL \to \Di$ which are solutions to  \eqref{once more1}. 
 In view of Theorem \ref{Second.thm} and  Corollary \ref{Dconverse}, this is equivalent to asking which simply connected domains are liquid with frozen boundary in the Lozenges model.

For this, we say that a rational map $R(z)$ is {\it univalent near } $\partial \Di$, if for some $\varepsilon >0$, $R$ is 
univalent in the annulus $\{ z: 1-\varepsilon < |z| < 1\}$.

\begin{Thm} \label{characterize}
Let $\LL \subset \C$ be a simply connected and bounded domain. Then there exists a continuous $W^{1,2}_{loc}$-solution  $f(z)$ to  \eqref{once more1} which is  proper as a  map  $f:\LL \to \Di$, if and only if  
\begin{equation} \label{rat.repre}
\partial \LL = R(\partial \Di),
\end{equation} 
where the rational function  $R(z)$ satisfies the following three conditions.

i) {\rm (univalence)} \;  $R$ is bounded and univalent near $\partial \Di$.

ii)  {\rm ($R'$ is self-reflective)} \; For some finite Blaschke product $B(z)$, 
\begin{equation} \label{selfreflect}
 R'(z) = \frac{B(z)}{z^2} \overline{ \,R'(1/ \bar z) \,}, \qquad z \in \C.
\end{equation} 
\, iii)  {\rm (bound on poles)}  \; $B'(z) \overline{ \,R(1/ \bar z) \,}$ is analytic on $\Di$.

\noindent Under these conditions,  $\partial \LL$ is the real locus of an algebraic curve. 

Furthermore, on the unit circle $R(z) = g(z)$, where $g$ is a homeomorphic solution to \eqref{eqLinB333}. The required  
solution to \eqref{once more1} is then given by $f = B \circ g^{-1}$ with 
$B(z)$ as above. 
\end{Thm} 
The properties of Theorem  \ref{characterize} are intrinsic and do not depend on the choice of the rational map $R(z)$. That is,
if i) - iii)  hold and $b(z)$ is a M\"obius transform preserving the unit disc, then $R \circ b$ satisfies the conditions i) - iii) with respect to the Blaschke product $B \circ b$. 
\begin{proof} 
Suppose first that we have a simply connected domain $\LL$ and a proper map $f:\LL \to \Di$ which is a solution to
\eqref{once more1}. Apply then Theorem \ref{propStoilow35} and factorisation \eqref{factor333} to find the homeomorphism  $g: \Di \to \LL$ for which
$$
\dbar g(z)=-B(z)\overline{\dv_z g(z)}, 	\qquad z \in \Di.
$$
The corresponding analytic factor $\gamma :=g+ B \,\overline{\,g\,}$ from \eqref{lemMapTeleoHolo1} is  rational
with the symmetry \eqref{oncemore}. 

We know that $g$ extends to $\partial \Di$ 
 with rational boundary values, implying  \eqref{rat.repre}. 
 Indeed,  by \eqref{bestsym2}, on the unit circle
$ g$ equals  
\begin{equation}\label{tied}
 R(z) := \gamma(z) -B(z) \frac{\gamma'(z)}{B'(z)}, 
\end{equation}
where one notes that the derivative of a finite Blaschke product does not vanish on the unit circle. 
The  curve $\partial \LL = R(\partial \Di)$ encircles  a simply connected domain.
As will be shown in Theorem  \ref{thm:geometry-curve},  the singularities of $\partial \LL $ are either first order interior cusps or tacnodes, with  the argument of the  tangent $ \frac{dR(e^{it})}{dt}$  increasing in $t$. Taking the  orientation   into account,  $R$ thus maps each interior normal of $\partial \Di$ to an interior normal of $\partial \LL$, c.f. also Figure \ref{atcusp}. This gives Condition i).
 
In addition, the rational functions $\gamma(z)$ and $R(z)$ become tied by the relation
\begin{equation}\label{rat.rel2}
\gamma(z) = R(z) + B(z)  \overline{ \,R(1/ \bar z) \,}, \qquad z \in \C,
\end{equation}
since the identity holds on the unit circle.
Finally, combine the identities  \eqref{tied} and \eqref{rat.rel2}. That shows first $ \overline{ \,R(1/ \bar z) \,} = \frac{\gamma'(z)}{B'(z)}$ and then, via a derivation,
$$R'(z) = - B(z) \left( \frac{\gamma'(z)}{B'(z)}\right)' = \frac{B(z)}{z^2}  \overline{ \,R'(1/ \bar z) \,}.
$$
This is  Condition ii). We also notice that $\gamma(z) =g(z)+ B(z)\,\overline{\,g(z)\,}$ is bounded in $\Di$, so that  $\gamma'(z) = B'(z) \overline{ \,R(1/ \bar z) \,}$ is analytic in $\Di$ which gives the remaining Condition iii).

We then need to prove the converse direction, that \eqref{rat.repre} with Conditions i) - iii) provide us a solution to \eqref{once more1} which is proper in $\LL$. For this, given the rational function $R(z)$ from the Conditions, define  $\gamma(z)$ via the identity  \eqref{rat.rel2}.
This makes $\gamma$ a   rational function which clearly satisfies the first of the symmetries \eqref{oncemore}.
Differentiating  \eqref{rat.rel2} and using Condition ii) gives us $\gamma'(z) = B'(z) \overline{ \,R(1/ \bar z) \,}$, which is analytic and bounded in $\Di$ by Condition iii). Thus also $\gamma(z)$ is analytic and bounded  in the unit disc.

In this setup let us  apply  the procedure of Proposition \ref{LemMapHolo}, i.e. \eqref{lemMapHoloTeleo1}, and use the pair  $(\gamma, B)$ to build a solution to 
\begin{align} \label{antilinear34}
\dbar g(z)=-B(z)\overline{\dv_z g(z)}, 	\qquad z \in \Di,
\end{align}
 smooth and locally quasiregular inside the unit disc. For points approaching  the boundary,  \eqref{bestsym2} shows that
$$g(z) \to \gamma(z) - \frac{B(z)}{B'(z)} \gamma'(z) = \gamma(z)  - B(z)  \overline{ \,R(1/ \bar z) \,} = R(z).$$ On the other hand, Stoilow factorisation represents $g = \psi \circ G$, where $G$ is a homeomorphism of $\Di$ and $\psi$ is analytic in $\Di$.  We can thus use  Condition i) and the argument principle to see that $\psi$ is univalent.  
Therefore $g $ is a homeomorphism in $\Di$, with 
 $\LL = g(\Di)$.

The required  map $f$ can now  be constructed by simply taking $f=B\circ g^{-1}$. This is a proper map from $\LL$ to $\Di$. Moreover, by \eqref{antilinear34} and \cite[p.34]{ATG} the inverse $G = g^{-1}$ satisfies the Beltrami equation
$$ \dbar G = (B \circ g^{-1}) \dv_z G.
$$
Since $f$ and $G$ have the same complex dilatation, we obtain \eqref{once more1}.  As discussed above, the last remaining claim, that  $\partial \LL = R(\partial \Di)$ is the real locus of an algebraic curve, will be shown in Theorem \ref{MultiplyLocus}. 
\end{proof}

\subsection{The holomorphic factor $\gamma(z)$}

Proposition \ref{LemMapHolo} and the auxiliary function $\gamma(z) = g(z)+B(z) \,\overline{g(z)}$ give a good starting point for the study of general teleomorphic functions. It is useful 
to analyse the factor $\gamma(z)$ in more detail in the case when $g$ is a  homeomorphic solution to  
\eqref{eqLinB333}, with $B(z)$  a finite Blaschke product. Via Theorem \ref{characterize}, Proposition \ref{propParaSolutionBlaschke} and \eqref{tied}, this also gives an explicit finite dimensional parametrisation of all simply connected liquid domains with frozen boundary, having 
 $\deg(B) -2$ cusps on $\partial \LL$. 

The Blaschke product $B(z)$ may have zeros at the origin, and therefore    
  it is convenient to write the product in the form
  $$ B(z) \;  = \;  \eta \, z^m \prod_{k=1}^{d-m}\bigg(\frac{z - z_k}{1-\overline{z_k}z}\bigg) \;   = \;  \eta \,  z^m  \frac{N(z)}{D(z)}, \qquad | \eta  | = 1,
  $$
  where  the polynomials
  \begin{align}
\label{blaschkefactors}
N(z) =  \eta \,  \prod_{k=1}^{d-m} ({z - z_k}) \quad  {\rm and } \quad D(z) = \prod_{k=1}^{d-m} (1-\overline{z_k}z), \qquad {\rm with} \; z_k \neq 0.
\end{align}
Also $z_k = 0$ for $d-m <  k \leq d$,  $ \{ z_k \}_{k=1}^d \subset \Di$ and $d = \deg(B)$ is the degree of the rational map.

\begin{Prop}
\label{propParaSolutionBlaschke}
Assume that $g:\Di\to \UU$ is a bounded and homeomorphic solution to \eqref{eqLinB333}. Then
 the holomorphic factor  $\gamma(z)=g(z)+B(z)\overline{g(z)}$ has the form
$$ \gamma(z) = \frac{P(z)}{D(z)},
$$
where $D(z)$ is as in \eqref{blaschkefactors} and where 
\begin{align}\label{pe12}
P(z)= \left\{
 \begin{array}{ll}
  \alpha \prod_{j=1}^d(z-e^{i\theta_j}), & \text{ if }\quad g(z) \neq   0 \; {\rm in } \; \Di \\ \\
   \alpha (z-z_0)(1-\overline{z_0}z) \prod_{j=1}^{d-2}(z-e^{i\theta_j}), & \text{ if } \; g(z_0) = 0, \; \;  z_0 \in \Di,
   \end{array} \right.
\end{align}
for some $\alpha \in \C$, $z_0\in \Di$ and $e^{i\theta_1},...,e^{i\theta_d}\in \Sii^1$. 
\end{Prop}

\noindent {\it Proof}. Suppose first that $g(z)$ has no zeros in $\Di$. Then 
from its definition we see that neither can $\gamma(z)$  have  any zeros in the unit disc. Combined with  the symmetry 
$\gamma(z) = B(z)  \overline{ \,\gamma(1/ \bar z) \,}$ this implies that $\gamma(z)$ is non-vanishing  also in $\C\backslash \overline{\Di}$. 
Hence, in this case all the zeros of $\gamma(z)$ are located on the unit circle $\Sii^1$. 

Furthermore, the symmetry also shows  that $\gamma(z)$ has the same poles of the same order as $B(z)$ has in $\overline{\C}\backslash \overline{\Di}$. In particular,  at $\infty$ the  function $\gamma(z)$ has a pole of order $m$. Therefore, $\gamma(z)$ must have $d$ zeros on $\Sii^1$, say at $e^{i\theta_1},...,e^{i\theta_d}$, and admits the representation 
\begin{align} \label{gammafactor}
\gamma(z)=\alpha\,  \frac{\prod_{j=1}^d(z-e^{i\theta_j})}{\prod_{j=1}^{d-m}(1-\overline{z_j}z)}
\end{align}
for some constant $\alpha \in \C$. Further, if $0 \notin \partial \LL$ then $\gamma'(z) = B'(z) \overline{ \,R(1/ \bar z) \,} \neq 0$ on $\partial \Di$.
\medskip

In the case where $g(z)$ vanishes at some point  $z_0 \in \Di$, also $\gamma(z_0) = g(z_0) + B(z_0) \overline{g(z_0)} = 0$ and by
\eqref{lemMapHoloTeleo1} this is the only zero of $\gamma(z)$ in $\Di$. On the other hand, $z_0$ must be a simple zero: We know that besides $g$ also the inverse $ g^{-1}: \UU \to \Di$ is $C^{\infty}$-smooth, as a solution to \eqref{Belt}. Therefore from chain rule, see e.g.  \cite[p.34]{ATG}, it follows that $g_z \neq 0$ in $\Di$. Differentiating the defining identity of $\gamma(z)$ gives 
$$ \gamma'(z_0)=g_z(z_0)+B'(z_0)\overline{g(z_0)}+B(z_0)\overline{g_{\overline{ z }}(z_0)} = g_z(z_0) (1- |B(z_0)|^2) \neq 0.
$$

Finally using the symmetry of $\gamma(z)$, we see that if $B(z_0)\neq 0$ then $\gamma(1/\overline{z_0})=0$  while if $B(z)$ has a zero of order $k$ at $z_0$, $1 \leq k \leq d$, then  $\gamma(z)$ has a pole of order $k-1$ at $1/\overline{z_0}$. Elsewhere in $\hat{\C}\backslash \overline{\Di}$ the function $\gamma(z)$ is non-vanishing and has the same poles of the same order as $B$, including $\infty$. It follows that 
\begin{align*}
\hspace{4cm} \gamma(z)=\alpha_2(z-z_0)(1-\overline{z_0}z)\frac{\prod_{j=1}^{d-2}(z-e^{i\theta_j})}{\prod_{j=1}^{d-m}(1-\overline{z_j}z)}. \hspace{4cm} \Box
\end{align*}

\subsection{Univalent polynomials} \label{unipoly}
As a particular example, let us  explore the case when the 
Blaschke product in Theorem \ref{propStoilow35} is $B(z)=z^d$ for some positive integer $d \geqslant 2$. (We will see from the arguments below that there are no solutions for $d=1$). From Proposition \ref{propParaSolutionBlaschke} it follows that the holomorphic factor $\gamma(z) $ is a polynomial. With a translation of the liquid region, we may normalise  the homeomorphism in Theorem \ref{propStoilow35} by $g(0)=0$. In this case, see Proposition \ref{propParaSolutionBlaschke},
 $\gamma$ is a polynomial of degree $d-1$ of the form
\begin{align} \label{poly3}
\gamma(z)=
 \alpha z \prod_{i=1}^{d-2}(z-e^{i\theta_i}) \quad \alpha \neq 0.
\end{align}
The rational parametrization of $\partial \LL$ from Theorem \ref{propStoilow35} is given this time by a polynomial of degree $d-1$, $p(z)=a_1 z + \cdots a_{d-1} z^{d-1}$. Indeed, \eqref{tied} shows that
\begin{equation}
\label{eq:polynomial-p}
p=\gamma-\frac{1}{d}z \gamma', \quad p(0)=g(0)= 0.
\end{equation}

The symmetry \eqref{oncemore} takes now the form 
$\gamma(z) = z^d \,  \overline{ \,\gamma(1/ \bar z)\,}$.
Such polynomials are called \emph{self-inversive}, see \cite[Chapter 7]{S}. From the symmetry and \eqref{eq:polynomial-p}, or from \eqref{selfreflect}, we find that 
\begin{equation}
\label{eq:coeff}
\bar a_1 = (d-1) a_{d-1}.
\end{equation}
All  roots of such self-inverse polynomials occur on $\dv \Di$ or as conjugate pairs relative to $\dv \Di$.

Let us then denote by $S^*_{d-1}$ the class of polynomials of degree $d-1$  that are univalent, i.e. injective,  in $\mathbb D$ and satisfy \eqref{eq:coeff}. 
Note that imposing the condition $g(0)=0$ was for mere convenience, it only affects the value $a_0=g(0)$ (but might change the degree of $\gamma$). In this setting Theorem \ref{characterize}  reads as follows.

 \begin{Cor} \label{polyhomeo}
The polynomial parametrization $p$ from \eqref{eq:polynomial-p} is univalent in $\Di$, that is $p \in S^*_{d-1}$. 

Vice versa, every $p \in S^*_{d-1}$ arises in this fashion. That is, for every $p \in S^*_{d-1}$, there exists a $d$-to-$1$ proper mapping $f$ from $\mathcal{L}=p(\mathbb{D})$ onto $\mathbb{D}$, such that $\dbar f=f \dv_z f$. 
Moreover, the Blaschke term in the decomposition \eqref{factor333} is $B(z)=z^d$.
\end{Cor}
\begin{proof}
Since $g$ is homeomorphism, $p(\mathbb{S}^1)$ is a non-selfcrossing curve bounding the region $\mathcal L=g(\mathbb{D})$. It follows from the argument principle that $p$ takes each value of $\LL$ exactly once and thus $p$ maps $\mathbb{D}$ univalently onto $\LL$.  Note, in comparison, that a general rational function in Theorem \ref{characterize} can have poles in $\mathbb{D}$, and indeed will  in general not be globally univalent in $\mathbb{D}$.

In order to prove converse direction we need to show that every $p \in S^*_{d-1}$ satisfies the Conditions i) - iii) from Theorem \ref{characterize} with respect to $B(z) = z^d$. Here the first and third conditions are clear, while for the second, 
it follows from univalence and condition \eqref{eq:coeff}, see  \cite[Lemma 2.6]{LM14}, that $p'$ is self-inversive with respect to degree $d-2$,
\begin{equation}
\label{eq:pprime0} 
p'(z)= z^{d-2} \overline{p'(1/\bar z)}.
\end{equation}
This is Condition ii). Thus  Theorem \ref{characterize} gives the required solution $f=(g^{-1})^d$.
\end{proof}

When considering the special case $B(z)=z^d$, we are thus led to univalent polynomials in the unit disk. For an overview of the subject, see \cite[7.4]{S}.
Since $p$ is univalent in $\Di$, $p'$ does not vanish in $\Di$ and neither in $\C\setminus \overline{\Di}$ because of \eqref{eq:pprime0}. Thus all the $d-2$ critical points of $p$ are forced to be on the unit circle $\Sii^1$. Furthermore, these all have to be simple from univalence of $p$ and thus geometrically correspond to an (inward pointing) cusps.

\begin{ex}
The simplest examples come from placing the critical points at the roots of unity. 
The choice of $\gamma(z)=\frac{d}{d-1}(z+z^{d-1})$ leads to $p=z+\frac{1}{d-1} z^{d-1}$, with $p'=1+z^{d-2}$.
When $d=2$ the liquid region $\LL$ is the unit disk, for $d=3$ it is a cardioid, and in general it is an epicycloid with $d-2$ cusps.  
\end{ex}

By placing two critical points sufficiently close to each other the boundary curve might develop a double point (a tacnode). Extremal examples of this phenomenon are known as Suffridge curves \cite{LM14,S} -- these curves have $d-3$ tacnodes in addition to the $d-2$ cusps.

\subsection{Geometry of the boundary} \label{bdrygeometry}  

We saw in the previous section that with $B(z)=z^d$, the frozen boundary $\partial \mathcal{L}\,$ of a liquid domain  is parametrized by a univalent polynomial and is locally convex except at $d-2$ cusps. In the general case with $B$ a degree $d$ Blaschke product, Theorem \ref{propStoilow35} gives the boundary a parametrization by a rational function $R$.
This rational map need not be injective in all of the unit disk, however, the geometry of $\partial \mathcal{L}\,$ remains much the same, as we will next see.

For a  simply connected domain $\LL \subset \C$ with piecewise smooth boundary we let  $\tau(\zeta)$, $\zeta \in \partial \LL$, denote the unit tangent vector of $\dv \LL$, with direction induced by the counter-clockwise orientation of  the boundary.

\begin{Thm}
\label{thm:geometry-curve}
Suppose $\LL$ is simply connected with $\partial \LL$ frozen  and $f: \LL \to \Di$  
 as in Theorem \ref{propStoilow35}. Let $\; d = \deg(f)$. \, 
 Then 

i) The tangent vectors $\tau(\zeta)$ and  the boundary values $f(\zeta) \in \partial \Di$ are  related via the identity
\begin{align} \label{tau2}
f(\zeta)=-\tau(\zeta)^2, \qquad \zeta \in \partial \LL \setminus \{{\rm cusps} \}.
\end{align} 
\vspace{-.7cm}

ii) \, $\partial \mathcal{L}$ is locally strictly convex and  smooth, except at precisely $d-2$ cusps. 

\quad
Thus for every 
$\zeta \in \partial \LL$ outside the cusps and  
tacnodes,  $ B(\zeta,\varepsilon) \cap \LL$ is strictly convex
 for $\varepsilon > 0$ 
 
 \quad small 
enough. At the tacnodes  $\zeta \in \partial \LL$  the set
$ B(\zeta,\varepsilon) \cap \LL$ 
has two components, both convex. 
\vspace{.1cm}

iii) The map $f \colon \mathcal{L} \to \mathbb{D}$ extends continuously to the closure $\overline{\LL}$.

\end{Thm}
\begin{proof}
Let us first record how the tangent vector changes along the boundary. We follow here \cite[Lemma 2.7]{LM14}, where the polynomial case was covered.
For the rational boundary parametrisation  of $\partial \LL$, we use the symmetry \eqref{selfreflect} to find that
\begin{equation}
\label{tangenteqn} 
 z^2 R'(z)^2 = |R'(z)|^2 B(z),\qquad z \in \partial \Di,
\end{equation}
where $B(z)$ is the finite Blaschke product from \eqref{selfreflect}.
Thus on the unit circle, $z R'(z) = A(z) \sqrt{B(z)}$ for some continuous function $A \colon \partial \Di \to \mathbb{R}$. Here, since $|B(z)|=1$ on $ \partial \Di $, the squareroot  $\sqrt{B(z)}$ for $z \in \partial \Di$ is defined by $e^{i \frac12 \arg B(z)}$ with some continuous branch of the argument.
By Theorem \ref{simplecuspthm} all critical points of $R$ on $\partial \Di$ are simple. Thus the function $A$ changes sign exactly at each critical point of $R$.

  On the other hand, for $z=e^{it} \in \partial \Di$, we can identify the  tangent to $\partial \LL$ via $ \frac{dR(e^{it})}{dt}=i z R'(z)$. 
Hence the unit tangent field of $\partial \LL$ takes the form 
\begin{equation}
\label{eq:unit-tangent} 
\tau(\zeta)= i\,  \mbox{sgn}(A(z)) \sqrt{B(z)}, \qquad \zeta = R(z), \; \; z \in \mathbb{S}^1\setminus \{ \text{critical points of $R$}\}.
\end{equation}
But from Theorem \ref{propStoilow35}, $f\circ R(z) = B(z)$ on $\mathbb{S}^1$. This proves the claim i).

Since the map $g:\Di \to \LL$ from Theorem \ref{propStoilow35} is a homeomorphism, the boundary curve $g(\mathbb{S}^1) = R(\mathbb{S}^1)$ is non-selfcrossing. 
Further, the identity  \eqref{tau2}  shows that as one moves along $\partial \LL$ in the counter-clockwise direction, the argument of the tangent $\tau(\zeta)$ is strictly increasing  on each smooth arc of the boundary. That gives the convexity properties in claim ii), and shows that the  singularities of $\partial \LL$ are all simple cusps or tacnodes (double points).

Similarly, at each 
 cusp  the unit tangent vector on $\partial \LL$ makes a half-turn backwards, corresponding to the sign changes of $A$ in  \eqref{eq:unit-tangent}. Since  $\partial \LL$ is non-selfcrossing,  the unit tangent vector $\tau(\zeta)$ turns around once as we go along the boundary, while the term  $\sqrt{B(z)}$ from \eqref{eq:unit-tangent} turns around $d/2$ times. Thus \eqref{eq:unit-tangent}  forces us to have exactly $d-2$ cusps, (simple) critical points of $R$ on $\mathbb{S}^1$.

Finally, the mapping $g^{-1}: \LL \to \Di$ develops a discontinuity on the boundary $\partial \LL$ at its double points. However, such points share the same tangent line, meaning that $\tau \circ g(z)= - \tau \circ g(z')$, for the two pre-images $g(z) = \zeta = g(z')$. From \eqref{eq:unit-tangent} we thus have $B(z)=B(z')$ so that the function $f=B \circ g^{-1}$ is continuous even at double points of $\partial \LL$.
\end{proof}

 Theorem \ref{thm:geometry-curve}  has a curious consequence on the uniqueness  of solutions to the  Beltrami equation \eqref{once more1}: 
   There {\it exists at most one} proper map $f:\LL \to \Di$ solving  the equation.
   Thus  the  map $f$  and its properties are  intrinsic to $\LL$ ! 

\begin{Cor} \label{uniqueSC}
Let $\LL \subset \C$ be a bounded simply connected domain and assume $f_1, \, f_2 : \LL \to \Di$ are continuous proper maps, both solutions to $f_{\overline{ z }}(z) = f(z) \, f_z(z)$ in $\LL$.  Then  $ \, f_1 = f_2$.
\end{Cor}
\begin{proof} Given a triangle $N = {\overline{co}}\{p_1,p_2,p_3\}$,   Theorem \ref{converse} gives us
a harmonic homeomorphism $U: \Di \to N^\circ$  
and  two Lipschitz functions $h_1$ and $h_2$  on $\LL$,  
such that 
\begin{equation} \label{h12} 
 \nabla h_k(z) =  U \circ f_k(z) \quad {\rm with } \quad   {\rm div} \, \big(\nabla \sigma(\nabla h_k)\big)= 0,  \qquad z \in \LL, \; \; k=1,2.
\end{equation}
But  Theorem \ref{thm:geometry-curve} i) tells that $f_1$ and $f_2$  agree on $\partial \LL$, thus $ \nabla h_1$ and $ \nabla h_2$ have the same boundary values, in particular same tangential derivatives. Hence up to an additive constant, the functions 
$h_k$ agree on $\partial \LL$. As these are both minimizers for  $\int_{\LL}\sigma(\nabla h)dx$ we have  $h_1 = h_2$. 
With \eqref{h12} this proves the claim.
\end{proof}

The result holds for multiply connected domains, too, but requires little more work, see Corollary \ref{uniquef}.

\subsection{Topology of liquid regions}
\label{subsect:topologyFiniteConn}

Let us then turn to the geometry  of general  liquid regions $\LL$, as described in Definition \ref{LiqFrozen}.  To start with, the definition itself does not  require any connectivity properties of $\LL$, c.f.  Figure \ref{figQuasifrozen1}, but when the boundary is frozen, it is easy to see that $\LL$
cannot have infinitely many components.

\begin{Lem} \label{lem:FiniteComp} Suppose $N$ and $\sigma$ are as in \eqref{Pst}, and that in a bounded domain $\LL$ there  is a  solution  $h$ to \eqref{eq:EL34}, such that  $\nabla h: \LL \to N^\circ \setminus \Gg$   is a continuous and proper map. Then $\LL$  has at most finitely many connected components. 
\end{Lem}
\begin{proof}
Let $\{\LL_k\}$ be the components of $\LL$. Then for each component $\nabla h: \LL_k \to N^\circ \setminus \Gg$   is continuous and
 proper, thus  surjective. Choose  $w \in N^\circ \setminus \Gg$. Then for each $k$ we have   a point $z_k \in \LL_k$ with $\nabla h(z_k) = w$. If $\LL$ has infinitely many components, then the $\{z_k\}$ have a accumulation point $z_0 \in \partial \LL$. But then 
 $\nabla h(z_k) \not\to \partial N$ and $\nabla h$ cannot be proper on $\LL$.
  \end{proof}

Next show that each connected component of the liquid region is finitely connected,  when $\partial \LL$ is frozen. For that we need a few auxiliary results.
\begin{Lem}
\label{lem:ProperHolo111}
Suppose $\LL \subset \C$ is a bounded domain and  $f:\LL \to \Di$  a proper holomorphic map. If \,$0<r<1$, let $W_r$ be a component of $f^{-1}(\{w \in \Di: r<\vert w \vert<1\})$. Then $$f: W_r\to \{w\in \Di: r<\vert w \vert<1\}$$ is a proper holomorphic map and hence surjective.
 In addition, for every $0<r<1$, every boundary component of $\LL$ is contained in the closure 
 of $f^{-1}(\{w\in \Di: r<\vert w \vert<1\})$. 
\end{Lem}

\begin{proof}
Let $\{z_k\}_k\subset W_r$ be a sequence  such that $\lim_{k\to \infty}z_k=z^*\in \dv W_r$. In case $z^*\in \LL$, we have 
$\vert f(z^*)\vert=r$, while $\vert f(z^*)\vert =1$ if $z^*\in \partial \LL$. The first claim follows from this; similarly 
 for the second claim, it follows from properness of $f(z)$, that $|f(z)| \to 1$ when $z \to \partial \LL$.
\end{proof}

\begin{Lem}
\label{lem:ProperHolo222}
 Suppose $\LL\subset \C$ is a bounded multiply connected domain. Assume that $f:\LL\to \Di$ is a proper holomorphic map and let  $\Gamma_1$,  $\Gamma_2$  be two different components of $\partial \LL$. 
 
 Then for $0<r<1$ large enough,  $\Gamma_1$ and $\Gamma_2$ cannot lie in the closure of the same component  $W_r$, where $W_r$ defined as in Lemma \ref{lem:ProperHolo111}. 
\end{Lem}

\begin{proof}
Assume the contrary. Then for all $0<r<1$ the exists a path $\gamma_r:(0,1)\to W_r$ which connects $\Gamma_1$ and $\Gamma_2$. On the other hand, we can separate $\Gamma_1$ and $\Gamma_2$ by a curve $\tilde{\gamma}$ in $\LL$ going through the point $w_0$, where $f(w_0)=0$. 

Consider next the hyperbolic metric $d_\LL$ of $\LL$. That is, by the uniformization theorem for arbitrary planar domains \cite{FHW} we may choose the unit disc as  the universal covering  of $\LL$,  and then equip $\LL$ with the push-forward of the Poincare metric of $\Di$ by the  covering map. 
In the hyperbolic metric the length of $\tilde{\gamma}$ is bounded, say less than $M<\infty$. Moreover, since $\tilde{\gamma}$ separates $B_1$ and $B_2$ in $\LL$,
 the paths $\gamma_r$ and $\tilde{\gamma}$ must intersect, so that $d_\LL(\gamma_r(t),w_0)\leq M$ for some parameter $t\in (0,1)$. 
 
 On the other hand by the Ahlfors-Schwarz-Pick Lemma \cite{A}, analytic mappings are contractions in the hyperbolic metric, which gives 
\begin{align*}
d_\Di(f\circ \gamma_r(t),0)\leq d_\LL(\gamma_r(t),w_0)\leq M.
\end{align*}
Thus if we take $r<1$ so close to 1 that the hyperbolic distance $d_\Di(r,0)>M$, necessarily $|f\circ \gamma_r(t)| < r$ 
at this special parameter point, and therefore the path $\gamma_r$ cannot be contained in the component $W_r$. 
\end{proof}

\begin{Prop}
\label{prop:ProperHoloFinite11}
Assume that $\LL\subset \C$ is a bounded domain and that $f:\LL\to \Di$ is a proper holomorphic map. Then $\LL$ is finitely connected. 
\end{Prop}

\begin{proof}
By Rado's theorem, \cite[Rado's theorem, p. 219]{R}, a holomorphic map $f: \LL \to V$ is proper if and only if $deg_wf$ is constant and finite. But Lemmas \ref{lem:ProperHolo111} - \ref{lem:ProperHolo222} \, show for a proper holomorphic map $f:\LL\to \Di$ that $\deg_w f$ gives an upper bound for the number of connected components of the boundary $\partial \LL$. This proves the claim.
\end{proof}

\begin{Prop}
\label{prop:FiniteConn}
Let $\LL\subset \C$ be a bounded  domain and assume that $f:\LL\to \Di$ is a proper map solving the Beltrami equation $\dbar f(z)=f(z)\dv f(z)$ for $z\in \LL$. Then $\LL$ is finitely connected. 
\end{Prop}

\begin{proof} As in the proof of Theorem \ref{propStoilow}, we give $\LL$ a complex structure ${\mathcal A}$ by requiring the analytic charts to satisfy \eqref{charts}. Then
$(\LL, {\mathcal A})$ is a planar Riemann surface, admitting a uniformization $G: \LL\to \LL'$ to some planar domain $\LL'$. This solves the Beltrami equation $\dbar G(z)=f(z)\dv G(z)$ and by classical
Stoilow's factorization theorem,
$f=b \circ G$ where  
$b:\LL'\to \Di$ is a proper holomorphic map. By Proposition \ref{prop:ProperHoloFinite11} the domain $\LL'$, and hence $\LL$, is finitely connected. 
\end{proof}

If the surface tension has gas points, the liquid domains with frozen boundaries are necessarily multiply connected, as shown  by Theorem \ref{thm:h:proper}. On the other hand, in lack of gas  this does not happen for  a simply connected  polygonal  domain $\Omega$ in the variational problem \eqref{ELbasic}.

\begin{Thm}\label{thm:SimplyConnect}
Let $\Omega\subset\R^2$ be a bounded simply connected Lipschitz domain, $h_0$ an admissible boundary value
and $h$ the minimizer of the variational problem (\ref{ELbasic}), where 
 $\sigma$ has no gas points. 
 
 If  $\LL_0$  is  a component of the liquid domain of $h$ with boundary $\partial \LL_0$ frozen, then $\LL_0$  is 
simply connected.
\end{Thm}

  \begin{proof} 
First note that by Lemma \ref{lem:FiniteComp} $\LL$, the liquid domain of $h$ has finitely many connected components and by Proposition \ref{prop:FiniteConn} each component is finitely connected. 

We now argue by contradiction and assume that the component $\LL_0$ of $\LL$  is not simply connected. Then $\dv \LL_0$ has a component $\Gamma \subset \partial \LL_0$, which the domain $\LL_0$ separates from $\partial \Omega$. By Corollary \ref{remo}  the continuum  $\Gamma$ cannot be a singleton.

If the minimizer $h$ is not $C^1$ on $\Gamma$, say at a point $z_0\in \Gamma$, then by Theorem 1.3 in \cite{DeS10} there exists a line segment $\ell\subset \overline{\Om}$ connecting $z_0$ to $\dv \Om$ and such that $h(z)=h(z_0)+\langle p,z-z_0\rangle$ for all $z\in \ell$ and for some $p\in \mathscr{P}\subset \dv N$. But then $\ell$ must intersect $\LL_0 \subset \LL$. However $\nabla h(z)\in N^\circ$ for all $z\in \LL$, a contradiction. Thus $h$ is $C^1$ on $\Gamma$. 

On the other hand, from  Theorem \ref{repone} we have  the representation formula 
\begin{align*}
\nabla h(z)=\sum_{j} p_j \, \omega_{\Di}(f(z),I_j), \qquad z \in \LL_0,
\end{align*}
where $f:\LL_0 \to \Di$ is a solution to the Beltrami equation \eqref{beltB1}  and also a proper map, since $\partial \LL_0$ is frozen for $h$, c.f. \eqref{repsecond}. As $h \in C^1(\Gamma)$, then the representation formula forces 
$\nabla h(z) = p$ for some $p\in \mathscr{P}$ and  for all $z\in \Gamma$.

Next,  via Corollary  \ref{cor:propStoilow} and \eqref{mu12} we have the factorisation $f= b \circ g^{-1}$, where $b:\mathcal{D}\to \Di$ is a proper holomorphic map and $g:\mathcal{D} \to \LL_0$ a homeomorphism, with  $\mathcal{D}$ a circle domain.  Here, like for any holomorphic and proper map from $\mathcal{D}$ to $\Di$,  the boundary values $b: \dv \mathcal{D}\to \partial \Di$ take each component of $ \dv \mathcal{D}$  surjectively to $\partial \Di$. In particular, that holds for the component $g^{-1}( \Gamma)$. However, in view of the above representation formula, this contradicts the fact that $\nabla h(z)=p$ for all $z\in \Gamma$.  
  \end{proof}

\subsection{Multiply connected liquid domains} \label{multiply.connected}

We then sum up the properties of   the nonlinear Beltrami equation 
\begin{align}
\label{eqNonLinB2}
\dbar f(z) =f(z) \dv_z f(z), \qquad z \in  \LL,
\end{align}
in general (bounded) planar domains which are not necessarily simply connected. This situation  arises  as soon as one has gas components but it appears naturally as well in many other settings, see  e.g. Figure \ref{fig.first}. 

We already know from the previous section that proper solutions $f \colon \LL \to \Di$  to \eqref{eqNonLinB2} exist only in finitely  connected domains. To describe then the boundaries $\partial \LL$ in detail we first need to understand their parametrisations. For this, recall from Theorem \ref{propStoilow} the representation
\begin{align} \label{factormulti}
f=b \circ g^{-1},  
\end{align}
where $b \colon \DD \to \Di$  is  analytic,   $\DD$ is a circle domain  and $g:\DD\to \LL$ is a homeomorphic $W^{1,2}_{loc}$-solution to  the linear Beltrami equation
\begin{align}
\label{eqLinB22}
\dbar g(z) = - b(z) \,\overline{\dv_z g(z)},  \qquad z \in  \DD.
\end{align}
 Assuming  $f \colon \LL \to \Di$  to be proper, this makes $b \colon \DD \to \Di$  an analytic proper map. 

The  boundary values of the homeomorphism $g(z)$ in \eqref{factormulti} - \eqref{eqLinB22} give now a parametrisation of $\partial \LL$. To analyse its properties, consider again the holomorphic factor $\gamma(z) := g(z) + b(z)\, \overline{g(z)}$. As in  the proof of Proposition \ref{surprise1} we see from \eqref{lemMapHoloTeleo1}, since $g$ is  bounded and $b \colon \DD \to \Di$  is proper, that
 \begin{equation} \label{gammafmla}
 \gamma(z)-b(z) \,\overline{\gamma(z)} \to 0 \qquad {\rm as} \qquad  z \to \dv \DD.
  \end{equation}
In view of  Lemma \ref{help2}, this means that $\gamma$  extends meromorphically beyond the boundary circles of $\dv \DD$, in particular it is continuous up to $\dv \DD$. 
 \medskip
 
 In fact, the best way to describe the analytic extension of $\gamma(z)$ (and of the other relevant functions here)  is in terms of the Schottky double $\widehat \DD$ of $\DD$, for details of this construction see e.g. \cite[Section 2.2]{SchifferSpencer}. In our case $\widehat \DD$ is a compact Riemann surface containing $\DD$, equipped with an anti-analytic reflection $j(z)$ that fixes $\partial \DD $ pointwise.
 In this terminology, the argument of Lemma \ref{help2} shows that $\gamma(z)$ and $b(z)$ extend meromorphically to $\widehat \DD$ with transformation rules
 \begin{equation} \label{invariance33}
 \gamma\bigl( j(z) \bigr) = \frac{\overline{\gamma(z)}}{\, \overline{b(z)}\, } \quad {\rm and} \quad  b\bigl( j(z) \bigr) = \frac{1}{\, \overline{b(z)}\, }, \qquad z \in \widehat \DD.
 \end{equation}
 Since $b \colon \DD \to \Di$  is proper, Schwarz reflection principle shows that $b(z)$ extends  analytically across the boundary circles of $\partial \DD$. For the extension and $\zeta \in \partial \DD$ one has $|b(\zeta)|=1$, and   $b \bigl(\DD \cap B(\zeta,\varepsilon)\bigr) \subset \Di$.   These mapping properties show  that  $b(z)$ has no critical points on $\partial \DD$. 
 \smallskip
 
   \begin{Lem} \label{multirep3} 
   The function 
   \begin{equation} \label{aux3}
   \Phi(z,w) = \frac{b(w)\, \gamma(z) - \gamma(w)\, b(z)}{b(w) - b(z)} 
   \end{equation}
   is meromorphic in $\widehat \DD \times \widehat \DD$. Moreover,  
  \begin{equation} \label{circleg}
   g(z) =   \Phi(z,j(z)), \qquad z \in \DD.
  \end{equation}
 Thus $g(z)$ has  meromorphic boundary values,
  $g(z) = \Phi(z,z)$ for $z \in \partial \DD$, and it extends real analytically to the Schottky double $\widehat \DD$.
   \end{Lem}
  \begin{proof} It is clear that $\Phi(z,w)$ is meromorphic outside the diagonal of $\widehat \DD \times \widehat \DD$, and from Taylor series of the local representatives we see that
     \begin{equation} \label{SchottkyFi}
      \lim_{w \to z}\Phi(z,w) = \gamma(z) - b(z) \frac{\gamma'(z)}{b'(z)} = \frac{(\gamma/b)'}{(1/b)'}, \qquad z \in \widehat \DD. 
  \end{equation}
   Here $R(z) := \Phi(z,z)$ is meromorphic on $\widehat \DD$, in fact the ratio of  differentials of two meromorphic functions on
   the surface.
  From the transformation rules \eqref{invariance33} we then have
     \begin{equation} \label{SchottkyFi12} 
      \Phi(z,j(z)) =  \frac{b(j(z))\, \gamma(z) - \gamma(j(z)) \,b(z)}{b(j(z)) - b(z)} = \frac{\gamma(z) -b(z)  \,\overline{\gamma(z)}}{1- |b(z)|^2} = g(z), \quad z  \in \DD,
    \end{equation}
  which proves the second claim. 
   \end{proof}

  In a similar fashion the symmetries \eqref{invariance33} show that 
   \begin{equation} \label{SchottkyFi2}
    \overline{ \Phi\bigl(j(w),j(z)\bigr) } = \frac{\gamma(w) - \gamma(z)}{b(w) - b(z)}, \qquad z \neq w \in \widehat \DD,
  \end{equation}
    so that $ R^*(z) :=  \lim_{w \to z}  \overline{\Phi\bigl(j(w),j(z)\bigr)} = \gamma'(z)/b'(z)$ defines a meromorphic function on $\widehat \DD$. 
  \medskip

  In the case of lozenges model, Kenyon and Okounkov  \cite[Theorem 2]{KeOk07} showed  for a large class of special boundary values (with the  liquid domain $\LL$ simply connected)  that the frozen boundary $\partial \LL$  is algebraic. The following result, in combination with Theorems \ref{key.connection} and \ref{thm:h:proper3}, generalizes algebraicity to all dimer models, and also to all natural polygonal domains with oriented boundary values. In particular, here we allow multiply connected domains. 

\begin{Thm}
\label{MultiplyLocus}
Suppose $ \LL  \subset \C$ is a bounded   domain, supporting a continuous and proper map $f \colon  \LL   \to \Di$, which is a $W^{1,2}_{loc}$-solution  to
 \begin{equation} \label{once more13}
  \partial_{\overline{ z }} f(z) = f(z)  \partial_z f(z).
  \end{equation}
Then the boundary $\partial \LL $ is  the real locus of an algebraic curve (minus the finite set of possible isolated points of the curve). 

In addition,  $\partial \LL $  does not have   degenerate boundary components. It has  finitely many singularities, and these are either simple cusps or tacnodes.
\end{Thm} 
\begin{proof} To start with, by Corollary \ref{remo} the boundary $\partial \LL $  does not have  any isolated points. With Lemma \ref{multirep3} we represent $\partial \LL = g(\partial \DD)$   as the image of the boundary of a finitely connected circle domain, under the function $R(z) = \Phi(z,z)$ meromorphic on $\widehat \DD$. Thus by Theorem \ref{simplecuspthm} the singularities of $\partial \LL $ are all either simple cusps or tacnodes. 

In addition,  \eqref{SchottkyFi} - \eqref{SchottkyFi2} give $R^*(z) = \overline{(R \circ j)(z)}$, so that  $R(z) = \overline{R^*(z)} = \Phi(z,z) = g(z)$ for points $z \in \partial \DD$.
Combining this with the identity $\gamma(z) := g(z) + b(z)\, \overline{g(z)}$ implies
\begin{equation} \label{multisym1}
  R(z) =  \gamma(z) - b(z) R^*(z), \qquad z \in \widehat \DD.
\end{equation}
As $R^*(z) = \gamma'(z)/b'(z)$,  for the differentials we then have
\begin{equation} \label{multisym2}
R'(z) = -b(z) (R^*)'(z), \qquad z \in \widehat \DD,
\end{equation}
an analogue of \eqref{selfreflect}.

 Next, any two meromorphic functions on a compact Riemann surface are polynomially related. 
Thus we can find a non-trivial and irreducible  polynomial $P(\zeta,\omega)$ of two complex variables, such that $P(R, R^*)=0$ on the Schottky dual  $\widehat \DD$.  In particular,  
\begin{equation} \label{rlocus}
R(\partial \DD) \subset \{\zeta \in \C:  P(\zeta, \overline{\zeta}) = 0 \}.
\end{equation}
 Our task is now to show that these two sets are equal,  up to the  set of isolated points of the real locus $ \{\zeta \in \C:  P(\zeta, \overline{\zeta}) = 0 \}$.
 For this, let
\begin{equation} \label{curve}
 {\mathcal C}:=  \{(\zeta,\omega) \in \C \times \C:  P(\zeta, \omega) = 0 \},
\end{equation}
 let $\alpha:{\mathcal C} \to \overline{\, \mathcal C \,}$ stand for the embedding to  its projective closure 
 in $\Pro_2(\C)$ and let $S \subset  \mathcal C$ with $\overline{\, S  \,} \subset \overline{\, \mathcal C \,}$ be the set of singularities of $ {\mathcal C}$ and $\overline{\, \mathcal C \,}$, respectively. The curve allows a resolution of singularities, a map 
$\pi: \widehat \DD \to \overline{\, \mathcal C \,}$ which is biholomorphic away from $\overline{\, S  \,}$. 

Mimicking   \cite{BNg}, consider the finite set
$$E = \{ z \in  \widehat \DD:  \; \; {\rm either } \; \; R(z) = \infty, \; R^*(z) = \infty  \; \; {\rm or } \; \; (R \times R^*)(z) \in S \}.$$
 In the complement,  for $z \in \widehat \DD \setminus E$,  
 one can define
 $  H(z) :=   \pi^{-1} \circ \alpha \circ (R \times R^*)(z)$. This is a meromorphic function, and since $R \times R^*$ takes 
 $ \widehat \DD \setminus E$ to  $ \, \mathcal C \subset \C \times \C$,
 we can write
\begin{equation} \label{self}
 (R \times R^*)(z)  = (\beta_1 \times \beta_2) \circ H(z), \qquad z \in  \widehat \DD \setminus E.
\end{equation}
According to \cite[Lemma 8]{BNg} both $\beta_1$, $\beta_2$  extend to meromorphic functions of $ \widehat \DD$, and similarly, $H$ extends  meromorphically to $H:  \widehat \DD \to  \widehat \DD$. Moreover,   $\beta_1 \times \beta_2$ is injective  outside $H(E)$.

A main step in our argument is now to show that $H$ is injective outside  a finite set. Here the key is to use the special symmetries of $R$ and $R^*$. 
Namely, upon differentiating \eqref{self} the symmetry \eqref{multisym2} gives for the differentials, 
$$ (\beta_1' \circ H)(z)  = - b(z) \, (\beta_2' \circ H)(z), \qquad z \in  \widehat \DD.
$$

Thus if $H(z) = H(w)$ for $z,w \in \widehat \DD \setminus E $ then $b(z) = b(w)$, outside the common critical points of $R$ and $R^*$. But $b(\DD) = \Di$ with  $|b\bigl( j(z) \bigr)| = 1/ |b(z)|$.  Therefore we have only three possibilities, either both  $z, w \in \DD$,  they both lie on the boundary $\partial \DD$ or as a third case,
both $j(z), j(w) \in \DD$.

If we consider the first case $z, w \in \DD$, combining \eqref{multisym1} and  \eqref{self}   gives $\gamma(z) = \gamma(w)$. Inserting this information with $b(z) = b(w)$ to \eqref{SchottkyFi12}  shows that $g(z) = g(w)$. But $g$ is a homeomorphism on $\DD$, so that 
$z = w$. The same argument shows that  $H$ is injective in $j(\DD)$, outside the union of $E$ and the common critical points of $R$ and $R^*$.
Last, on  $\partial \DD$ the boundary values  of $ g$  (i.e. $R$)  is injective outside the possible tacnodes.
With \eqref{self}  the same holds for $H$.

All in all, we have shown that $H: \widehat \DD \to \widehat \DD$ is meromorphic and injective outside a finite set, which implies that $H$ is biholomorphic. 

To conclude the theorem, suppose $\zeta \in \C$ lies on the real locus of $P$, i.e.  $P(\zeta, \overline{\zeta}) = 0$. 
Then $ \zeta = R(z_0)$ and  $\overline \zeta = R^*(z_0) $ for some $z_0 \in \overline \C$, with
$$R(z_0) = \overline{ R^*(z_0)} = R\bigl(j( z_0)\bigr) \quad {\rm and} \quad  R^*(z_0) = \overline{ R( z_0)} = R^*\bigl(j( z_0)\bigr).
$$
We know that outside a finite set, $R \times R^*$ is injective. Thus outside this finite set we obtain $z_0 = j( z_0)$, meaning that $z_0 \in \partial \DD$
and $\zeta = R(z_0) \in R(\partial \DD)$. 
\end{proof}
\smallskip

\begin{rem} \label{rem:cloud-curves}
The  identities \eqref{multisym1} - \eqref{multisym2} imply that the dual curve of \eqref{curve} is parametrised by $\bigl( 1/\gamma(z), 1/\gamma^*(z) \bigr)$, where $\gamma(z) = g(z) + b(z)\, \overline{g(z)}$ is the associated analytic factor. In case the curve has genus zero, i.e. $\LL$ is simply connected and $R(z)$ is rational,   the representation \eqref{gammafactor}  
 indicates that the real locus of the dual curve is a winding curve, in the sense of \cite{KeOk07}.
\end{rem}

The geometric properties of liquid boundaries we found earlier generalise quickly to the multiply connected case. 
In  a bounded and  multiply connected domain $\LL$ with piecewise smooth boundary we let $\tau$ denote the tangent field on $\partial \LL$,
with direction induced from the orientation of $\LL$. Thus $\tau(\zeta)$ has counter-clockwise direction on the outer component of $\partial \LL$ and clockwise direction on the interior components.

\begin{Thm}
\label{MultiplyConect2}
Suppose  $ \LL  \subset \C$ is a bounded domain supporting a solution $f$ to \eqref{once more13}, such that $f \colon  \LL   \to \Di$ is a proper map.
 Then 
 
a) $f(z)$ is real analytic  in $ \LL  $ and extends continuously up to the boundary $\partial  \LL  $.  

b) If $\tau$ is the tangent field on $\partial \LL$, then $f(\zeta) = - \tau(\zeta)^2, \; \zeta \in \partial \LL \setminus \{ {\rm cusps} \}$. 

c) In the complement of the cusps and tacnodes, $\partial \LL$ is locally strictly convex in the sense of  

\quad Theorem \ref{thm:geometry-curve}  ii). 

d) If $\, \partial \LL_k\, $, $k =1, \dots,m$ 
 are the components 
of  $\partial \LL$, then  each $\, \partial \LL_k$ has $| d_k - 2|$ 
cusps, where 

\quad $d_k$ is the degree of  $f: \partial \LL_k  \to \partial \Di$.

\end{Thm} 
\begin{proof} From Theorem \ref{propStoilow} we have $f = b \circ g^{-1}$, where $b: \DD \to \Di$ is analytic and proper,
$\DD = \Di \setminus \cup_{k=1}^\ell D(z_k, \delta_k)$ is a circle domain and the homeomorphism $g:\DD \to \LL$ solves \eqref{eqLinB22}.
 The argument is now basically the same as in Theorem \ref{thm:geometry-curve}. If $S_k = \partial D(z_k, \delta_k)$ and we are given a component $\partial \LL_k = g(S_k)$ of $\partial \LL$,  instead of  $\widehat \DD$  it is here convenient to use the reflection across $S_k$,
$$j_k(z) = z_k +  \frac{r_k^2}{\, \overline{ z } - \overline{ z_k } \,}, $$
 and  define $R^*(z) = \overline{(R \circ j_k)(z)}$ with $ b\bigl( j_k(z) \bigr) = 1 / \overline{b(z)}$. This again leads to $R'(z) = -b(z) (R^*)'(z)$ near $S_k$.  Arguing then as in Theorem \ref{thm:geometry-curve} proves the claims. We leave the details to the reader.
\end{proof}

In  particular, Theorem \ref{MultiplyConect2} tells that the bounded components of the complement $\C \setminus {\overline { \LL }}$ each have three or more outward  cusps.

As for simply connected domains in the earlier subsection, 
 Theorem \ref{MultiplyConect2}  (with Theorem \ref{MultiplyLocus})
  leads to the uniqueness of  proper maps $f:\LL \to \Di$ solving  the Beltrami equation \eqref{once more13}.

\begin{Cor} \label{uniquef} 
Let $\LL \subset \C$ be a bounded domain and assume $f_1, \, f_2 : \LL \to \Di$ are continuous proper maps, both solving the Beltrami equation $f_{\overline{ z }}(z) = f(z) \, f_z(z)$ in $\LL$.  Then  $ \, f_1 = f_2$.
\end{Cor}
\begin{proof} First  by Theorem  \ref{MultiplyConect2} b), the maps $f_1$ and $f_2$ agree on the boundary
$\partial \LL$. Second, Theorem \ref{propStoilow} gives us the factorisations $f_j = b_j \circ g_j^{-1}, $ where 
$g_j :\DD_j \to \LL$ is a homeomorphic solution to $\, \dbar g=- \, b_j(z) \overline{\, \dv_z g \,} \,$  and the
$b_j$ are analytic with $|b_j(z)| < 1$ in $\DD_j$. And third, the homeomorphisms have meromorphic boundary values,
$$  {g_j}_{\big|_{\partial \DD_j}} = {R_j}_{\big|_{\partial \DD_j}}, \qquad R_j \; \; {\rm meromorphic \; on\;  \;} \DD_j, \;\; j=1,2. 
$$

As in \eqref{rlocus} we have polynomials $P_j(z,w)$ with $P_j(R_j, R_j^*) = 0$, such that $\partial \LL = R_j(\partial \DD_j)$ is the real locus of $P_j$. Since  on the boundary, for  $z \in \partial \DD_j$, we have $P_k(R_j(z), R_j^*(z)) = 0$, the identity holds in all of $\C$ and we can take $P_1 = P_2 =: P$. 

The  curve $ {\mathcal C}$ in \eqref{curve} has now two parametrizations,  by $R_1 \times R_1^*$ and by $R_2 \times R_2^*$. In addition, we see from the proof of Theorem \ref{MultiplyLocus} that both are proper, i.e. injective  outside a finite set on the respective Schottky doubles $\widehat{\DD}_j$. As proper parametrizations are unique up to an automorphism, we see that $R_1 = R_2 \circ \psi$ for some conformal automorphism $\psi: \DD_1 \to \DD_2$.

On the other hand, the analytic coefficients $b_j \, ( = f_j \circ   g_j) $ have now the same boundary values up to the automorphism $\psi$, thus $b_1 = b_2 \circ \psi$.  Similarly, the analytic factors 
 $ \gamma_j(z) = R_j(z) + b_j(z) R_j^*(z)$ 
 in  \eqref{multisym1} are equal up to $\psi$, i.e. $\gamma_1 = \gamma_2 \circ \psi$. Finally, inserting these to the last identity in 
 \eqref{SchottkyFi12}  we see that $g_1 = g_2 \circ \psi$, and finally from the factorization that $f_1 = f_2$.
 \end{proof}

Since for any solution $f$ to \eqref{multimyy} the expression $f_0 =  \mu_{\sigma}(f) $ satisfies the universal Beltrami equation \eqref{once more13}, Corollary  \ref{uniquef} also implies

\begin{Cor} \label{uniqauto}
 The proper maps $f:\LL \to \DD$ solving  $f_{\overline{ z }} = \mu_{\sigma}(f) \, f_z$ are (if they exist)  unique  up to a conformal automorphism  
 preserving $\mu_{\sigma}$. 
\end{Cor}

\subsubsection{Proof of Theorem \ref{First.thm}.}

If $h$ is a $C^1$-solution to $\text{div}\, \big(\nabla \sigma(\nabla h)\big)= 0$ in a bounded domain $\LL$, and
$\, \nabla h : \LL \to N^\circ \setminus \Gg$ is  a proper map, then Theorem \ref{key.connection}
constructs from
$\nabla h$ a solution $f$ to the Beltrami equation \eqref{once more13}, such that $f:\LL \to \Di$ is proper. Thus the
 first two claims, a) and b) of Theorem \ref{First.thm}, follow from Theorems \ref{propStoilow35} and   \ref{MultiplyLocus}. The claim c) follows from  Theorem \ref{MultiplyConect2}. \hfill  $\Box$

\subsection{Boundary regularity of solutions and the Pokrovsky - Talapov Law} \label{sect:boundaryRegu}

We conclude this section with the boundary regularity properties for solutions to the Beltrami equations \eqref{multimyy} and \eqref{once more13}, 
as well as the for the minimizers of \eqref{ELbasic}, under the appropriate properness assumptions. Since now the question is about local properties, it is most convenient to discuss  the results in the setting of partially frozen boundaries.

Thus let $\LL$ be a bounded  domain,   and $f: \LL \to \DD$ a continuous $W^{1,2}_{loc}$-solution either to  the universal equation \eqref{once more13} (so that 
$\DD = \Di$ is the unit disc)  
or to  \eqref{multimyy}, with $\DD$ a general circle domain. 
Furthermore, we assume that there is a connected part  $\Gamma \subset \dv \LL$ with a smooth crosscut $\gamma$ of  $\LL$, such that $\gamma \cup \Gamma$  is the boundary of a simply connected domain contained in $\LL$. That is,  $\gamma$  is a simple arc of  $\LL$, with end points on  $\partial \LL$. In addition, of the map we assume that it is proper on $\Gamma$, which in this connection means that
\begin{equation} \label{partialfro}
f(z) \to \partial \DD \quad {\rm as } \quad  z \to \Gamma \quad {\rm in } \quad \LL.
\end{equation}

Furthermore, recall from Remark \ref{locfrozen} that if, for instance,  we have a solution to the equation
 $\, {\rm div} \, \big(\nabla \sigma(\nabla h)\big)= 0$ with 
 $$\nabla h(z) \to \partial N \cup \mathscr G  \qquad {as} \quad z \to \Gamma, \; z \in\LL,$$
 then \eqref{partialfro} holds for the map $f =  (\mathcal{H}_\sigma' \circ L_\sigma  \,) \circ \nabla h\, $, a solution  of \eqref{once more13}.

A prototype of \eqref{partialfro} is, naturally, the situation of Theorems \ref{MultiplyLocus} and \ref{MultiplyConect2}, where $f: \LL \to \Di$ is a proper map and $\Gamma \subset \partial \LL$ is a subarc of a component of $\partial \LL$. The general situation does not differ much from this. We first discuss the case where the target domain is the unit disc.

\begin{Thm}
\label{thm:localboundary}
Let  $f: \LL \to \Di$ be a continuous $W^{1,2}_{loc}$-solution the equation $f_{\overline{ z }} = f \, f_z$, 
and assume that  \eqref{partialfro} holds for an arc $\Gamma \subset \dv \LL$ as above. 

 Then $\Gamma$ is an analytic curve with at most finitely many singularities  $\{ \zeta_j \}$. The (possible) singularities  are all either first order cusps or tacnodes.   
 
 Furthermore, for $z_0 \in \Gamma \setminus \{ \zeta_j \}$ outside the cusps, $f \in C^{1/2}\bigl( B_\varepsilon(z_0) \cap \LL\bigr)$ where the H\"older exponent is optimal.
\end{Thm}
\begin{proof}

Call $\, \widetilde \LL \,$ the domain bounded by  $\gamma \cup \Gamma$. 
Then from Theorem \ref{propStoilow} we obtain  a homeomorphism $g: \Di \to \widetilde \LL\,$  and a holomorphic function $b:\Di \to \Di$ such that $f=b\circ g^{-1}$ in $\widetilde \LL$. In particular, $g$ is teleomorphic solving
\begin{align}\label{againeqn}
\dbar g(z)=-b(z)\, \overline{\dv g(z)}, \quad z\in \Di. 
\end{align}
If, say, in the sense of prime ends $I := g^{-1}(\Gamma) \subset \partial \Di$ then  \eqref{partialfro} tells that 
\begin{equation} \label{limit43}
|b(z)| \to 1, \quad {\rm  when} \; z \to I \quad {\rm  in} \; \Di.
\end{equation}
Consequently,  Theorems \ref{surprise2} and \ref{simplecuspthm} show that $\Gamma  = R(I)$ for a function $R$ analytic and locally injective on 
the interval. 

The  cusps of $\Gamma$ are simple by Theorem \ref{simplecuspthm}.  Moreover, the symmetries \eqref{symmetries} can be used, as in the proof of Theorem \ref{thm:geometry-curve}, to prove the local convexity of 
$\,\Gamma$, i.e.  that  for
$\zeta \in \Gamma$ outside the cusps and  
tacnodes,  $ B(\zeta,\varepsilon) \cap \LL$ is convex
 for $\varepsilon > 0$ small. This proves the above claims on the geometry of  $\Gamma$. 
 
 For the regularity of the map,  by  Proposition \ref{surprise1} the coefficient $b$ extends analytically  across $I$, and has no critical points on the interval. Thus  outside the critical points of $R$, for $z_0 \in \Gamma \setminus \{ \zeta_j \}$, the  factorisation $f = b \circ g^{-1}$  in $\widetilde \LL$ combined with  Proposition  \ref{bdryreg1}  shows that for $r> 0$ small enough
 $$ |f(z) - f(w)| \leq C|z-w|^{1/2}, \qquad z, w \in B_r(z_0) \cap \LL,
 $$ 
 where the H\"older exponent $1/2$  optimal.
 Finally, the continuity at possible tacnodes follows as in Theorem \ref{thm:geometry-curve}.\end{proof}
 
 \begin{Cor}
 \label{onesdided}
 Let  $f: \LL \to \DD$ be a solution to the equation 
 $f_{\overline{ z }} = \mu_\sigma(f) \, f_z$, 
and let the domain $\LL $ and $\Gamma \subset \dv \LL$  be  as in Theorem \ref{thm:localboundary}, so that \eqref{partialfro} holds. 

 Then all conclusions of Theorem \ref{thm:localboundary} remain true, except that at the possible tacnodes $f$ has  one-sided limits, but retains the H\"older continuity with exponent $1/2$.
 \end{Cor}
 
 \begin{proof}
For solutions $f$ to  \eqref{multimyy}, Theorem \ref{thm:localboundary} holds for  $\widehat f := \mu_\sigma \circ f$, which solves \eqref{once more13}. Since as an analytic and proper map $\mu_\sigma: \DD \to \Di$ does not have critical points on $\partial \DD$, we have $f \in C^{1/2}\bigl( B_\varepsilon(z_0) \cap \LL\bigr)$ outside the cusps and tacnodes. However, at possible tacnodes the argument gives only one-sided limits.
\end{proof}
\smallskip

\begin{rem} Under the assumptions of  Theorem \ref{thm:localboundary} we  avoid  a possible accumulation of gas components (or  other components of $\C \setminus \LL$) 
on $\Gamma $. Whether such an accumulation can really happen 
we leave as an open question. 

On the other hand, if $\Gamma$ is not connected, it is possible to have  accumulation points of cusps in the reference polygon $\Om$, even such  that for any ball $B_r(z_0)$  centered at an accumulation point $z_0$, $\mathcal{H}^1(B_r(z_0)\cap \dv \LL_F)=+\infty$. That such situations do occur can be seen from constructions in \cite{Duse15a}. 
\end{rem}
 
It can also happen that instead of \eqref{limit43} one has 
$I =\{z_0\}\in \dv \Di$  a singleton; for a concrete example of this see e.g.  \cite{Duse15a}, Proposition 4.7.

At the singularities the boundary behaviour is a little more complicated. Applying the argument above together  with Corollary \ref{gmap} one obtains the following explicit description.
\smallskip

\begin{Cor}
\label{bdryRegu3} Suppose the mapping $f: \LL \to \Di$,   domain $\LL $ and the connected part $\Gamma \subset \dv \LL$ are as in Theorem \ref{thm:localboundary}, with \eqref{partialfro} holding.
If $\, \zeta_j \in \Gamma$ is a cusp of $\Gamma$, then $f \in C^{1/3}\bigl( B_\varepsilon(\zeta_j) \cap \LL\bigr)$, where the H\"older exponent is optimal.   Moreover, 

a) There is a line $\ell$ transversal to the cusp at $\zeta_j$,  such that  
$$f \in C^{1/3}(\ell \cap B(\zeta_j,\varepsilon) ).$$

b)  However,  in the direction $\tau$ of the cusp, $f \in C^{1/2}(\tau \cap B(\zeta_j,\varepsilon) )$. 

\quad In particular, if $f: \LL \to \Di$ is a proper map solving Equation \eqref{once more13}, then $f \in C^{1/3}( \, \overline{ \LL \, }\, )$.
\end{Cor} 

 In view of  the  factorisation $f = b \circ g^{-1}$ as  in proof of Theorem \ref{thm:localboundary},  one sees from Figure  \ref{atcusp} how $f$ behaves at $\zeta_j$ in other than the cusp  direction.

\subsubsection{Proofs of  Theorems \ref{thm:main2}, \ref{GAS}, \ref{main.belt.regularity} and \ref{thm:localregboundary}} \label{proofPTalapov}

 \noindent{\it Proof of Theorem \ref{main.belt.regularity}}. On collecting the previous results,   claims  a) - c)  follow from 
Theorems \ref{MultiplyLocus} and  \ref{MultiplyConect2}. Claims d) and f) are given by Corollary \ref{bdryRegu3}. The last remaining claim e) is a consequence of   Theorem \ref{thm:localboundary}. \hfill $\Box$ 
\medskip

Theorem \ref{thm:localregboundary} is a part of Theorem \ref{thm:localboundary}.  However, Theorems \ref{GAS} and \ref{thm:main2}  require their details.
\medskip

\noindent{\it Proof of Theorem \ref{GAS}}.
Assume that the surface tension $\sigma$ has gas components, as defined  in (2.1). We  know that the  map  $  f \equiv L_\sigma \circ \nabla h :\LL \to {\rm Dom}(\mathcal H)$ is  proper,  by our assumption of $\partial \LL$ being frozen. But in addition, as a solution to the Beltrami equation \eqref{multipartic} $f$ is also a composition of an analytic function and a homeomorphism, c.f. Corollary  \ref{cor:propStoilow}.  
 Since the Lewy-transform is a homeomorphism, thus also $\nabla h:\LL \to N^\circ \setminus \mathscr{G}$ is  a proper, discrete and open mapping. 

Therefore, given any of the gas points $q_k \in \mathscr G \subset N^\circ$, there is a component of $\partial \LL$,
say $\partial \LL_k$, for which $\nabla h(z) \to q_k$ when $z \to \partial \LL_k$ in $\LL$.
   
   On the other hand, given a  bounded domain $W \subset \Omega$, if the  minimizer of \eqref{ELbasic} is affine on $\partial W$, then the minimizer is affine in all of $W$, since the surface tension $\sigma$ is convex. In particular, the above $\partial \LL_k$ must be the boundary of a bounded component $U_{q_k}$ of $\C \setminus {\overline{\LL}}$. Since $\Omega$ is bounded and simply connected we see that
   $U_q \subset \Omega$ and moreover, that $U_q$ has the  properties of a gas domain required by Theorem \ref{GAS}. \hfill $\Box$
   \medskip

\noindent{\it Proof of Theorem \ref{thm:main2}}.  Assume in a bounded  domain $\LL \subset \C$ the Euler-Lagrange equation \eqref{eq:EL34}
admits a solution $h$, 
$$ \text{div}\, \big(\nabla \sigma(\nabla h)\big)= 0  \; \; {\rm in} \; \; \LL, \quad {\rm such \; \;that } \quad  \nabla h(z) \to \partial N \cup \mathscr G  \; \; {\rm as} \; \; z \to \partial \LL.$$
 Here $\sigma$ is any surface tension as in \eqref{Pst2}. 

 Applying Remark \ref{myygen} we find then a circle domain $ \mathcal{D} = \Di \setminus \cup_{k=1}^\ell D(z_k, \delta_k)$ and a  proper map   
$f: \LL \to \DD$ solving  $f_{\overline{ z }} = \mu_{\sigma} (f) f_z$,
 tied together with $h$ by the relation
\begin{equation} \label{hrepre34}
\nabla h(z)=\sum_{j=1}^m \, p_j \, \omega_{\mathcal{D}}\bigl(\, f(z); I_j\bigr) +  \sum_{k=1}^\ell \, q_k \,  \omega_{\mathcal{D}}\bigl(\, f(z); S_k\bigr), \qquad z \in \LL.
\end{equation}
Here $p_j$ are the corners and quasifrozen points and $q_k$ the gas points in the gradient constraint $N$. Further, $S_k = \partial D(z_k, \delta_k)$, and the $I_j \subset \partial \Di$ are disjoint open arcs with union of their closures covering $\partial \Di$.

From this representation, and much of the previous results, Theorem \ref{thm:main2} readily follows. Indeed, for the
 first claim,  if  ${\mathcal I}_{j,\ell} \subset \partial \LL$  is a component of  $f^{-1}(I_j )$, then \eqref{hrepre34} gives 
   $$
\nabla h(z) \to p_j \quad {\rm as } \quad z \to z_0 \in {\mathcal I}_{j,\ell},
$$
while  at an endpoint of the arc ${\mathcal I}_{j,\ell}$ the gradient $\nabla h(z)$ fails to be continuous.

Similarly, for any component  ${\mathcal J}_{k,r} \subset \partial \LL $ of $f^{-1}( S_k )$, we see that  $\nabla h(z) \to q_k$ as $z \to z_0 \in {\mathcal J}_{k,r}$.
The union of the ${\mathcal I}_{j,\ell}$ and ${\mathcal J}_{k,r}$ covers $\partial \LL$ up to a finite set, which gives claim a) in
 Theorem \ref{thm:main2}.

Concerning  claim b), we give the proof only for corners and quasifrozen points $p_0 = p_j \in \mathscr{P}\bigcup\mathscr{Q} $; the case of gas points $p_0 = q_k \in \mathscr{G}$ is completely analogous, with Lemma \ref{subdiff12} replacing the role of Lemma \ref{subdiff}. 

For corners and quasifrozen points, first use Lemma \ref{subdiff} to see that  the conformal uniformisation $\psi: \DD \to  \text{Dom}(\mathcal{H}) $  extends analytically 
across each arc $I_j$. Next, since the solution $f$ in \eqref{hrepre34} satisfies 
\begin{equation} \label{decomp} 
L_\sigma \circ \nabla h = \psi \circ f,
\end{equation}
  Corollary \ref{onesdided} then shows that   $(\nabla \sigma \circ \nabla h)(z) $ has a definite limit, contained on the boundary of the subdifferential $\partial \sigma(p_j)$, when $z \to z_0 \in {\mathcal I}_{j,\ell}$. 
 Finally, according to 
 Theorem \ref{sigma1},
\begin{equation} \label{sigmalimit}
 \lim_{z\to z_0,z\in\LL}\nabla \sigma\big(\nabla h(z)\big)=\lim_{\tau\to 0^+}
\nabla\sigma\big(p_j+\tau(\hat p-p_j)\big),
\end{equation}
for some $\hat p \in N^\circ \setminus \Gg$. 

For b) it hence remains to identify the limiting direction $\hat p-p_j$ of the gradient $\nabla h(z)$ in  \eqref{sigmalimit}. Here 
use Lemma \ref{subdiff} to choose a smooth arc-length parametrisation $\zeta(t)$, $-\infty < t < \infty$, of the boundary of $\, \partial \sigma(p_j)$.    In particular, if $\widehat{\nabla}\sigma(p_j; \hat p-p_j) = \zeta(t_0)$, then according to  Proposition \ref{prop:sigma:property4} the tangent
    $\zeta'(t_0)$ is orthogonal to $\hat p-p_j$.
    
Write now the  identities \eqref{345kolme} in the form  ${\mathcal H} \left( {\overline{p + \nabla \sigma(p) }} \right) = p - \nabla \sigma(p)$. Using this at $p = p_j+\tau(\hat p-p_j)$ and  taking the limit with Theorem \ref{sigma1}, we  arrive at  $ {\mathcal H} \left( {\overline{ p_j + \zeta(t)}}\right) = p_j -  \zeta(t)$ for $  t \in \R.$ 
Differentiating this leads to 
    $$ {\mathcal H}' \left( {\overline{ p_j +\zeta(t)}}\right)  {\overline{ \zeta'(t)}}\,  = \, - \zeta'(t) \quad \Rightarrow \quad    {\mathcal H}'\left( {\overline{ p_j + \zeta(t)}}\right) \, = - \zeta'(t)^2.$$
As the last step, use \eqref{decomp} and recall that  
the map 
$$\widehat f := \mu_{\sigma}(f) = {\mathcal H}_\sigma'  \circ L_{\sigma}  \circ \nabla h,$$ is a solution to \eqref{once more13}. 
      Thus if  now  $\, z \to z_0 \in {\partial \mathcal L}$, $\, \nabla h(z) \to p_j\, $  
       and $\,  \nabla \sigma \circ \nabla h(z) \to \zeta(t_0)  \in  \partial \sigma(p_j)$,
   we see that $$\widehat f (z_0) =  {\mathcal H}'  \left(\, \overline{ p_j + \zeta(t_0)}\, \right) = - \zeta'(t_0)^2.$$  
   But from Theorem \ref{main.belt.regularity}
   we also know that  $\widehat f (z_0) = - \tau(z_0)^2$, where $\tau(z_0)$ is the tangent to the boundary of the liquid domain $\partial \LL$ at $z_0$. Thus the two tangents, $\zeta'(t_0)$ and $\tau(z_0)$, are parallel, so that 
  $\, \hat p - p_j $, the asymptotic direction of $\nabla h(z)$ as $z \to z_0$ in $\LL$, is orthogonal to both.
With this    claim b) in  Theorem \ref{thm:main2} follows. 

For c),  the boundary regularity of $h$ and the Pokrovsky-Talapov law, suppose $z_0 \in {\mathcal I}_{j,\ell} \subset  \partial \LL$  with $n_{z_0}$ the inner normal to $\LL$ at $z_0$.  Then
$$ h(z_0 + \delta \, n_{z_0}) - h(z_0) -    \delta \,  \langle  n_{z_0}, p_j \rangle \, = \int_0^\delta \langle n_{z_0},  \nabla h(z_0 + t n_{z_0})  - p_j \rangle dt.
$$ 
To estimate the term under the integral sign, first via \eqref{decomp}, 
$\nabla h = L_\sigma^{-1} \circ \psi \circ f  = U \circ f$, where  $U$ is a harmonic homeomorphism, c.f. \eqref{Lewyinverse}. 
On the other hand, Corollary \ref{cor:propStoilow} tells that $f = \eta \circ g^{-1}$ where $\eta$ is a proper analytic map between two circle domains $\DD'$ and $\DD$, thus has a non-vanishing derivative on the boundary. Then arguing as in \eqref{radial}, we see that the normal derivative $\partial_n (U \circ \eta) (g^{-1}(z_0) ) \neq 0$.

 Last, 
with the symmetry \eqref{circleg} the teleomorphic homeomorphism $g: \DD' \to  \LL$  has vanishing normal derivatives on $\partial \DD'$. 
As in Proposition \ref{bdryreg1}, the derivates $\partial_{\zeta} g(z_0)$ in  other directions are non-zero, however,  they all become  tangential on the boundary. Therefore $g^{-1}$ preserves the normal direction, it even  takes  a cone of directions around $n_{z_0}$ asymptotically  to the normal direction at $w_0 \in \partial \DD'$, $g(w_0) = z_0$. 
Thus as in Proposition \ref {bdryreg1} one has  $|g^{-1}(z_0 +t n_{z_0}) - g^{-1}(z_0)| \simeq t^{1/2}$. Combining the estimates  completes the proof of Theorem \ref{thm:main2}.  \hfill $\Box$

Combining with Corollary \ref{bdryRegu3} one can describe the behaviour of $\nabla h(z)$ also at a cusp of $\partial \LL$ (unless that happens to be an endpoints of some ${\mathcal I}_{j,\ell}$).

\addtocontents{toc}{\vspace{-4pt}}
\section{Frozen Extensions}\label{sect:properness}

We then turn to the  polygonal domains and the boundary values on them, that are  the natural  candidates
 for frozen phenomena, for the minimizers of \eqref{ELbasic}. As discussed in Subsection \ref{subsubsect:distinguished},
the simulation images suggest that one should consider the natural domains $\Omega\subset \R^2$,  as in Definition \ref{Def:naturaldomain}, 
and the natural boundary values $h_0$ on them, from  Definition \ref{def:extremal}.

The problem in the study of \eqref{ELbasic} is that with singular surface tensions $\sigma$ such as ours, there are no  general boundary extension methods available for the corresponding minimizers, and thus new approaches are needed. For this purpose, the goal of this section is to cover the basic properties of the {\it frozen extensions} as given in Definition \ref{oldold}. With these tools, Theorem   \ref{thm:h:proper3}  gives then  frozen boundaries for the corresponding liquid domains.

Let us start with the simplest case.

\begin{Lem}\label{thm:frozenextension1}
Let $\Omega\subset\R^2$ be a  natural domain for a closed convex polygon $N$, as in Definition \ref{Def:naturaldomain}, and assume in addition that $\Omega$ is convex. If $h_0$ is a natural boundary value, then $(\Omega, h_0)$ admits a frozen extension at any point $z_0\in \partial \Omega$.
\end{Lem}
\begin{proof} Let $\Omega$ be a convex natural domain with $d$ vertices $\{ z_1,...,z_d\}$, set $z_{d+1}=z_1$ and let $h_0$ a natural boundary value.

Suppose either $z_0 \in \partial \Omega$ is a corner, $z_0= z_j$, or it lies on the open interval $(z_j,z_{j+1})$.
Choose then  a parallelogram  $\mathcal{P}$ such that $z_0\in \mathcal P^\circ$  and  that one pair of the sides of $\mathcal{P}$  is parallel to $z_j-z_{j-1}$ and the other is parallel  to $z_{j+1}-z_j$, c.f. Figure \ref{fig:ExtL2} below. We choose $\mathcal{P}$  so small that it intersects only one side of $\partial \Omega$ if $z_0\in (z_j,z_{j+1})$ and two sides if $z_0=z_j$. 
This allows us to define a new domain $\widehat{\Om}=P \cup \Om$. 

Since $\Omega$ is natural and also $h_0$ is natural, condition (\ref{zjpn}) tells that
\[
 h_0(z)=\langle p_n, z-z_j\rangle+h_0(z_j)
\]
for all $z$ in the line segments $[z_{j-1},z_{j}]$ and $[z_{j},z_{j+1}]$. Here 
$p_n, n=1,\dots, k$, is the vertex such that (\ref{ND1})  or (\ref{ND3})  holds.
Now we define a boundary height function $\widehat{h}_0$ on $\partial \widehat{\Omega}$ as follows: Set $\widehat{h}_0(z)=h_0(z)$ when $z\in \partial \widehat{\Omega} \cap \partial\Omega$, and 
\begin{equation}\label{Hhat}
\widehat{h}_0(z)=\langle p_n, z-z_j\rangle+h_0(z_j)
\end{equation}
when $z\in \partial \widehat{\Omega}\setminus\partial \Omega$.
Our task is then to show that the upper and lower obstacles in the enlarged domain $\widehat{\Omega}$ satisfy
\begin{equation}\label{Mhat}
\widehat{M}(z)=\widehat{m}(z)=\langle p_n, z-z_j\rangle+h_0(z_j), \qquad \forall \;  z\in \widehat{\Omega}\setminus \overline{\Omega}.
\end{equation}
The proof of (\ref{Mhat}) is easy. Indeed, 
consider a point $z\in  \widehat{\Om}\setminus  \Om$. Then  $z$ intersects a line $\ell_1$ parallel to $z_{j+1}-z_j$ (or $\ell_2$ parallel to $z_{j-1}-z_j$) that intersects $\dv \widehat{\Om}\setminus \dv \Om$ at two points $w_1$ and $w_2$, respectively. See Figure \ref{fig:ExtL2}.
\begin{figure}[H]
\centering
\begin{tikzpicture}[xscale=2,yscale=2]
\draw (-1,0)--(2,0);
\draw[xshift=-1.75cm,yshift=-0.5cm] (0,0)--(1,0)--(1.5,1)--(0.5,1)--(0,0);
\draw(2,0) node{$\bullet$};
\draw(-1,0) node{$\bullet$};
\draw (-1,0)--(-1.5,-1);
\draw(-1.5,-1) node{$\bullet$};
\draw(-1.7,-1) node{$z_{j-1}$};
\draw(2.2,0) node{$z_{j+1}$};
\draw(-0.95,-0.1) node{$z_{j}$};
\draw(-0.8,0) node{$\bullet$};
\draw(-0.8,0.1) node{$z_0$};
\draw(1,-0.5) node{$ \Om$};
\draw(-1.2,0.25) node{$\bullet$};
\draw(-1.1,0.2) node{$z$};
\draw[thick,dashed] (-2,0.25)--(1,0.25);
\draw[thick,dashed] (-1.7,-0.75)--(-0.95,0.75);
\draw(1,0.4) node{$ \ell_1$};
\draw(-0.85,0.75) node{$ \ell_2$};
\draw(-1.4,0.55) node{$\mathcal{P}$};
\end{tikzpicture}
\caption{Frozen extension}
\label{fig:ExtL2}
\end{figure}
Now either (\ref{ND1}) or (\ref{ND3}) holds. In each case
 we see from \eqref{struct.supportfcn}, that  either
\begin{equation}\label{fun1}
h_N(z-w_1)=\langle p_n, z-w_1\rangle, \quad h_N(w_2-z)=\langle p_n, w_2-z\rangle,
\end{equation}
or 
\begin{equation}\label{fun2}
h_N(z-w_2)=\langle p_n, z-w_2\rangle, \quad h_N(w_1-z)=\langle p_n, w_1-z\rangle,
\end{equation}
If (\ref{fun1}) holds, we have that
\[
\widehat{M}(z)=\min_{w\in \partial \widehat{\Omega}}\big(h_N(z-w)+\widehat{h}_0(w)\big)\leq h_N(z-w_1)+\widehat{h}_0(w_1)
=\langle p_{n}, z-z_j\rangle+h_0(z_j),
\]
and 
\[
\widehat{m}(z)=\max_{w\in \partial \widehat{\Omega}  }\big(-h_N(w-z)+\widehat{h}_0(w)\big)\geq -h_N(w_2-z)+\widehat{h}_0(w_2)=\langle p_{n}, z-z_j\rangle+h_0(z_j).
\]
This proves the equalities (\ref{Mhat}), and the proof in case (\ref{fun2}) is analogous. 
\end{proof}

If the natural domain $\Omega$ has concave corners, establishing frozen extensions  becomes more complicated.
Namely, assume that $z_0=z_j$ is a concave corner with respect to $\Om$
for an index $j=1,2,...,d$, and assume  (\ref{ND1}) holds for the corner $p_n$ of $N$. The other possibility (\ref{ND3}) can be analysed similarly.

Under these assumptions 
we have four different possible geometric configurations for the sides of $N$ and $\Omega$ at $p_n$ and $z_j$, respectively, as illustrated  in the following Figures.

{\begin{figure}[H]
  \centering
  \begin{minipage}[b]{0.4\textwidth}

\begin{tikzpicture}[xscale=1.5,yscale=1.5]
\draw[thick] (0,0)--(3,0);
\draw[thick] (0,0)--(2,2);
\filldraw(0,0)  circle (0.6pt);
\filldraw(3,0)  circle (0.6pt);
\filldraw(2,2)  circle (0.6pt);
\draw[thick] (1.5,2)--(1.5,0.5);
\draw[thick] (1.5,0.5)--(0.5,1.5);
\filldraw(1.5,2)  circle (0.6pt);
\filldraw(1.5,0.5)  circle (0.6pt);
\filldraw(0.5,1.5)  circle (0.6pt);
\draw(2.3,2) node{$z_{j-1}$};
\draw(0,-0.2) node{$z_j$};
\draw(3,-0.2) node{$z_{j+1}$};
\draw(-0.5,0) node{$\Omega$};
\draw(1.5,2.2) node{$p_{n+1}$};
\draw(1.5,0.3)node{$p_{n}$};
\draw(0.2,1.6)node{$p_{n-1}$};
\end{tikzpicture}

    \caption*{$(i)$} \label{one1}
  \end{minipage}
  \hfill
  \begin{minipage}[b]{0.4\textwidth}
 \begin{tikzpicture}[xscale=1.5,yscale=1.5]
\draw[thick] (0,1)--(3,1);
\draw[thick] (3,1)--(1,-1);
\filldraw(0,1)  circle (0.6pt);
\filldraw(3,1)  circle (0.6pt);
\filldraw(1,-1)  circle (0.6pt);
\draw[thick] (1.5,1.8)--(1.5,0.5);
\draw[thick] (1.5,0.5)--(0.5,1.5);
\filldraw(1.5,1.8)  circle (0.6pt);
\filldraw(1.5,0.5)  circle (0.6pt);
\filldraw(0.5,1.5)  circle (0.6pt);
\draw(-0.4,1) node{$z_{j+1}$};
\draw(3.3,1) node{$z_{j}$};
\draw(1,-1.2) node{$z_{j-1}$};
\draw(2.5,1.5) node{$\Omega$};
\draw(1.5,2) node{$p_{n+1}$};
\draw(1.5,0.3)node{$p_{n}$};
\draw(0.2,1.6)node{$p_{n-1}$};
\end{tikzpicture}
    \caption*{$(ii)$}
  \end{minipage}
  
   \end{figure}
  
   \begin{figure}[H]
  \vfill
  \begin{minipage}[b]{0.4\textwidth}
  
\begin{tikzpicture}[xscale=1.5,yscale=1.5]
\draw[thick] (0,0)--(3,0);
\draw[thick] (0,0)--(-1,-1);
\filldraw(0,0)  circle (0.6pt);
\filldraw(3,0)  circle (0.6pt);
\filldraw(-1,-1)  circle (0.6pt);
\draw[thick] (1.5,1)--(1.5,-0.5);
\draw[thick] (1.5,-0.5)--(0.5,0.5);
\filldraw(1.5,1)  circle (0.6pt);
\filldraw(1.5,-0.5)  circle (0.6pt);
\filldraw(0.5,0.5)  circle (0.6pt);
\draw(-0.4,0) node{$z_{j}$};
\draw(3.3,0) node{$z_{j+1}$};
\draw(-1,-1.2) node{$z_{j-1}$};
\draw(-1.5,-0.5) node{$\Omega$};
\draw(1.5,1.2) node{$p_{n+1}$};
\draw(1.5,-0.7)node{$p_{n}$};
\draw(0.2,0.6)node{$p_{n-1}$};
\end{tikzpicture}
    \caption*{$(iii)$}
  \end{minipage}
  \hfill
  \begin{minipage}[b]{0.4\textwidth}
  
\begin{tikzpicture}[xscale=1.5,yscale=1.5]
\draw[thick] (0,1)--(3,1);
\draw[thick] (3,1)--(4,2);
\filldraw(0,1)  circle (0.6pt);
\filldraw(3,1)  circle (0.6pt);
\filldraw(4,2)  circle (0.6pt);
\draw[thick] (1.5,2)--(1.5,0.5);
\draw[thick] (1.5,0.5)--(0.5,1.5);
\filldraw(1.5,2)  circle (0.6pt);
\filldraw(1.5,0.5)  circle (0.6pt);
\filldraw(0.5,1.5)  circle (0.6pt);
\draw(-0.4,1) node{$z_{j+1}$};
\draw(3.3,1) node{$z_{j}$};
\draw(4.2,2.2) node{$z_{j-1}$};
\draw(2,0.5) node{$\Omega$};
\draw(1.5,2.2) node{$p_{n+1}$};
\draw(1.5,0.3)node{$p_{n}$};
\draw(0.2,1.6)node{$p_{n-1}$};
\end{tikzpicture}
    \caption*{$(iv)$}
  \end{minipage}

  \end{figure} }

A version of the argument of Lemma \ref{thm:frozenextension1} works in concave corners only if the gradient constraint and the sides of $\Omega$ have the same ``orientation'', meaning 
that only the cases (i) and (ii) in the above figure can occur. A formal definition is  as follows.

\begin{Def}\label{def:extremalorient}
Let $\Omega$ be a natural domain with
$d$ vertices $\{ z_1,...,z_d\}$.
We say that a natural boundary value $h_0:\partial\Om\to \mathbb R$ is oriented if for any non-convex corner $z_j$ of $\dv \Om$ either we have 
\begin{equation}\label{def:extorient1}
\frac{(p_{n-1}-p_{n})\wedge (z_{j-1}-z_{j})}{\vert (p_{n-1}-p_{n})\wedge (z_{j-1}-z_{j})\vert}=\frac{(p_{n+1}-p_{n})\wedge (z_{j+1}-z_{j})}{\vert (p_{n+1}-p_{n})\wedge (z_{j+1}-z_{j})\vert},
\end{equation}
where $p_n$ is the vertex such that (\ref{ND1}) and (\ref{zjpn}) hold, or we have that
\begin{equation}\label{def:extorient2}
\frac{(p_{n-1}-p_{n})\wedge (z_{j+1}-z_{j})}{\vert (p_{n-1}-p_{n})\wedge (z_{j+1}-z_{j})\vert}=\frac{(p_{n+1}-p_{n})\wedge (z_{j-1}-z_{j})}{\vert (p_{n+1}-p_{n})\wedge (z_{j-1}-z_{j})\vert}.
\end{equation}
when (\ref{ND3}) and (\ref{zjpn}) hold. Here $\wedge$ denotes the exterior product in the exterior algebra $\Lambda \R^2$. 
\end{Def}

We posed the  above restrictions on the natural boundary values for technical reasons, for our argument to work, but it is curious to notice that  also in the simulations
oriented and non-oriented boundary values have qualitatively different behaviours, in the neighbourhood of the concave corners. For this see below in Figure \ref{figConcaveDifferent} two simulations of random domino tilings on polygonal domains, where the left simulation has two concave corners which are {\it not} oriented, while in the right figure concave corners  are oriented, in the sense of Definition \ref{def:extremalorient}.  In particular one notices that the topology of the respective liquid domains appears different.  

\begin{figure}[H] 
\centering{}
\includegraphics[scale=0.55]{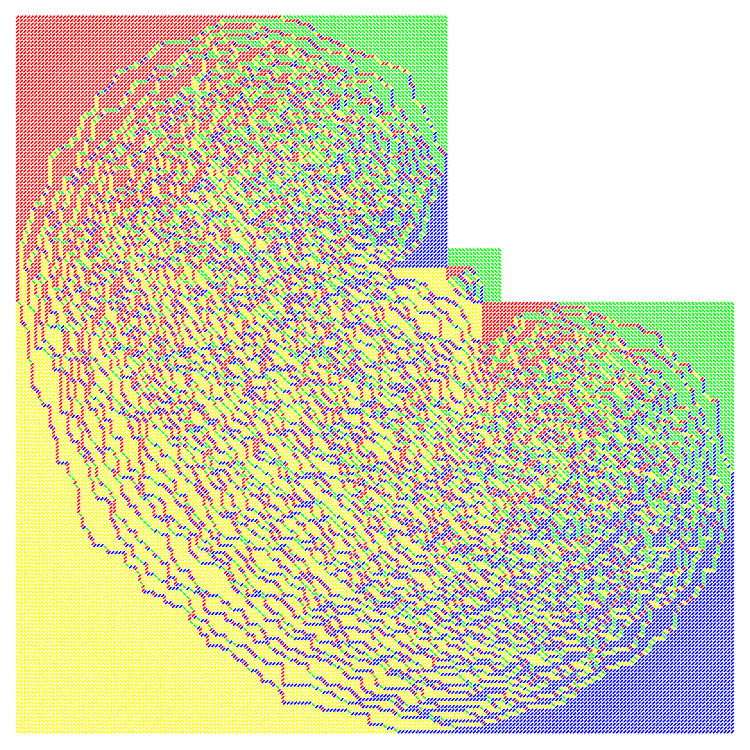} 
\includegraphics[scale=0.55]{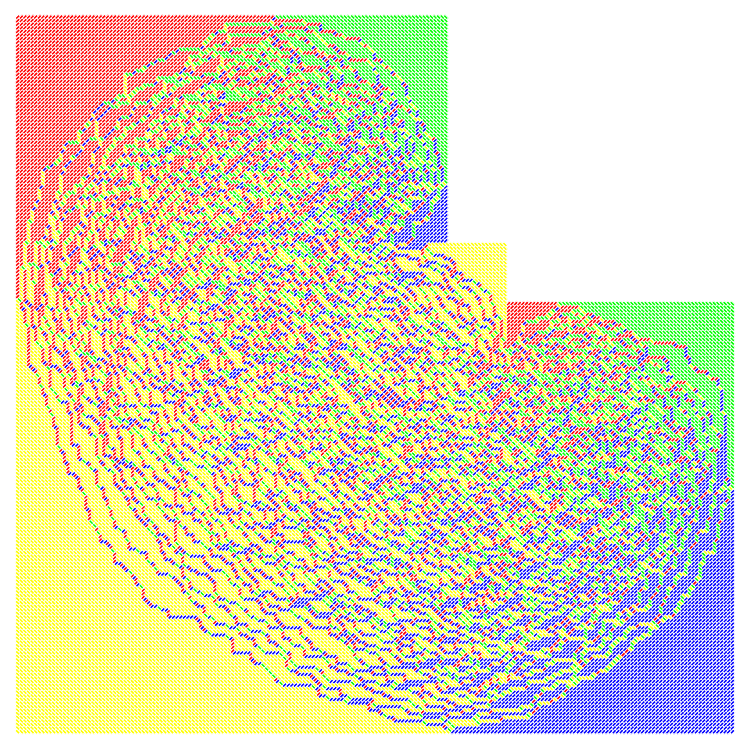}
\caption{Simulation  of random domino tilings on polygonal domains; on left two non-oriented concave corners,   on  right  all corners are oriented.
 Simulation image courtesy of Sunil Chhita.}
\label{figConcaveDifferent}
\end{figure}

 One observes already on the 
microscopic level the different behaviours of the oriented and non-oriented boundary values at concave corners, in constructing them as in  Example \ref{microdominoes} from tilings with the simplest discrete boundary values. For an illustration see  Figure
 \ref{kolme} below.

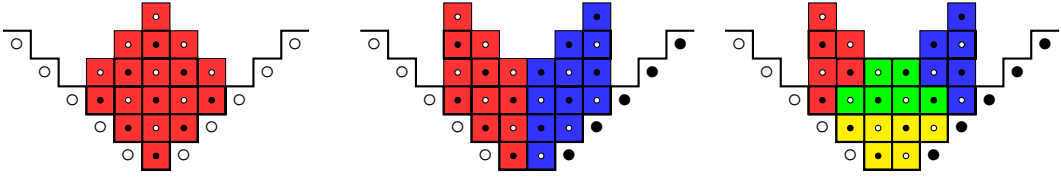
\begin{figure}[H]
  \centering
  \begin{minipage}[b]{0.4\textwidth}
    \begin{tikzpicture}[xscale=0.37,yscale=0.37]
\draw[thick] (0,0)--(1,0)--(1,-1)--(2,-1)--(2,-2)--(3,-2)--(3,-3)--(4,-3)--(4,-4)--(5,-4)--(5,-5)--(6,-5)
--(6,-4)--(7,-4)--(7,-3)--(8,-3)--(8,-2)--(9,-2)--(9,-1)--(10,-1)--(10,0)--(11,0);

\filldraw[xshift=5cm,yshift=-5cm,fill=red!80!white] (0,0)--(1,0)--(1,2)--(0,2)--(0,0);
\draw[xshift=5cm,yshift=-5cm,thick] (0,0)--(1,0)--(1,1)--(0,1)--(0,0)--(1,0)--(1,1);
\filldraw[xshift=5cm,yshift=-5cm,fill=black] (0.5,0.5) circle [radius=0.1cm];
\filldraw[xshift=5cm,yshift=-5cm,fill=white] (0.5,1.5) circle [radius=0.1cm];

\filldraw[xshift=4cm,yshift=-4cm,fill=red!80!white] (0,0)--(1,0)--(1,2)--(0,2)--(0,0);
\draw[xshift=4cm,yshift=-4cm,thick] (0,0)--(1,0)--(1,1)--(0,1)--(0,0)--(1,0)--(1,1);
\filldraw[xshift=4cm,yshift=-4cm,fill=black] (0.5,0.5) circle [radius=0.1cm];
\filldraw[xshift=4cm,yshift=-4cm,fill=white] (0.5,1.5) circle [radius=0.1cm];

\filldraw[xshift=6cm,yshift=-4cm,fill=red!80!white] (0,0)--(1,0)--(1,2)--(0,2)--(0,0);
\draw[xshift=6cm,yshift=-4cm,thick] (0,0)--(1,0)--(1,1)--(0,1)--(0,0)--(1,0)--(1,1);
\filldraw[xshift=6cm,yshift=-4cm,fill=black] (0.5,0.5) circle [radius=0.1cm];
\filldraw[xshift=6cm,yshift=-4cm,fill=white] (0.5,1.5) circle [radius=0.1cm];

\filldraw[xshift=3cm,yshift=-3cm,fill=red!80!white] (0,0)--(1,0)--(1,2)--(0,2)--(0,0);
\draw[xshift=3cm,yshift=-3cm,thick] (0,0)--(1,0)--(1,1)--(0,1)--(0,0)--(1,0)--(1,1);
\filldraw[xshift=3cm,yshift=-3cm,fill=black] (0.5,0.5) circle [radius=0.1cm];
\filldraw[xshift=3cm,yshift=-3cm,fill=white] (0.5,1.5) circle [radius=0.1cm];

\filldraw[xshift=7cm,yshift=-3cm,fill=red!80!white] (0,0)--(1,0)--(1,2)--(0,2)--(0,0);
\draw[xshift=7cm,yshift=-3cm,thick] (0,0)--(1,0)--(1,1)--(0,1)--(0,0)--(1,0)--(1,1);
\filldraw[xshift=7cm,yshift=-3cm,fill=black] (0.5,0.5) circle [radius=0.1cm];
\filldraw[xshift=7cm,yshift=-3cm,fill=white] (0.5,1.5) circle [radius=0.1cm];

\filldraw[xshift=5cm,yshift=-3cm,fill=red!80!white] (0,0)--(1,0)--(1,2)--(0,2)--(0,0);
\draw[xshift=5cm,yshift=-3cm,thick] (0,0)--(1,0)--(1,1)--(0,1)--(0,0)--(1,0)--(1,1);
\filldraw[xshift=5cm,yshift=-3cm,fill=black] (0.5,0.5) circle [radius=0.1cm];
\filldraw[xshift=5cm,yshift=-3cm,fill=white] (0.5,1.5) circle [radius=0.1cm];

\filldraw[xshift=4cm,yshift=-2cm,fill=red!80!white] (0,0)--(1,0)--(1,2)--(0,2)--(0,0);
\draw[xshift=4cm,yshift=-2cm,thick] (0,0)--(1,0)--(1,1)--(0,1)--(0,0)--(1,0)--(1,1);
\filldraw[xshift=4cm,yshift=-2cm,fill=black] (0.5,0.5) circle [radius=0.1cm];
\filldraw[xshift=4cm,yshift=-2cm,fill=white] (0.5,1.5) circle [radius=0.1cm];

\filldraw[xshift=6cm,yshift=-2cm,fill=red!80!white] (0,0)--(1,0)--(1,2)--(0,2)--(0,0);
\draw[xshift=6cm,yshift=-2cm,thick] (0,0)--(1,0)--(1,1)--(0,1)--(0,0)--(1,0)--(1,1);
\filldraw[xshift=6cm,yshift=-2cm,fill=black] (0.5,0.5) circle [radius=0.1cm];
\filldraw[xshift=6cm,yshift=-2cm,fill=white] (0.5,1.5) circle [radius=0.1cm];

\filldraw[xshift=5cm,yshift=-1cm,fill=red!80!white] (0,0)--(1,0)--(1,2)--(0,2)--(0,0);
\draw[xshift=5cm,yshift=-1cm,thick] (0,0)--(1,0)--(1,1)--(0,1)--(0,0)--(1,0)--(1,1);
\filldraw[xshift=5cm,yshift=-1cm,fill=black] (0.5,0.5) circle [radius=0.1cm];
\filldraw[xshift=5cm,yshift=-1cm,fill=white] (0.5,1.5) circle [radius=0.1cm];

\draw (0.5,-0.5) node {$\circ$};
\draw (1.5,-1.5) node {$\circ$};
\draw (2.5,-2.5) node {$\circ$};
\draw (3.5,-3.5) node {$\circ$};
\draw (4.5,-4.5) node {$\circ$};
\draw (6.5,-4.5) node {$\circ$};
\draw (7.5,-3.5) node {$\circ$};
\draw (8.5,-2.5) node {$\circ$};
\draw (9.5,-1.5) node {$\circ$};
\draw (10.5,-0.5) node {$\circ$};

\end{tikzpicture}

  \end{minipage}

 \hspace{-1.9cm}
  \begin{minipage}[b]{0.4\textwidth}
  \begin{tikzpicture}[xscale=0.37,yscale=0.37]
\draw[thick] (0,0)--(1,0)--(1,-1)--(2,-1)--(2,-2)--(3,-2)--(3,-3)--(4,-3)--(4,-4)--(5,-4)--(5,-5)--(7,-5)
--(7,-4)--(8,-4)--(8,-3)--(9,-3)--(9,-2)--(10,-2)--(10,-1)--(11,-1)--(11,0)--(12,0);
\draw (0.5,-0.5) node {$\circ$};
\draw (1.5,-1.5) node {$\circ$};
\draw (2.5,-2.5) node {$\circ$};
\draw (3.5,-3.5) node {$\circ$};
\draw (4.5,-4.5) node {$\circ$};
\draw (7.5,-4.5) node {$\bullet$};
\draw (8.5,-3.5) node {$\bullet$};
\draw (9.5,-2.5) node {$\bullet$};
\draw (10.5,-1.5) node {$\bullet$};
\draw (11.5,-0.5) node {$\bullet$};


\filldraw[xshift=5cm,yshift=-5cm,fill=red!80!white] (0,0)--(1,0)--(1,2)--(0,2)--(0,0);
\draw[xshift=5cm,yshift=-5cm,thick] (0,0)--(1,0)--(1,1)--(0,1)--(0,0)--(1,0)--(1,1);
\filldraw[xshift=5cm,yshift=-5cm,fill=black] (0.5,0.5) circle [radius=0.1cm];
\filldraw[xshift=5cm,yshift=-5cm,fill=white] (0.5,1.5) circle [radius=0.1cm];

\filldraw[xshift=4cm,yshift=-4cm,fill=red!80!white] (0,0)--(1,0)--(1,2)--(0,2)--(0,0);
\draw[xshift=4cm,yshift=-4cm,thick] (0,0)--(1,0)--(1,1)--(0,1)--(0,0)--(1,0)--(1,1);
\filldraw[xshift=4cm,yshift=-4cm,fill=black] (0.5,0.5) circle [radius=0.1cm];
\filldraw[xshift=4cm,yshift=-4cm,fill=white] (0.5,1.5) circle [radius=0.1cm];

\filldraw[xshift=3cm,yshift=-3cm,fill=red!80!white] (0,0)--(1,0)--(1,2)--(0,2)--(0,0);
\draw[xshift=3cm,yshift=-3cm,thick] (0,0)--(1,0)--(1,1)--(0,1)--(0,0)--(1,0)--(1,1);
\filldraw[xshift=3cm,yshift=-3cm,fill=black] (0.5,0.5) circle [radius=0.1cm];
\filldraw[xshift=3cm,yshift=-3cm,fill=white] (0.5,1.5) circle [radius=0.1cm];

\filldraw[xshift=3cm,yshift=-1cm,fill=red!80!white] (0,0)--(1,0)--(1,2)--(0,2)--(0,0);
\draw[xshift=3cm,yshift=-1cm,thick] (0,0)--(1,0)--(1,1)--(0,1)--(0,0)--(1,0)--(1,1);
\filldraw[xshift=3cm,yshift=-1cm,fill=black] (0.5,0.5) circle [radius=0.1cm];
\filldraw[xshift=3cm,yshift=-1cm,fill=white] (0.5,1.5) circle [radius=0.1cm];

\filldraw[xshift=4cm,yshift=-2cm,fill=red!80!white] (0,0)--(1,0)--(1,2)--(0,2)--(0,0);
\draw[xshift=4cm,yshift=-2cm,thick] (0,0)--(1,0)--(1,1)--(0,1)--(0,0)--(1,0)--(1,1);
\filldraw[xshift=4cm,yshift=-2cm,fill=black] (0.5,0.5) circle [radius=0.1cm];
\filldraw[xshift=4cm,yshift=-2cm,fill=white] (0.5,1.5) circle [radius=0.1cm];

\filldraw[xshift=5cm,yshift=-3cm,fill=red!80!white] (0,0)--(1,0)--(1,2)--(0,2)--(0,0);
\draw[xshift=5cm,yshift=-3cm,thick] (0,0)--(1,0)--(1,1)--(0,1)--(0,0)--(1,0)--(1,1);
\filldraw[xshift=5cm,yshift=-3cm,fill=black] (0.5,0.5) circle [radius=0.1cm];
\filldraw[xshift=5cm,yshift=-3cm,fill=white] (0.5,1.5) circle [radius=0.1cm];


\filldraw[xshift=6cm,yshift=-5cm,fill=blue!80!white] (0,0)--(1,0)--(1,2)--(0,2)--(0,0);
\draw[xshift=6cm,yshift=-5cm,thick] (0,0)--(1,0)--(1,1)--(0,1)--(0,0)--(1,0)--(1,1);
\filldraw[xshift=6cm,yshift=-5cm,fill=white] (0.5,0.5) circle [radius=0.1cm];
\filldraw[xshift=6cm,yshift=-5cm,fill=black] (0.5,1.5) circle [radius=0.1cm];

\filldraw[xshift=7cm,yshift=-4cm,fill=blue!80!white] (0,0)--(1,0)--(1,2)--(0,2)--(0,0);
\draw[xshift=7cm,yshift=-4cm,thick] (0,0)--(1,0)--(1,1)--(0,1)--(0,0)--(1,0)--(1,1);
\filldraw[xshift=7cm,yshift=-4cm,fill=white] (0.5,0.5) circle [radius=0.1cm];
\filldraw[xshift=7cm,yshift=-4cm,fill=black] (0.5,1.5) circle [radius=0.1cm];

\filldraw[xshift=8cm,yshift=-3cm,fill=blue!80!white] (0,0)--(1,0)--(1,2)--(0,2)--(0,0);
\draw[xshift=8cm,yshift=-3cm,thick] (0,0)--(1,0)--(1,1)--(0,1)--(0,0)--(1,0)--(1,1);
\filldraw[xshift=8cm,yshift=-3cm,fill=white] (0.5,0.5) circle [radius=0.1cm];
\filldraw[xshift=8cm,yshift=-3cm,fill=black] (0.5,1.5) circle [radius=0.1cm];

\filldraw[xshift=6cm,yshift=-3cm,fill=blue!80!white] (0,0)--(1,0)--(1,2)--(0,2)--(0,0);
\draw[xshift=6cm,yshift=-3cm,thick] (0,0)--(1,0)--(1,1)--(0,1)--(0,0)--(1,0)--(1,1);
\filldraw[xshift=6cm,yshift=-3cm,fill=white] (0.5,0.5) circle [radius=0.1cm];
\filldraw[xshift=6cm,yshift=-3cm,fill=black] (0.5,1.5) circle [radius=0.1cm];

\filldraw[xshift=7cm,yshift=-2cm,fill=blue!80!white] (0,0)--(1,0)--(1,2)--(0,2)--(0,0);
\draw[xshift=7cm,yshift=-2cm,thick] (0,0)--(1,0)--(1,1)--(0,1)--(0,0)--(1,0)--(1,1);
\filldraw[xshift=7cm,yshift=-2cm,fill=white] (0.5,0.5) circle [radius=0.1cm];
\filldraw[xshift=7cm,yshift=-2cm,fill=black] (0.5,1.5) circle [radius=0.1cm];

\filldraw[xshift=8cm,yshift=-1cm,fill=blue!80!white] (0,0)--(1,0)--(1,2)--(0,2)--(0,0);
\draw[xshift=8cm,yshift=-1cm,thick] (0,0)--(1,0)--(1,1)--(0,1)--(0,0)--(1,0)--(1,1);
\filldraw[xshift=8cm,yshift=-1cm,fill=white] (0.5,0.5) circle [radius=0.1cm];
\filldraw[xshift=8cm,yshift=-1cm,fill=black] (0.5,1.5) circle [radius=0.1cm];
\end{tikzpicture}

  \end{minipage}
   \hspace{-1.8cm}
  \begin{minipage}[b]{0.4\textwidth}
    \begin{tikzpicture}[xscale=0.37,yscale=0.37]
\draw[thick] (0,0)--(1,0)--(1,-1)--(2,-1)--(2,-2)--(3,-2)--(3,-3)--(4,-3)--(4,-4)--(5,-4)--(5,-5)--(7,-5)
--(7,-4)--(8,-4)--(8,-3)--(9,-3)--(9,-2)--(10,-2)--(10,-1)--(11,-1)--(11,0)--(12,0);
\draw (0.5,-0.5) node {$\circ$};
\draw (1.5,-1.5) node {$\circ$};
\draw (2.5,-2.5) node {$\circ$};
\draw (3.5,-3.5) node {$\circ$};
\draw (4.5,-4.5) node {$\circ$};
\draw (7.5,-4.5) node {$\bullet$};
\draw (8.5,-3.5) node {$\bullet$};
\draw (9.5,-2.5) node {$\bullet$};
\draw (10.5,-1.5) node {$\bullet$};
\draw (11.5,-0.5) node {$\bullet$};


\filldraw[xshift=3cm,yshift=-1cm,fill=red!80!white] (0,0)--(1,0)--(1,2)--(0,2)--(0,0);
\draw[xshift=3cm,yshift=-1cm,thick] (0,0)--(1,0)--(1,1)--(0,1)--(0,0)--(1,0)--(1,1);
\filldraw[xshift=3cm,yshift=-1cm,fill=black] (0.5,0.5) circle [radius=0.1cm];
\filldraw[xshift=3cm,yshift=-1cm,fill=white] (0.5,1.5) circle [radius=0.1cm];

\filldraw[xshift=4cm,yshift=-2cm,fill=red!80!white] (0,0)--(1,0)--(1,2)--(0,2)--(0,0);
\draw[xshift=4cm,yshift=-2cm,thick] (0,0)--(1,0)--(1,1)--(0,1)--(0,0)--(1,0)--(1,1);
\filldraw[xshift=4cm,yshift=-2cm,fill=black] (0.5,0.5) circle [radius=0.1cm];
\filldraw[xshift=4cm,yshift=-2cm,fill=white] (0.5,1.5) circle [radius=0.1cm];

\filldraw[xshift=3cm,yshift=-3cm,fill=red!80!white] (0,0)--(1,0)--(1,2)--(0,2)--(0,0);
\draw[xshift=3cm,yshift=-3cm,thick] (0,0)--(1,0)--(1,1)--(0,1)--(0,0)--(1,0)--(1,1);
\filldraw[xshift=3cm,yshift=-3cm,fill=black] (0.5,0.5) circle [radius=0.1cm];
\filldraw[xshift=3cm,yshift=-3cm,fill=white] (0.5,1.5) circle [radius=0.1cm];


\filldraw[xshift=8cm,yshift=-3cm,fill=blue!80!white] (0,0)--(1,0)--(1,2)--(0,2)--(0,0);
\draw[xshift=8cm,yshift=-3cm,thick] (0,0)--(1,0)--(1,1)--(0,1)--(0,0)--(1,0)--(1,1);
\filldraw[xshift=8cm,yshift=-3cm,fill=white] (0.5,0.5) circle [radius=0.1cm];
\filldraw[xshift=8cm,yshift=-3cm,fill=black] (0.5,1.5) circle [radius=0.1cm];

\filldraw[xshift=7cm,yshift=-2cm,fill=blue!80!white] (0,0)--(1,0)--(1,2)--(0,2)--(0,0);
\draw[xshift=7cm,yshift=-2cm,thick] (0,0)--(1,0)--(1,1)--(0,1)--(0,0)--(1,0)--(1,1);
\filldraw[xshift=7cm,yshift=-2cm,fill=white] (0.5,0.5) circle [radius=0.1cm];
\filldraw[xshift=7cm,yshift=-2cm,fill=black] (0.5,1.5) circle [radius=0.1cm];

\filldraw[xshift=8cm,yshift=-1cm,fill=blue!80!white] (0,0)--(1,0)--(1,2)--(0,2)--(0,0);
\draw[xshift=8cm,yshift=-1cm,thick] (0,0)--(1,0)--(1,1)--(0,1)--(0,0)--(1,0)--(1,1);
\filldraw[xshift=8cm,yshift=-1cm,fill=white] (0.5,0.5) circle [radius=0.1cm];
\filldraw[xshift=8cm,yshift=-1cm,fill=black] (0.5,1.5) circle [radius=0.1cm];


\filldraw[xshift=5cm,yshift=-5cm,fill=yellow] (0,0)--(2,0)--(2,1)--(0,1)--(0,0);
\draw[xshift=5cm,yshift=-5cm,thick] (0,0)--(2,0)--(2,1)--(0,1)--(0,0)--(1,0)--(1,1);
\filldraw[xshift=5cm,yshift=-5cm,fill=black] (0.5,0.5) circle [radius=0.1cm];
\filldraw[xshift=5cm,yshift=-5cm,fill=white] (1.5,0.5) circle [radius=0.1cm];

\filldraw[xshift=4cm,yshift=-4cm,fill=yellow] (0,0)--(2,0)--(2,1)--(0,1)--(0,0);
\draw[xshift=4cm,yshift=-4cm,thick] (0,0)--(2,0)--(2,1)--(0,1)--(0,0)--(1,0)--(1,1);
\filldraw[xshift=4cm,yshift=-4cm,fill=black] (0.5,0.5) circle [radius=0.1cm];
\filldraw[xshift=4cm,yshift=-4cm,fill=white] (1.5,0.5) circle [radius=0.1cm];

\filldraw[xshift=6cm,yshift=-4cm,fill=yellow] (0,0)--(2,0)--(2,1)--(0,1)--(0,0);
\draw[xshift=6cm,yshift=-4cm,thick] (0,0)--(2,0)--(2,1)--(0,1)--(0,0)--(1,0)--(1,1);
\filldraw[xshift=6cm,yshift=-4cm,fill=black] (0.5,0.5) circle [radius=0.1cm];
\filldraw[xshift=6cm,yshift=-4cm,fill=white] (1.5,0.5) circle [radius=0.1cm];


\filldraw[xshift=5cm,yshift=-2cm,fill=green] (0,0)--(2,0)--(2,1)--(0,1)--(0,0);
\draw[xshift=5cm,yshift=-2cm,thick] (0,0)--(2,0)--(2,1)--(0,1)--(0,0)--(1,0)--(1,1);
\filldraw[xshift=5cm,yshift=-2cm,fill=white] (0.5,0.5) circle [radius=0.1cm];
\filldraw[xshift=5cm,yshift=-2cm,fill=black] (1.5,0.5) circle [radius=0.1cm];

\filldraw[xshift=4cm,yshift=-3cm,fill=green] (0,0)--(2,0)--(2,1)--(0,1)--(0,0);
\draw[xshift=4cm,yshift=-3cm,thick] (0,0)--(2,0)--(2,1)--(0,1)--(0,0)--(1,0)--(1,1);
\filldraw[xshift=4cm,yshift=-3cm,fill=white] (0.5,0.5) circle [radius=0.1cm];
\filldraw[xshift=4cm,yshift=-3cm,fill=black] (1.5,0.5) circle [radius=0.1cm];

\filldraw[xshift=6cm,yshift=-3cm,fill=green] (0,0)--(2,0)--(2,1)--(0,1)--(0,0);
\draw[xshift=6cm,yshift=-3cm,thick] (0,0)--(2,0)--(2,1)--(0,1)--(0,0)--(1,0)--(1,1);
\filldraw[xshift=6cm,yshift=-3cm,fill=white] (0.5,0.5) circle [radius=0.1cm];
\filldraw[xshift=6cm,yshift=-3cm,fill=black] (1.5,0.5) circle [radius=0.1cm];
\end{tikzpicture}
   \end{minipage}
   \caption{Left: Gives oriented boundary values; \quad Middle and Right: Unoriented ones} \label{kolme}
\label{fig:OrientedCorner}
  \end{figure}

For a general convex gradient constraint $N$ and a natural domain $\Omega$, once we have oriented and natural boundary values $h_0$
it is not difficult to construct the frozen extensions.

\begin{Thm}\label{thm:frozenextension2}
Let $\Omega\subset\R^2$ be a natural domain for a closed convex polygon $N$,  and   $h_0$  an oriented  natural boundary value on $\partial \Omega$. Then $(\Omega, h_0)$ admits a frozen extension at every $z_0 \in \partial \Omega$.
\end{Thm}
\begin{proof}
If $z_0 \in \partial \Omega$ is not a corner point, or $z_0$ is a convex corner of $\Omega$, the argument of Lemma \ref{thm:frozenextension1} applies. Thus we only need to consider concave corners of $ \partial \Omega$.

Of the four different possible corner configurations, as discussed with the figures before Definition \ref{def:extremalorient},  at the corners of $\Omega$ oriented boundary values allow only the cases i) and ii). 
Assuming this, we add a parallelogram $\mathcal{P}$ with corner point $z_j$ in its interior and sides parallel to $z_{j+1}-z_j$ and $z_{j}-z_{j-1}$ as shown in Figure \ref{fig:GeoCase(ii)}. 
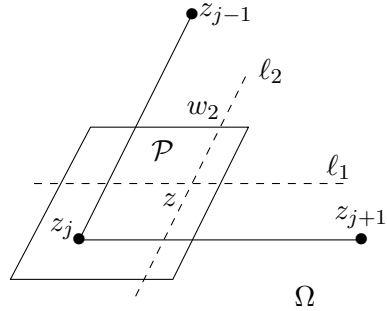
\begin{figure}[H]
\centering
\begin{tikzpicture}[xscale=1.5,yscale=1.5]
\draw (-1,0)--(1.5,0);
\draw(1.5,0) node{$\bullet$};
\draw(-1,0) node{$\bullet$};
\draw(-1,0)--(0,2);
\draw(0,2) node{$\bullet$};
\draw(0.3,2) node{$z_{j-1}$};
\draw(1.5,0.2) node{$z_{j+1}$};
\draw(-.2,0.36) node{$z$};
\draw(0.1,1.15) node{$w_2$};
\draw(-1.15,0.1) node{$z_{j}$};
\draw(1,-0.5) node{$ \Om$};
\draw(-1.61,-.35)--(-.9,1)--(0.5,1)--(-.17,-.35)--(-1.61,-.35);
\draw[dashed] (-0.5,-0.5)--(-0.25,0)--(0.25,1)--(0.5,1.5);
\draw[dashed] (-1.4,0.5)--(1.35,0.5);
\draw(1.3,0.65) node{$\ell_1$};
\draw(0.7,1.5) node{$\ell_2$};
\draw (-0.26,0.8) node{$\mathcal{P}$};
\end{tikzpicture}
\caption{Frozen extension}
\label{fig:GeoCase(ii)}
\end{figure}
Consider then the new domain $\widehat{\Om}= \mathcal{P} \cup \Om $, and a new boundary value on $\dv \widehat{\Om}$ as
in \eqref{Hhat}.  We claim that (\ref{Mhat}) holds for all $z\in \widehat{\Omega}\setminus\Omega$. Indeed, 
given a point $z\in \widehat{\Om}\setminus  \Om$,
there are two  lines $\ell_1$ and $\ell_2$ passing through $z$, where $\ell_1$ is parallel to $z_{j+1}-z_j$ and $\ell_2$ is parallel to $z_j-z_{j-1}$. Here $\ell_1$ intersects $\dv \widehat{\Om}\setminus \dv \Om$ at a unique point $w_1$ and $\ell_2$ at another one $w_2$. In the configuration i),
\[ 
h_N(z-w_2)=\langle p_n, z-w_2\rangle, \quad h_N(w_1-z)=\langle p_n, w_1-z\rangle, 
\]
while in case ii)  the roles of  $w_1$ and $w_2$ are interchanged in the above identity.   Arguing then as in the proof of Lemma \ref{thm:frozenextension1}  gives (\ref{Mhat}). Thus   $\Omega$ admits a frozen extension also at the corners. 
\end{proof}

The above argument fails for general natural boundary values. However, in case the convex gradient constraint $N$ is  {\it a triangle}, it is possible to modify the approach and find frozen extensions for any natural boundary value $h_0$, at any boundary point $z_0$, even if $N$ would have gas points or quasi-frozen points.

\begin{Thm}\label{thm:frozenextension2tri}
Assume that $N$ is a triangle. If  $\Omega\subset\R^2$ is a natural domain for $N$ and $h_0$ is a natural boundary value, then 
$\Omega$ admits a frozen extension with respect to $h_0$.
\end{Thm}
\begin{proof} It suffices to find a frozen extension at those concave corners $z_j$ of $\Omega$ where, among the four different possible corner configurations discussed before Definition \ref{def:extremalorient}, either of 
 the cases  iii) or iv) occurs. Further, $\Omega$ is a natural domain, thus at $z_j$ satisfies   either (\ref{ND1})  or (\ref{ND3}), for a vertex $p_n \in N$.  It is enough to consider the situation where we have iii) with (\ref{ND1}), all other cases work in an analogous way.
 
  Rather than using the previous extensions, we apply a different one, which works only when $N$ is a triangle.
Namely,  the two sides of $\Omega$ at $z_j$ are orthogonal to $p_n - p_{n-1}$ and to $p_n - p_{n+1}$, respectively. Let $\nu$ be the direction orthogonal to the remaining side of $N$, to $p_{n+1} - p_{n-1}$, and  
      define two closed parallelograms $\mathcal{P}_1$ and $\mathcal{P}_2$ outside of $\Om$ as in Figure \ref{figFrozenExtConcave},
such that $z_j$ is the corner point of both parallelograms. They are chosen so small that
they do not intersect with $\Om$.
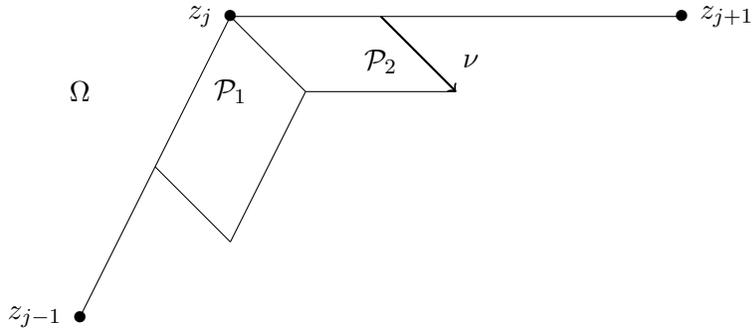
\begin{figure}[H]
\centering
\begin{tikzpicture}[xscale=2,yscale=2]
\draw (-1,0)--(2,0);
\draw(2,0) node{$\bullet$};
\draw(-1,0) node{$\bullet$};
\draw (-1,0)--(-2,-2);
\draw(-2,-2) node{$\bullet$};
\draw(-2.3,-2) node{$z_{j-1}$};
\draw(2.3,0) node{$z_{j+1}$};
\draw(-1.2,0) node{$z_{j}$};
\draw(-2,-0.5) node{$ \Om$};
\draw (0,0)--(0.5,-0.5);
\draw (-1,0)--(-0.5,-0.5);
\draw (-0.5,-0.5)--(0.5,-0.5);
\draw (-0.5,-0.5)--(-1,-1.5);
\draw (-1,-1.5)--(-1.5,-1);
\draw(-1,-0.5) node{$\mathcal{P}_1$};
\draw(0,-0.3) node{$\mathcal{P}_2$};
\draw[thick,->] (0,0)--(0.5,-0.5);
\draw(0.6,-0.3) node{$\nu$};
\end{tikzpicture}
\caption{Concave corner}
\label{figFrozenExtConcave}
\end{figure}
The sides of $\mathcal{P}_1$ are parallel to $z_j-z_{j-1}$ and $\nu$, and 
those of $\mathcal{P}_2$ parallel to $z_{j+1}-z_j$ and $\nu$, respectively.
Set then $\widehat\Om=({\Om}\cup {\mathcal{P}_1}\cup {\mathcal{P}_2})^\circ$ and define a boundary value $\widehat{h}_0:\partial\widehat{\Om}\to \R$ by 
\[ \widehat{h}_0(z)=\langle p_{n-1}, z-z_{j}\rangle+h_0(z_j)\]
when $z\in  \partial \mathcal{P}_1\setminus \widehat{\Om}$, by
\[ \hat{h}_0(z)=\langle p_{n+1}, z-z_{j}\rangle+h_0({z_j})\]
when $z\in  \partial \mathcal{P}_2\setminus \widehat{\Om}$, and  by $\widehat{h}_0(z)=h_0(z)$ when $z\in \partial \Omega\cap\partial \widehat{\Omega}$. 
It suffices to show that 
\begin{equation}\label{Mhat2}
\widehat{M}(z)=\widehat{m}(z)=\langle p_{n-1}, z-z_j\rangle+h_0(z_j)
\end{equation}
for all $z\in \mathcal{P}_1$, and that
\begin{equation}\label{Mhat3}
\widehat{M}(z)=\widehat{m}(z)=\langle p_{n+1}, z-z_j\rangle+h_0(z_j)
\end{equation}
for all $z\in \mathcal{P}_2$. We only prove (\ref{Mhat2}); the proof of (\ref{Mhat3}) is similar. For (\ref{Mhat2}), 
 fix a point $z\in \mathcal{P}_1$ and draw two lines $\ell_1$ and $\ell_2$ passing through $z$ such that $\ell_1$ is parallel to $\nu$ and $\ell_2$ is parallel to $z_j-z_{j-1}$. Let $w_1$ be the intersection point of $\ell_1$ and $\partial \widehat{\Om}\cap\partial \mathcal{P}_1$, and $w_2$ be that of  $\ell_2$ and $\partial\widehat{\Om}\cap\partial\mathcal{P}_1$. 
Under these conditions,
\begin{equation}\label{fun3}
h_N(w_2-z)=\langle p_{n}, w_2-z\rangle = \langle p_{n-1}, w_2-z\rangle
\end{equation}
and 
\begin{equation}\label{fun4}
h_N(z-w_1)= \langle p_{n-1}, z-w_1\rangle.
\end{equation}
Note here that (\ref{fun3}) holds for a general convex polygon $N$ while for
(\ref{fun4}) one needs $N$ to be a triangle.  
Now, as in the proof of Lemma \ref{thm:frozenextension1}, the claims (\ref{Mhat2}) - (\ref{Mhat3}) follow  from (\ref{fun3}) and 
(\ref{fun4}). 
In addition, $\widehat{M}=h_0$ on $\partial \Omega$. This completes the proof of theorem.
\end{proof}

Once a minimizer $h$ admits a frozen extension, the properness of $\nabla h$ follows readily.

\begin{Thm} \label{thm:h:proper3} Let $\Omega\subset \R^2$ be a natural domain and $h_0$  a natural boundary height function. Assume that either the gradient constraint $N$ is a triangle or, for a general gradient constraint $N$, that the boundary value $h_0$ is oriented. 

If the minimizer  $h$ to the variational problem \eqref{ELbasic} has a non-empty liquid domain $\LL \neq \emptyset$, then 
$$ \nabla h: \LL\to N^\circ\setminus
\mathscr G \quad \mbox{is proper},$$
that is, the boundary $\partial \LL$ is frozen.
\end{Thm}
\begin{proof}
Let $\{z_j\}_{j=1}^\infty\subset \LL$ be a sequence of points converging to a point $z_0\in \partial \LL$. 
We are then to show that 
\begin{equation}\label{fun5}
\lim_{j\to \infty} \textrm {dist}\,\big(\nabla h(z_j), \, \partial N^\circ \cup \mathscr{G} \big)=0.
\end{equation}
There are two cases: $z_0\in \Omega$ or $z_0\in \partial\Omega$. In the first case $z_0\in \Omega$, suppose that (\ref{fun5}) is not true. Hence
there is a subsequence, still denoted by $\{z_j\}_{j=1}^\infty$,  such that $\nabla h(z_j)\to p\in N^\circ\setminus \mathscr{G}$. 

Here recall  Theorem \ref{thm:DeS22}, which states that $\Gamma\circ\nabla h:\Omega\to \mathbb{S}^2$ is continuous, for a mapping  $\Gamma:N\to \mathbb{S}^2$ continuous in $N$, homeomorphic between $N^\circ$ and $\mathbb{S}^2 \setminus \{ \xi \}$, 
and taking $\partial N$ to the point $\{\xi\}\subset \mathbb{S}^2$. 
In particular,
\begin{equation}\label{fun6}
\Gamma\big(\nabla h(z_0)\big)= \Gamma (p)\neq \xi,  
\end{equation}
so that $\nabla h$ is continuous at $z_0$ and $\nabla h(z_0)=p \in N^\circ\setminus \mathscr{G}$. This means
that $z_0\in \LL$ by the definition of the liquid region, and contradicts  the fact that $z_0\in \partial\LL$.  

In the second case, we have $z_0\in \partial \Omega$. Now  we apply the frozen extension theorems, 
Theorem \ref{thm:frozenextension2}
  in case $h_0$ is oriented and Theorem \ref{thm:frozenextension2tri} in case $N$ is a triangle. 
With these  $(\Omega, h_0)$ admits a frozen extension  at any boundary point. Therefore, 
there is a domain $\widehat\Om\supset\Omega$ such that $z_0\in \widehat\Om$ is an interior point, and  a there is a
  boundary value $\widehat{h}_0$ on $\partial\widehat{\Omega}$ such that the upper and lower obstacles 
 $\widehat M (z) = \widehat m(z)$ on $ {\overline {\widehat \Omega}} \setminus \Omega$ and  such that on $\partial\Omega$ these obstacles agree with $h_0(z)$.
 
Let us  then consider the variational  problem \eqref{ELbasic} in $\widehat \Omega$, among the admissible class 
$\mathscr{A}_N(\widehat{\Omega},\widehat{h}_0)$. 
We know that the minimizer $\widehat{h}$ coincides with $h$ in $\Omega$ and 
necessarily they have the same liquid region. But now $z_0\in \widehat{\Omega}$ is an interior point. Thus we are in the first case, and  \eqref{fun5} follows.   \end{proof}
\smallskip

We conjecture that,  for any  natural domain $\Omega\subset \R^2$ and for any  natural boundary height function $h_0$   on $\partial \Omega$, that if the  minimizer $h$ to the corresponding variational problem \eqref{ELbasic} has a non-empty liquid domain $\LL \neq \emptyset$,  then the  whole boundary  $\partial \LL$ is always frozen.

\addtocontents{toc}{\vspace{-4pt}}
\section{Global Structure of the Limit Shape}\label{sec:topology}

Natural domains $\Omega$ and natural oriented  boundary values $h_0$ are the natural candidates for situations where one should expect frozen phenomena. Indeed, as we saw in Theorem \ref{thm:h:proper3}, if there  at all is a non-empty liquid domain $\LL$ for such a pair $(\Omega,h_0)$, then all of the boundary $\partial \LL$  is necessarily frozen. On the other hand, even within this class there are simple examples - such as the one below - where the minimizer $h$ is piecewise affine with $\nabla h$ taking values only in the corners of the constraint $N$, and thus trivially there is no liquid domain. But  Theorem \ref{thm:liquid:affine1} and Proposition \ref{thm:liquid:affine2} below show that such examples are, in fact,  the only obstructions. Thus the question of the existence of a liquid domain is a combinatorial one.

\begin{ex}\label{noliquid}
Let us present an example in  the uniform lozenge model, of a natural domain $\Omega$ and a natural boundary value $h_0$, such that $h_0$ is piecewise affine but  \emph{not affine} and yet there is no liquid domain.

For this, let $N$ be an even sided triangle shown on left in Figure \ref{figPA-bv1} and let $\Om$ be the regular hexagon, natural with respect to $N$ with natural boundary values. The corresponding upper McShane extension $M$ is  shown in middle of Figure \ref{figPA-bv1}.

\begin{figure}[H]
\centering
\begin{tikzpicture}[xscale=2,yscale=2]
\draw[thick,yshift=1cm] (-2,0)--(-1.5,{-sqrt(3)/2})--(-1,0)--(-2,0);
\draw[yshift=1cm] (-1.5,-0.25) node{$N$};
\draw[yshift=1cm] (-1.5,-1) node{$p_1$};
\draw[yshift=1cm] (-2.2,0) node{$p_3$};
\draw[yshift=1cm] (-0.8,0) node{$p_2$};
\draw[thick] (0,0)--({sqrt(3)/2},-1/2)--({sqrt(3)},0)--({sqrt(3)},1)--({sqrt(3)/2},1.5)--(0,1)--(0,0);
\draw[thick,blue] ({sqrt(3)/2},-1/2)--({sqrt(3)/2},1/2);
\draw[thick,blue] ({sqrt(3)/2},1/2)--(0,1);
\draw[thick,blue] ({sqrt(3)/2},1/2)--({sqrt(3)},1);
\draw[xshift=4cm] ({sqrt(3)/2-3.3},1.3) node{$\Om$};
\draw(1,1) node{$p_1$};
\draw (0.5,0.2) node{$p_2$};
\draw (1.2,0.2) node{$p_3$};
\draw[thick] ({sqrt(3)/2+2.5},-1/2)--({sqrt(3)+2.5},0)--({sqrt(3)+2.5},1)--({sqrt(3)/2+2.5},1.5)--(2.5,1);
\draw[thick] ({sqrt(3)/2+2.5},-1/2)--({sqrt(3)/2+2.5},1/2);
\draw[thick] ({sqrt(3)/2+2.5},1/2)--(2.5,1);
\draw[xshift=2.5cm] ({sqrt(3)/2+0.7},1.3) node{$\Om'$};
\end{tikzpicture}
\caption{Natural boundary value allowing no liquid domain}
\label{figPA-bv1}
\end{figure}
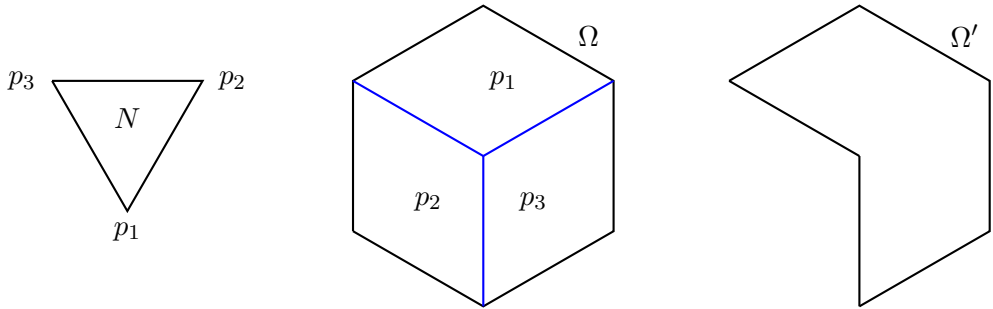

Consider the new domain $\Omega'$ as shown on right in  Figure \ref{figPA-bv1}, with boundary value $h_0'=M\vert_{\dv \Om'}$. Then $\Om'$ is a natural domain and $h_0'$ is a natural boundary value. By similar reasoning as in the previous section, 
we see that the upper and lower McShane extensions for $\dv \Om'$ and $h_0'$ satisfy $M'(x)=m'(x)$ for all $x\in \Om'$. Moreover, the boundary value $h_0'$ is only piece-wise affine. \hfill $\Box$
\end{ex}
\smallskip

Let us then turn to the question of the structure  of the minimizers outside the liquid regions, for general domains and boundary values.
Recall here the coincidence set  $\Lambda = \{x\in \overline{\Om}: h(x)=m(x) \; {\rm or } \; h(x)=M(x)\} $, where $M(x)$ and $m(x)$ are 
the upper and lower obstacles, respectively, as in  \eqref{def:Mm}.

\begin{Thm}\label{thm:liquid:affine1}
Let $\Omega\subset\R^2$ be a bounded Lipschitz domain, $h_0$ an admissible boundary value
and $h$ the minimizer of variational problem (\ref{ELbasic}). Then  $h$ is countably piecewise affine in 
$\Omega\setminus(\Lambda\cup\overline{\LL})$
with gradient having values in $ \mathscr{P}\cup\mathscr{Q}\cup\mathscr{G}$.
\end{Thm}

\begin{proof}
By Theorem 
\ref{thm:DeS12} of De Silva and Savin \cite{DeS10}, we know that 
$h\in C^1(\Omega\setminus \Lambda)$.  Moreover,  
outside  $\Lambda\cup \overline{{\LL}}$ 
the gradient $\nabla h$ takes values in $ \partial N \cup \mathscr{G}$, so that the set
\[ E= \left\{ z\in \Omega\setminus(\Lambda\cup\overline{\LL}): \nabla h(z)\in \partial N\setminus (\mathscr{P}\cup\mathscr{Q})\right\}\]
is open. We claim that
\begin{equation}\label{claim:E}
E =\emptyset.
\end{equation}
This means that the image of the open set $\Omega\setminus(\Lambda\cup\overline{\LL})$ under the continuous map $\nabla h$ is contained in   the finite set 
$\mathscr{P}\cup\mathscr{Q}\cup\mathscr{G}$,  
which implies that $h$ is countably piecewise affine.  

Thus it suffices to show \eqref{claim:E}. We may assume that $0\in N^\circ$, and argue by contraction. Suppose that $E$ is not empty. Since it is open, there is a ball 
$B_r(z_0)\subset E$.  Let us then construct a function $u\in \Af$ such that $\nabla u(z)=0$ for all $z\in B_{r^\prime}(z_0)$,
where $0<r^\prime<r$. For this, note that by definition $m(z)<M(z)$ for all 
$z\in \Om \setminus \Lambda$. Hence for some $0<r^\prime<r$  small enough we find $c\in \R$ such that $m(z)<c<M(z)$ whenever 
$|z-z_0| < {r^\prime}$.

In the following we denote $B _{r^\prime}(z_0)$ by $B$.
Consider next the constant function on $B$ with value $c$, and let $M_{c}(z)$ be its upper McShane extension from $B$ to $\R^2$,
\[ M_c(z)=\min_{w\in B} \big(h_N(z-w)+c\big), \quad z\in \R^2.\]
Then we define a function $u:\Omega\to \R$ by
\[ u(z)=\min \big( M(z), M_c(z)\big).\]
 Since
$c< M$ in $B$, we have $u=M_c=c$ in $B$, and hence that
$\nabla u=0\in N^\circ$ in $B$. Since  $\nabla M(z)\in N$ and $M_c(z)\in N$ for a.e. $z\in \Omega$, we have that
$\nabla u(z)\in N$ for a.e. $z\in \Omega$. And since $m(z)<c$ for $z\in B$, we have $M_c(z)\ge m(z)$ for $z\in \overline{\Omega}$. In particular, 
$M_c(z)\ge m(z)=h_0(z)$ for $z\in \partial \Omega$. This shows that
$u=M=h_0$ on $\partial \Omega$. Thus $u\in \Af$, as required. 

Continuing with the proof of the Theorem, let us start from the Gâteaux derivative inequality $dI_\sigma[h; u-h]\geq 0$ in \eqref{gateaux1}, valid since $h$ is the minimizer of the variational problem (\ref{ELbasic}). In concrete terms,
\begin{equation}\label{minimzernonnegtive}
\int_{\Omega\setminus B}d\sigma(\nabla h(z); \nabla u(z)-\nabla h(z))dz+\int_{B}d\sigma(\nabla h(z), \nabla u(z)-\nabla h(z))dz\ge 0.
\end{equation} 
Since the Gâteaux derivative of $\sigma$ is bounded from above, see \eqref{GatDerivative}, 
we have
\begin{align*}
\int_{\Omega\backslash B}d\sigma(\nabla h(x), \nabla u(x)-\nabla h(x))dx \le 2 \vert \Omega\setminus B\vert \max_N |\sigma| <\infty.
\end{align*}
We know that $ d\sigma(p,q-p)=-\infty$ for    $p\in \dv N\backslash (\Pp\cup \Qq)$ and $q\in N^\circ$. Since for all $z\in B$, $\nabla u(z)=0\in N^\circ$ and $\nabla h(z)\in
\dv N\backslash (\Pp\cup \Qq)$, thus
\begin{align*}
\int_{B}d\sigma(\nabla h(z), \nabla u(z)-\nabla h(z))dz=-\infty.
\end{align*}
But that is in  a contradiction with (\ref{minimzernonnegtive}). This proves the claim (\ref{claim:E}), and hence the Theorem.
\end{proof}

Since Theorem 
\ref{thm:DeS12} requires  the minimizer $h$ of the variational problem (\ref{ELbasic}) to coincide with one of the obstacles $M(x)$ or $m(x)$ whenever it is not $C^1$-smooth, 
 Theorem  \ref{thm:liquid:affine1} gives strong global rigidity for the structure of the minimizer. 
 Also, for natural boundary values one can control the minimizers on the coincidence set.

\begin{Prop}\label{thm:liquid:affine2}
Let $\Omega\subset\R^2$ be a natural domain and $h_0$ a natural boundary value. Then the upper obstacle $M$ and the lower obstacle $m$ are both piecewise affine functions with gradient having values in $\mathscr{P}$.
\end{Prop}
\begin{proof} We give the proof only  for the upper obstacle $M$. Let $\Omega$ ba a natural domain with $d$ vertices $\{ z_1,...,z_d\}$ 
(setting $z_{d+1}=z_1$). Let $h_0$ be a natural boundary value and $M$ the upper obstacle defined as in (\ref{def:Mm}). 
For each $j=1,..,d$, we consider the McShane extension from the line segment $[z_j, z_{j+1}]$ to $\R^2$
\[ M_j(z)=\min_{w\in [z_j,z_{j+1}]} \, \big( h_N(z-w)+h_0(w)\big). \]    
Then it is easy to prove that 
\[ M(z)=\min_{1\le j\le d} M_j(z), \quad z\in \R^2.\]
 
 For the Proposition it is thus enough to show  that each $M_j$ is  piecewise affine  in $\R^2$,  with gradient having values in $\mathscr{P}$.  To prove this, we note that since $\Omega$ is natural there is a vertex $p_n \in N$ with
\[ \langle z_{j+1}-z_j, p_{n+1}-p_n\rangle =0,\]
 and since $h_0$ is natural,
\begin{equation}\label{h0pn} h_0(w)=\langle w-z_j, p_n\rangle + h_0(z_j), \qquad  w\in [z_j,z_{j+1}].
\end{equation}
 Moreover, 
 we have either 
\begin{equation}\label{case1:zj} h_N(z_{j+1}-z_j)= \langle z_{j+1}-z_j, p_n\rangle,  \qquad {\rm or } \qquad h_N(z_{j}-z_{j+1})=\langle z_{j}-z_{j+1}, p_n\rangle.
\end{equation}

It is enough to only  consider the first option in (\ref{case1:zj}). In this case,  consider the McShane extension of $h_0$ from point $z_j$ to $\R^2$,
\[ \overline{ M}_j(z)= h_N(z-z_j)+h_0(z_j), \quad z\in \R^2.\]
From \eqref{supportf} we see that $\overline{M}_j(z)$ is piecewise affine 
and its gradient has values in $\mathscr{P}$.  Thus it suffices to prove that for all $z\in \R^2$,
\begin{equation}\label{Mjtilte}
M_j(z)=\overline{M}_j(z).
\end{equation}
Clearly, we have $M_j(z)\le \overline{M}_j(z)$ for  $z\in \R^2$. To prove $M_j(z)\ge \overline{M}_j(z)$, we let $w$ be a point in the line segment $[z_j,z_{j+1}]$. First, with the triangle inequality for the support function $h_N$,
\[ h_N(z-w)+h_N(w-z_j)\ge h_N(z-z_j)\]
for any $z\in \R^2$. Second, because of the first option in (\ref{case1:zj}), we have for any $w\in [z_j,z_{j+1}]$ that
\[ h_N(w-z_j)=\langle w-z_j, p_n\rangle .\]
It then follows from the above two inequalities and (\ref{h0pn})  that
\[ h_N(z-w)+h_0(w)\ge  h_N(z-z_j)+\langle z_j - w, p_n\rangle + h_0(w) = h_N(z-z_j)+h_0(z_j) \]
for any $z\in \R^2$. This shows that
$M_j(z)\ge \overline{M}_j(z)$ for $z\in \R^2$, and completes the proof.  
\end{proof}

\subsubsection{Proofs of Theorems \ref{thm:liquid:affine}  and \ref{thm:h:proper}. }

Theorem \ref{thm:liquid:affine} is an immediate consequence of  Theorem \ref{thm:liquid:affine1} and Proposition \ref{thm:liquid:affine2}. \hfill $\Box$

For Theorem \ref{thm:h:proper}, in $\Omega\setminus \overline{\LL}$ the minimizer is piecewise affine by Theorem \ref{thm:liquid:affine},
and if $\LL \neq \emptyset$, with gradient constraint $N$ a triangle the boundary $\partial \LL$ is frozen by   Theorem \ref{thm:h:proper3}. 
Similarly, for gas domains the argument of Theorem \ref{GAS} applies. Indeed, since 
$h_0$ is a natural boundary value and $\partial \Omega$ is polygonal,  we have a gas component $U_q \subset \Omega$ with $\nabla h \equiv q \in \Gg$ on $\partial \LL \cap \partial U_q$, even if $\Omega$ is not simply connected. This completes the proof of  Theorem \ref{thm:h:proper}. \hfill $\Box$
\smallskip

In addition to the above Theorems, in the case of a general gradient constraint $N$ and surface tension $\sigma$, we can use Theorem \ref{thm:h:proper3} to show that $\nabla h: \LL\to N^\circ\setminus
\mathscr G $ is a proper map,  if $h_0$ is oriented and the liquid domain $\LL \neq \emptyset$. Thus combining with Theorem \ref{thm:liquid:affine} this gives

\begin{Thm}\label{thm:h:proper2}
Let $\Omega$ be a natural domain, and $h_0$  a natural and oriented boundary value. Suppose that $h$ is the minimizer of variational problem (\ref{ELbasic})
among the class $\mathscr{A}_N(\Omega,h_0)$. 

Then, as in Theorem \ref {thm:h:proper}, either $h$ is piecewise affine in $\Omega$, or 
else $h$ has a liquid domain  with frozen boundary. If $q \in  \Gg$ is a gas point for $\sigma$, then there is a non-empty gas domain $U_q \subset \Omega$ with $\nabla h \equiv q$ in  $U_q$. 
\end{Thm}

 We conjecture that Theorem \ref{thm:h:proper2} remains true even without assuming   $h_0$  being  oriented.

\addtocontents{toc}{\vspace{-4pt}}
\section{Minimality}
\label{section:minimality}

\subsection{Minimality of a function}\label{subsection:minimality}

In this subsection we present a  method to show that a given function $h$ from the admissible class $\mathscr{A}_{N}(\Omega,h_0)$ is actually the minimizer of the variational problem (\ref{ELbasic}).
We start with a   general proposition which gives a sufficient condition for $h$ to be the minimizer. 
\begin{Prop}\label{prop:divergefree}
Function $h$ is a minimizer of (\ref{ELbasic}) among the class $\mathscr{A}_{N}(\Omega,h_0)$, provided  there is a vector field $\Phi\in 
L^1(\Omega; \mathbb{R}^2)$ such that $\Phi(z)\in \partial \sigma\big(\nabla h(z)\big)$
for almost all $z\in \Omega$ and 
\[ \rm{div}\, \Phi=0 \quad \text{in }\Omega\]
in the sense of distributions.
\end{Prop} 
\begin{proof} Let $u\in\mathscr{A}_{N}(\Omega,h_0)$.   With \eqref{defsubdif34} we have 
\[ d\sigma\big(\nabla h(z); \nabla u(z)-\nabla h(z)\big)\ge \langle \Phi(z), \nabla u(z)-\nabla h(z)\rangle\]
for almost all $z\in \Omega$, since $\Phi(z)\in \partial \sigma\big(\nabla h(z)\big)$. 
Thus,
\[\int_\Omega d\sigma\big(\nabla h(z); \nabla u(z)-\nabla h(z)\big)\, dz
\ge \int_\Omega \langle \Phi(z), \nabla u(z)-\nabla h(z)\rangle\, dz=0.\]
The last equality follows by a simple approximation argument from the facts that $\Phi$ is  divergence free,  that
$\Phi\in L^1(\Omega;\mathbb{R}^2)$ and  that both $u$ and $h$ are Lipschitz continuous
in $\Omega$ with the same boundary value.  This proves the claim, since the above inequality holds for all $u\in\Af$. \end{proof}

Assume then the  surface tension $\sigma$ is as in \eqref{Pst2}. Let $h \in \mathscr{A}_N(\Om,h_0)$, and assume the function has a liquid domain $\LL \subset \Om$ with $\partial \LL$ frozen, so that \eqref{def:liquiddomain2} - \eqref{frozen.def} hold. Thus
\begin{equation}\label{equ:hstar}
\text{div}\, \big(\nabla\sigma(\nabla h)\big)=0 \quad \text{in }\mathcal L.
\end{equation}
Our  task is now  to find explicit conditions that guarantee  that $h$ is the actual minimizer of the variational problem (\ref{ELbasic}) in the original domain $\Omega$. For this we use Proposition \ref{prop:divergefree} and  assume that $\sigma$ has no gas points. In this case our construction  gives a vector field $\Phi$ which is continuous, except on a finite set of line segments where $\vert\Phi\vert$ is not bounded. In the gas regions  the vector field $\Phi$ will, in general, be discontinuous. We plan to discuss this issue in a  future work. 

By our assumptions $\nabla h: \mathcal L\to N^\circ\setminus \mathscr G$ is a proper map. 
Then thanks to Theorem \ref{thm:main2}, we know that there is a finite set of
 singular points $\{z_j\}_{j=1}^n \subset \partial \mathcal L$ such that for any $1 \leq j \leq n$ and
for any point $z_0$ in the 
arc of $\partial \LL$ joining $z_j$ and $z_{j+1}$, one has
\begin{equation}\label{equ:z0}
\nabla h(z_0)=\lim_{z\to z_0, z\in \LL} \nabla h(z)=p_0,
\end{equation}
where $p_0\in \mathscr P\cup\mathscr Q$ is the same for all points on the given boundary arc. 
Moreover, 
\begin{equation}\label{equ:z00}
 (\nabla \sigma \circ \nabla h)(z_0)= \lim_{z\to z_0, z\in \LL}
\nabla \sigma(\nabla h(z))=\widehat{\nabla} \sigma(p_0; p-p_0),
\end{equation}
where $p\in N^\circ\setminus \mathscr G$ such that the unit normal  $\nu(z_0)$ to $ \LL$ at $z_0$ is parallel to the vector $p-p_0$. 

 We further assume that the candidate minimizer $h$ is a (finitely) piecewise affine function in
$\Omega\setminus \overline{ \LL \,}$ with $\nabla h \in \mathscr P\cup \mathscr Q$. 
In addition, in each component of $\Omega\setminus  \overline{ \LL \,}$
we assume that
$h$ enjoys the {\it frozen star ray} property, 
defined as follows.

\begin{Def}\label{def:ray} Suppose $h \in \mathscr{A}_N(\Om,h_0)$, and \eqref{equ:hstar} - \eqref{equ:z00} hold in a subdomain $\LL \subset \Om$. 
We  say that  $h$ 
has the {\rm frozen star ray property}, if it is finitely piecewise affine in $\Omega\setminus  \overline{ \LL \,}$ and 
in each component of 
$\Omega\setminus  \overline{ \LL \,}$ there is a family of rays such that the following holds: 

i) Each ray starts from a point $z_0\in \partial\LL$ and  is tangent to $\partial \LL$ at $z_0$; 

ii) For each $z_0\in \partial \LL\setminus \{ z_1,z_2,...,z_n\}$, if $\, \nabla h(z_0) = p_0 \in  \mathscr P\cup \mathscr Q$ in (\ref{equ:z0}), then the intersection

\vspace{-.17cm}
\quad of this component of $\Omega\setminus \overline{ \LL \,}$ and the ray starting from $z_0$ 
 lies in 
\begin{equation}\label{def:fp0}
\mathcal F_{p_0}=\{ z\in \Omega\setminus \overline{ \LL \,}: \nabla h(z)=p_0\}.
\end{equation}

iii) The rays in the family do not intersect each other  inside 
 this component of $\Omega\setminus \overline{ \LL \,}$;

iv) The union of the rays in the family covers this component of $\Omega\setminus \overline{\LL \,}$.
\end{Def}
 For an illustration see Figure \ref{figNonIntersect1}.

\begin{figure} 
\centering
\begin{tikzpicture}[xscale=2,yscale=2]
\draw[thick] (-3,-0)--(0,0)--(0,3);
\draw[thick,color=blue,rotate=-45,domain=0:3.1,shift={(0.7,0)}]   plot (-\x+1,{(0.4*\x)^1.5}); 
\draw[thick,color=blue,rotate=-45,domain=0:3.1,shift={(0.7,0)}]   plot (-\x+1,{-(0.4*\x)^1.5});  
\filldraw (-2,0)  circle (0.8pt);
\filldraw (0,2)  circle (0.8pt);
\draw[thick,color=blue] plot [smooth , tension=0.7] coordinates {(-2,0)(-2.5,-0.1)(-3,-0.3)}; 
\draw[thick,color=blue] plot [smooth , tension=0.7] coordinates {(0,2)(0.1,2.5)(0.3,3)}; 
\draw (-2,0.2) node{$z_{j+1}$};
\draw (0.4,2) node{$z_{j}$};
\filldraw (0,-0.55)  circle (0.8pt);
\draw (0,-0.7) node{$z_{0}$};
\draw[dashed] (0,-0.55)--(-1.3,0);
\draw[dashed] (0,0)--(1,-1);
\draw[dashed] (0,1.3)--(0.55,0);
\draw (-1,-0.4) node {$\dv \LL$};
\draw (1.1,0.8) node {$ \LL$};
\end{tikzpicture}
\caption{}
\label{figNonIntersect1}
\end{figure}

Suppose that all of the above assumptions hold. Then we claim that 
$h$ is a minimizer of variational problem (\ref{ELbasic}), with its boundary value $h_0$. Indeed, we define 
$\Phi:\Omega\to \mathbb R^2$ as follows
\[ \Phi(z)=\begin{cases}  (\nabla \sigma \circ \nabla h)(z), \quad &\text{if } z\in \LL;\\
 (\nabla \sigma \circ \nabla h)(z_0), & \text{if } z \text{ lies 
on the ray starting at } z_0, \text{ tangent to }\LL, 
\end{cases}
\]
where $ (\nabla \sigma \circ \nabla h)(z_0)$ is defined as in (\ref{equ:z00}).
By Proposition \ref{prop:divergefree}, we only need to prove that
the vector fields $\Phi$ defined as above satisfy that $\Phi\in L^1(\Omega;\mathbb R^2)$, $\Phi(z)\in \partial \sigma(\nabla h(z))$ for 
almost all $z\in \Omega$ and
that it is  divergence free. 

First, from Theorem \ref{repone}, 
$\nabla h = U_\sigma \circ f$ where $f: \LL \to \DD$ is proper and solves $f_{\overline{ z }} = \mu_{\sigma} (f) f_z$.
Thus $\nabla \sigma(\nabla h) = (\nabla  \sigma  \circ U_\sigma)(f)$. This leads, c.f. Corollary \ref{sigmarep}, to
$$ |\Phi(z)| =  |(\nabla \sigma \circ \nabla h)(z)| \leq C \sum_{j=1}^n \bigl | \log|f(z) - f(z_j)| \bigr|, \qquad z\in \LL.
$$
From  Theorem \ref{propStoilow35}, Proposition \ref{bdryreg1} and Corollary \ref{gmap} we thus see that $\Phi\in L^1(\Omega; \mathbb R^2)$. 
Actually, $\Phi$ is continuous in $\Omega$, except on  the  rays starting at  
the singular points $z_1,z_2,..., z_n$, where 
$\vert \Phi\vert$ is not bounded. Second, due to (\ref{equ:z0}), (\ref{equ:z00}), iii) in Definition \ref{def:ray} and Theorem \ref{sigma1}, we know that
$\Phi(z)\in \partial \sigma(\nabla h(z))$ for all $z\in \Omega$, except on a  finite number of rays.

Therefore the main point of the argument is to  show that $\Phi$ is  divergence free in $\Omega$. 
Clearly, in $\LL$ this is true by Equation (\ref{equ:hstar}). Since 
$\Phi$ is continuous  on $ \partial \LL\setminus \{ z_1,z_2,...,z_n\}$,
we only need to show that it is  divergence free in each component 
$W$ of $\Omega\setminus \overline{\LL \,}$. We first show this in 
$W_0=W\cap \mathcal{F}_{p_0}$ for any $p_0\in   \mathscr P\cup \mathscr Q$
such that $W_0\neq \emptyset$.  
Notice that for $z\in W_0$ lying on the ray starting at $z_0$, 
\[ \Phi(z)=  (\nabla \sigma \circ \nabla h)(z_0)= \widehat{\nabla} \sigma(p_0; p-p_0),\]
where $p\in N^\circ\setminus \mathscr G$ with $p-p_0$ parallel to the unit normal $\nu(z_0)$ to $\LL$ at $z_0$.
To simplify the notation, let us write  $\Psi (s,t)= (\Psi_1(s,t), 
\Psi_2(s,t))=
\widehat{\nabla} \sigma(p_0; p-p_0)$ where $(s,t)=p-p_0$. By  \eqref{GeneralizedGradient3a}, 
$\Psi$ is homogeneous of degree zero, that is, for $0<\lambda<1$
\begin{equation}\label{equ:Psi0} 
\Psi(\lambda s, \lambda t)=\Psi(s,t).
\end{equation}
Since $\Psi$ is homogeneous, we may extend the domain  of $\Psi$
so that the above equality holds for all $\lambda\in \mathbb{R}_+$.  Also, by rotating $N$ when necessary, and  transforming  $\sigma$ and $h$ accordingly, we may assume that the $s-$axis intersects  the domain  of $\Psi$ only at the origin.   This gives us
\[ \Psi(s,t)=\Psi(s/t,1).\]
By Theorem \ref{sigma1} and Corollary \ref{subdiff}  we know that $\Psi$ is smooth outside the origin.
Further, the differential matrix $D\Psi$ is  symmetric.  That, with a differentitation of (\ref{equ:Psi0}) with respect to $\lambda$, shows
\begin{equation}\label{equ:Psi} 
s\partial_s\Psi_1(s,t)+t\partial_s\Psi_2(s,t)=s\partial_s\Psi_1(s,t)+
t\partial_t \Psi_1(s,t)  = 0.
\end{equation}
With these notations, we rewrite 
$\Phi(z)=\Psi(\varphi(z),1)$, 
where $\varphi:\mathcal F_{p_0} \cap W \to \mathbb R$ is given by 
$\varphi(z)=\varphi(z_0) = s/t$, for $(s,t)=p-p_0$ such that
(\ref{equ:z00}) holds and for $z$ on the ray of Definition \ref{def:ray},  starting at $z_0$ and tangent there to 
$\LL$.
 Note that
\begin{equation}\label{equ:varphi0}
\varphi(z_0)=\nu_1(z_0)/\nu_2(z_0),
\end{equation}
where $\nu(z_0)=(\nu_1(z_0), \nu_2(z_0))$ is the inward normal unit of $\LL$ at $z_0$.

In this notation
\[\text{div}\, \Phi(z)=\partial_s \Psi_1(\varphi(z),1)\partial_x\varphi(z)
+\partial_s \Psi_2(\varphi(z),1)\partial_y\varphi(z), \quad z=(x,y).\]
Thus by (\ref{equ:Psi}), $\Phi$ being divergence free in $\mathcal F_{p_0}$ is equivalent  
to the real Burgers equation
\begin{equation}\label{equ:realburgers}
\partial_x\varphi(z)  - \varphi(z)\partial_y\varphi(z) =0, \quad z=(x,y)\in \mathcal{F}_{p_0}.
\end{equation}
Now it is easy to show that the function $\varphi$ defined in above is actually
a solution to equation (\ref{equ:realburgers}), due to (\ref{equ:varphi0}) and  the fact that $\varphi(z)=\varphi(z_0)$ for all $z$ that lie in the ray tangent to $\LL$ at $z_0$. 
Indeed, if $\varphi$ is a continuous solution 
to equation (\ref{equ:realburgers}), it follows 
by the method of characteristic curves that $\varphi(z)=\varphi(z_0)$ on each characteristic line
\[ y-y_0=-\varphi(z_0) (x-x_0), \quad z=(x,y), \quad z_0=(x_0,y_0).
\]
From (\ref{equ:varphi0}) one sees that the above characteristic line is tangent to $\LL$ at $z_0$, and is the ray starting from $z_0$  in Definition \ref{def:ray}. Also,  the rays in Definition \ref{def:ray} do not intersect in $\Omega\setminus \overline{\LL \,}$. Thus $\varphi$ is a solution to equation (\ref{equ:realburgers}), and  $\Phi$ is divergence free in $W_0 = W \cap \mathcal F_{p_0}$. 

Next, we show that $\Phi$ is of divergence free in the component $W$ of $\Omega \setminus \overline{\LL}$. The proof is now easy and  follows from  Proposition \ref{prop:sigma:property}. We only need to consider the case that $W$ has two parts: $W_0=W\cap \mathcal{F}_{p_0}$ and $W_1=W\cap \mathcal{F}_{p_1}$, which are separated by the ray starting from a singular point $z_j$ on $\partial \LL$. Here $p_0, p_1
\in \mathscr P\cup \mathscr Q$ are two neighbouring points, c.f \eqref{general}. 

Let us  here consider 
two points $z_0^\varepsilon, z_1^\varepsilon\in \partial \LL$ such that
they both approach  $z_j$ as $\varepsilon \to 0$. We assume $L_0^\varepsilon \subset \mathcal{F}_{p_0}$ and 
$L_1^\epsilon \subset \mathcal{F}_{p_1}$ for the corresponding rays starting from $z_0^\varepsilon$ and $z_1^\varepsilon$,
respectively. Denote by $W_\varepsilon\subset W$  the small region between
$L_0^\varepsilon$ and $L_1^\varepsilon$, see Figure \ref{figNonIntersect} below.

\begin{figure} 
\centering
\begin{tikzpicture}[xscale=2,yscale=2]
\draw[thick,color=blue,rotate=-45,domain=0:3.1,shift={(0.7,0)}]   plot (-\x+1,{(0.4*\x)^1.5}); 
\draw[thick,color=blue,rotate=-45,domain=0:3.1,shift={(0.7,0)}]   plot (-\x+1,{-(0.4*\x)^1.5}); 
\draw (0.6,-1) node {$z_{0}^\eps$};
\draw (1.0,-0.6) node {$z_{1}^\eps$};
\filldraw (0.6,-0.82)  circle (0.6pt);
\filldraw (0.8,-0.6)  circle (0.6pt);
\draw[red] (0.6,-0.82)--(-1.6,0.6);
\draw[red] (0.8,-0.6)--(-0.6,1.7);
\draw (-1.7,0.7) node {$L_0^\eps$};
\draw (-0.8,1.7) node {$L_1^\eps$};
\draw (-1.1,1.2) node {$W_\eps$};
\draw (1.3,-1.3) node{$z_j$};
\draw (-1,1)--(0,0)--(1,-1)--(1.2,-1.2);
\draw (-1,-0.4) node {$\dv \LL^\star$};
\draw (1,0.5) node {$ \LL^\star$};
\filldraw (-1.4,{-0.2*1.42/2.12+0.6})  circle (0.6pt);
\draw[->] (-1.4,{-0.2*1.42/2.12+0.6})--(-1.2,{-0.2*1.42/2.12+0.6+0.2*2.12/1.42});
\filldraw (-0.4,{1.7-0.2*2.3/1.4})  circle (0.6pt);
\draw[->] (-0.4,{1.7-0.2*2.3/1.4})--(-0.6,{1.7-0.2*2.3/1.4-0.2*1.4/2.3});
\draw (-1.4,{-0.2*1.42/2.12+0.6-0.1}) node {$z$};
\draw ({-0.4+0.15},{1.7-0.2*2.3/1.4}) node {$z$};
\draw ({-1.4+0.3},{-0.2*1.42/2.12+0.6+0.1}) node {$\nu(z)$};
\draw ({-0.4-0.05},{1.7-0.2*2.3/1.4-0.3}) node {$\nu(z)$};
\draw (-1.7,0.2) node {$\mathcal{F}_{p_0}$};
\draw (-0.2,1.7) node {$\mathcal{F}_{p_1}$};

\end{tikzpicture}
\caption{}
\label{figNonIntersect}
\end{figure}

Now let $\varphi\in C^\infty_0(W)$ be a cut-off function. We have
\begin{equation}\label{Uepsilon1}
\int_W  \langle \Phi(z), \nabla \varphi(z)\rangle\, dz=\lim_{\epsilon\to 0} \int_{W\setminus W_\epsilon} \langle \Phi(z), \nabla \varphi(z)\rangle\, dz,
\end{equation}
since $\Phi\in L^1(\Omega;\R^2)$. Integration by parts, with  the fact that $\Phi$ is  divergence free in 
$W\cap \mathcal{F}_{p_0}$ and $W\cap \mathcal{F}_{p_1}$, gives us 
\begin{equation}\label{Uepsilon2}
\int_{W\setminus W_\epsilon} \langle \Phi(z), \nabla \varphi(z)\rangle\, dz=
\int_{L_0^\varepsilon}\langle \Phi(z), \nu(z)\rangle \varphi(z)\, dS(z)
+\int_{L_1^\varepsilon}\langle \Phi(z), \nu(z)\rangle\varphi(z)\, dS(z),
\end{equation}
where $\nu(z)$ is the inward normal unit vector of $W_\varepsilon$ at $z$ and the integrals on the right hand side are with respect to the one dimensional Lebesgue measure.

Next, let us use \eqref{equ:z00} at the point $z_0^\varepsilon$, that
\[ \Phi(z)= (\nabla \sigma \circ \nabla h)(z_0^\varepsilon)=
\lim_{z\to z_0^\varepsilon, z\in \LL}
\nabla \sigma(\nabla h(z))=\widehat{\nabla} \sigma(p_0; p_\varepsilon-p_0), \qquad z \in L_0^\varepsilon.
\]
 Since $L_0^\varepsilon$ is tangent to $\partial \LL$, we know from Theorem \ref{thm:main2} that
$p_\varepsilon-p_0$ is parallel to $\nu(z)$ for all $z\in L_0^\varepsilon$. 
As $\varepsilon$ vanishes, we may assume that $p_\varepsilon$ approaches to
a point $p$ that lies in the line segment $(p_0,p_1)$ on $\partial N$. 
We may also assume that
$\nu(z)=(p_\varepsilon-p_0)/\vert p_\varepsilon-p_0\vert$ for all $z\in 
L_0^\varepsilon$. The other possibility
$\nu(z)=-(p_\varepsilon-p_0)/\vert p_\varepsilon-p_0\vert$ can be proved similarly. Then we have for all $z\in L_0^\varepsilon$
\[ \langle \Phi(z), v(z)\rangle=\langle \widehat{\nabla} \sigma(p_0; p_\varepsilon-p_0), (p_\varepsilon-p_0)\rangle/\vert p_\varepsilon-p_0\vert.
\]
Thus by Corollary \ref{ortoasymptotes} we have 
\begin{equation}\label{Uepsilon3}
\lim_{\varepsilon\to 0} \int_{L_0^\varepsilon}\langle \Phi(z), \nu(z)\rangle \varphi(z)\, dS(z)= \int_L \varphi(z)(\sigma(p)-\sigma(p_0))/\vert p-p_0\vert \, dS(z),
\end{equation}
where $L$ is the ray starting from $z_j$. 
Similarly, we have 
\begin{equation}\label{Uepsilon4}
\lim_{\varepsilon\to 0} \int_{L_1^\varepsilon}\langle \Phi(z), \nu(z)\rangle \varphi(z)\, dS(z)= \int_L \varphi(z)(\sigma(p)-\sigma(p_1))/\vert p-p_1\vert \, dS(z).
\end{equation}
Since $\sigma$ is an affine  function in the line segment $(p_0,p_1)$ and $p\in (p_0,p_1)$, we know that
\begin{equation}\label{Uepsilon5}
(\sigma(p)-\sigma(p_0))/\vert p-p_0\vert +(\sigma(p)-\sigma(p_1))/\vert p-p_1\vert=0.
\end{equation}
Thus (\ref{Uepsilon1}) and (\ref{Uepsilon2}) give 
\begin{equation}\label{Uepsilon6}
 \int_W \langle \Phi(z), \nabla \varphi(z)\rangle\, dz \, = \; \lim_{\varepsilon\to 0}\left(\int_{L_0^\varepsilon}\langle \Phi(z), \nu(z)\rangle \varphi(z)\, dS(z)+ \int_{L_1^\varepsilon}\langle \Phi(z), \nu(z)\rangle \varphi(z)\, dS(z)\right)=0.
\end{equation}
Thus $\Phi$ is divergence free in $W$. In conclusion, we have now proven that under the above criteria $h$ is the minimizer.

\begin{Thm}\label{thm:hstar:minimizer}
Suppose $\Omega\subset \mathbb{R}^2$ a bounded Lipschitz domain, $\sigma$ 
a convex function satisfying (\ref{Pst2}) and  $h\in \mathscr{A}_N(\Omega, h_0)$. 

Assume $h$ has a liquid region $\LL$ with frozen boundary, i.e. it is a solution of the Euler-Lagrange equation \eqref{equ:hstar}  and \eqref{def:liquiddomain2} - \eqref{frozen.def} hold for $h$  in a subdomain $\LL \subset \Om$. Assume also  that $h$ enjoys the 
frozen star ray property in each component of $\Omega\setminus \overline{\LL}$. Then $h$ is the minimizer of variational problem
(\ref{ELbasic}) among the class $\mathscr{A}_N(\Omega, h_0)$.
\end{Thm}
  
The above argument and  proof of Theorem \ref{thm:hstar:minimizer}
 makes it very suggestive that there is a strong connection between our approach to Theorem \ref{thm:hstar:minimizer}, and the so called tangent method,  see e.g. in \cite{CS} and \cite{DiFG}.

\subsubsection{Proof of Theorem \ref{Second.thm}.}\label{subsection:proofoftheoremmain3}

Given a polygonal natural domain $\Om$ and natural boundary values $h_0$ on $\partial \Om$, one approach
 to understand the corresponding minimizer of \eqref{ELbasic} and the liquid domain $\LL$ it possibly creates, 
 is to find 
 a good ``candidate minimizer''  $h \in \mathscr{A}_N(\Omega, h_0)$. The difficult task is then  to show that this candidate does indeed minimize the integral in \eqref{ELbasic}, among the admissible functions in $\Om$.  In many concrete situations our
Theorem \ref{thm:hstar:minimizer} appears here a flexible tool.

On the other hand,  one can use Theorem \ref{thm:hstar:minimizer} to a converse direction. 

\begin{Thm}\label{Nminimizer} Suppose $N \subset \R^2$ is a convex polygon and $\sigma$ a surface tension as in \eqref{Pst2}, without gas or quasifrozen points in $N^\circ$. 
Assume also that $\, h^\star$ is a Lipschitz solution of the Euler-Lagrange equation \eqref{equ:hstar} in a bounded Jordan domain $\mathcal U \subset \C$. 

If $\partial \mathcal U$ is frozen for $h^\star$, i.e. if $\, \nabla h^\star: \mathcal U \to N^\circ$ is a proper map, then there is a natural polygonal domain $\Omega\supset \mathcal U$ and a natural boundary value $h_0$ on $\partial \Omega$ such that $\mathcal U = \LL$,  the liquid domain  of the minimizer $h$ for the variational problem \eqref{ELbasic} among $\mathscr{A}_N(\Omega,h_0)$. Also, $h_{| \mathcal U} = h^\star$.
\end{Thm}

\begin{proof} We use the representation of Theorem \ref{repone}, 
 \begin{eqnarray} \label{general2}
\nabla h^* (z) = \sum_{j=1}^k \, p_j \, \omega_{\Di}\bigl( f(z); I_j \bigr), \qquad  z \in \mathcal U,
      \end{eqnarray} 
      where $f:\mathcal U \to \Di$ is a proper map solving $f_{\overline{ z }} = \mu_{\sigma} (f) f_z$. Here  $\mu_\sigma: \Di \to \Di$ is a proper analytic map, thus in fact a finite Blaschke product.
    The $\{ p_j \}$ are the corners of $N$, and as usual, $I_j \subset \partial \Di$ are open intervals with closures  covering the unit circle, and $\omega_{\Di}(\zeta; I_j)$ is the harmonic measure of $I_j$ in the unit disc.
    In particular, Theorem \ref{MultiplyLocus} applied to $\widehat f =  \mu_{\sigma} (f)$
tells that $\partial \mathcal U$ is the real locus of an algebraic curve, thus real analytic outside its cusps  and tacnodes. 
Note that tacnodes are ruled out by our assumption that $\mathcal U$ is a Jordan domain.

By Theorem \ref{First.thm},  we have $f\in C(\overline{\mathcal U})$. 
The inverse image under $f$ of each interval $I_j \subset \partial \Di$  has finitely many  connected components, pairwise disjoint open arcs $\Gamma_{j,\ell} \subset \partial \mathcal U$,    
   \begin{eqnarray} \label{arcs.on.bdry}
    f^{-1}(I_j ) = \bigcup_{\ell} \Gamma_{j,\ell}, \qquad {\rm with} \quad \nabla h(z) \to p_j \quad {\rm as } \quad z \to \Gamma_{j,\ell}.
             \end{eqnarray} 
If $z_0 \in \partial \mathcal U$ is a common endpoint of  arcs $\Gamma_{j,n}$ and $\Gamma_{j+1,\ell}$, then
 the tangent of $\partial \mathcal U$ at $z_0$ is orthogonal to $p_{j+1} - p_j$.

If, say, $\Gamma = \Gamma_{j,k}$ is one of the connected components,  consider a subdivision of the arc into a finite number of points $\{z_1,z_2,..., z_{m}\}\subset \Gamma$.
Let us then connect  outside of $\mathcal U$ the point $z_n$ to $z_{n+1}$, $n=1,2,...,m-1$,
by a line segment orthogonal to $p_j - p_{j-1}$ followed by  a  segment orthogonal to $p_{j+1}- p_j$. In this way, we have a zig-zag curve joining the endpoints of $\Gamma$, which will define part of the boundary of natural domain $\Omega$, as in Figure \ref{zigzag4} below. For the figure we use the standard lozenge triangle $N$ with vertices $\{(0,0),(0,1),(1,0)\}$.

\begin{figure}[H] 
\label{zigzag4} 
\centering{}
\includegraphics[scale=0.34]{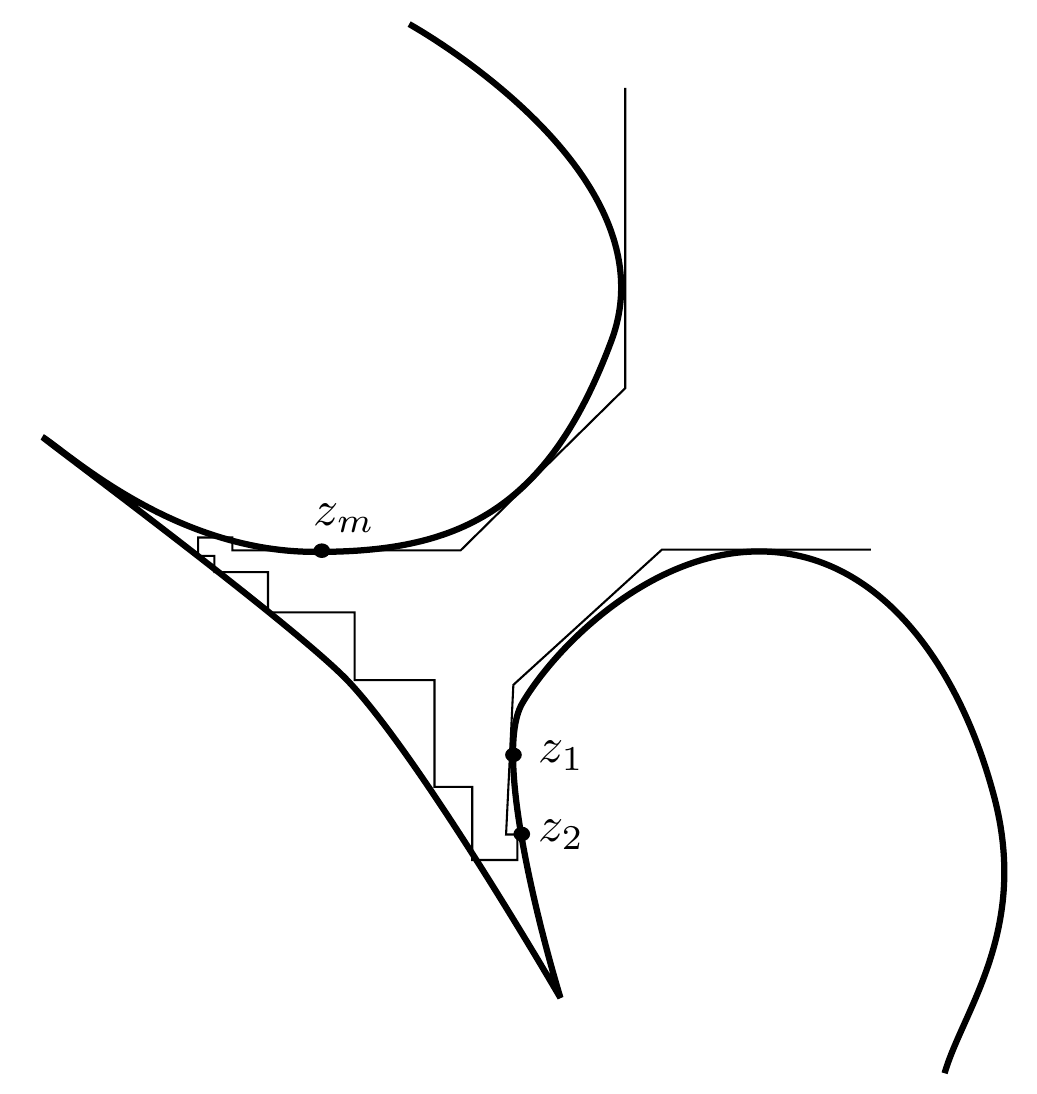}
\caption{Construction of natural domain $\Omega$}
 \label{zigzag4}
\end{figure}

In the extended regions, we define $h^\star$ to be the affine function with $\nabla h^\star=p_j$,
so  that it coincides with $h^\star$ on $\Gamma$.  
In the above procedure, we are allowed to choose the subdivision. We will avoid cusps and make sure that there is at most one cusp
in each part of $\Gamma$ joining $z_n$ and $z_{n+1}$, $n=1,2,...,m-1$. 
We do the above construction for all connected components 
of $  f^{-1}(I_j )$, $j=1,\dots,k$, with the obvious modifications. We also define $h^\star$
in the extended areas by continuity and $\nabla h^\star=p_j$.  

Let us first assume that the end-points $z_1$ and $z_m$ are not cusp points. By choosing the sudivision $\{ z_1, z_2,..., z_m\}$ fine enough, we can make sure that $\partial \Omega$ has no self-intersection. If, say, the end-point $z_m$ is a cusp point, that is, the tangent at the cusp is orthogonal to $p_{j+1}-p_j$, then the last zig-zag segment $J$ in the construction is tangential to the cusp. When we proceed to the next facet corresponding to $p_{j+1}$, the first zig-zag segment in the construction is also tangential to the cusp and can be arranged to coincide with the segment $J$. In this case we remove the double segment $J$ from the construction of $\partial \Omega$ to ensure that $\partial \Omega$ is a polygon. After doing the construction for all the arcs $\Gamma_{j,\ell}$, we thus have a natural domain $\Omega$, $h^\star$ well-defined in $\Omega$ and $h^\star |_{\partial \Omega}$ defining natural boundary values.

After this extension operation we see that $h^\star$ enjoys the properties i), ii) and iii) in Definition \ref{def:ray} in each component of $\Omega\setminus \overline{\mathcal U}$.
The condition iv) is not satisfied but in the areas, disjoint triangles, which the rays do not cover we extend $\Phi$ continuously, to be a constant vector in each.
This gives us a decomposition $\Om =  \Om_0 \cup T_1 \cup \cdots \cup T_\ell$, a disjoint union where $\Om_0$ has the star ray property, with conditions i)-iv), and $\Phi$ is constant   in each  triangle $T_n$. Also, $\Phi$ is continuous across the common boundary points of $T_n$ and $\Om_0$.

Clearly, the vector field $\Phi$ thus defined satisfies the assumptions in Proposition \ref{prop:divergefree}.
Thus $h^\star$ is the minimizer of variational problem (\ref{ELbasic}) among the class $\mathscr{A}_N(\Omega, h_0)$ with $h_0=h^\star$ on $\partial \Omega$. By the construction, the liquid region of $h^\star$ is $\mathcal U$. This proves the theorem.
\end{proof}

\begin{rem} In case of tacnodes, the argument of Theorem \ref{Nminimizer} gives a polygonal domain $\Om  \supset \LL$ (with boundary not necessarily natural) and an extension $h^\star$ with piecewise affine boundary values $h_0$,  minimizing (\ref{ELbasic}) within $\mathscr{A}_N(\Omega, h_0)$.
\end{rem}
\smallskip

 In conclusion, combining the above arguments with our previous results on the universal Beltrami equation \eqref{beltrami2345} give us 

{\bf Proof of Theorem \ref{Second.thm}}: Suppose $\LL$ is a bounded Jordan domain and $\sigma$ is some surface tension as in \eqref{Pst}, with $\Gg = \emptyset$.

 If  $h \in C^1(\LL)$ is a solution to the Euler-Lagrange equation 
\, $  \rm{div} \big(\nabla \sigma(\nabla h)\big)= 0$ in $\LL$, such that $\nabla h: \LL \to N^\circ$ is a proper map,  then Theorem \ref{key.connection} shows that the domain $\LL$ also supports  a proper map  $f:\LL \to \Di$ which is  a solution  to the Beltrami equation $\partial_{\overline{ z }} f(z) = f(z)  \partial_z f(z)$. 

In this setup, since $\LL$ is simply connected by assumption, combining Corollary \ref{Dconverse} with  Remark \ref{domH} - see also  the discussion after
Remark \ref{domH} - now gives us in $\LL$ a $C^1$-solution $h_0$ to 
$$  \rm{div} \big(\nabla \sigma_{_{Lo}}(\nabla h_0)\big)= 0,$$
where $\sigma_{_{Lo}}$ is  the surface tension \eqref{gradsigma23} for the standard lozenges model and $\nabla h_0 : \LL \to N_{_{Lo}}$ is  a proper map. Finally, since the lozenges surface tension $ \sigma_{_{Lo}}$ has no gas points, combining this with Theorem  \ref{Nminimizer} concludes the proof of Theorem \ref{Second.thm}.

\begin{rem}\label{loc.reg} 
Suppose we have a solution $h$ to \eqref{eq:EL34} in a domain $\LL$, such that a part  $\Gamma \subset \partial \LL$ is frozen for $h$, i.e. $\nabla h(z) \to \partial N \cup \Gg$ as $z \to \Gamma$. Then the argument  above extends  $h$ locally  across $\Gamma$ to a partially polygonal domain $\Om \supset \LL$. The extension $h^\star$ is piece-wise affine in $\Om \setminus \LL$ and again minimizes the integral $\int_{\Om}\sigma(\nabla h)dx$ among  $\mathscr{A}_N(\Omega, h_0)$.

 In this way, arguing as in Theorem  \ref{Second.thm} and using Remark \ref{locfrozen} in the case where properness holds only for a part of the boundary,  we have the universality in the geometry of locally frozen boundaries: Locally any frozen boundary in any dimer model, with or without gas or quasifrozen points, is the locally frozen boundary in the lozenges model. In particular, here there are no connectivity restrictions for the liquid domain in question.
\end{rem}

\addtocontents{toc}{\vspace{-4pt}}
\section{Universal Edge Fluctuation Conjecture} \label{conjectures}

 Dimer models are so called \emph{free fermionic} or \emph{determinantal models}. Unfortunately, to give the definition here would require a lengthy detour, and the reader is instead referred to the papers \cite{G19,JohHav,Ke} for an introduction to this topic. It has been shown in many one and two dimensional statistical models that determinantal point process, once suitably rescaled, converges to so called universal statistical processes. Again, to give a precise meaning of this would be rather technical and the readers are again referred to the papers \cite{G19,JohHav}. In view of the classification of the local regularity of the minimizers at the frozen boundary of dimer models in Theorem \ref{thm:main2} 
 one expects that this imposes a strong rigidity of which local universal random scaling limits that can occur at the frozen boundary. Indeed, we have the six different cases which, in view of Theorem \ref{thm:main2}, cover all possible situations:

\begin{itemize}\label{cases}
\item[(1)] $z_0$ is a smooth point of $\dv \LL$ and $\nabla h$ is Hölder continuous with exponent $1/2$ up to the boundary from inside the liquid domain.  
\item[(2)] $z_0$ is a first order cusp and $\nabla h$ is Hölder continuous with exponent $1/3$ up to the boundary from inside the liquid domain. 
\item[(3)] $z_0$ is a first order tacnode and $\nabla h$ is Hölder continuous with exponent $1/2$ up to the boundary from inside the liquid domain. 
\item[(4)]  $z_0$ is a smooth point of $\dv \LL$ and $\nabla h$ is {\bf not continuous} up to the boundary from inside the liquid domain.    
\item[(5)]  $z_0$ is a first order cusp and $\nabla h$ is {\bf not continuous} up to the boundary from inside the liquid domain. 
\item[(6)] $z_0$ is a first order tacnode and $\nabla h$ is {\bf not continuous} up to the boundary from inside the liquid domain. 
\end{itemize}

It has been shown for special natural domains with suitable assumptions that scaling limits in the respective cases are
\begin{itemize}\label{scaling}
\item[(i)] The extended Airy process. See for example the papers \cite{J13,Duse15c,Pet12, AH21}.
\item[(ii)] The Pearcey process. See for example the papers \cite{OR07,ACvM06, HYZ}.
\item[(iii)] The Tacnode process. See for example the paper \cite{AJM14}.
\item[(iv)] The GUE corner process. See for example the papers \cite{JN06,OR06, AG}.
\item[(v)] The Cusp-Airy process. See for example the paper \cite{Duse15d}.
\item[(vi)] The discrete Tacnode process. See for example the papers \cite{AJvM1,AJvM2}. 
\end{itemize}

The cases $(1)-(3)$ corresponds to so called \emph{continuous-continuous} point processes, whereas cases $(4)-(6)$ corresponds to so called \emph{continuous-discrete} point processes, see the discussion in \cite{JohHav}. More importantly, the cases (1)-(3) are \emph{stable} in the sense one expects that they do not depend on the convergence rate of the underlying discrete hexagon graph to the natural domain. On the other hand, the cases $(4)-(6)$ require much stronger convergence assumptions. 

We now state our conjecture, see \cite{AG, AH21, HYZ} for the recent  progresses for this conjecture.

\begin{Con}[Universal Edge Fluctuation Conjecture]
Let $\Om$ be a natural domain with natural boundary values for a dimer model, without gas and quasifrozen points. Then under suitable assumptions of the convergence of the underlying  bipartite graph $G_n$, converging to $\Om$ in a suitable sense, all scaling limits of the determinantal point processes at the frozen boundary are given by the above list (i) -(vi) 
corresponding to the local regularity cases $(1)-(6)$ in Section \ref{cases} and interpolation processes of these. 
\end{Con}

kari.astala@helsinki.fi, \quad erik.duse@kth.se, \quad istvan.prause@abo.fi, \quad xiao.zhong@helsinki.fi

\end{document}